\documentclass{amsart}
\usepackage{amssymb,mathrsfs}
\usepackage{color,umoline}
\usepackage[dvipsnames]{xcolor}
\usepackage{graphicx}
\usepackage{enumitem}
\usepackage[font=small,labelfont=bf]{caption}
\usepackage{tikz, caption}
\usepackage{relsize}
\usepackage[mathscr]{eucal}
\usepackage{dsfont}
\usepackage{verbatim}
\usepackage[draft]{todonotes}
\usepackage[
 citecolor=black,  
  hypertexnames=false]{hyperref}
\usepackage[dvipsnames]{xcolor}
\usetikzlibrary{arrows.meta}
\usepackage{subfig}

\def\wt#1{\widetilde{#1}}

\newcommand{\N}{\mathbb{N}}
\newcommand{\R}{\mathbb{R}}

\newcommand{\x}{\xi}

\newcommand{\supp}{\operatorname{supp}}
\newcommand{\dist}{\operatorname{dist}}

\newtheorem{thm}{Theorem}[section]
\newtheorem{defn}[thm]{Definition}
\newtheorem{prop}[thm]{Proposition}
\newtheorem{cor}[thm]{Corollary}
\newtheorem{lem}[thm]{Lemma}
\numberwithin{equation}{section}
\newtheorem{rem}{Remark}

\newcommand{\inp}[2]{\langle #1, #2\rangle}

\newcommand{\mmf}{(\mu\mu')^\frac14}

\newcommand{\mmt}{\sqrt{\mu\mu'}}
\newcommand{\mm}{{\mu\mu'}}
\newcommand{\Sa}{S_\ast(x,y)}
\newcommand{\D}{\mathcal D}

\newcommand{\sxy}{S_c(x,y)}
\newcommand{\sxyl}{2^{-l}s+\sxy}

\newcommand{\vepc}{\varepsilon_{\!\circ}}

\newcommand{\ol}[1]{\overline{#1}}

\newcommand{\smu}{\sqrt\mu}
\newcommand{\smp}{\sqrt{\mu'}}
\newcommand{\mom}{\big({\mu'}/\mu\big)}

\newcommand{\mup}{(\mu')}

\newcommand{\pnchi}{\chi_{k}^{+,\nu}}

\newcommand{\pxkl}{\chi_{k}^{\pm,\nu}}
\newcommand{\spxkl}{\chi_{k}^{\sigma,\nu}}

\newcommand{\spxklp}{\chi_{k'}^{\sigma', \nu}}
\newcommand{\pxklp}{\chi_{k'}^{\pm, \nu}}
\newcommand{\fP}{\mathfrak P}
\newcommand{\fq}{\mathfrak q}
\newcommand{\tr}{\intercal}
\newcommand{\fc}{\mathbf c}
\renewcommand{\sxy}{S_c}
\renewcommand{\sxyl}{S_c^l}

\newcommand{\ksim}{k\sim_\nu k'}

\newcommand{\cP}{\mathcal P}
\newcommand{\dpq}{\delta(p,q)}

\newcommand{\che}{\chi_{\sqrt{\lambda}E}}

\newcommand{\ppqpair}{$(1/p, 1/q)$}

\newcommand{\ppq}{(1/p,1/q)}
\newcommand{\Be}{\begin{equation}}
\newcommand{\Ee}{\end{equation}}
\newcommand{\scalep}{\frac{d-1}2\dpq-\frac d2}

\newcommand{\scaleppp}{\frac{d-1}2\dpp-\frac d2}
\newcommand{\scalepq}[2]{\frac{d-1}2\delta(#1,#2)-\frac d2}

\newcommand{\mU}{\mathbf U}
\newcommand{\cQ}{\mathcal Q}
\newcommand{\cD}{\mathcal D}
\newcommand{\schi}{\widetilde \chi_{\mathsmaller {\sqrt{\lambda}}}}

\newcommand{\ddk}{\partial_s^k \cP}
\newcommand{\fA}{\mathbb  A}

\newcommand{\bpqp}{{\mathrm B}_{p,p'}(\lambda,\mu, \mu’)}

\newcommand{\cB}{\mathcal B}
\newcommand{\dpp}{\delta(p,p')}
\newcommand{\sxyls}{\sxyl(x,y,s)} 

\newcommand{\eval}{2\mathbb N_0+d}
\newcommand{\delc}{\delta_\circ}
\newcommand{\rp}{1/p}
\renewcommand{\rq}{1/q}

\newcommand{\cI}{\mathcal I}
\newcommand{\sxylt}{\sxyl(x,y, t)}

\newcommand{\su}{\sin{S_c}}
\newcommand{\cu}{\cos S_c}
\newcommand{\sul}{\sin\sxyl}
\newcommand{\cul}{\cos\sxyl}
\newcommand{\sut}{\sin^2 S_c}
\newcommand{\sult}{\sin^2\sxyl}
\newcommand{\ba}{\mathbf a}
\newcommand{\bc}{\mathbf c}

\newcommand{\oo}[1]{O(#1)}

\newcommand{\chpl}{\chi^\pm_{\lambda, \mu}}
\newcommand{\chplpp}{\chi^\pm_{\lambda, \mu'}}
\newcommand{\chplp}{\chi^+_{\lambda, \mu}}
\newcommand{\chpln}{\chi^-_{\lambda, \mu}}

\newcommand{\aas}{A_\mu^\sigma\times A_{\mu'}^{\sigma'}}
\newcommand{\aask}{A_k^{\sigma, \nu}\times A_{k'}^{\sigma', \nu}}
\newcommand{\aack}{A_k^{\circ, \nu}\times A_{k'}^{\circ, \nu}}
\newcommand{\spp}{{\sigma'}}

\newcommand{\diam}{{\rm diam }} 
\newcommand{\Ac}{A^\circ}

\newcommand{\bp}{\cB_{p,j}}

\newcommand{\ps}{p_\ast}

\newcommand{\mmfpj}{\chi_{\lambda,\mu}^e\Pi_\lambda [\psi_j]\chi_{\lambda,\mu'}^e}
\newcommand{\mmfp}{\chi_{\lambda,\mu}^e\Pi_{\lambda'}\chi_{\lambda,\mu'}^e}

\newcommand{\nuc}{{\nu_\circ}}
\newcommand{\proj}{{\mathlarger \wp}_k}
\newcommand{\wproj}{\widetilde{\mathlarger \wp}_k}
\begin{document}

\author[Jeong]{Eunhee Jeong}
\address[Jeong]{Department of Mathematics Education and  Institute of Pure and Applied Mathematics, Jeonbuk National University, Jeonju 54896, Republic of Korea}
\email{eunhee@jbnu.ac.kr}

\author[Lee]{Sanghyuk Lee}
\address[Lee]{Department of Mathematical Sciences and RIM, Seoul National University, Seoul 08826, Republic of  Korea}
\email{shklee@snu.ac.kr}

\author[Ryu]{Jaehyeon Ryu}
\address[Ryu]{Department of Mathematics and Institute of Pure and Applied Mathematics, Jeonbuk National University, Jeonju 54896, Republic of Korea}
\email{jhryu67@jbnu.ac.kr}

\title{
Hermite  spectral  projection operator}

\begin{abstract} 
We study $L^p$--$L^q$ estimate for the spectral projection operator $\Pi_\lambda$  associated to the Hermite operator $H=|x|^2-\Delta$ in $\mathbb R^d$. 
Here $\Pi_\lambda$ denotes the projection to the subspace spanned by the Hermite functions  which are the eigenfunctions of $H$ with eigenvalue $\lambda$. Such estimates were previously available only for $q=p'$, equivalently with $p=2$ or $q=2$ (by $TT^*$ argument) except for the estimates which are straightforward consequences of interpolation between those estimates. As shown in the works of Karadzhov, Thangavelu, and Koch and Tataru, the local and global estimates for $\Pi_\lambda$ are of different nature. Especially, $\Pi_\lambda$ exhibits complicated behaviors near the set $\sqrt\lambda\mathbb  S^{d-1}$. Compared with the spectral projection operator associated to the Laplacian, $L^p$--$L^q$ estimate for $\Pi_\lambda$ is not so well understood up to now for  general $p,q$. In this paper we consider  $L^p$--$L^q$ estimate for $\Pi_\lambda$ in a general framework including the local and global estimates with $1\le p\le 2\le q\le \infty$ and undertake the work of characterizing  the sharp bounds on $\Pi_\lambda$. We establish various new sharp estimates in extended ranges of $p,q$. First of all, we  provide a complete characterization of the local estimate for $\Pi_\lambda$ which was first considered by Thangavelu. Secondly, for $d\ge5$, we prove the endpoint $L^2$--$L^{2(d+3)/(d+1)}$ estimate for $\Pi_\lambda$ which has been left open since the work of Koch and Tataru. Thirdly, we extend the range of $p,q$ for which the operator $\Pi_\lambda$ is uniformly bounded from $L^p$ to $L^q$.
\end{abstract}

\maketitle

\section{Introduction}  
Let $H$ denote the Hermite operator $-\Delta+|x|^2$ on $\mathbb R^d$. The operator $H$ has a discrete spectrum $\lambda\in 2\mathbb N_0+d$, where $\mathbb N_0=\mathbb N\cup \{0\}$.  For $\alpha\in \mathbb N_0^d$, let $\Phi_\alpha$ be the $L^2$--normalized Hermite function
which is given as a tensor product of the Hermite functions in $\mathbb R$. The set  $\{\Phi_\alpha: \alpha\in \mathbb N_0^d\}$  forms an orthonormal basis in $L^2$ and  $\Phi_\alpha$ is   
an  eigenfunction of $H$ with eigenvalue $2|\alpha|+d$  where $|\alpha|=\alpha_1+\cdots+\alpha_d$.  We denote by $\Pi_\lambda$  the spectral  projection operator to the vector space  spanned by eigenfunctions with  the eigenvalue $\lambda$, which   is given by
\begin{align*}
    \Pi_\lambda f = \underset{\alpha:d+2|\alpha|=\lambda}{\sum}\langle f,\Phi_\alpha \rangle \Phi_\alpha. 
\end{align*}
Then, for $f\in L^2$ we  have  $f=\sum_{\lambda\in 2\mathbb N_0+d}\Pi_\lambda f$, which is in fact  the Hermite expansion.

For an operator $T$ which maps any function in the Schwartz class $\mathcal S(\mathbb R^d)$ to a  measurable function in $\mathbb R^d$, we define 
\[ \|T\|_{p\to q}=\sup\{ \|Tf\|_q : \|f\|_p\le 1, f\in \mathcal S(\mathbb R^d) \}. \]  
Let $E, F\subset \mathbb R^d$ be measurable subsets. 
In this paper, we are concerned with  the following form of estimate: 
\begin{align}\label{rstest}
    \|\chi_E \Pi_\lambda   \chi_F \|_{p\to q}\le B .
\end{align}
Here, $B$ is a constant  which may depend on $d,p,q,E,F$ and $\lambda$.  The sharp bound $B$ in terms of $\lambda$ has  been of interest in connection with Bochner-Riesz summability of the Hermite operator $H$. 
See, Askey and Wainger \cite{AW},  Karadzhov \cite{K94}, and Thangavelu \cite{Th98}.  In particular, the uniform estimate $\eqref{rstest}$ with $E=F=\mathbb R^d$ and $B$ independent of $\lambda$ has applications to unique continuation problem for the parabolic operator.  In the work of Escauriaza \cite{e00},  the estimate due to Karadzhov \cite{K94} was used  as a main tool  to prove  Carleman estimates  for the heat operator.

We may compare the operator $\Pi_\lambda$ with the spectral projection  associated to $-\Delta$ of which boundedness is better understood. 
The operator 
\Be 
\label{proj-laplace}
{\mathlarger \wp}_k f= \frac1{(2\pi)^d}\int_{k-1\le |\xi|^2<k}e^{ix\cdot\xi}\,\widehat f(\xi)\, d\xi 
\Ee  can be regarded as a natural counterpart of $\Pi_\lambda$ where $k\ge 1$. We have rather complete characterization of the bound $\|\proj\|_{p\to q}$ for $(p,q)\in [1,2]\times [2, \infty]$ (see Corollary \ref{laplacian}) thanks to  the well-understood $L^p$--$L^q$ bound on  the restriction-extension operator $R^*\!Rf:=( \widehat f\, |_{ \mathbb S^{d-1}})^\vee$\footnote{Here, $R$ denotes the restriction operator to the sphere, i.e.,  $g\to \widehat g|_{\mathbb S^{d-1}}$ and $R^*$ is the adjoint operator of $R$.}(see Theorem \ref{rr*}). In fact,  using the spherical coordinates, the boundedness of 
${\mathlarger \wp}_k $ can be deduced from that of  $R^*\!R$ (see Section \ref{transplantation} for more details).  
However, the estimate \eqref{rstest}  has mainly been studied so far with $p=2$,  equivalently, $q=p'$ or $q=2$ (by  $TT^*$ argument and duality) except the estimates which follow by interpolation between those estimates. The estimates in  the works of Karadzhov \cite{K94}, Thangavelu \cite{Th98}, and 
Koch and Tataru \cite {T05} all fall in this category and  the estimate \eqref{rstest} has not been considered with general $(p,q)$.
In this paper we investigate the bound \eqref{rstest} in a general framework considering $(p,q)\in [1,2]\times [2, \infty]$.

Estimates for spectral projections (to  each eigenspace or  spectral window of fixed length) and, more generally,  spectral measures which are associated to specific differential operators have been widely studied.  The most typical example is the estimate for the operator $R^*\!R$,  
\begin{align}\label{est:euclid}
\Big  \| \int_{\mathbb S^{d-1} } e^{ix\cdot \omega} \widehat f(\omega) d\sigma(\omega) \Big \|_q\lesssim \|f\|_p.
\end{align} 
When $(p,q)=\big(\frac{2(d+1)}{d+3}, \frac{2(d+1)}{d-1}\big)$, the estimate is equivalent to the Fourier restriction estimate for the sphere known as the Stein-Tomas theorem \cite{tomas}.
The estimates for the spectral projection operators which are given by the spherical harmonics were obtained by  Sogge \cite{sogge-1986} (also, see \cite{KL18})
and his result was later  extended to the projection operators defined by the eigenfunctions of  the elliptic operator on compact 
Riemannian manifold \cite{sogge-1988} which forms a basis for $L^2$. These estimates may be regarded as a variable coefficient version of the restriction theorem (see  \cite{so93} for a detailed discussion), 
and have a broad range of applications to related problems  such as Bochner--Riesz summability of eigenfunction expansions \cite{sogge-1987} and  unique continuation problems  \cite{sogge-1990, T01}.  See \cite{GHS}  for results concerning   operators  on non-compact manifold with asymptotic assumption and applications of spectral projection estimates to the study of spectral multipliers \cite{SS, CQSY}.   There are also results on the spectral projection for the twisted Laplacian \cite{SZ, KR}.

\subsection*{Global estimate}  We first consider the global estimate, that is to say, \eqref{rstest} with $E=F=\mathbb R^d$. 
When $d=1$,  the sharp bound on $\|\Pi_\lambda\|_{p\to q}$  is easy to obtain by making use of  the bound on $L^p$ norm of   the 
Hermite function in $\mathbb R$ (see Lemma \ref{S4-bd_Her}) and sharpness of these bounds  follows from duality and Lemma \ref{S4-bd_Her}.  
Hence,  throughout this paper, we assume 
 $d\ge 2$ unless it is mentioned otherwise. In higher dimensions  the problem is no longer trivial except some special cases. 
 In particular, there are easy bounds such as $\|\Pi_\lambda\|_{2\to 2}\le 1$ and 
 \begin{equation}  
 \label{plambda-100} 
 \|\Pi_\lambda \|_{1\to \infty} \lesssim \lambda^{\frac {d-2}2}
  \end{equation}
that follows from  \cite[Lemma 3.2.2]{Th93}. Bounds for the other cases are not straightforward and require different approach. Compared with the 
 projection operator $\proj$, the boundedness  properties  of the projection operator $\Pi_\lambda$ exhibits more complicated behavior.

\subsubsection*{Uniform bound for $\Pi_\lambda$} In \cite{K94}, Karadzhov showed that there exists a uniform constant $\mathcal B$ independent of  $\lambda$ such that
\begin{align}\label{karad-global}
    \|\Pi_\lambda \|_{2\to \frac{2d}{d-2}}\le \mathcal B. 
    \end{align}
      Making use of  \eqref{karad-global}, he showed the sharp $L^p$--Bochner-Riesz summability of the Hermite expansion for $p\ge \frac{2d}{d-2}$ and $p\le \frac{2d}{d+2}$.  Furthermore,  the estimate \eqref{karad-global} was recently  utilized by Chen, the second named author, Sikora, and Yan  \cite{CLSY} to prove  $L^p$ boundedness of the maximal  Bochner-Riesz means for the Hermite expansion,  and  by  Chen, Li, Ward, and Yan \cite{CLWY} to obtain  the endpoint weak type estimate at the critical summability index.
     As mentioned in the above, this kind of  uniform estimate has  been used to study  unique continuation problem for  parabolic equations (\cite{e00, ev01}). 
By duality and $TT^*$ argument,  the estimate \eqref{karad-global}  is equivalent to the uniform bounds  $\|\Pi_\lambda \|_{ {2d}/(d+2)\to 2}\le \mathcal B$ and 
$\|\Pi_\lambda\|_{2d/(d+2)\rightarrow 2d/(d-2)}\le \mathcal B^2$. Interpolation gives  $\|\Pi_\lambda \|_{p\to q}\le C$ with a uniform constant $C$ for $p,q$ satisfying $\frac{2d}{d+2}\le p\le 2\le q\le \frac{2d}{d-2}$. Escauriaza and Vega \cite{ev01} made use of those uniform estimates to prove the Carleman estimate for the heat 
operator on  an extended collection of  mixed norm spaces. We refer the reader to \cite{ev01,ef03,Fe03,T09} for further related developments in this direction. 
Our first  result extends  the previously known range of uniform boundedness.  

\begin{defn} If $X=(a, b)\in [1/2,1]\times[0,1/2]$, we define $X':=(1-b, 1-a)$.  Likewise, we also define
$\mathfrak Z':=\{ X': X\in \mathfrak Z\}$ for  a set $\mathfrak Z\subset [1/2,1]\times[0,1/2]$.  For $X, Y\in [1/2,1]\times[0,1/2]$, $X\neq  Y,$ let $[X,Y]$, $(X,Y)$  denote the closed,  open line segments connecting the points $X$ and $Y,$ respectively. 
The half open  line segments  $(X, Y]$, $[X, Y)$ are similarly defined.   
Let us set 
\begin{align*}
\mathfrak E=  \mathfrak E(d):=\left(\frac{d+2}{2d},\frac 12\right),\; 
\quad 
\mathfrak F=        \mathfrak F(d):=\left(\frac{d^2+2d-4}{2d(d-1)},\frac{d-2}{2(d-1)}\right).
\end{align*}
\end{defn}

\begin{thm}\label{glo-est}  Let  $d\ge 3$. 
Suppose $(1/p,1/q)$ is contained in the closed pentagon $\mathscr P$ with vertices $ \mathfrak E,$  $\mathfrak E'$,  $\mathfrak F$, $\mathfrak F'$, and $(1/2,1/2)$ from which the vertices $\mathfrak F$, $\mathfrak F'$ are removed. Then, there is a constant $C$,  independent of $\lambda$,  such that  the estimate 
\Be
\label{unif}
    \|\Pi_\lambda \|_{p\to q}\le C
\Ee
holds. Moreover,  the estimate $\|\Pi_\lambda f\|_{q,\infty}\le C\|f\|_{p,1}$ holds
with $C$ depending only on $d$ if $(1/p,1/q)= \mathfrak F$, $\mathfrak F'$.  When $d=2$,  the uniform estimate \eqref{unif} holds for all $1\le p\le 2\le q\le \infty$.
\end{thm}

\begin{figure}[t]
\centering
\begin{tikzpicture}[scale=0.7]
\filldraw[fill=gray!35] (0,0)--(6,0)--(6,6)--(0,6);
\draw[<->] (0,7) node[left]{$\frac1q$}--(0,0)--(6.8,0)node[right]{$\frac1p$};
\draw (0,0) rectangle (6,6);
\node[left] at (0,0) {$0$};\node[below] at (0,0) {$\frac12$};
\node[left] at (0,6) {$\frac12$};\node[below] at (6,-0.1) {$1$};
\end{tikzpicture}
\quad \quad 
\begin{tikzpicture}[scale=0.7]
\draw[<->] (0,7) node[left]{$\frac1q$}--(0,0)--(6.8,0)node[right]{$\frac1p$};
\draw (0,0) rectangle (6,6);
\filldraw[fill=gray!35](0,6)--(0,4.5)--(6/7,27/7) --(15/7,36/7)--(1.5,6);
\draw [line width=0.2mm]  (0,6)--(0,4.5)--(6/7,27/7) --(15/7,36/7)--(1.5,6)--(0,6); 
\node[left] at (0,4.5) {$\mathfrak E'$}; \node[above] at (1.5,6) {$ \mathfrak E$}; 
\draw (6/7,27/7) circle [radius=0.07]; \node[below] at (6.5/7,27/7) {$\mathfrak F'$};
\draw (15/7,36/7) circle [radius=0.07]; \node[right] at (15/7,36/7) {$\mathfrak F$}; 
\node[left] at (0,0) {$0$};\node[below] at (0,0) {$\frac12$};
\node[left] at (0,6) {$\frac12$};\node[below] at (6,-0.1) {$1$};
\end{tikzpicture}
\caption{The range of $p,q$ for which \eqref{unif} holds:  $d=2$ (left)  and $d\ge 3$ (right). }
\label{fig:uniform}
\end{figure}
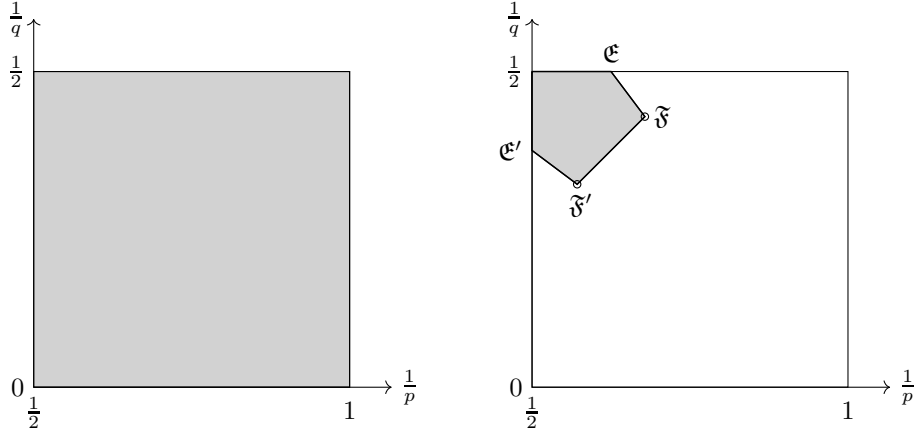

Here,  $\|\cdot\|_{p,r}$ denotes the  norm of  the  Lorentz space $L^{p,r}(\mathbb R^d)$ (for example, see \cite{St71}). To extend the range of uniform boundedness we devise a kind of $TT^*$  argument which combines the typical $TT^\ast$ argument (\cite{keel-tao, fos, vil})  with a representation formula for the projection operator $\Pi_\lambda$ (see \eqref{proj-op}). 
The extension of the range  $p,q$ for \eqref{unif}  to the off-diagonal $p,q$  is in analogue with  enlargement of the  range of pairs $(p,s)$, $(q,r)$ for which  the $L_t^sL_x^p$--$L_t^rL_x^q$ inhomogeneous Strichartz estimate holds. In particular,  the admissible range of the inhomogeneous Strichartz estimate for the Schr\"odinger equation was extended by  Foschi \cite{fos} and Vilela \cite{vil} beyond the range given by  the admissible pairs of the homogeneous Strichartz estimate (see  Keel and Tao \cite{keel-tao}).

When $d=2$, the estimate  \eqref{unif}  is easy to obtain using  duality and the uniform $L^1$--$L^\infty$ estimate \eqref{plambda-100} for $\Pi_\lambda$.   
In higher dimensions, there is a gap between the  range in 
Theorem \ref{glo-est} and the necessary condition on $p,q$ which we currently have. In fact, 
 \eqref{unif} holds  true only  if  $(1/p, 1/q)\in \widetilde {\mathscr P}:=\{(a,b)\in  [1/2,1]\times[0,1/2]:  a-b\le 2/d,   (d-1)/d\le a+b\le (d+1)/d\}$ as can  be easily seen  from 
 duality and 
 the lower bounds \eqref{local2} and \eqref{local3}  in Section 3.  
 So, uniform boundedness of $\|\Pi_\lambda\|_{p\to q}$  remains open when  $(1/p, 1/q)\in   \widetilde {\mathscr P}\setminus \mathscr P$, $d\ge 3$. 
 The  current situation is similar to that of the inhomogeneous Strichartz estimate for the Schr\"odinger equation whose optimal range of boundedness remains open  for $d\ge 3$ (see for example 
 \cite{fos, vil}).

\subsubsection*{$L^2$--$L^q$ estimate} For $(p,q)\in [1,2]\times [2,\infty]$, we set 
 \[ \delta(p,q):=\frac1p-\frac1q.\]
  A systematic study on how  $\|\Pi_\lambda \|_{2\to q}$  (with $2\le q\le \infty$) depends on $\lambda$ was carried out by Koch and Tataru \cite{T05}. They obtained almost  complete  characterization of the $L^2$--$L^q$ estimate.   For $d\ge 2$, it was shown  that 
\[ 
\|\Pi_\lambda \|_{2\to q}\sim \begin{cases}   \lambda^{-\frac12\delta(2,q)},    \ \,  & 2\le q< \frac{2(d+3)}{d+1}, 
                                                                                   \\   \lambda^{-\frac16+\frac d6\delta(2,q)},  \    &  \frac{2(d+3)}{d+1} <q\le \frac{2d}{d-2},  
                                                                                    \\  \lambda^{-\frac12+\frac d2\delta(2,q)} , \  \  &  \frac{2d}{d-2}\le q\le \infty.   \end{cases}   
\]
The estimate for the missing case $q={2(d+3)}/(d+1)$ is  important  in that other estimates can be recovered by interpolation  between the missing estimate,  the trivial $L^2$ bound, 
 and relatively easy \eqref{karad-global} and  $L^2$--$L^\infty$(\eqref{plambda-100})  estimates. 
By means of the estimate for $\Pi_\lambda$ over annulus $A_{\lambda,\mu}^\pm$ in \cite{T05} (see \eqref{kt-nu} and Theorem 
\ref{thm-annest}) and summation over the dyadic annuli, one can easily see  (\cite{T07}) that 
\[
    \|\Pi_\lambda \|_{2\to \frac{2(d+3)}{d+1}}\le C\lambda^{-\frac{1}{2(d+3)}}(\log\lambda)^\frac{d+1}{2(d+3)}.
    \]
    We refer the reader forward to \eqref{def-annalus} for the definition of $A_{\lambda,\mu}^\pm$. 
When $d=1$, using Lemma \ref{S4-bd_Her},  it is easy to see that  such estimate  fails to hold without extra logarithmic factor. However,  it was conjectured  in \cite{T05} that the natural $L^2$--$L^{{2(d+3)}/{(d+1)}}$ estimate  continues to hold without the logarithmic factor  for $d\ge 2$ since it is improbable that there is an eigenfunction concentrating  on every annulus $A_{\lambda,\mu}^\pm$. 

We show the conjecture is true for $d\ge 5$.

\begin{thm}\label{endpoint}
For $d\ge 5$ the  estimate 
\Be\label{intro-conj}
 \|\Pi_\lambda \|_{2\to \frac{2(d+3)}{d+1}}\le C\lambda^{-\frac{1}{2(d+3)}}
\Ee
 holds with $C$ independent of $\lambda$. 
\end{thm}

It is likely that the estimate \eqref{intro-conj} holds for  $2\le d\le 4$ but 
our argument in this paper is not enough to recover the  endpoint estimates in those cases.  
Besides \eqref{intro-conj} we also obtain new sharp global estimates for $p,q$ on a certain range. Those estimates (and interpolation between them) significantly extend 
 the range of $p,q$ for which the sharp estimate holds. However, there  still remain  regions of $p,q$  where we do not have the sharp estimate. 
At the present time  it seems  that  complete characterization of sharp bounds on $\|\Pi_\lambda\|_{p\to q}$  is a difficult problem.

Let $\mu', \mu \in \{ 2^k:k\in \mathbb Z\}$ and  set 
\Be  
\label{def-annalus}
A^{\pm}_{\mu}:=\Big\{x:   \pm(1- |x|)\in [\mu, 2\mu]\Big\}, \quad   A^{\pm}_{\lambda, \mu}:=\Big\{x:   \lambda^{-\frac12} x\in A^{\pm}_{\mu} \Big\}.
\Ee
We also denote
\[ \chi^\pm_\mu = \chi_{A^{\pm}_\mu}, \quad \chi^\pm_{\lambda, \mu}=\chi_{A^{\pm}_{\lambda,\mu}}.\]
The main idea for the proof of Theorem \ref{endpoint} is to consider estimate for the operator 
$\chpl \Pi_\lambda \chplpp$ with $\mu\ge \mu'$   and to establish  \emph{unbalanced}  improvement to its bound: 
\Be
\label{unbalanced}    \|  \chi_{\lambda,\mu}^\sigma\Pi_\lambda  \chi^{\sigma'}_{\lambda, \mu'}  \|_{q'\to q} \lesssim  
\lambda^{-\delta(q',q)}(\mu\mu')^{\frac14-\frac{d+3}{4}\delta(q',q)}  \mom^{c}, \quad \sigma, \sigma'\in \{+, -\} 
\Ee
for some $c>0$ if $2\le q<2(d+1)/(d-1)$  (see Theorem \ref{improvedbound} and the rescaled  equivalent estimate \eqref{eq:main}).   Apart from the additional factor  $(\mu'/\mu)^{c}$ 
the estimate \eqref{unbalanced} follows from the known bound  \eqref{kt-nu} via duality which is again equivalent to the estimate with $\mu=\mu'$ (also, see Theorem \ref{thm-annest}). Indeed, note that  $\|  \chi_{\lambda,\mu}^\sigma\Pi_\lambda  \chi^{\sigma'}_{\lambda, \mu'}  \|_{q'\to q}\le 
\|  \chi_{\lambda,\mu}^\sigma\Pi_\lambda \|_{2\to q}\|  \chi_{\lambda, \mu'}^{\sigma'}\Pi_\lambda   \|_{2\to q}$ because $\Pi_\lambda=\Pi_\lambda^2$.  
When $\mu\sim\mu'$, we can not expect any improvement of the bound since the estimate \eqref{unbalanced} is optimal. Nonetheless, \eqref{unbalanced} tells that  the bound significantly improves
if $\mu\gg \mu'$. This is a new phenomenon which has not been observed before.   
However, as expected,  it  is much  more involved to obtain the bound  \eqref{unbalanced}
with the additional improvement $(\mu'/\mu)^{c}$.   In Section \ref{sec:endpoint}  and Section \ref{proof-end},  a great deal of detailed analysis is devoted to establishing  the estimate \eqref{unbalanced}. 
We provide a brief explanation about our approach near the end of the introduction.

\subsection*{Local  $L^p$--$L^q$ estimate} 
In \cite{Th98}, Thangavelu  considered the estimate for $\Pi_\lambda$ over compact set and  he showed  that such local estimate  has a better bound than that of  the uniform estimate  \eqref{karad-global}.  More precisely, it was shown in \cite{Th98}  that, for any compact set $E\subset\mathbb R^d$, there is a constant $C=C(p,q,d,E)$ such that
\Be  
\label{local-thanga}
\|\chi_ E \Pi_\lambda\chi_E \|_{\frac{2(d+1)}{d+3}\to \frac{2(d+1)}{d-1}}\le C\lambda^{-\frac{1}{d+1}}.\footnote{By $TT^*$ argument, this is equivalent to $ \|\Pi_\lambda\chi_E \|_{\frac{2(d+1)}{d+3} \to 2 }\le C^\frac12\lambda^{-\frac{1}{2(d+1)}}$.}
 \Ee
As to be seen later, the difference between the local and global bounds arises from the behavior of $\Pi_\lambda$ near the set $\sqrt \lambda \mathbb S^{d-1}$. 
We generalize this type of local estimate to any $p,q$ satisfying $1\le p\le 2\le q\le  \infty$, and  obtain   a complete characterization of such estimate  under the condition 
\eqref{det-con} below.  In order to state our results we need to introduce some notations.

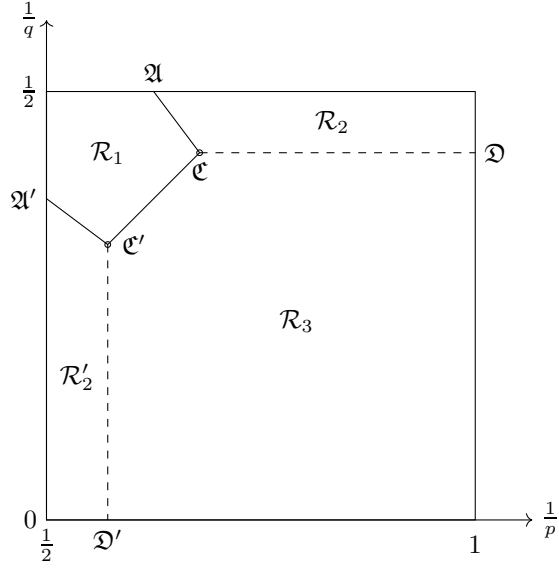
\begin{figure}
\begin{tikzpicture}[scale=0.95]
\draw[<->] (0,7) node[left]{$\frac1q$}--(0,0)--(6.8,0)node[right]{$\frac1p$};
\draw (0,0) rectangle (6,6);
\draw (0,6)--(0,4.5)--(6/7,27/7) --(15/7,36/7)--(1.5,6);
\draw[dashed] (6/7,27/7)--(6/7,0); \draw[dashed] (15/7,36/7)--(6,36/7); 
\node[left] at (0,4.5) {$\mathfrak A'$}; \node[above] at (1.5,6) {$ \mathfrak A$}; 
\node[below] at (6/7,0) {$\mathfrak D'$}; \node[right] at (6,36/7) {$\mathfrak D$}; 
\node[] at (6/7,36/7) {$\mathcal R_1$};
\node[] at (3/7,14/7) {$\mathcal R_2'$}; \node[] at (28/7,39/7) {$\mathcal R_2$};
\node[above] at (3.5,2.5) {$\mathcal R_3$};
\draw (6/7,27/7) circle [radius=0.04]; \node[right] at (6.5/7,27/7) {$\mathfrak C'$};
\draw (15/7,36/7) circle [radius=0.04]; \node[below] at (15/7,36/7) {$\mathfrak C$}; 
\node[left] at (0,0) {$0$};\node[below] at (0,0) {$\frac12$};
\node[left] at (0,6) {$\frac12$};\node[below] at (6,-0.1) {$1$};
\end{tikzpicture}
\caption{The points $\mathfrak A,$ $\mathfrak C $, $\mathfrak D$, and the regions $\mathcal R_1$, $\mathcal R_2,$ $\mathcal R_3$. }
\label{local-type}
\end{figure}

\begin{defn}
\label{points-abcd}
Let $\mathfrak A=\mathfrak A(d),$ $ \mathfrak C=\mathfrak C(d)$, and $\mathfrak D=\mathfrak D(d)\in [1/2,1]\times[0,1/2]$ be  given    by
\begin{align*}
    &\mathfrak A=\left( \dfrac{d+3}{2(d+1)}, \  \dfrac12\right), 
    \  \mathfrak C=\left(\dfrac{d^2+4d-1}{2d(d+1)},\dfrac{d-1}{2d}\right),
   \ \mathfrak D=\left(1,\frac{d-1}{2d}\right).
\end{align*}
Let $\mathcal R_1$ denote the closed pentagon with vertices $(\frac12, \frac12),  \mathfrak A,  \mathfrak C,  \mathfrak C', \mathfrak A'$ from which  
two points $\mathfrak C$ and $\mathfrak C'$ are removed,  and   $\mathcal R_2$ be the closed trapezoid with vertices $\mathfrak A , (1,1/2), $ $\mathfrak D,$ $ \mathfrak C$ 
from which the closed line segment  $[\mathfrak C, \mathfrak D]$ is removed. 
Let $\mathcal R_3$ be  the closed pentagon with vertices  $\mathfrak C, \mathfrak D, (1,0), \mathfrak D',$ and $\mathfrak C'$ from which the line segments $[\mathfrak C, \mathfrak D]$ and 
$[\mathfrak C', \mathfrak D']$ are removed. $($See Figure \ref{local-type} below$)$.
\end{defn}

For $(p,q)\in [1,2]\times [2,\infty]$,  we define the exponent  $\beta(p,q)$ by setting 
\Be
\label{def-beta} 
\beta(p,q)=
         \begin{cases} 
                 \  -\frac12\delta(p,q),    & \ \big(\frac1p,\frac1q\big)\in\mathcal R_1,  
                                              \\  
                 \ \frac{d}{2}\big(\frac1p+\frac1q\big) -\frac{d+1}{2}, & \ \big(\frac1p,\frac1q\big)\in\mathcal R_2, 
                    \\
               \     \frac{d-1}{2}-\frac{d}{2}\big(\frac1p+\frac1q\big),& \ \big(\frac1p,\frac1q\big)\in\mathcal R_2',
                   \\
               \    \frac{d}{2}\delta(p,q) -1,           & \ \big(\frac1p,\frac1q\big)\in \ol{\mathcal R}_3 .
\end{cases} 
\Ee
Here, $\overline {\mathcal R}_3$ denotes the closure of $\mathcal R_3$.  Clearly, $\beta(p,q)$ is well defined. Indeed, note that $ \beta(p,q)=\max \big(-\tfrac12\delta(p,q),  -1+\tfrac{d}{2}\delta(p,q), -\tfrac{d+1}{2}+\tfrac{d}{2}\big(\frac1p+\frac1q\big), \tfrac{d-1}{2}-\tfrac{d}{2}\big(\frac1p+\frac1q\big)\big). $
For $x,y\in \mathbb R^d$, we define 
\[\mathcal D(x,y):=1+\inp xy^2-|x|^2-|y|^2.\]
We are now  ready to state our result concerning the local estimate.

\begin{thm}\label{thm-locest}  Let $d\ge 2$ and $(p,q)\in [1,2]\times [2,\infty]$.  
Suppose that $E \subset \mathbb R^d$ is a measurable set satisfying 
\Be
\label{det-con}
\inf\{  |\mathcal D(x,y)|: x,y\in E \}\ge c
\Ee
for some $c>0$. 
Then, if  $(1/p,1/q)\not\in [\mathfrak C,\mathfrak  D]\cup[\mathfrak C', \mathfrak  D']$,  we have the estimate 
\begin{align}\label{est-loc}    \|\chi_{\sqrt \lambda E}\Pi_\lambda\chi_{\sqrt{\lambda} E}\|_{p\to q}\lesssim    \lambda^{\beta(p,q)}
\end{align}
 where $\sqrt \lambda E=\{\sqrt\lambda x: x\in E\}$. 
Additionally, $(i)$  if $(1/p,1/q)\in (\mathfrak C, \mathfrak D]$ we have the weak type estimate 
$  \|\chi_{\sqrt{\lambda}E}\Pi_\lambda\chi_{\sqrt{\lambda}E}\|_{L^p\to L^{q,\infty}}\lesssim \lambda^{\beta(p,q)} $, 
and  $(ii)$   if $(1/p,1/q)=\mathfrak C$ or $\mathfrak C'$, we  have the restricted weak type estimate $
    \|\chi_{\sqrt{\lambda}E}\Pi_\lambda\chi_{\sqrt{\lambda}E}\|_{ L^{p,1}\to L^{q,\infty}} \lesssim   \lambda^{\beta(p,q)}
   $.  
\end{thm}

 Here $\| \cdot \|_{L^{p,r}\to L^{q,s}}$ denotes the operator norm from $L^{p,r}$ to $L^{q,s}$. There is a strong resemblance between Theorem \ref{thm-locest} and Corollary \ref{laplacian}.  
From  Lemma \ref{implication}  and  Lemma \ref{equivalence} we see  that  the local estimates in Theorem \ref{thm-locest}  imply  those  in Corollary \ref{laplacian} in some cases.  There is a  similar implication from the  local  uniform  $L^p$  bound on Bochner-Riesz means  of the Hermite expansion to the $L^p$ bound on the classical  Bochner-Riesz means \cite{Th87, KST} (also, see \cite{LR}
for recent progress concerning  Bochner-Riesz means  of the Hermite expansion). 
 
 Remarkably,  the estimates in Theorem \ref{thm-locest} are  sharp, and the range of  $p,q$  is also optimal. 
The most typical example of the set $E$ which satisfies \eqref{det-con} is the ball $B(0,r):=\{x: |x|<r \} $ with $r<1/\sqrt 2$.  If $E=B(0,r)$ and $r<1/\sqrt 2$,  from Theorem \ref{thm-locest} and Proposition \ref{prop:lolo}  we see
\[ \|\chi_{\sqrt \lambda E}\Pi_\lambda\chi_{\sqrt{\lambda} E}\|_{p\to q}\sim \lambda^{\beta(p,q)} \]
if $(1/p,1/q)\in[1/2,1]\times[0,1/2]\setminus([\mathfrak C,\mathfrak  D]\cup[\mathfrak C', \mathfrak  D'])$. 
Indeed, the upper bound  follows from Theorem \ref{thm-locest} and 
 the lower bound is shown in Proposition \ref{prop:lolo}.  If  $(1/p,1/q)\in [\mathfrak C,\mathfrak  D]\cup[\mathfrak C', \mathfrak  D']$  the estimate \eqref{est-loc} fails, see 
Lemma \ref{implication} and Theorem \ref{rr*}.   Thus, Theorem \ref{thm-locest} provides a complete characterization of $(p,q)$ for which \eqref{est-loc} holds under the assumption \eqref{det-con}.

In some cases, especially, $p=2$, $q=p'$, and $q=2$,  these sharp  (local) bounds  coincide  with the previously known estimates  \cite{K94, Th98, T05}.  However, 
due to change of the regimes  (see Figure \ref{local-type})  the other sharp estimates are not accessible  by interpolation between those known estimates.  The implication in Lemma \ref{implication} remains valid with Lorentz spaces  as long as $q>1$. So, it is not possible to strengthen the weak type estimates  in Theorem \ref{thm-locest}  ($(i)$)  by replacing 
$L^{q,\infty}$ with the smaller space $L^{q,r}$, $r<\infty$, because the same  is true for  the restriction-extension operator $R^*\!R$ (see Theorem \ref{rr*}). 
Only the problem of determining whether the restricted weak type estimate ($(ii)$ in Theorem \ref{thm-locest})  can be improved to the weak type estimate remains open, not to mention the 
same problem for $R^*\!R$.  We  refer the reader forward to Section 
\ref{transplantation} regarding  boundedness of the operator $R^*\!R$.

The quantity $\cD(x,y)$ is closely related  to the behavior of  the kernel of $\Pi_\lambda$.  
If $(x(t),\xi(t))$ is a Hamiltonian flow of $H-\lambda$ such that $(x(0),\xi(0))=(x,\xi)\in B_{2d}(0,\lambda^\frac12)$.
Then a computation shows $\mathcal D(\lambda^{-\frac12}x,\lambda^{-\frac12}x(t)) = \lambda^{-2}(\inp{\xi}{\xi(t)})^2\ge 0$. 
Thus, the curve $x(t)$ passing through $x$ is contained in the set $\{y:\mathcal D(\lambda^{-\frac12}x,\lambda^{-\frac12}y)\ge 0\}$. 
Since the eigenfunctions of $H$ with eigenvalue $\lambda$ are highly concentrated near the Hamilton flows in the phase space,  one may expect that 
 the kernel $\Pi_\lambda(x,\cdot)$ is essentially supported in the region $\{y:\mathcal D(\lambda^{-\frac12}x,\lambda^{-\frac12}y)\ge 0\}.$
Indeed, this heuristics can be justified rigorously (see  Section 4 and 5).
The quantity $\cD(x,y)$ controls  the estimate for the oscillatory integral $I_j$ (see $\eqref{S3-ptkernel}$ for its definition) of which phase $\cP$   is given by \eqref{p-def}.
The integral $I_j$  appears in the decomposition \eqref{decomp-proj} of a kernel expression of the projection operator $\Pi_\lambda$.   
Since the derivative of the phase function $\cP$  is given by \eqref{ph-d},  the discriminant  $\cD(x,y)$ of the quadratic equation $\cQ(x,y,\tau)=0$  controls the distance between the critical points of the phase function, so  the condition \eqref{det-con} guarantees that the second order derivative of the phase function is not vanishing. This allows us to have a better decay estimate for $I_j$.  A similar observation played an important role in proving the local estimate  \eqref{local-thanga}  due to Thangavelu \cite{Th98}.  
We elaborate his observation to get the optimal local estimate for the projection operator $\Pi_\lambda$.

\subsubsection*{Local estimates  over $A_{\lambda,\mu}^+$ } 
If the condition \eqref{det-con} is not satisfied any more, the boundedness of  $\Pi_\lambda$  becomes more complicated. 
Estimate over the region near the sphere  $\sqrt\lambda \mathbb S^{d-1}:=\{x:|x|=\sqrt \lambda\} $ is of special interest since the kernel of  $\Pi_\lambda$ exhibits different behaviors across  $\sqrt\lambda \mathbb S^{d-1}$.  For more about  the estimates over the intermediate region, see Corollary \ref{local-r}. Koch and Tataru \cite{T05} considered the estimate for the localized projection operator $\chpl\Pi_\lambda$ and obtained   the sharp $L^2$--$L^q$ estimates with $q\ge 2$. More precisely,  the following were  shown (see \cite[Theorem 3]{T05}): 
\begin{equation}  
\label{kt-nu}
 \|\chplp\Pi_\lambda \|_{2\to q}\sim 
                                      \begin{cases}
\lambda^{-\frac12\delta(2,q)}\mu^{\frac14-\frac{d+3}{4}\delta(2,q)}, &   \ 2\le q\le \frac{2(d+1)}{d-1},
\\
(\lambda\mu)^{-\frac12+\frac d2\delta(2,q)}, &   \   \frac{2(d+1)}{d-1} \le q\le \infty 
                                      \end{cases}
\end{equation}
for $\lambda^{-\frac23}\le\mu\le 1/4$. Actually, a slightly stronger  weighted $L^2$ estimate  was shown.
However,  the result in  \cite{T05} is essentially equivalent to \eqref{kt-nu} except for the cases of  critical exponents, which  can also be recovered by  
making use of the estimate \eqref{unbalanced}.  If $\chplp$ in the \eqref{kt-nu} is replaced by $\chpln$, similar but stronger estimates for $p,q$ in the full range $1\le p\le 2\le q\le \infty$  are easier to show (see \eqref{est:ext--} in Proposition \ref{prop:estext}).  
 By duality, \eqref{kt-nu} gives the equivalent bounds on  $\|\chplp \Pi_\lambda \chplp\|_{q'\to q}$.  As mentioned before, mere interpolation between these  estimates and those in \eqref{kt-nu}  does not yield  the sharp estimate. 
 
In Section \ref{off-d} we obtain  the optimal bound on $\|\chplp \Pi_\lambda \chplp\|_{p\rightarrow q}$ for an extended range of $p,q$.  
More precisely,  it is natural to expect the following estimates hold for $1\le p\le 2\le q\le \infty$, possibly except for some endpoint cases:   
\begin{align}\label{est-ann}
    \|\chi_{\lambda,\mu}^+ \Pi_\lambda\chi_{\lambda,\mu}^+\|_{p\to q}\le C\lambda^{\beta(p,q)}\mu^{\gamma(p,q)},
\end{align}
where the exponent $\beta(p,q)$ is given by \eqref{def-beta} and  
\begin{align*}
    \gamma(p,q)=\begin{cases} 
                 \  \frac12-\frac{d+3}{4}\delta(p,q),    & \  \big(\frac1p,\frac1q\big)\in\mathcal R_1,  
                                              \\  
                 \ d\big(\frac{1}{2p}+\frac1q\big)-\frac{3d+1}{4}, & \ \big(\frac1p,\frac1q\big)\in\mathcal R_2, 
                    \\
               \     \frac{3d-1}{4}-d\big(\frac1p+\frac{1}{2q}\big),& \  \big(\frac1p,\frac1q\big)\in\mathcal R_2',
                   \\
               \    \frac{d}{2}\delta(p,q) -1,           &  \ \big(\frac1p,\frac1q\big)\in\overline {\mathcal R}_3.
\end{cases}
\end{align*}
The estimate \eqref{est-ann} is a natural off-diagonal extension of \eqref{kt-nu} and  it is not difficult to show that the exponent in \eqref{est-ann} can not be improved to any better one, 
see Proposition \ref{lower-mu}.
We verify \eqref{est-ann} on a larger range $p,q$.  Our result subsumes  what is obtained in \cite{T05} and significantly extends  the range of the sharp bound. 
The result is summarized in Theorem \ref{thm-annest}.  
Among others, worth mentioning is  the weak type $(1,2d/(d-1))$ estimate  which corresponds to the point $\mathfrak D$ in Figure \ref{local-type} (see  Theorem \ref{thm-annest}). 
 Though   the range where the sharp bounds are available is considerably extended, there are regions where  the optimal bound is  left unknown.

\subsection*{Our approach} Before closing the introduction, we briefly discuss our approach. We make use of an explicit representation  for the projection operator $\Pi_\lambda$ which 
is given by an oscillatory integral. The local estimates in Theorem \ref{thm-locest} are obtained by obtaining the estimates for the oscillatory integral which we combine with a form of $TT^*$.   
However, compared with the abovementioned local estimates, the estimates \eqref{unbalanced} and \eqref{est-ann} are much more difficult to prove.  The kernel of  $\Pi_k$  is related to  the phase function $\mathcal P(x,y,s)$ (see \eqref{p-def} for its definition and  see, also,  Lemma \ref{kernel} and \eqref{S3-ptkernel}).  Control of associated oscillatory integral becomes more complicated  over $A_{\lambda, \mu}^+\times A_{\lambda, \mu'}^+ $ since the quantity  $|\cD(x,y)|$ may vanish. If the zeros of $\partial_s \mathcal P$ and $\partial_s^2 \mathcal P$ are well separated,  the lower bound for each derivative can be exploited to get the desirable bound. However,  to make  the problem even worse,   zeros of  $\partial_s \mathcal P$ and $\partial_s^2 \mathcal P$   get close  to each other as $\cD(x,y)\to 0$ (see \eqref{qcos}).  This naturally gives a rise to an Airy type integral for which  we can only expect decay of $\lambda^{-\frac13}$  at best. 
Thus,  we can not get a good estimate which allows us  to recover the sharp bound
To get around this difficulty, we perform additional decomposition away from the zero of the equation $\partial_s^2 \mathcal P=0$,  which salvages the correct decay $\lambda^{-\frac12}$. 
We elaborate this idea in Section \ref{proof-end} and  Section \ref{off-d} to prove the sharp estimates. This involves the decomposition \eqref{decomp} which leads our to consider 
various different cases.

\subsubsection*{\bf Organization} In Section \ref{sec:tt-st} we obtain an explicit representation  for the projection operator $\Pi_\lambda$ and formalize a form of $TT^*$ argument which plays a crucial role throughout the paper and we then use it  to prove Theorem \ref{glo-est}.  Section 3 is devoted to proving the local estimates.  In Section 4, we consider the estimate over the annulus $A_{\lambda, \mu}^\pm$ and prove various preparatory estimates based on  the sectorial decomposition of the annuli, which we utilize  to show Theorem \ref{endpoint}. In Section \ref{sec:endpoint}, we begin the proof of  the estimate  \eqref{unbalanced}  (with unbalanced improvement) which actually proves the endpoint estimate \eqref{intro-conj}. By establishing  estimates  for various specific cases  we manage to  single out the most difficult case, i.e., Proposition \ref{critical-muj} which we do not prove  until Section \ref{proof-end}.  In  Section \ref{off-d},  we obtain new sharp off-diagonal estimates by making use of the estimates obtained  in the previous sections and the $TT^*$ argument in  Section \ref{sec:tt-st}.

\subsubsection*{\bf Notations} Throughout the paper  $d\ge2$ unless it is mentioned otherwise.
\vspace{-5pt}
\begin{enumerate} 
[leftmargin=.5cm, labelsep=0.3 cm, topsep=0pt]
\item[$\bullet$]  For  nonnegative quantities   $A$ and $B$,    $B \lesssim A$ means that there is a constant $C$, depending only on dimensions such that 
$B\le CA $.  Likewise,   $A\sim B$ if and only if   $B \lesssim A$ and $A \lesssim B$. Slightly abusing the conventional notation, by $D=O(A)$  we denote 
$|D|\lesssim A$. 
\item[$\bullet$] Additionally, we denote $A\gg B$ if there is a large constant $C$ depending only on the dimensions  such that $A\ge CB$. 
\item[$\bullet$]   We occasionally write $x=(x_1, \ol x)\in \mathbb R\times \mathbb R^{d-1}$,  $x=(x_1, x_2, \widetilde x)\in \mathbb R\times\mathbb R\times  \mathbb R^{d-2}$.
\item[$\bullet$]    $B_d(x,r):=\{y\in \mathbb R^d: |y-x|<r\}. $
\item[$\bullet$]  In addition to the sets $A_\mu^\pm$ and $A_{\lambda, \mu}^\pm$, we also define  $A_{\mu}^\circ=\{x:   |1- |x|| \le  2\mu\}$ and  $A_{\lambda, \mu}^\circ=\{x: \lambda^{-1/2}x\in  A_{\mu}\}$.
Then, we set  $\chi_\mu^\circ=\chi_{A_\mu^\circ}$, and $\chi_{\lambda, \mu}^\circ=\chi_{A_{\lambda, \mu}^\circ}$. 
\item[$\bullet$]   
For a given operator $T$ we denote by $T(x,y)$ the kernel of $T$.  

\item[$\bullet$] Let $\psi$ be a  function in $C^\infty_c([\frac14,1])$ such that $\sum \psi(2^j t) =1$ for $t>0$. We also denote   $\widetilde\psi(t)=\psi(|t|)$.
\item[$\bullet$]   If it is not mentioned otherwise,  $\psi_j$ always denotes  a smooth function which is supported in $[2^{-2-j},2^{-j}]$  and satsifies $|\frac{d^l}{dt^l}\psi_j|\le C2^{j l}$ for $l=0,1,\dots.$

\item[$\bullet$]   The constants $c$ and $\vepc$ are small constants depending only on the dimensions.

\item[$\bullet$]    Let $\partial_x:= ( \partial_{x_1}, \dots, \partial_{x_d})^\intercal$, $\partial_x^\intercal : = ( \partial_{x_1}, \dots,  
\partial_{x_d})$,  and $\partial_x \partial_y^\intercal=(\partial_{x_i}\partial_{y_j})_{1\le i,j\le d}$. 

\item[$\bullet$]  
We  set  
$ \mathfrak a(t) :=(2\pi i\:\sin t)^{-\frac d2}{e^{i \pi d/4}}.$

\item[$\bullet$]  By $\mathbf I_d$ we denote the $d\times d$  identity matrix.
\end{enumerate}

\section{The projection operator $\Pi_\lambda$ and $TT^*$ argument}
\label{sec:tt-st}

In this section, we obtain an explicit expression  for the kernel of  the Hermite projection operator $\Pi_\lambda$ and deduce some useful properties which are to be used later. We also formulate a $TT^*$ argument which is adapted to the projection operator $\Pi_\lambda$. 

\subsection{Representation of the projection operator} We begin with noting that  the  Hermite-Schr\"odinger propagator $e^{-itH} f$ is given  by 
\Be 
\label{schro-op}
 e^{-itH} f=\sum_{\lambda\in \eval} e^{-it\lambda} \Pi_\lambda f
  \Ee
 for  $f\in \mathcal S(\mathbb R^d)$.
Since $\Pi_\lambda f$ decays rapidly in $\lambda$ (see Lemma \ref{rapid-decay} below) for  $f\in \mathcal S(\mathbb R^d)$, the sum 
 $ \sum_{\lambda\in \eval} e^{-it\lambda} \Pi_\lambda f$ converges uniformly. Furthermore,  it is clear that $e^{-itH} f$ is smooth on $
\mathbb R^d\times \mathbb R$. Orthogonality of the Hermite functions yields 
\begin{equation}
\label{iso}
\|e^{-itH} f\|_2=\|f\|_2,  \  \  \forall t\in \mathbb R.
\end{equation} 

\begin{lem}\label{rapid-decay} Let $f\in \mathcal S(\mathbb R^d)$. Then, for any $N$,  there is a constant $C=C(N, f)$ such that $|\inp {f} {\Phi_\alpha}|\le C (1+|\alpha|)^{-N} $ and 
$
\|\Pi_\lambda f\|_\infty \le C \lambda^{-N}.
$ 
\end{lem} 
Since the eigenvalues $\lambda, \lambda'\in  2\mathbb N_0+d$,   $\lambda-\lambda'\in 2\mathbb Z$. So, $\frac1{2\pi}\int_{-\pi}^\pi e^{i\frac t2(\lambda-\lambda')}  dt $ $=\delta(\lambda-\lambda')$. By Lemma \ref{rapid-decay}, it follows  that 
\[  \Pi_\lambda f=\sum_{\lambda'\in\eval}  \frac1{2\pi} \int_{-\pi}^\pi  e^{i\frac t2(\lambda-\lambda')}  dt\, \Pi_{\lambda'} f=  \frac1{2\pi} \int_{-\pi}^\pi    \sum_{\lambda'\in\eval}  
e^{i\frac t2(\lambda-\lambda')} \Pi_{\lambda'} f dt\] 
for all $f\in \mathcal S(\mathbb R^d)$ because the series converges uniformly. Hence, combining this with \eqref{schro-op} we get 
  \begin{align}
 \Pi_\lambda f = \frac{1}{2\pi}\int_{-\pi}^{\pi} e^{i\frac t2(\lambda-H)}fdt, \quad \forall f\in \mathcal S(\mathbb R^d).
 \label{proj-op}
 \end{align}
This observation together with the explicit kernel representation of the operator $e^{-itH}$ serves as an important tool throughout the paper.

\begin{proof}[Proof of Lemma \ref{rapid-decay}] Since $\Phi_\alpha$ is an eigenfunction with its eigenvalue $2|\alpha|+d$,  $H^N \Phi_\alpha=(d+2|\alpha|)^N \Phi_\alpha$. Integration by parts  gives $ \inp {f}{\Phi_\alpha}=(d+2|\alpha|)^{-N} \inp {H^N f}{\Phi_\alpha}.$
Since  $\|\Phi_\alpha\|_1\le C|\alpha|^{{d}/4}$ with $C$ independent of $\alpha$ (see Lemma \ref{S4-bd_Her}) and 
$H^N f\in \mathcal S(\mathbb R^d)$, we have 
$
 | \inp {f}{\Phi_\alpha}|\le C (d+2|\alpha|)^{-N+d/4}
 $
for any $N$.  From Lemma \ref{S4-bd_Her}, it follows that $|\Phi_\alpha|\lesssim 1$. 
A better bound  is possible but this is enough for our purpose. 
Thus, we have
\[ |\Pi_\lambda f|\le \sum_{\alpha: d+2|\alpha|=\lambda}  | \inp {f}{\Phi_\alpha}|\, \|\Phi_\alpha\|_\infty
\lesssim \sum_{\alpha: d+2|\alpha|=\lambda}  ( d+2|\alpha|)^{-N+\frac {d}4}\lesssim \lambda^{-N+\frac {5d-4}4}. \] 
Since $N$ can be taken to be arbitrarily large, we get the desired bounds. 
\end{proof}

Now we recall that  the operator $e^{-itH}$ can be  expressed as follows: 
\begin{align}\label{Schrodinger}
    e^{-itH}f=\mathfrak{a}(2t) 
    \int  e^{ip(x,y,t)}f(y)\,dy,
\end{align}
where $p(x,y,t) = \frac{|x|^2+|y|^2}{2}\, \cot 2t -\inp x y\,\csc2t$ and $ \mathfrak a(t) =(2\pi i\:\sin t)^{-\frac d2}{e^{i \pi d/4}}.$
See S\"ojgren and Torrea \cite{SoTo},  and  Thangavelu \cite[p.11]{Th87}  for a detailed discussion.  Combining \eqref{proj-op} and \eqref{Schrodinger},
we obtain the following representation of 
the projection operator  $\Pi_\lambda$.

\begin{lem}
\label{kernel}
Let $\lambda\in 2\mathbb N_0+d$ and set 
\begin{align*}
  {\phi}_\lambda:=\phi_\lambda (x,y,t):=\frac{\lambda t}{2}+\frac{|x|^2+|y|^2}{2}\cot t-  \inp xy \csc t.
\end{align*}Then, for all $f\in \mathcal S(\mathbb R^d)$,  we have 
\[
    \Pi_\lambda f =   \frac{1}{2\pi}  \int_{-\pi}^{\pi}   \mathfrak a(t) \int   e^{i\phi_\lambda(x,y,t)}   f(y)\, dy\,dt. \footnote{A branch of $(2\pi i\:\sin t)^\frac d2$ should be chosen suitably.} 
   \]
\end{lem}

We note that  $\phi_\lambda(\mU x,\mU y,t)=\phi_\lambda(x,y,t)$ for all $\mathbf U\in \mathrm O(d)$ where $\mathrm O(d)$ is the orthonormal group in $\mathbb R^d$.  In viewpoint of the kernel, the operator is invariant under the map $(x,y)\to (\mU x, \mU y)$, that is to say, $\Pi_\lambda (\mU x, \mU y)=\Pi_\lambda(x,y) $.  Thus,  it follows that 
\Be
\label{rot}
  (\Pi_\lambda f)\circ \mU=\Pi_\lambda (f\circ \mU). 
\Ee
Kochneff \cite{kochneff} made the same observation  in a different way which was based on the properties of the Hermite functions. 

Since  $\sin t$ vanishes at $t=0$ and $t=\pi, -\pi$, we decompose the operator away from those points to avoid the singularities.  We denote by $\psi$  the function in $C^\infty_c([\frac14,1])$ such that $\sum \psi(2^j t) =1$ for $t>0$. Then, let 
$ \psi^0$ be the smooth function which  satisfies,   for $t\in(-\pi,\pi)\setminus\{0\}$,
\begin{equation}
\label{j-decomp}
 \psi^0+ \sum_{j\ge 4} \big( \psi(2^j t)+ \psi(-2^j t)  +\psi(2^j(t+\pi))+\psi(2^j(\pi-t))\big)=1.
\end{equation}
So, $ \psi^0$ is supported in the interval $[-\pi, \pi]$ and vanishes near  $0, \pi,$  and $-\pi$.  
  For a bounded continuous  function $\eta$ supported in $[-\pi, \pi]$ and  any $\lambda\in \mathbb R$, we define 
   \Be \label{proj-pi}\Pi_\lambda[\eta]:= \int \eta(t)e^{i\frac t2(\lambda-H)}dt.\Ee
  Clearly, the definition  of $\Pi_\lambda[\eta]$ makes sense for any real number $\lambda$. From the isometry  \eqref{iso} it follows that 
  \Be
  \label{easy-l2} 
  \|\Pi_\lambda[\eta]\|_{2\to 2} \le \|\eta\|_1.
  \Ee
For simplicity   let us denote 
\[ \psi_j^\pm(t) := \psi(\pm 2^j t), \qquad \psi_j^{\pm \pi} (t):=\psi(2^j(\pi\mp t)), 
\qquad \psi_j^0:= 
\begin{cases} \psi^0,  &j=4 
\\
\ 0,  &{ j>4}
\end{cases} .\]
For the last  one we adopt  a notational convention which facilitates statements of our results.
By \eqref{j-decomp} we can decompose the operator $\Pi_\lambda$ as follows:  
\begin{equation}
\label{decomp-proj}
 \Pi_\lambda f =  \Pi_\lambda[ \psi^0] f +  \sum_{j\ge 4}  \Big( \Pi_\lambda [\psi_j^+]f  + \Pi_\lambda  [\psi_j^-]f  + \Pi_\lambda  [\psi_j^{+\pi}]f  + \Pi_\lambda  [\psi_j^{-\pi}] f\Big)
 \end{equation}
for $f\in \mathcal S(\mathbb R^d)$ and  $\lambda\in \eval$.  This decomposition of $\Pi_\lambda$ is clearly valid because the right hand side of \eqref{decomp-proj} 
converges to $\Pi_\lambda$ as a  bounded operator on $L^2$. This is easy to show using \eqref{easy-l2}. Indeed, by \eqref{easy-l2}  we have 
the estimate 
$\| \Pi_\lambda [\psi_j^\kappa ]\|_{2\to 2}\lesssim 2^{-j}$, $\kappa= \pm, \pm \pi.$ Thus the convergence is clear. 

We prove estimates for $\Pi_\lambda$ by obtaining  those estimates for the individual operators  appearing in \eqref{decomp-proj}. 
Once we have the estimates for $ \Pi_\lambda [\psi_j^+]$, 
by making use of symmetric property of $e^{-itH} f$  (\eqref{symmetric})  we can handle  the operators $\Pi_\lambda  [\psi_j^-]$, $\Pi_\lambda  [\psi_j^{+\pi}]$, and $\Pi_\lambda  [\psi_j^{-\pi}]$ via simple changes of variables. So, we mainly work with  $ \Pi_\lambda [\psi_j^+]$.  The estimate for $\Pi_\lambda[ \psi^0]$ is relatively easier than that for   $ \Pi_\lambda [\psi_j^+]$.  
We  observe the following symmetric properties  of the phase function $\phi_\lambda$:  
 \[ \phi_\lambda(x,y,-t)=-\phi_\lambda(x,y,t), \ \  \  \phi_\lambda(x,y,\pm (\pi-t))=\pm\lambda \frac \pi 2 \mp\phi_\lambda(x,-y,t) . \]
Thus, by a simple change of variables  it follows that
\Be 
\label{symmetric} 
\begin{aligned}
\Pi_\lambda[\psi_j^{-}](x,y)&=C(d) \ol{ \Pi_{\lambda}[\psi_j^{+}](x,y)}, 
\\ 
\Pi_\lambda[\psi_j^{\pm \pi}](x,y)&= C(d,\lambda, \pm) \Pi_{\lambda} [\psi_j^{\mp}](x,-y)
\end{aligned}
\Ee
with  $C(d)$, $C(d,\lambda, \pm)$ satisfying $|C(d)|=|C(d,\lambda, \pm)|=1$. This  clearly implies $\|\Pi_\lambda[\psi_j^{+}]\|_{p\to q}= \|\Pi_{\lambda} [\psi_j^{\kappa}]\|_{p\to q}$ for $\kappa=-, \pm \pi$, 
so we only need to obtain the bounds on $\|\Pi_\lambda[\psi_j^{+}]\|_{p\to q}$ and $\|\Pi_\lambda[\psi^0]\|_{p\to q}.$

We define the rescaled operators $\fP_\lambda$,  $\fP_\lambda[\eta]$ whose kernels are given by  
\Be 
\label{eq:scaled-op}
\begin{aligned}
\fP_\lambda(x,y)&=\Pi_\lambda(\sqrt \lambda x, \sqrt \lambda y), 
\\
 \fP_\lambda[\eta](x,y)&= \Pi_\lambda[\eta](\sqrt \lambda x, \sqrt \lambda y).
\end{aligned}
\Ee
These operators, instead of those defined  by $\Pi_\lambda$, are sometimes more convenient to work with. 
By scaling it is clear that 
\begin{equation}
\label{eq:norm-scaled}
  \|  \chi_{E} \mathfrak \fP_\lambda [\eta]  \chi_{F} \|_{p\to q} =  \lambda^{\frac d2(\frac1p-\frac1q-1)}\|  \chi_{\sqrt\lambda E} \Pi_\lambda [\eta] \chi_{\sqrt\lambda F} \|_{p\to q}.  
   \end{equation}

\subsection{A $TT^\ast $ argument}
The argument below  allows us to deduce off-diagonal estimates from  $L^1$--$L^\infty$ and $L^2$ estimates.   
  The following is a variant of the typical $TT^\ast$ argument  (see \cite{keel-tao, fos, vil}) . 

\begin{lem}
\label{tt-st} Let $E$ be a measurable subset of $\mathbb R^d$ and $b> 0$.
Let  $\mathfrak Q(b)\subset [1/2,1]\times [0,1/2]$ denote the closed quadrangle with vertices $(\frac12, \frac12)$,  $(\frac12, \frac{b}{2(b+1)})$, 
$(\frac{b+2}{2(b+1)}, \frac12)$, and  $(1,0)$. Suppose that, for $j\ge 0$, 
\begin{equation}
\label{1-infty}
\| \chi_{E} \Pi_\lambda [\eta_j]\chi_{E}\|_{1\to \infty}\lesssim    \beta 2^{bj}
\end{equation}
whenever $\eta_j$ is a smooth function supported in $[2^{-2-j},2^{-j}]$ $($or $[-2^{-j}, -2^{-2-j}])$ and $|\frac{d^l}{dt^l}\eta_j|\le C2^{j l}$ for $l=0,1.$
Then, if  $\ppq$ is contained in  $\mathfrak Q(b)$, for $j\ge 2$ we have the estimate
\begin{equation} 
\label{bpjb} 
\| \chi_{E} \Pi_\lambda [ \psi_j ] \chi_{E}\|_{p\to q}\lesssim     \beta^{\dpq} 2^{j(-1+(b+1)\dpq)}.
\end{equation}
\end{lem}

\begin{proof} 
By \eqref{easy-l2}, it is clear that $\| \chi_{E} \Pi_\lambda [ \eta_j ] \chi_{E}\|_{2\to 2}\lesssim    2^{-j}.$
Interpolating this with the estimate \eqref{1-infty},   for  $1\le p\le 2$  we have 
\begin{equation}
\label{bpjb2} 
\| \chi_{E} \Pi_\lambda [ \eta_j ] \chi_{E}\|_{p\to p'}\lesssim    \beta^{(\frac2p-1)} 2^{j(-1+(b+1)(\frac2p-1))}.
\end{equation}
Clearly, by \eqref{symmetric} the same estimate also holds if $\eta_j$ is a smooth function supported in $[-2^{-j}, -2^{-2-j}]$ and $|\frac{d^l}{dt^l}\eta_j|\le C2^{j l}$ for $l=0,1.$  
 Thus,  \eqref{bpjb2} remains valid  for  $\eta_j$ supported in $[-2^{-j}, -2^{-2-j}]\cup [2^{-2-j},2^{-j}]$ as long as $|\frac{d^l}{dt^l}\eta_j|\le C2^{j l}$ for $l=0,1.$

 It is sufficient to show  \eqref{bpjb} with $q=2$ and $p={2(b+1)}/(b+2)$ because the other estimates follow from duality and interpolation between this  estimate and \eqref{bpjb2}.  
We claim that 
\begin{equation}
\label{2p}
  \| \Pi_\lambda [ \psi_j ] \chi_{E}  f\|_{2}\lesssim 2^{-\frac j2} \beta^{\frac1p-\frac12}
    \|f\|_p
  \end{equation}  
holds with $p={2(b+1)}/(b+2)$. This clearly implies the desired estimate \eqref{bpjb} with $q=2$ and $p={2(b+1)}/(b+2)$.

Let $\widetilde\psi(t)=\psi(|t|)$. 
For $k\ge j-2$, we set
\begin{align}
\label{pk}    (\Pi_\lambda [ \psi_j ] ^*\Pi_\lambda [ \psi_j ] )_k&=\iint\wt\psi(2^k(t-s))\psi_j(t)\psi_j(s)e^{i\frac{t-s}2(\lambda-H)}dsdt\,.
\end{align}
By \eqref{proj-op}  we note 
$%
  \| \Pi_\lambda [ \psi_j ] \chi_{E} f\|_{2}^2 =\langle   \chi_{E}  \iint \psi_j(t)\psi_j(s)e^{i\frac{t-s}2(\lambda-H)}   \chi_{E}  f dsdt, f\rangle.
$ 
Thus   we have \Be
 \label{ttsum}
   \| \Pi_\lambda [ \psi_j ] \chi_{E} f\|_{2}^2 =  \sum_{k\ge j-2}  \left\langle   \chi_{E} (\Pi_\lambda [ \psi_j ] ^*\Pi_\lambda [ \psi_j ] )_k  \chi_{E} f, f\right\rangle.
\Ee
After a simple change of variables we observe that 
\Be
\label{bk}
 \chi_{E}(\Pi_\lambda [ \psi_j ] ^*\Pi_\lambda [ \psi_j ] )_k  \chi_{E} f= \int\psi_j(s) \chi_{E}  \Pi_\lambda[\widetilde \psi(2^k\cdot)\psi_j(\cdot+s)]\chi_{E} f ds
\Ee
 Since $k\ge j-2$, we apply the estimate \eqref{bpjb2}  to $\chi_{E}  \Pi_\lambda[\widetilde \psi(2^k\cdot)\psi_j(\cdot+s)]\chi_{E}$  and  get 
\[  |  \left\langle   \chi_{E} (\Pi_\lambda [ \psi_j ] ^*\Pi_\lambda [ \psi_j ] )_k  \chi_{E} f, g\right\rangle|\lesssim    \beta^{(\frac2p-1)}  2^{-j} 2^{k(-1+(b+1)(\frac2p-1))} \|f\|_p\|g\|_p\] 
for $1\le p\le 2$. Now, using \eqref{ttsum},  summation along $k$ yields, for ${2(b+1)}/(b+2)<p\le 2$,
\[  \| \Pi_\lambda [ \psi_j ] \chi_{E} f\|_{2}\lesssim    \beta^{(\frac1p-\frac12)}  2^{j(-1+(b+1)(\frac1p-\frac12))} \|f\|_p.\] 
Hence, we have \eqref{bpjb} for $q=2$ and ${2(b+1)}/(b+2)<p\le 2$, and 
by duality we additionally  have  \eqref{bpjb} for $p=2$ and $2\le q<\frac{2b+2}{b}$, so interpolation between these estimates and the estimate 
\eqref{bpjb2} gives the desired estimate \eqref{bpjb}  for  $p,q$ such that  $\ppq$ is contained in $\mathfrak Q(b)$ but  not in the line segments 
$[(\frac12,\frac b{2(b+1)}), (1,0))$, $[(\frac {b+2}{2(b+1)},\frac12), (1,0))$.  

We now use those estimates to obtain the estimate \eqref{2p} with $p={2(b+1)}/(b+2)$. In fact, using \eqref{bk} and the estimate \eqref{bpjb} obtained  in the above  we get 
\begin{align} \label{bibi}
  | \left\langle   \chi_{E} (\Pi_\lambda [ \psi_j ] ^*\Pi_\lambda [ \psi_j ] )_k  \chi_{E} f, g\right\rangle|\lesssim    \beta^{(\frac1p-\frac1q)}  2^{-j} 2^{k(-1+(b+1)(\frac1p-\frac1q))} \|f\|_p\|g\|_{q'}
  \end{align}
for  $p,q$ such that  $(1/p, 1/q)$ is contained in the interior of $\mathfrak Q(b)$. This allows us to apply the bilinear interpolation argument 
(e.g.,  Keel and Tao \cite{keel-tao}). Thus, we obtain
\[ \sum_{k\ge j-2}  | \left\langle   \chi_{E} (\Pi_\lambda [ \psi_j ] ^*\Pi_\lambda [ \psi_j ] )_k  \chi_{E} f, g\right\rangle|\lesssim    \beta^{(\frac1p-\frac1q)}  2^{-j} \|f\|_p\|g\|_{q'}\] 
provided that $(1/p, 1/q)$ is contained in the interior of $\mathfrak Q(b)$ and $(b+1)(1/p-1/q)=1$. Therefore, taking $q'=p$ in the above estimate,  by  \eqref{ttsum}  we get the desired estimate \eqref{2p} with $p=\frac{2b+2}{b+2}$.
 \end{proof}

We occasionally use the following  simple lemma which is useful for getting the estimates of critical cases, if  combined with real interpolation. 

\begin{lem}\label{s-trick}  Let $1\le p_0, p_1,  q_0, q_1\le \infty$ and $\epsilon_0, \epsilon_1>0$.   Let $T_j$, $j\in \mathbb Z $ be sublinear operators defined from $L^{p_k}\to  L^{q_k}$ with $\|T_j\|_{p_k\to q_k}\le B_k 2^{ j (-1)^k \epsilon_k }$ for $k=0,1$.  Let us set $\theta=\frac {\epsilon_0}{\epsilon_0+\epsilon_1}$, 
$\frac1p_\ast=\frac \theta{p_1}+ \frac {1-\theta}{p_0},$ and $\frac1q_\ast=\frac \theta{q_1}+\frac {1-\theta}{q_0}. $  
Then we have the following\,$:$
\vspace{-6pt}
\begin{enumerate} 
[leftmargin=.65cm, labelsep=0.15 cm, topsep=0pt]
\item[$(a)$] If $p_0=p_1=p$ and $ q_0\neq q_1$, then   $\|\sum_j T_j f  \|_{{q_\ast,\infty}}\lesssim B_0^{1-\theta}B_1^\theta \|f\|_{p}$. 
\item[$(b)$] If $q_0=q_1=q$ and $ p_0\neq p_1$, then $\|\sum_j T_j f  \|_{{q}}\lesssim B_0^{1-\theta}B_1^\theta \|f\|_{p_\ast,1}$.
\item[$(c)$] If $p_0\neq p_1$ and $ q_0\neq q_1$, then  $\|\sum_j T_j f  \|_{{q_\ast,\infty}}\lesssim B_0^{1-\theta}B_1^\theta \|f\|_{p_\ast, 1}$.
\end{enumerate}
\end{lem} 

The third $(c)$ assertion is known as \emph{`Bourgain's summation trick'}.  See \cite[Section 6.2]{Car99} for a formulation in abstract setting.  The first $(a)$ and the second $(b)$ statements give  a little  better estimate than the restricted weak type estimate if $p=p_0=p_1$ or $q=q_0=q_1$. As far as the authors are aware, this observation first appeared in  Bak \cite{bak}. A slight modification of the  argument  in \cite{bak} (also see  \cite[Lemma 2.3]{LS}) shows $(b)$ and $(c)$.

\begin{rem}  Thanks to Lemma \ref{s-trick},  a more elementary approach to the estimate \eqref{2p} with $p=\frac{2b+2}{b+2}$  
 is possible  since \eqref{bibi} is equivalent to
 \[\|\chi_{E} (\Pi_\lambda [ \psi_j ] ^*\Pi_\lambda [ \psi_j ] )_k  \chi_{E} f\|_q\lesssim \beta^{(\frac1p-\frac1q)}  2^{-j} 2^{k(-1+(b+1)(\frac1p-\frac1q))} \|f\|_p.\] 
 Using  this and $(c)$ in Lemma \ref{s-trick}, for $p,q$ satisfying that  $\ppq$ is contained in the interior of $\mathfrak Q(b)$ and $(b+1)(\frac1p-\frac1q)=1$    we have $\|\sum_{k}\chi_{E} (\Pi_\lambda [ \psi_j ] ^*\Pi_\lambda [ \psi_j ] )_k  \chi_{E} f\|_{q,\infty}$ $ \lesssim  \beta^{(\frac1p-\frac1q)}  2^{-j}  \|f\|_{p,1}$.   Real interpolation between these restricted weak type estimates gives, in particular,  the estimate \eqref{bpjb} with  $p={2(b+1)}/(b+2)$ and $q={2(b+1)}/b$, which  yields \eqref{2p} with $p={2(b+1)}/(b+2)$.
\end{rem}

As  an application of Lemma \ref{tt-st} we prove Theorem \ref{glo-est}.

\subsection{Proof of Theorem \ref{glo-est}}  
By  interpolation and duality it suffices to show the assertion regarding the restricted weak type $(p,q)$  estimate with $\ppq=\mathfrak F, \mathfrak F'$  in Theorem \ref{glo-est}. 
By duality  it is enough to show 
\Be\label{restricted-weak}
\|\Pi_\lambda f \|_{\frac{2d(d-1)}{(d-2)^2},\infty}\le C\|f\|_{\frac{2(d-1)}{d},1}.
\Ee
Indeed, once  we have the estimate \eqref{restricted-weak} by duality and interpolation  it follows that the estimate \eqref{unif} holds for $p,q$ with 
$(\rp,\rq)\in (\mathfrak F,\mathfrak F')$. In particular, since 
$\Pi_\lambda^* \Pi_\lambda =\Pi_\lambda$, the estimate \eqref{unif}  with $p={2d}/(d+2)$ and $q={2d}/(d-2)$ gives the estimate \eqref{unif} for $q=2$, $p={2d}/(d+2)$, which is 
equivalent to \eqref{karad-global}. We interpolate those estimates  with 
the  $L^2$ estimate  $\|\Pi_\lambda f\|_{2}\le \|f\|_2$  and get the desired estimate for the rest of $\ppq\in \mathscr P$ (see Figure \ref{fig:uniform}). 

By Lemma \ref{kernel}  we have \eqref{1-infty} with $E=\mathbb R^d$, $\beta=1$,  and $b=\frac{d-2}{2}$. 
Using this and  Lemma \ref{tt-st}, we now have, for $j\ge 4$, 
\Be
\label{eq:strong}
\|   \Pi_\lambda [\psi_j^{\kappa} ] \|_{p\to q}\lesssim   2^{j(\frac{d}{2}(\frac1p-\frac1q)-1)},  \quad \kappa=\pm,  \;  \pm \pi,
\Ee
whenever $(\rp, \rq)$ is contained in the closed  quadrangle  $\mathfrak  Q(\tfrac{d-2}2)$ which has  vertices $(\frac12, \frac12)$,  $(\frac12, \frac{d-2}{2d})$, 
$(\frac{d+2}{2d}, \frac12)$, and  $(1,0)$.  It is sufficient to show \eqref{eq:strong} with $\kappa=+$ since  we may use \eqref{symmetric} to get the same estimates for the other cases. By  repeating the argument in the proof of  Lemma \ref{tt-st}  it is easy to see  that 
\begin{equation}
\label{psi0}
 \|   \Pi_\lambda [ \psi^0 ] \|_{p\to q}\lesssim 1\end{equation}
for the same range of $p, q$ as in the above. Thus, combining the estimates \eqref{eq:strong} and \eqref{psi0} and taking sum over $j$, from 
\eqref{decomp-proj}
 we get \eqref{unif} 
if  $(\rp, \rq)$ is contained in $\mathfrak Q(\tfrac{d-2}2)$ 
and $\frac1p-\frac1q<\frac2d$. The estimates in the borderline case ($1/p-1/q=2/d$) can be obtained  by using the summation trick ($(c)$ in Lemma \ref{s-trick}) as before. 
Indeed,  using the above two estimates we get restricted weak type $(p,q)$ estimate for $p,q$ satisfying 
$(\frac1p, \frac1q)\in\mathfrak  Q(\tfrac{d-2}2)$ and $\frac1p-\frac1q=\frac 2d$. We particularly have   
\[
 \|\sum_j \Pi_\lambda [\psi_j^{\kappa}] f\|_{\frac{2d(d-1)}{(d-2)^2},\infty}\le C\|f\|_{\frac{2(d-1)}{d},1},  \quad \kappa=\pm,  \;  \pm \pi.
 \]
Consequently, combining this with \eqref{psi0} we get \eqref{restricted-weak}.
\qed

\subsection{Estimate for $\Pi_\lambda[\psi_j](\sqrt \lambda x,\sqrt \lambda y)$
}\label{sec:local} 
We investigate the behavior of the kernel  of $\Pi_\lambda[\psi_j]$.
We obtain a bound on $\Pi_\lambda[\psi_j](\sqrt \lambda x,\sqrt \lambda y)$ in terms of $|\cD(x,y)|$ (see Lemma \ref{oscillatory}), which  is to be used to  prove  the local estimates in Theorem \ref{thm-locest}. The estimate for the kernel alone  is not enough to obtain the sharp  bound on $\chi_E\Pi_\lambda \chi_F$ as $|\cD(x,y)|$ gets close to $0$ over $E\times F$. 
However, when we handle the  estimates over the region near $\sqrt \lambda \mathbb S^{d-1}$, the results (Lemma \ref{oscillatory} and Corollary \ref{pj-det}) in this subsection become instrumental in controlling  the 
minor part, thus  those results lead us to distinguish the major part (see Section \ref{sec:endpoint}).

Let us set 
\begin{align}
\label{p-def} 
\cP &=\cP(x,y,s):=\frac{s}{2}+\frac{|x|^2+|y|^2}{2}\cot s-  \inp xy \csc s,
\\
\label{q-def} 
 \quad \cQ&=  \mathcal{Q}(x,y, \tau)
:=(\tau- \inp xy)^2-\mathcal D(x,y).
\end{align}
Note  that $\cP(x,y,s)=\phi_1(x,y,s)$. 
Also, for $j=0,1,\dots, $ we set 
\begin{align}\label{S3-ptkernel}
    I_j(x,y)=\fP_\lambda[\mathfrak a^{-1}\psi_j] (x,y)
    :=\int\psi_j(s) e^{i\lambda \mathcal{P}(x,y,s)}ds.   
    \end{align}
    Here $\psi_0\in C_c(0, \pi)$. 
Clearly,  $\Pi_\lambda[\mathfrak a^{-1}\psi_j](\sqrt \lambda x, \sqrt \lambda y)= I_j(x,y). $ We now note that  
\begin{align}
\label{ph-d}
   \partial_s \mathcal P(x,y,s) 
           =-\frac{\mathcal{Q}(x,y,\cos s)}{2\sin^2s},  
\end{align}
and the stationary point of $\mathcal P(x,y,\cdot)$ is  given by the zeros of $\cQ(x,y,\cos \cdot)$. 
 As $\cD(x,y)$ is the discriminant of the quadratic equation $\cQ(x,y,\tau)=0$,  the estimate 
 for $I_j$ is  controlled by the value of  $\cD(x,y)$ which regulates  the nature of stationary point of the phase function $\mathcal P(x,y,\cdot)$.  
 Under the condition  \eqref{det-con} these points are well separated, so we can obtain an improved bound when considering the operator $\Pi_\lambda$ over the set $\sqrt \lambda E$.  
 
 The following estimate is to be crucial  in proving the sharp local estimate (Theorem \ref{thm-locest}) and  the endpoint estimate \eqref{intro-conj}. 

\begin{lem}
\label{oscillatory} Let $\psi_0\in C_c^\infty (0, \pi)$, and for $j\ge 1$  let $\psi_j$ be a smooth function supported in $[2^{-2-j},2^{-j}]$ $($or $[-2^{-j},-2^{-2-j}]$$)$ and $|\frac{d^l}{dt^l}\psi_j|\le C2^{j l}$ for $l=0,1$.
 Then, if $|\mathcal D(x,y)|\neq 0$ and $|x|, |y|\le 2$,   we have
\begin{equation}
\label{osci-est}
 |I_j(x,y)|\lesssim  2^{-\frac{j}{2}}\lambda^{-\frac12}|\mathcal D(x,y)|^{-\frac14}
 \end{equation}
 for $j\ge 0$.
\end{lem}
The following is rather a straightforward  consequence of the above lemma. 
\begin{cor} 
\label{pj-det}
Let $E, F\subset B(0,2)$ be measurable.  Suppose 
that 
$
|\cD(x,y)|\ge \rho^2 
$
for  $( x,  y)\in E\times F$.
Then, for $p,q$ satisfying $1\le p\le 2$ and $1/p+1/q=1$, we have 
\begin{equation}   
\label{apjb} \| \chi_{\sqrt \lambda E} \Pi_\lambda  [ \psi_j ] \chi_{\sqrt \lambda F} f\|_{q}\lesssim  \lambda^{-\frac12(\frac1p-\frac1q)} 2^{( \frac{d+1}2(\frac1p-\frac1q)-1)j}   \rho^{-\frac12(\frac1p-\frac1q)}\|f\|_p
\end{equation}
 for $j\ge 0$.
\end{cor}

This is easy to show. By rescaling $(x,y)\to (\sqrt \lambda x,  \sqrt \lambda y)$, Lemma \ref{oscillatory} and \eqref{S3-ptkernel},  it follows that 
the kernel of $ \chi_{\sqrt \lambda E} \Pi_\lambda  [ \psi_j ] \chi_{\sqrt \lambda F}$ is bounded by $ C2^{\frac{d-1}{2}j}\lambda^{-\frac12} \rho^{-\frac12}$ since 
the same estimate holds for $\chi_{E} \fP_\lambda [ \psi_j ] \chi_F$. Thus,  we have $ \| \chi_{\sqrt \lambda E} \Pi_\lambda  [ \psi_j ] \chi_{\sqrt \lambda F} f\|_{\infty}\lesssim  2^{\frac{d-1}{2}j}\lambda^{-\frac12} \rho^{-\frac12} \|f\|_1$. 
We get  \eqref{apjb} by interpolating the above estimate and  the estimate 
$ \|  \chi_{\sqrt \lambda E} \Pi_\lambda  [ \psi_j ] \chi_{\sqrt \lambda F} f\|_{2}\lesssim 2^{-j}\|f\|_2$ which follows from \eqref{easy-l2}.

The rest of this section is devoted to the proof of Lemma \ref{oscillatory}.  We use  the following  which is  known as the Van der Corput lemma, for example see 
\cite[pp.\,332--334]{St93}.

\begin{lem}\label{osc-lemma} Let $\lambda\ge 1$ and let $\phi$ be a smooth function on an interval $I=[a,b]$ and $A\in \mathrm C^1$ supported in $[a,b]$.
Suppose that $| \frac{d^k}{ds^k}\phi|\ge L$ on $I$ and suppose additionally  that $\phi^\prime$ is monotone on $I$ when $k=1$. 
Then, the following  estimates hold with $C$ independent of $a,b$, $\phi$, $A$\,$:$
\begin{align*}
 \Big|\int_a^b e^{i\lambda\phi(s)} ds\Big|
  & \le{C}(\lambda L)^{-1/k},
  \\
   \Big|\int e^{i\lambda\phi(s)}A(s)ds\Big|
    &\le{C}(\lambda L)^{-1/k}\|A'\|_1.
\end{align*}
\end{lem}
In Lemma \ref{osc-lemma} the monotonicity assumption on $\phi^\prime$  is not essential. Since the first estimate holds independently of interval $[a,b]$,  this uniform estimate remains valid as long as the interval can be divided into a fixed number of  intervals on which $\phi'$ is monotone.
The second  estimate is a straightforward consequence of 
the first estimate  and integration by parts. See \cite[pp.\,332--334]{St93} for the detail. 

In addition, for later use we  record  the following which can be shown by repeated use of  integration by parts.

\begin{lem}\label{S3-osc} Let $0<\mu\le 1$ and let $\lambda\ge 1$. Suppose $A$ is a smooth function supported in an interval of length $\sim \mu$, and $\phi$ is smooth on the support of $A$. If  $|\phi'|\gtrsim  L$, $|\frac{d^k}{ds^k} \phi |\lesssim  L \mu^{1-k}$, $k= 1,\dots$, and 
$|\frac{d^k}{ds^k} A| \lesssim  \mu^{-k}$,  $k=0, 1,\dots$  on the support of $A$, then    
we have 
\Be
\label{high-order}
\Big|\int e^{i\lambda\phi(s)} A(s)ds\Big|
    \le{C}  \mu(1+\lambda\mu L)^{-N}
\Ee 
for  any $N>0$.
\end{lem}

Indeed, this can  be  shown by making change of variables $s\to \mu s+s_0$ where $s_0\in \supp A$. Setting  
$\tilde \phi(s)=(\mu L)^{-1}\phi(\mu s+s_0)$ and $\tilde A(s)=A(\mu s+s_0)$, we see that   $|\tilde\phi'|\gtrsim  1$, $|\frac{d^k}{ds^k} \tilde A| \lesssim  1$, and $|\frac{d^k}{ds^k} \tilde \phi |\lesssim  1$, $k=0, 1,\dots$  on the support of $\tilde A$. Then the integral is equal to $\mu\int e^{i\lambda\mu L\tilde \phi(s)}\tilde A(s)ds$. Thus, routine integration by parts gives the estimate
\eqref{high-order}.

\subsection{Proof of Lemma~\ref{oscillatory}}  
The proof we give below  is rather elementary and it is based on Lemma \ref{osc-lemma} and decomposition away from degeneracy points. However to get the estimate \eqref{osci-est}
we need to consider various cases, separately. 
Since  $\mathcal D(x,y)=\mathcal D(-x,y)$, we may assume \[ \inp xy\ge 0\] because  otherwise we need only to replace $x$ with $-x$. To show  Lemma~\ref{oscillatory} we consider the cases $\mathcal D(x,y) >0$ and $\mathcal D(x,y)\le 0$, separately. The second case is easier to handle.    
We first show  \eqref{osci-est} for  $j\ge 1$ and then  \eqref{osci-est} with $j=0$, which  is much easier.

\subsubsection*{\bf When $\mathcal D(x,y)\le 0$} Clearly, from \eqref{q-def} we have $\mathcal Q(x,y, \tau)\ge |\mathcal D(x,y)|$. 
We additionally  assume $\inp xy\le 1$ for the moment.  We handle the other case $\inp xy>1$ later.

Let  $S_\ast(x,y)\in [0, \pi/2]$ be the number such that 
\[
\cos \Sa=\inp xy.\]
We  distinguish the  two cases:
\begin{align*}
 {\rm I}: \Sa \not \in [2^{-3-j}, 2^{1-j}], \\ 
  {\rm I\!I}: \Sa \in [2^{-3-j}, 2^{1-j}].
   \end{align*}
We first consider the case ${\rm I}$. From \eqref{ph-d} and \eqref{q-def} we see  $|\partial_s \mathcal P(x,y,s)|\gtrsim 2^{-2j}+2^{2j} |\mathcal D|\gtrsim  |\mathcal D|^\frac12$ on the support of $\psi_j$. By Lemma {\ref{osc-lemma}} we get 
$ |I_j|\lesssim  \min(\lambda^{-1}|\mathcal D|^{-\frac12}, 2^{-j})$ 
and  \eqref{osci-est} follows.  
For the second case  ${\rm I\!I}$  we make a decomposition of the integral $I_j$ away from $\Sa$ as follows:
\[ I_j=\sum_{k\ge j-2} {I}_{j,k}:= \sum_{k\ge j-2} \int\widetilde\psi(2^k(s-\Sa))\psi_j(s)e^{i\lambda\mathcal P(x,y,s)}ds.\] 
If $s$ is contained in the support of the cutoff function $\psi(2^k(\cdot-\Sa))\psi_j$, by \eqref{ph-d} and \eqref{q-def} we have 
\[ |\partial_s \mathcal P(x,y,s)|\gtrsim 2^{2j}\Big( \big(\int_{\Sa}^s \sin \tau d\tau\big)^2 +|\mathcal D|\Big)  
\gtrsim 2^{2j} \Big( 2^{-2j}2^{-2k} +|\mathcal D|\Big) . \] 
Again, by Lemma {\ref{osc-lemma}} we see $|{I}_{j,k}| \lesssim \min(\lambda^{-1}2^{2k},  2^{-k}) $ if $2^k\le  2^{-j} |\mathcal D|^{-\frac12}$, and 
$|{I}_{j,k}| \lesssim  \min(\lambda^{-1}2^{-2j}|\cD|^{-1}, 2^{-k})$  if  $2^{k}>    2^{-j} |\mathcal D|^{-\frac12}$.
Hence, splitting the sum into two cases and  taking geometric mean of the two bounds in each case, we have 
\begin{align*}
   |{I}_{j}|   
\lesssim  \sum_{k:2^k \le 2^{-j} |\mathcal D|^{-\frac12}}  \lambda^{-\frac12} 2^{\frac k2}           \ +\sum_{k: 2^k >2^{-j} |\mathcal D|^{-\frac12}} 2^{-j}  \lambda^{-\frac12}   |\mathcal D|^{-\frac12}   2^{-\frac k2}   
  \lesssim 2^{-\frac{j}{2}}\lambda^{-\frac12}|\mathcal D|^{-\frac14} .        
\end{align*}
We now consider the case $ \inp xy\ge 1$. Clearly, we have 
 $(\cos s-\inp xy)^2
\ge  (1-\cos s)^2$. 
So, from \eqref{ph-d} and \eqref{q-def} we have $|\partial_s \mathcal P(x,y,s)|\gtrsim (2^{-2j}+ 2^{2j} |\mathcal D|)
\gtrsim \sqrt{|\mathcal D|}$ on the support of $\psi_j$. Thus,
by Lemma \ref{osc-lemma} we get $ |I_j|\lesssim  \min(\lambda^{-1}|\mathcal D|^{-\frac12}, 2^{-j}),$
which gives the desired bound.

{\bf When $\mathcal D(x,y)>0$.} Recalling  that  we are assuming  $\inp xy\ge 0$, we handle the case $\inp xy\le 1$ first.   Note that 
$\cQ(x,y,1)=|x-y|^2,$  $\cQ(x,y,-1)=|x+y|^2$. Since $0\le \inp xy\le 1$ and $\cD(x,y)>0$, we see that  $\cQ(x,y,\tau)=0$ has two distinct roots  $r_1> r_2$  which are separated by $2\sqrt{\cD}$ and $r_1, r_2\in [-1,1]$.  To show \eqref{osci-est}  we also assume \[ \cQ(x,y,0)\ge 0\] 
 for the moment. Thus $r_1, r_2$ are positive.  We handle  the remaining cases  $\cQ(x,y,0)<0$ and $\inp xy>1$ later.

Let $s_1$, $s_2\in [0, \pi]$ be positive numbers such that $\cos s_i =r_i$, $i=1,2$,  and $s_1< s_2$.  We distinguish the following four cases: 
\begin{align*}
    { \mathrm A}\!&: \ s_2< 2^{-3-j}, 
    \\
       \ {\mathrm  B}\!&:  \  2^{1-j}< s_1,  
       \\
        {\mathrm  C}\!&:  \   s_1< 2^{-3-j},    2^{1-j} <s_2,   
        \\  {\mathrm D}\!&:  \ s_1 \text{ or } s_2\in [2^{-3-j} , 2^{1-j}].
        \end{align*} 

{\bf Case $\mathrm A$.}  Since $0\le s_1<s_2$,  the distances between $s$ and $s_1$, and between $s$ and $s_2$ are $\sim 2^{-j}$ for $s\in \supp \psi_j$. We note 
\Be
\label{eq:phase} 
 |\partial_s \mathcal P(x,y,s)|\gtrsim   \frac{|\int_{s_1}^s \sin \tau d\tau \int_{s_2}^s \sin \tau d\tau|}{2\text{sin}^2s}.
 \Ee
 Thus, it follows that   $|\partial_s \mathcal P(x,y,s)|\gtrsim  2^{-2j}$ for $s\in \supp \psi_j$ and $2\sqrt{\cD}=r_1-r_2=\cos s_1-\cos s_2\lesssim 2^{-2j}$, so by Lemma \ref{osc-lemma}  we get 
$ |I_j|\lesssim  \min(\lambda^{-1} 2^{2j}, 2^{-j})\le \lambda^{-\frac12} 2^{\frac j2}\lesssim   \lambda^{-\frac12} 2^{-\frac j2} \mathcal D^{-\frac14}.$

{\bf Case $\mathrm B$.} Since the distance between $s$ and $s_1$ is $\gtrsim 2^{-j}$ and $ \cos s-  \cos s_2\ge  \cos s_1-\cos s_2 = 2\sqrt{ \mathcal D}$ for $s\in \supp \psi_j$,  
from \eqref{eq:phase} it follows that 
$ |\partial_s \mathcal P(x,y,s)|\gtrsim  \sqrt{\mathcal D}$ for $s\in \supp \psi_j$. 
By Lemma {\ref{osc-lemma}} and the trivial estimate we have 
$ |I_j|\lesssim   \min(\lambda^{-1}  \mathcal D^{-\frac12}, 2^{-j})$
and we get the desired estimate \eqref{osci-est}.

{\bf Case $\mathrm C$.}  As in the {\bf Case $\mathrm B$}, it is sufficient to show $|\partial_s \mathcal P(x,y,s)|\gtrsim  \sqrt{\mathcal D}$ for $s\in \supp \psi_j$.   
Since $s_1< 2^{-3-j}$ and $ 2^{1-j} <s_2$ we clearly have $2\sqrt{\cD}=\cos s_1-\cos s_2\gtrsim 2^{-2j}$. If $\sqrt{\cD}\sim 2^{-2j}$, 
from \eqref{eq:phase} we see $|\partial_s \mathcal P(x,y,s)|\gtrsim  \sqrt{\mathcal D}$ for $s\in \supp \psi_j$. If $\sqrt{\cD}\gg 2^{-2j}$,  $\cos s-\cos s_2\sim  \sqrt{\cD}$ because 
$\cos s_1-\cos s\sim 2^{-2j}$ for $s\in \supp \psi_j$.   Hence,  from \eqref{eq:phase} we see $|\partial_s \mathcal P(x,y,s)|\gtrsim  \sqrt{\mathcal D}$ for $s\in \supp \psi_j$.

{\bf Case $\mathrm D$.} In this case $\partial_s \mathcal P$ may vanish on the support of $\psi_j$, so we decompose the integral  $I_j$ away  from the zeros of $\partial_s \mathcal P$. 
We further  distinguish the two cases  
\[ 2^{-j}< 2^{-4} \mathcal D^\frac14, \quad  2^{-j}\ge  2^{-4}\mathcal D^\frac14.\]
If $2^{-j}< 2^{-4} \mathcal D^\frac14 $,  only  $s_1$  can be contained in $[2^{-4-j} , 2^{2-j}]$ because $s_1$ is positive. We 
  decompose $I_j$ as follows:
  \[    I_j=\sum_{k\ge j-2} \widetilde {I}_{j,k}:=\sum_{k\ge j-2}  \int\widetilde\psi(2^k(s-s_1))\psi_j(s)e^{i\lambda\mathcal P(x,y,s)}ds.\] 
 From \eqref{eq:phase} we see that
$ |\partial_s \mathcal P(x,y,s)| \gtrsim 2^j2^{-k}|\cos s-\cos s_2| \gtrsim 2^{j-k}\mathcal D^\frac12$ on the support of $\psi(2^k(\cdot-s_1))\psi_j$ because $|\cos s-\cos s_2|\ge |\cos s_1-\cos s_2|-|\cos s-\cos s_1|\gtrsim \mathcal D^\frac12 $ since we are assuming $2^{-j}< 2^{-4} \mathcal D^\frac14$. 
 Thus by Lemma {\ref{osc-lemma}} we get $  |\widetilde{I}_{j,k}|\lesssim  \min( \lambda^{-1} 2^{k-j} \mathcal D^{-\frac12}, 2^{-k} ).$
Summation over $k$ yields the  desired estimate \eqref{osci-est}.

We turn to the case $2^{-j}\ge  2^{-4}\mathcal D^\frac14$ where  the interval $[2^{-4-j} , 2^{2-j}]$ may  contain both $s_1$ and $s_2$. 
Taking it into account that  $|s_1-s_2|\sim 2^j\sqrt \cD$, 
we break the integral away from  these two vanishing points. 
Let us set 
\[
 \psi_\ast(s)= \widetilde \psi\Big( 2^{8-j}\cD^{-\frac12}\Big(s-\frac{s_1+s_2}2\Big)\Big)
 \]
and 
\begin{align*} 
{I}_{j, \ast} &=\int \psi_j(s)\psi_\ast(s)e^{i\lambda\mathcal P(x,y,s)}ds, 
\\
{I}_{j,k}^1&=\int\widetilde\psi(2^k(s-s_1))\psi_j(s)(1- \psi_\ast(s))e^{i\lambda\mathcal P(x,y,s)}ds,
\\
{I}_{j,k}^2&=\int\widetilde\psi(2^k(s-s_2))\psi_j(s) (1- \psi_\ast(s))e^{i\lambda\mathcal P(x,y,s)}ds.
\end{align*}

Then it follows that \[   I_j= {I}_{j, \ast}+ \sum_{k} {I}_{j,k}^1+ \sum_{k} {I}_{j,k}^2.\]
If $s$ is  contained in the support of $\psi_\ast$,  we see  $ |\partial_s \mathcal P(x,y,s)|\gtrsim   2^{2j}\cD$ using \eqref{eq:phase}. Thus, by Lemma \ref{osc-lemma} we get 
\[ |I_{j, \ast}|\lesssim \min(\lambda^{-1}2^{-2j} \mathcal D^{-1},   \mathcal D^\frac12 2^j)\lesssim \lambda^{-\frac12}2^{-\frac j2}|\cD|^{-\frac14}.
\] 
 From \eqref{eq:phase} we note that 
$
 |\partial_s \mathcal P(x,y,s)|
\gtrsim  2^j2^{-k}\cD^{\frac12}
 $
 on the support of $\wt\psi(2^k(\cdot-s_1))\psi_j(1-\psi_\ast)$.
By Lemma \ref{osc-lemma}  we have $  |{I}_{j,k}^1|\lesssim    \min( \lambda^{-1} 2^{k-j} \mathcal D^{-\frac12}, 2^{-k} ). $ 
Summation along $k$ yields $|\sum_{k}{I}_{j,k}^1|\lesssim   \lambda^{-\frac12} 2^{-\frac j2} \D^{-\frac14}$.
We can similarly handle  ${I}_{j,k}^2$. Note that $ |\partial_s \mathcal P(x,y,s)| \gtrsim   2^j2^{-k}\cD^{\frac12}$ on the support of $\wt\psi(2^k(\cdot-s_2))\psi_j(1-\psi_\ast)$ and, hence,  $| {I}_{j,k}^2|\lesssim    \min( \lambda^{-1} 2^{k-j} \mathcal D^{-\frac12}, 2^{-k} ).$ 
Summation over $k$ gives  $|\sum_k{I}_{j,k}^2|\lesssim  \lambda^{-\frac12} 2^{-\frac j2} \D^{-\frac14}.$
Combining the estimates for $|I_{j, \ast}|$, $|\sum_k{I}_{j,k}^1|$, and $|\sum_k{I}_{j,k}^2|$, we get the bound \eqref{osci-est}.

We now show \eqref{osci-est} for the cases $\cQ(x,y,0)<0$ and  $\inp xy>1$. 

{\bf When $\cQ(x,y,0)<0.$} In this case we have two roots $r_1$ and $r_2$ which  are  of different signs. 
This implies $s_2>\pi/2>s_1>0.$ We  consider the three cases: 
\[s_1\le 2^{-3-j}, \ \    2^{-3-j}\le s_1\le 2^{1-j}, \ \ s_1\ge 2^{1-j}. \] 
 If $s_1\le 2^{-3-j}$, $\mathcal D$ should be $\sim 1$ and $|\partial_s \mathcal P(x,y,s)|\gtrsim 1$ on the support of $\psi_j$. 
Thus, by Lemma \ref{osc-lemma} we get  $|I_j|\lesssim  \min (\lambda^{-1}, 2^{-j})$ and the desired estimate \eqref{osci-est}.  If  $ 2^{-3-j}\le s_1\le 2^{1-j}$, we have $\mathcal D\sim 1$ and  we make additional decomposition $I_j=\sum_k I_{j,k}$ where
\[{I}_{j,k}=\int\wt\psi(2^k(s-s_1))\psi_j(s)e^{i\lambda\mathcal P(x,y,s)}ds.\] 
Using \eqref{eq:phase} we note that $|\partial_s \mathcal P(x,y,s)|\gtrsim 2^{j} 2^{-k}$ on the support of $\wt\psi(2^k(\cdot-s_1))\psi_j$.   By Lemma \ref{osc-lemma} we have 
$|{I}_{j,k}|\lesssim \min(\lambda^{-1}2^{-j}2^k, 2^{-k})$. Summation over $k$ gives the desired estimate since $\mathcal D\sim 1$. Finally, we deal with the case $s_1\ge 2^{1-j}$ in which $|\partial_s \mathcal P(x,y,s)|\gtrsim  |\cD|^\frac12$ because $\cos s-\cos s_2 \ge \cos s_1-\cos s_2=2\sqrt \cD$ and $\cos s-\cos s_1 \gtrsim 2^{-2j}$.  Therefore, 
$|I_j|\lesssim \min(\lambda^{-1} |\cD|^{-\frac12}, 2^{-j})$ and we get the estimate \eqref{osci-est}.

{\bf When $\inp xy>1.$}  Since $\cQ(x,y,1)=|x-y|^2$, both roots $r_1, r_2$ are bigger than $1$.  Thus, we have
\[ |\partial_s \mathcal P(x,y,s)|\gtrsim 2^{2j} \big((r_1 -1)+(1-\cos s)\big)\big((r_2-r_1)+(r_1-\cos s)\big )
\gtrsim |\cD|^\frac12 \]
on the support of $\psi_j$. By Lemma  \ref{osc-lemma}  $| {I}_{j}|\lesssim  \min( \lambda^{-1} \mathcal D^{-\frac12}, 2^{-j} ).$  
This gives the desired estimate  \eqref{osci-est} and this completes the proof when $j\ge 1$. 

Finally we consider \eqref{osci-est} with $j=0$.  Since  $\psi_0$ is supported away from $\pi,$ $0$, this case  is much easier to show because there is no singularity in $s$ which is related to the  $\sin s$.  As before, one can prove the estimate \eqref{osci-est} considering the cases $\cD\le 0$ and $\cD>0$. The first case is easy to show. For the latter  we may assume $\cQ(x,y,0)\ge 0$ and $\inp xy\le 1$ since otherwise the previous argument works without modification. The rest of the argument is the same as before, so we omit the details. 
\qed

\section{Local estimate:  Proof of Theorem \ref{thm-locest}}
In this section we prove Theorem \ref{thm-locest} and show  the lower bounds on the $L^p$--$L^q$ operator norm of  $\chi_{\sqrt \lambda E}\Pi_\lambda\chi_{\sqrt{\lambda} E}$ including 
the failure of the estimate \eqref{est-loc} for $(1/p,1/q)\in [\mathfrak C,\mathfrak  D]\cup[\mathfrak C', \mathfrak  D']$. This shows 
that the bounds in Theorem \ref{thm-locest} are sharp. 
  In addition, when $q=p'$, we  obtain some estimates (Proposition \ref{away-sphere}) which strengthen the previous local estimates. 
We need those estimate  for the proof  of Theorem \ref{endpoint}. 

\subsection{Proof of Theorem \ref{thm-locest}}  In order to prove  Theorem \ref{thm-locest} we  mainly use  the decomposition \eqref{decomp-proj},  the estimates in Lemma  \ref{oscillatory},  and Lemma \ref{tt-st}. The following provides the sharp $L^1$--$L^2$ estimate for $ \Pi_\lambda [\psi_j]$.

\begin{lem}\label{l-infty} Let $\psi_j\in C_c^\infty(-\pi, \pi)$ be a smooth function with its support contained in an interval of length $\sim 2^{-j}$. 
Suppose  $ 2^j\lesssim \lambda$ and $|\frac{d^l}{dt^l}\psi_j|\le C2^{j l}$ for $l=0, 1,2,\dots$, then  we have
\begin{align}
\label{ho12}
 &\| \Pi_\lambda [\psi_j] f\|_2 \lesssim  2^{-\frac{j}{2}}\lambda^{\frac{d-2}{4}}\|f\|_1,    
 \\
 \label{ho10} &  \| \Pi_\lambda [\psi_j] f\|_\infty \lesssim  \lambda^{\frac{d-2}{2}}\|f\|_1. 
\end{align}
\end{lem}

\begin{proof}
From  \eqref{proj-pi} and  \eqref{schro-op}  we have
\Be
\label{spectral-proj}
\Pi_\lambda [\phi]   f=  \sum_{\lambda'} \widehat {\phi} ( 2^{-1}(\lambda'-\lambda))\Pi_{\lambda'} f
\Ee
for any smooth function $\phi\in C_c^\infty(-\pi,\pi)$. Orthogonality between the Hermite projections  gives 
\[ \| \Pi_\lambda [\psi_j] f\|_2^2
                     =   \sum_{\lambda' }   | \widehat {\psi_j} ( 2^{-1}(\lambda'-\lambda))|^2\|\Pi_{\lambda'} f\|_2^2
                                   = \sum_{\lambda' }   | \widehat {\psi_j} ( 2^{-1}(\lambda'-\lambda))|^2\inp{\Pi_{\lambda'} f}{f} \] 
                                    since $\Pi_\lambda^2= \Pi_\lambda$. 
Using $|\widehat \psi_j(\tau)|\lesssim 2^{-j}(1+ 2^{-j}|\tau|)^{-N}$ for any $N$ and  the estimate \eqref{plambda-100},  we get 
\begin{align*}
         \| \Pi_\lambda [\psi_j] f\|_2^2
                    \lesssim     2^{-2j} \sum_{\lambda' } \big( 1+ 2^{-j}|\lambda-\lambda'| \big)^{-2N}  (\lambda')^{\frac {d-2}2} \|f\|_1^2  \lesssim     2^{-j}\lambda^{\frac {d-2}2} \|f\|_1^2
\end{align*}
since  $2^j\lesssim  \lambda$.
This gives the estimate \eqref{ho12}. The proof of  \eqref{ho10} is similar. Using \eqref{spectral-proj} as before and 
\eqref{plambda-100},  we get 
$\| \Pi_\lambda [ \psi_j ]  \|_{1\to \infty}\lesssim  
2^{-j} \sum_{\lambda'}  (1+2^{-j}|\lambda-\lambda'|)^{-N} (\lambda')^\frac{d-2}2\lesssim  \lambda^\frac{d-2}2 $ 
because $2^j\lesssim \lambda$. 
\end{proof}

\begin{proof}[Proof of Theorem \ref{thm-locest}]  
We begin with recalling  \eqref{decomp-proj}. Since $|\mathcal D(x,y)|>c$ for $x,y\in E$,  from  Lemma \ref{oscillatory} and \eqref{symmetric}  we have 
\[    
\|   \che \Pi_\lambda [\psi_j^{\kappa} ] \che \|_{1\to \infty} \lesssim  \lambda^{-\frac12} 2^{\frac{d-1}{2}j}  
\]
for $j\ge 4$  and $\kappa=0,  \pm,  \pm \pi. $
Applying Lemma \ref{tt-st} and \eqref{symmetric} immediately yields 
\begin{equation} 
\label{basic} 
\begin{aligned} 
   \|   \che \Pi_\lambda [\psi_j^{\kappa} ] \che \|_{p\to q}&\lesssim  \lambda^{-\frac12\dpq} 2^{(\frac{d+1}{2}\dpq-1)j}, \  \  \  \kappa=0,  \pm,  \pm \pi
\end{aligned} \end{equation}
provided that \ppqpair\,  is contained in the close quadrangle $ \mathfrak Q(\tfrac{d-1}2)$ which has vertices 
$(\tfrac 12, \tfrac12)$,  $\mathfrak A=\big( \frac{d+3}{2(d+1)},\frac12\big)$, $\mathfrak A'$ and $(1,0).$

The bounds on  $ \| \che \Pi_\lambda [\psi^{0} ] \che\|_{p\to q}$ are easy to show. For $\ppq\in  \mathfrak Q(\tfrac{d-1}2)$, which contains $\mathcal R_1$,  the desired bound $ \| \che \Pi_\lambda [\psi^{0} ] \che\|_{p\to q} \lesssim \lambda^{-\frac12 \dpq} $  follows from  \eqref{basic}. We have the estimate \eqref{basic} with $\ppq=\mathfrak C, (1,0)$ and $\kappa=0$,  and  $\| \che \Pi_\lambda [\psi^{0} ] \che\|_{1\to 2} \lesssim \lambda^{(d-2)/4} $ from the estimate 
 \eqref{ho12}. Interpolation between these estimates gives the better bound $\|\che \Pi_\lambda [\psi^{0} ] \che \|_{p\to q}\lesssim \lambda^\gamma$ with $\gamma< {\beta(p,q)} $ provided that $\ppq$ is in  the closed triangle with vertices $(1, \frac12), \mathfrak C,\mathfrak (1,0)$ which contains $ (\mathfrak C, \mathfrak D].$  Thus, by  interpolation and  duality  we have
 \Be 
 \label{pi0}
 \|\che \Pi_\lambda [\psi^{0} ] \che \|_{p\to q}\lesssim \lambda^{\beta(p,q)}
 \Ee
  for $1\le p\le 2\le q\le \infty.$   

Consequently, we  have only to consider the estimates  for  $ \sum_{j\ge 4}   \che \Pi_\lambda [\psi_j^{\kappa} ] \che$, $\kappa= \pm,$ $ \pm \pi$, and show 
 the same bound on those operators  as in Theorem \ref{thm-locest}, that is to say, 
\Be
\label{eq:j-sum}
\| \sum_{j\ge 4}   \che \Pi_\lambda [\psi_j^{\kappa} ] \che \|_{p\to q }\lesssim  \lambda^{\beta(p,q)}, \  \  \  \kappa= \pm, \  \pm \pi.
\Ee
 Using \eqref{basic}, summation over $j$ gives   
the estimate  \eqref{eq:j-sum}  
provided that $\ppq\in \mathcal R_1\setminus[\mathfrak C, \mathfrak C']$.  Furthermore, by  $(c)$ in Lemma \ref{s-trick}  we have 
\begin{equation}
\label{rest-weak}
     \| \sum_{j\ge 4}   \che \Pi_\lambda [\psi_j^{\kappa} ] \che \|_{L^{p,1} \to L^{q,\infty} }\lesssim  \lambda^{-\frac12\dpq}, \  \  \  \kappa= \pm, \  \pm \pi    
     \end{equation} 
     for  $\ppq=\mathfrak C  $ and $\ppq= \mathfrak C'$,  which satisfy $\dpq=2/(d+1)$. Interpolation between the estimates in \eqref{rest-weak}  establishes 
the estimate \eqref{eq:j-sum} for $\ppq\in \mathcal R_1$ since we already have \eqref{eq:j-sum}  for $\ppq\in \mathcal R_1\setminus[\mathfrak C, \mathfrak C']$. 
This also proves the assertion $(ii)$ in Theorem  \ref{thm-locest}.

We have established the estimate \eqref{est-loc} with  $(p,q)=(2,2)$, $(1,\infty),$ and $(\frac{2(d+1)}{d+3}, 2)$ and the restricted weak type 
$(p,q)$ estimate for $\ppq=\mathfrak C  $   by the estimates in the above. Thus, to complete  the proof we need only to show \eqref{eq:j-sum} with $(p,q)=(1,2)$ and 
the weak type  estimate 
\Be 
\label{weak-pq} \| \sum_{j\ge 4}   \che \Pi_\lambda [\psi_j^{\kappa} ] \che \|_{L^p\rightarrow L^{\frac{2d}{d-1},\infty}} \lesssim   \lambda^{ \frac d2(\frac1p-\frac{d-1}{2d})- 1}, \  \  \  \kappa= \pm, \ \pm \pi, 
\Ee
for $1\le p<\frac{2d(d+1)}{d^2+4d-1}$,   which shows the assertion $(i)$  because we have already obtained \eqref{pi0}. 
Duality and interpolation provide
 all the $L^p$--$L^q$ estimates  stated in Theorem \ref{thm-locest}  (see Figure \ref{local-type}).  We also note  $\beta(1,2)=\frac {d-2}4$, $\beta(p,q)=\frac d2(\frac1p-\frac{d-1}{2d})- 1$ if $(1/p,1/q)\in(\mathfrak C,\mathfrak D]$, and $\beta(p,q)=-\frac{1}{d+1}$ if $(1/p,1/q)=\mathfrak C$.

Since $ \sum_{j: 2^j> \lambda }\psi_j^\kappa=  \widetilde {\psi}^\kappa(\lambda\cdot)$ $a.e.$ for some $\widetilde {\psi}^\kappa\in C_c^\infty(-\pi,\pi),$ we write
\[\sum_{j\ge 4  }  \Pi_\lambda [\psi_j^\kappa] 
 = \sum_{j: 2^j\le \lambda }  \Pi_\lambda [\psi_j^\kappa]  +  \Pi_\lambda [\widetilde {\psi}^\kappa(\lambda\,\cdot)] . \]
From Lemma \ref{l-infty} we have 
\begin{align}
  \| \Pi_\lambda [\widetilde \psi^\kappa(\lambda \cdot)] f\|_2& \lesssim  \lambda^{-\frac12}\lambda^{\frac{d-2}{4}}\|f\|_1,  \  \  \  \kappa= \pm, \ \pm \pi, 
   \nonumber
 \\
   \| \Pi_\lambda [\psi_j^\kappa] f\|_2 &\lesssim  2^{-\frac{j}{2}}\lambda^{\frac{d-2}{4}}\|f\|_1, \  \  \  \kappa= \pm, \ \pm \pi,  \ \ \ 2^j\lesssim \lambda.
   \label{12j}
\end{align}
Then,  \eqref{eq:j-sum} with $(p,q)=(1,2)$ is clear from these two estimates.

It now remains to show \eqref{weak-pq}. 
We first deal with $\sum_{j: 2^j\le \lambda }  \Pi_\lambda [\psi_j^\kappa] $. 
Interpolation between the estimates  \eqref{12j} and   \eqref{basic} with $\ppq=(1,0),$ $\mathfrak C$ gives 
\begin{align*}
    \| \che \Pi_\lambda [\psi_j^{\kappa} ] \che\|_{p\rightarrow q}\lesssim2^{jd(\frac{d-1}{2d}-\frac1q)}\lambda^{ \frac d2(\frac1p+\frac1q)- \frac{d+1}{2}}, \  \  \  \kappa= \pm, \ \pm \pi
\end{align*}
if  $(1/p,1/q)$ is in the closed triangle  $\mathscr T$ with vertices $(1, \frac12), \mathfrak C,(1,0)$. 
Once we have this estimate, fixing $p\in [1, {2d(d+1)}/{(d^2+4d-1)})$ and choosing two $q_0, q_1$ such that $ q_0<2d/(d-1)<q_1$ and $(1/p, 1/q_0), (1/p, 1/q_1)\in \mathscr T$,  we have two estimates from the previous estimates with $p=p_0=p_1$ and $q=q_0, q_1$. Then, we  apply the summation trick ($(a)$ in Lemma \ref{s-trick}) to  these two estimates to get the weak type estimate for 
$p,q$  such that  $\ppq\in (\mathfrak C,\mathfrak D]$.\footnote{Figure \ref{local-type} is  helpful here.} Hence, for $1\le p<\frac{2d(d+1)}{d^2+4d-1}$, we establish  the estimate
\begin{align*}
    \| \sum_{j: 2^j\le \lambda }  \che \Pi_\lambda [\psi_j^{\kappa} ] \che\|_{L^p\rightarrow L^{\frac{2d}{d-1},\infty}}\lesssim  \lambda^{ \frac d2(\frac1p-\frac{d-1}{2d})- 1}, \  \  \  \kappa= \pm, \ \pm \pi. 
\end{align*}
Now we handle  $\Pi_\lambda [\widetilde {\psi}^\kappa(\lambda\cdot)]$. From \eqref{ho12} and \eqref{ho10} we have  
\begin{align*}
\| \che \Pi_\lambda [\widetilde \psi^\kappa(\lambda\cdot)]  \che \|_{1\to 2} \lesssim \lambda^{\frac{d-4}4}, 
\\ \| \che \Pi_\lambda [\widetilde \psi^\kappa(\lambda\cdot)]  \che \|_{1\to \infty} \lesssim \lambda^{\frac{d-2}2} 
\end{align*} 
 for $\kappa= \pm, \ \pm \pi$. 
Since $ \che \Pi_\lambda [\widetilde \psi^\kappa(\lambda\cdot)]  \che= \sum_{j: 2^j> \lambda } \che \Pi_\lambda [ \psi_j^\kappa] \che$, the estimate \eqref{basic} and Lemma \ref{s-trick} give the restricted weak type $(p,q)$ estimate when $\ppq=\mathfrak C$ with bound $\lambda^{-\frac12\dpq} $. Real interpolation among those estimates gives 
\[
    \big\|\che \Pi_\lambda [\widetilde \psi^\kappa(\lambda\cdot)]  \che \big\|_{p\to q}\lesssim  \lambda^{ \frac d2(\frac1p-\frac1q)-1} , \  \  \  \kappa= \pm, \ \pm \pi
\]
provided  that  $\ppq$ is in  the closed triangle with vertices $(1, \frac12), \mathfrak C,(1,0)$ from which  $ \mathfrak C$ is excluded. This clearly 
yields  the desired estimate  for the operator  $\che \Pi_\lambda [\widetilde \psi^\kappa(\lambda\cdot)]  \che$ if $\ppq\in (\mathfrak C,\mathfrak D]$. This  completes the proof of 
\eqref{weak-pq}.
 \end{proof}

\subsection{Lower bound for $\|\chi_{ \sqrt\lambda E} \Pi_\lambda \chi_{\sqrt\lambda E}\|_{p\to q}$ } 
We prove the upper bound  \eqref{est-loc} in Theorem \ref{thm-locest} is sharp. Recalling that 
$\cD(x,y)\gtrsim 1$ for $x,y \in B(0,1/2)$, 
 for the purpose it is enough to show the following.

\begin{prop}\label{prop:lolo} Let $\widetilde \chi_{\sqrt \lambda}=\chi_{ B(0,\sqrt \lambda/2)}$. 
Let $d\ge2$, $\lambda\gg 1$, and $1\le p, q\le \infty.$ Then, we have 
\begin{align}
\| \widetilde \chi_{\mathsmaller {\sqrt{\lambda}}}\, \Pi_\lambda \widetilde \chi_{\mathsmaller {\sqrt{\lambda}}}\|_{p\to q} &\gtrsim \lambda ^{\frac{d-1}2- \frac d2(\frac1p+\frac1q)},  
                   \label{local2}
\\
\| \widetilde \chi_{\mathsmaller{\sqrt{\lambda}}}\,  \Pi_\lambda \widetilde \chi_{\mathsmaller {\sqrt{\lambda}}}\|_{p\to q} &\gtrsim \lambda ^{-\frac12(\frac1p-\frac1q)},    
                  \qquad\, 1\le p<4, 
                      \label{local1}
\\
\| \widetilde \chi_{\mathsmaller{\sqrt{\lambda}}}\, \Pi_\lambda \widetilde \chi_{\mathsmaller {\sqrt{\lambda}}}\|_{p\to q}  &\gtrsim \lambda ^{-1+\frac d2(\frac1p-\frac1q)} .  
\label{local3}
\end{align}
\end{prop}

In order to prove Proposition \ref{prop:lolo}, we use the asymptotic properties of the Hermite functions which are the eigenfunctions of the one-dimensional Hermite operator. Let $h_k(t)$ denote the $L^2$-normalized $k$-th Hermite function  of which eigenvalue is  $2k+1$. We make use of the following lemma from \cite{T05}.   Also see \cite{AAR} and \cite{F93}.

\begin{lem}\cite[Lemma 5.1]{T05}\label{S4-asym} Let $\mu=\sqrt{2k+1}$. We set 
\[
    s^-_\mu(t)=\int_0^t\sqrt{|\tau^2-\mu^2|}d\tau\quad \text{and}\quad s^+_\mu(t)=\int^t_\mu\sqrt{|\tau^2-\mu^2|}d\tau. 
\]
Then,  the following hold\,$: $
\[  \  \  \  h_{2k}(t)
                =\begin{cases}
                       a_{2k}^-(\mu^2-t^2)^{-\frac14} \big(\!\cos s^-_\mu(t)+\mathcal E\big),     & \ \qquad\qquad\quad \ |t|<\mu-\mu^{-\frac13}, 
                        \\
                    \qquad \qquad O(\mu^{-\frac16}),                                           & \ \ \, \mu-\mu^{-\frac13}<|t|<\mu+\mu^{-\frac13}, 
                     \\
                      a^+_{2k}e^{-s^+_\mu(|t|)}(t^2-\mu^2)^{-\frac14}\big(\!1+\mathcal E\big),     &\ \,  \ \mu+\mu^{-\frac13}<  |t|,
                \end{cases}
\]
\vspace{4pt}
\[    h_{2k+1}(t)
              =\begin{cases}
               a_{2k+1}^-(\mu^2-t^2)^{-\frac14}\big(\! \sin\,s^-_\mu(t)+\,\mathcal E\big),   & \qquad\qquad \  \, \ |t|<\mu-\mu^{-\frac13},  
               \\
              \qquad \qquad  O(\mu^{-\frac16}),      & \mu-\mu^{-\frac13}<|t|<\mu+\mu^{-\frac13},  
              \\
             a^+_{2k+1}e^{-s^+_\mu(|t|)}(t^2-\mu^2)^{-\frac14}\big(\!1+\mathcal E\big),    &\mu+\mu^{-\frac13}<|t|,
             \end{cases}
             \]
where
$    |a^{\pm}_k|\sim 1$  and $ \mathcal E=O\big(|t^2-\mu^2|^{-\frac12}\left||t|-\mu\right|^{-1}\big) .$
\end{lem}

We also make use of the following lemma concerning the asymptotic behavior of  the $L^p$ norm of $h_k$.  
A more precise result  can be found in \cite{Cohen}.  In fact, the following  is not difficult to show  using Lemma \ref{S4-asym}. 

 \begin{lem}\cite[Lemma 1.5.2]{Th93}\label{S4-bd_Her} Let $k\ge 2$. Then, we have  
  \begin{equation}\label{bd1} \|h_k\|_{L^p(\R)}\sim  	\begin{cases}k^{\frac1{2p}-\frac14},&1\le p<4,  \\     	k^{-\frac18}\log{k},&\quad p=4,  \\     	k^{-\frac1{6p}-\frac1{12}},& 4<p\le\infty.  	\end{cases}\end{equation} \end{lem}

To prove Proposition \ref{prop:lolo}, we choose suitable functions $f$  which yield  the particular lower bounds on the operator norms of $\schi  \Pi_\lambda\schi$. Some of the constructions of such function $f$ are inspired by those in Koch and Tataru \cite{T05}. 
 Compared with the examples in  \cite{T05} where only $f\in L^2$ are considered,  we have to handle functions in  $L^p$.  This gives a rise to additional difficulty in obtaining the precise lower bounds by exploiting orthogonality among  involved functions which are no longer in $L^2$.

\subsubsection*{Proof of \eqref{local2}}
Let $\lambda$ be the $N$-th eigenvalue,  $\lambda=2N +d$, and let   $\ell={(2\sqrt{d})^{-1}}{\sqrt{\lambda}}$. Set
\begin{align*}
    J=\big\{\alpha\in\mathbb N^d \,: \,|\alpha|=N,\; {N}/(16d) \le\alpha_j\le {N}/(8d),  \ \  2\le j\le d\big\},
\end{align*}
and $Q_{\ell}=\{x\in\mathbb{R}^d : |x_j|\le\ell\}$. For each $\alpha\in J$,  let  $c_\alpha\in \{-1,1\}$,  which is  to be chosen later. We define a function $f$ by
\[
    f = \chi_{Q_\ell}(x)\sum_{\alpha\in J}c_\alpha\Phi_{\alpha}(x).
\]
Then, $
    \|f\|_{L^p}\le\big\|\sum_{\alpha\in J}c_\alpha\Phi_{\alpha}(x)\big\|_{L^2}|Q_\ell|^{\frac1p-\frac12}\lesssim\lambda^{\frac{d-1}{2}+\frac{d}{2}(\frac1p-\frac12)}
$ by H\"older's inequality and orthonormality of $ \{\Phi_{\alpha} \}$. Hence, 
for \eqref{local2} it is enough to show that   
\begin{equation}\label{co4_low}
  \|\schi \Pi_\lambda(\schi f)\|_{L^q(\mathbb{R}^d)}\gtrsim\lambda^{\frac{3d}4-\frac{d}{2q}-1},\quad \lambda\gg 1.
\end{equation}
We  
set 
\begin{align*}
 a_{u,v}&:=\int_{-\ell}^{\ell}h_u(t)h_v(t)dt,
\\
A_{\alpha,\beta}&:=\int_{Q_\ell}\Phi_{\alpha}(x)\Phi_{\beta}(x)dx.
\end{align*} 
Then, $A_{\alpha,\beta} = \prod_{i=1}^d a_{\alpha_{i},\beta_{i}}$
and  
$\schi \Pi_{\lambda}(\schi f)(x) 
    = \sum_{\alpha\in J}  \sum_{\beta: |\beta|=N} c_\alpha A_{\alpha, \beta} \Phi_\beta(x)$  for $x\in B(0, \sqrt\lambda/2)$. So, we may write
\begin{align}
\label{sumsum}
    \schi &\Pi_{\lambda}(\schi f)(x) 
       =  \sum_{\alpha\in J}c_{\alpha}A_{\alpha, \alpha}\Phi_{\alpha}(x)+ \sum_{\alpha\in J} \  \ \sum_{\beta: |\beta|=N, \beta\not= \alpha} c_{\alpha} A_{\alpha, \beta} \Phi_{\beta}(x) .
      \end{align}
Using Lemma~\ref{S4-asym}, it is easy to see that $a_{u,u}\sim 1$ if $u\sim N$.  Consequently, it follows that 
$A_{\alpha, \alpha}\sim 1 \ \text{ if } \ \alpha\in J. $
However, if $u\neq v$, we can show $a_{u,v}$  is exponentially decaying, thus we  can disregard the terms $A_{\alpha,\beta}$ with $\alpha\neq \beta$.  
More precisely, we claim 
\Be\label{auv-final}
A_{\alpha, \beta}\lesssim e^{-c\lambda}
\Ee 
if $\alpha\in J$, $|\beta|=N$, and $\alpha\neq \beta$.  Assuming this for the moment, we proceed to show 
\eqref{co4_low}. 
Since there are as many as $O(N^{2d-2})$  of $(\alpha, \beta)$ in the summation,  from \eqref{sumsum}  and \eqref{auv-final}  it follows that 
\Be  
\label{ha4}
\schi \Pi_{\lambda}(\schi f)(x) = \sum_{\alpha\in J}c_{\alpha}A_{\alpha, \alpha}\Phi_{\alpha}(x)+O(N^C e^{-cN})
\Ee
for some $C, c>0$.   Let $\alpha\in J$. Since $\alpha_j\sim N$ for each $j$ and $t\sqrt{\mu^2-1} \le  s^-_{\mu}(t)\le t\mu $ if $0\le t\le 1$ with $\mu=\sqrt{2\alpha_j+1}$, 
we can choose a constant $c>0$ such that  $\sin s^-_{\mu}(t)\sim 1$ and  $ \cos s^-_{\mu}(t)\sim  1$    if $ t\in (c/\sqrt{\lambda}, 2c/\sqrt{\lambda})$. Therefore,  each $h_{\alpha_j}(x_j)$ has the same sign with its absolute value $\sim \lambda^{-\frac14}$ if $ x_j\in (c/\sqrt{\lambda}, 2c/\sqrt{\lambda})$. By Lemma \ref{S4-asym}, we can choose $c_\alpha\in\{-1,1\}$ such that $c_\alpha\Phi_{\alpha}(x)$ has positive value $\sim \lambda^{-\frac d4}$  if  $ x_j\in (c/\sqrt{\lambda}, 2c/\sqrt{\lambda})$, $j=1, \dots, d$. We note that $|J|\sim\lambda^{d-1}$. Thus, by our  choice of $c_\alpha$  we have 
\[\|\sum_{\alpha\in J}c_{\alpha}A_{\alpha, \alpha}\Phi_{\alpha}\|_{L^q\big((c/\sqrt{\lambda}, 2c/\sqrt{\lambda})^d\big)}\gtrsim \lambda^{\frac{3d}4-\frac{d}{2q}-1}. \]
 By this and \eqref{ha4}  we get the desired \eqref{co4_low} provided that $\lambda$ is large enough.

In order to complete the proof of \eqref{local2}
it remains to show \eqref{auv-final}. 
Recalling the identity $2(u-v)h_uh_v =h_uh_v''-h_u''h_v$  \cite[pp 1-3, Chapter 1]{Th93}, we have 
\[a_{u,v}
                          =\frac{1}{2(u-v)}\int_{-\ell}^{\ell} h_u(s)h^{''}_v(s)-h_u^{''}(s)h_v(s)ds.\]
                          Thus, by  integration by parts we see that  
\begin{align*}
    a_{u,v}
                =\frac{1}{2(u-v)} \Big(h_u(\ell)h^{'}_v(\ell)-h_u(-\ell)h^{'}_v(-\ell)-h_u^{'}(\ell)h_v(\ell)+h_u^{'}(-\ell)h_v(-\ell)\Big)
\end{align*}
 if $u\neq v  $. 
Since $h_u$ is odd if $u$ is odd and $h_u$ is even otherwise,  $h_u h^{'}_v$ is even if $u+v$ is odd and $h_u h^{'}_v$ is odd otherwise. 
Thus, $a_{u,v}=0$ if $u+v$ is odd.  Using the identity $h_u'(s)=sh_u(s)-\sqrt{2u+2}\,h_{u+1}(s) $(\cite[pp. 1--3, Chapter 1]{Th93}), we may rewrite the last displayed  as follows: 
\Be\label{auv}
   a_{u,v}
 =\frac{1+(-1)^{u+v}}{\sqrt{2}(u-v)}\big(\sqrt{u+1}\,h_{u+1}(\ell)h_v(\ell)-\sqrt{v+1}\,h_u(\ell)h_{v+1}(\ell)\big).
\Ee
Since $\ell-\sqrt{2u+1}\gtrsim \sqrt\lambda\sim \ell  $ if $2^{-4}N/d\le u\le 2^{-3}N/d$,
we have $s^+_{(2u+1)^{1/2}}(\ell)\gtrsim \lambda$. Hence,  by Lemma \ref{S4-asym}  it follows that 
\[|h_u(\ell)|\lesssim e^{-s^+_{\sqrt{2u+1}}(\ell)} \lesssim e^{-c\lambda} \]
if $2^{-4}N/d\le u\le  2^{-3}N/d$.
Combining this with \eqref{auv} and \eqref{bd1},  we have 
\[
    |a_{u,v}| \lesssim e^{-c\lambda}
\]
for $2^{-4}N/d\le u\le  2^{-3}N/d$ and $v\le N$ if $u\neq v$.  Note that  $ |a_{u,v}| \lesssim 1$ for any $u,v$. Clearly,  there is at least one $j\in \{2, \dots, d\}$ such that  $\alpha_j\neq \beta_j$ if $\alpha\neq \beta$, $\alpha\in J$, and $|\beta|=N$. Therefore, using the above estimate,  we 
get \eqref{auv-final} because $A_{\alpha,\beta} = \prod_{i=1}^d a_{\alpha_{i},\beta_{i}}$. 
  \qed

\subsubsection*{{Proof of \eqref{local1}}}  
We write  $x=(x_1,\ol x)\in \R\times \R^{d-1}$ and set
\[    f(x):=\Phi_{\alpha_\circ}(x)\chi_{Q_\ell}(x)=h_N(x_1)e^{-\frac{|\ol x|^2}{2}}\chi_{Q_\ell}(x)\]
with $\alpha_\circ=(N,0,\cdots,0)$.  Here $Q_\ell$ is the same  as in the proof of \eqref{local2}.
 From Lemma \ref{S4-bd_Her} and Lemma \ref{S4-asym}, it is clear that
$\|f\|_{L^p}\lesssim\lambda^{\frac1{2p}-\frac14}$
for $1\le p<4$ since $\lambda =2N+d.$ Moreover,  by the same argument as before  (see \eqref{sumsum})  we easily see 
\[  | \schi \Pi_{\lambda}(\schi f)(x) |\gtrsim | \Phi_{\alpha_\circ}(x)|-  O(\lambda^C e^{-c\lambda})
\]
because only the diagonal term has significant contribution. So, using Lemma  \ref{S4-asym} we get
$ \| \Phi_{\alpha_\circ}\|_{L^q(\{x :|x|<\sqrt {\lambda }/2\})}\gtrsim  \lambda^{\frac1{2q}-\frac14}$. Therefore,   \eqref{local1} follows if  $\lambda$ is large enough. 
\qed

\subsubsection*{Proof of \eqref{local3}} 
We begin with claiming  that, for any $0\le u\le N$,
\begin{align}
\label{hlowerbd}
    |h_u(t)|\gtrsim (2u+1)^{1/4}\lambda^{-\frac12}, \quad t\in [2^{-4}\lambda^{-\frac12}, 2^{-3} \lambda^{-\frac12}].
\end{align}
To see this, we set $\mu=\sqrt{2u+1}$ for $u\le N$. Then, $\mu \le\sqrt \lambda$. By a simple calculation, we have
$    |t|\mu/{\sqrt 2}\le s^-_\mu(t)\le |t|\mu $ for $|t|\le\mu/2$ where $s^-_{\mu}(t)$ is given  in Lemma \ref{S4-asym}. Since 
$ \cos{s^-_\mu(t)}\sim 1$ and $\sin{s^-_\mu(t)}\sim \mu t$  if $ t\in [2^{-4}\lambda^{-\frac12}, 2^{-3} \lambda^{-\frac12}],$
 the claim \eqref{hlowerbd} follows from Lemma \ref{S4-asym}.

Let $f=\chi_{q_{\lambda}}$ where $q_\lambda=\{x\in\mathbb{R}^d: 2^{-4}/\sqrt\lambda <|x_j|<2^{-3}/\sqrt\lambda, \ j=1,\dots, d\}.$ Then we have
\begin{align*}
    \|\schi \Pi_\lambda(\schi f)\|_{L^q}
    &\ge \Big\|\sum_{\alpha:|\alpha|=N}\Phi_\alpha\int_{q_\lambda}\Phi_\alpha(y) dy\Big\|_{L^q(q_\lambda)}.
\end{align*}
By \eqref{hlowerbd} it follows that  $ \Phi_\alpha(x)\Phi_\alpha(y)\gtrsim\lambda^{-d} \prod_{j=1}^d(2\alpha_j+1)^{1/2}$ for  $x,y\in q_\lambda$. Thus, 
the right hand side of the above inequality is bounded below by 
\begin{align*}
         C\lambda^{-\frac {3d}2-\frac{d}{2q}}\sum_{\alpha:|\alpha|=N}\prod_{j=1}^d(2\alpha_j+1)^{1/2} \gtrsim\lambda^{-1-\frac {d}{2q}}. 
\end{align*} This gives 
$\|\schi \Pi_\lambda(\schi f)\|_{L^q}\gtrsim\lambda^{-1-\frac {d}{2q}}$.
Since $\|f\|_{L^p}\sim \lambda^{-\frac{d}{2p}}$, we get \eqref{local3}.  \qed

\subsection{Boundedness of  the operator  $\wp_k$ and a transplantation result 
}\label{transplantation}
We consider the estimate
\begin{equation}\label{lplq} \|\chi_{\sqrt{\lambda} B}\Pi_{\lambda}\chi_{\sqrt{\lambda}B}\|_{p\to q}\lesssim\lambda^{\frac d2\delta(p,q)-1},
\end{equation}
which holds for $\ppq\in\mathcal R_3$ and discuss how the estimate \eqref{lplq} is related to its counterpart of the Laplacian, that is to say, the estimate for $\proj$. In order to put our discussion in a proper context,  let us consider the following estimate: 
\begin{align}\label{scaled-euclid} \big\|\proj   \big\|_{p\to q}\lesssim k^{-1+\frac{d}{2}\dpq} 
\end{align}
for $k\ge 1$ where $\proj$ is defined by \eqref{proj-laplace}. The following shows that  the estimate \eqref{scaled-euclid}  is equivalent to \eqref{est:euclid}.
\begin{lem} 
\label{equivalence}
The  estimate  \eqref{scaled-euclid} holds for all $k\ge 1$  if and only if the estimate  \eqref{est:euclid} for $R^*\! R$ holds.  
The equivalence remains valid with  $L^p$, $L^q$ spaces replaced, respectively,  by the Lorentz spaces  $L^{p,r}$, $L^{q,s}$ if $1<p,q\le \infty$ and $1\le r, s\le \infty$.  
\end{lem}

Indeed,  let us set
  \[ \wproj  f=\frac k{(2\pi)^d}\int_{\sqrt{1-1/k} \,\le |\xi| <1}e^{ix\cdot\xi}\,\widehat f(\xi) d\xi.\]  
 Since $\| \proj\|_{p\to q} =k^{\frac d2 \dpq-1}\| \wproj \|_{p\to q}$ by scaling,  the estimate  \eqref{scaled-euclid}  is  equivalent to 
 \begin{align}\label{est:saturated} 
 \big\| \wproj \big\|_{p\to q} \lesssim   1. 
 \end{align}
  Letting $k\to \infty$, we get  \eqref{est:euclid}. Conversely, making use of the spherical coordinates\footnote{ We write $\int_{\sqrt{1-1/k}\le |\xi|<1}\,e^{ix\cdot\xi}\,\widehat f(\xi) d\xi =C_d\int^{1}_{\sqrt{1-1/k}}\int e^{ix\cdot r\omega}\,\widehat f(r\omega) d\omega\, r^{d-1} dr$.} and Minkowski's inequality it is easy to see that \eqref{est:euclid} 
 implies \eqref{est:saturated} and then \eqref{scaled-euclid} via scaling.  Extension to the Lorentz spaces is clear since $L^{p,r}$ is a Banach space if $1<p \le \infty$ and $1\le r\le \infty$.

 The operator $R^*\!R$ is imbedded in a family of operators which are called the Bochner-Riesz operators of negative order defined by 
\[ S^{\alpha}\!f(x)=(2\pi)^{-d}\int e^{ix\cdot\xi}\,\frac{(1-|\xi|^2)_+^\alpha}{\Gamma(\alpha+1)}\widehat f(\xi) d\xi, \quad \alpha>-1, \] 
where $\Gamma$ is the gamma function. 
For $\alpha\le -1$ the definition is extended by analytic continuation. In fact, $S^{-1}=R^*\!R$. The $L^p$--$L^q$ boundedness  problem of $S^\alpha$ was studied by some authors, B\"orjeson \cite{B86},  Carbery and Soria \cite{CS88}, Sogge \cite{sogge-1986}, Bak,  McMichael, and  Oberlin \cite{bmo},   Bak \cite{bak}, Guti\'{e}rrez \cite{Gu99}, Cho, Kim, Lee, and Shim \cite{CKLS}.  
The complete set of the necessary  conditions for $L^p$--$L^q$ boundedness was shown by B\"orjeson \cite{B86}.  The problem  was settled by  Bak  \cite{bak} for $d= 2$  but  it remains open  for $d\ge 3$.  For the most recent development see Kwon and Lee \cite{KL19}.  In particular, concerning  $R^\ast\! R$ we have a complete characterization of $L^p$--$L^q$ boundedness.

\begin{thm}[\cite{B86, bmo, bak,  Gu99}]  \label{rr*} The operator $R^\ast\! R$  is bounded from $L^p$  to $L^q$ if and only if $\ppq\in \mathcal R_3$. 
Furthermore,  we have  $\|R^\ast\! R f\|_{q,\infty}\lesssim \|f\|_p$ if $\ppq\in (\mathfrak C, \mathfrak D]$,\footnote{By duality we also have the estimate  $\|R^\ast\! R f\|_{q}\lesssim \|f\|_{p,1}$ if $\ppq\in (\mathfrak C', \mathfrak D']$.} 
 and 
 $\|R^\ast\! R f\|_{q,\infty}\lesssim \|f\|_{p,1}$  if $\ppq\in \mathfrak C, \mathfrak C'$. 
\end{thm}

 Using  Theorem \ref{rr*},  Lemma \ref{equivalence}, and the Stein-Tomas theorem, we  can obtain the sharp  $L^p$--$L^q$ estimate for $\proj$.
 In fact, we prove Corollary \ref{laplacian} below.

 Combining Theorem \ref{rr*}  and Lemma \ref{equivalence} gives the estimate \eqref{scaled-euclid} with  $\ppq\in \mathcal R_3$ including the weak and restricted type estimates 
 for $\ppq\in (\mathfrak C, \mathfrak D]\cup (\mathfrak C', \mathfrak D']$. 
 On the other hand, we have $\|\proj\|_{2\to q}\le  k^{\frac{d-2}4-\frac d{2q}}$  for $q={2(d+1)}/(d-1)$, which can  be shown  similarly as before,  using  the spherical coordinates, the Stein-Tomas theorem, and Plancherel's theorem.  The estimate 
 $\|\proj\|_{2\to \infty}\lesssim   k^{\frac{d-2}4}$ follows from the Cauchy-Schwarz  inequality and Plancherel's theorem. These two estimates respectively correspond to the points   
$\mathfrak A'$ and $(1/2,0)$ in  Figure \ref{local-type} and then  duality gives the estimates for $\ppq=\mathfrak A$, and $(1,1/2)$. 
Together with $\|\proj  \|_{2\to 2}\lesssim  1$, interpolation between all those estimates and the estimate \eqref{scaled-euclid}  with $\ppq\in \mathcal R_3$ yields  the bound 
\[\|\proj\|_{p\to q} \lesssim  k^{\beta(p,q)}\] for $\ppq\in [1/2,1]\times [0,1/2]\setminus([\mathfrak C,\mathfrak  D]\cup[\mathfrak C', \mathfrak  D'])$ where $\beta(p,q)$ is given by 
\eqref{def-beta}. The opposite inequality can be shown without difficulty. In fact, from duality and scaling   we need only to show 
  \[\|\wproj\|_{p\to q} \gtrsim 
  \max( k^{1-\frac{d+1}2 \dpq}, 1, k^{  \frac{d}{q} -\frac{d-1}{2}}).\] 
  The second lower bound is clear. Since the multiplier of the operator $\proj$ is radial and supported in $O(k^{-1})$-neighborhood of the sphere $\mathbb S^{d-1}$,  the first and the third lower bounds  can be shown  by using, respectively,  the Knapp type example and the asymptotic expansion of the Bessel function (for example, see \cite{B86}).  Therefore, we obtain

\begin{cor}  
\label{laplacian} 
Let $d\ge 2$ and $(1/p,1/q)\in[1/2,1]\times[0,1/2]$. If  $(1/p,1/q)\not\in [\mathfrak C,\mathfrak  D]\cup[\mathfrak C', \mathfrak  D']$,  we have 
 \[ \|\proj \|_{p\to q}\sim k^{\beta(p,q)}. \]
Additionally, $(i)$   we have  
$  \|\proj\|_{L^p\to L^{q,\infty}}\le C k^{\beta(p,q)} $ for $(1/p,1/q)\in (\mathfrak C, \mathfrak D]$, 
and  $(ii)$ we  have $
    \|\proj\|_{ L^{p,1}\to L^{q,\infty}} \le C  k^{\beta(p,q)}
   $  if $(1/p,1/q)=\mathfrak C$ or $\mathfrak C'$.
   \end{cor}

 In what follows we show that the estimate \eqref{lplq} implies  \eqref{est:euclid}. 
\begin{lem} \label{implication} Let $B$ be a ball with small radius centered at the origin.  
 Suppose the estimate 
\Be\label{lplq-weaker}
\|\chi_{B}\Pi_{\lambda}\chi_{B}\|_{p\to q}\lesssim\lambda^{\frac d2\delta(p,q)-1}
\Ee 
holds. Then we have the estimate \eqref{est:euclid}.\end{lem}

Combining this with Theorem \ref{rr*}, we see that \eqref{lplq} holds if and only if $(1/p,1/q)$ lies in $\mathcal R_3$.  Transplantation of $L^p$ bounds for differential operators was shown by Kenig, Stanton, and Tomas \cite{KST}. Our argument below is similar to the one  in \cite{KST} since it also relies on a scaling argument. Unlike $L^p$ bound,    $L^p$--$L^q$ estimate is not scaling invariant, so the argument is not so simple as  in \cite{KST}. The particular form of bound in \eqref{lplq}  becomes crucial.   Furthermore,  our argument extends to general second order elliptic operators without difficulty as long as the associated spectral projection operator satisfies the same form of  bound as in \eqref{lplq}. 

In order to prove Lemma \ref{implication}, we recall the following which is a special case of the celebrated theorem due to H\"ormander \cite[Theorem 5.1]{H68}. 

\begin{thm}\label{thm:hor}Let $P$ be a self-adjoint elliptic differential operator of order 2 with $C^{\infty}$-coefficients
on $\mathbb R^d$  and $p$ be its principal part. Then, for $x,y$ in a compact subset and sufficiently close to each other, we have \begin{align*}
 \Big| e(x,y,\lambda)-(2\pi)^{-d}\int_{p(y,\xi)<\lambda}e^{i\psi(x,y,\xi)}d\xi\Big|\le C(1+|\lambda|)^{\frac{d-1}{2}},
 \end{align*}
 with $C$ independent of $\lambda$ 
 where $e(x,y,\lambda)$ is the spectral function of $P$, i.e., the kernel of the spectral projection operator $\Pi_{[0,\lambda]}$\footnote{Here, the operator $\Pi_{[0,\lambda]}$ is defined by  the typical spectral resolution.}   and $\psi$ is the function which is homogeneous in $\xi$ of degree 1 and satisfies $ p(x,\nabla_x\psi)=p(y,\xi)$ and 
 \begin{align}
 \label{spectral-phase} \psi(x,y,\xi)=\left\langle x-y,\xi\right\rangle+O(|x-y|^2|\xi|).
 \end{align}
 \end{thm}

\begin{proof}[Proof of Lemma \ref{implication}]
 Let $k, \nu$ be large positive integers. Then, we consider an auxiliary  projection operator $ \widetilde \Pi$ given by 
\[ \widetilde \Pi=\sum_{k\nu< \lambda \le(k+1)\nu}\Pi_\lambda .\] 
 By the triangle inequality and the assumption \eqref{lplq-weaker} we have 
\[
\| \chi_{ B} \widetilde \Pi \chi_{B}\|_{p\rightarrow q}\le\sum_{k\nu\le\lambda\le(k+1)\nu}\|\chi_{ B}\Pi_\lambda\chi_{B}\|_{p\to q}\lesssim k^{-1}(k\nu)^{\frac{d}{2}(\frac1p-\frac1q)}.
\] 
Let $f,g$ be nontrivial functions in $C^{\infty}_c(\mathbb{R}^d)$ such that $\supp f$, $\supp g$ are contained in $B$.
Since $\inp{\widetilde \Pi f}{g}=\inp{\chi_{B} \widetilde \Pi\chi_{ B}f}{g}$, we have
\begin{align*}
 \Big|\iint \widetilde \Pi(x,y)f(x)g(y)dxdy\Big|\lesssim k^{-1}(k\nu)^{\frac{d}{2}(\frac1p-\frac1q)}\|f\|_{p}\|g\|_{q'}
 \end{align*}
  Rescaling $(x,y)\to (\nu^{-\frac12}x,\nu^{-\frac12}y)$ gives the equivalent estimate 
 \begin{align}\label{est:trans} 
 \Big|\nu^{-\frac{d}{2}}\iint \widetilde \Pi(\nu^{-\frac12}x,\nu^{-\frac12}y)F(x) G(y) dxdy\Big|\lesssim k^{-1+\frac{d}{2}(\frac1p-\frac1q)} \|F\|_p \|G\|_{q'}
 \end{align}
 provided that $F$ and $G$ are supported in $\sqrt\nu B$.  Taking the radius of $B$ small enough, we may apply Theorem \ref{thm:hor} because 
 $\nu^{-\frac12}x$ and $\nu^{-\frac12}y$ are close enough.  Since $\widetilde \Pi(\nu^{-\frac12}x,\nu^{-\frac12}y) 
       = e(\nu^{-\frac12}x,\nu^{-\frac12}y,(k+1)\nu)-e(\nu^{-\frac12}x,\nu^{-\frac12}y, k\nu)$, 
 making use of Theorem \ref{thm:hor} and changing variables $\xi\to \nu^{\frac12}\xi$, we observe that 
 \begin{align*} 
 \widetilde \Pi(\nu^{-\frac12}x,\nu^{-\frac12}y) 
       &=(2\pi)^{-d}\int_{k\nu\le |\xi|^2<\nu(k+1)}e^{i\psi(\nu^{-\frac12}x,\nu^{-\frac12}y,\xi)}d\xi+R_{\nu,k}(x,y) \\ 
       &=(2\pi)^{-d}\nu^\frac d2\int_{k\le  |\xi|^2 <(k+1)}e^{i\psi(\nu^{-\frac12}(x,y),\nu^\frac12 \xi)}d\xi+R_{\nu,k}(x,y),
 \end{align*}
 where $R_{\nu,k}(x,y)=O(|k\nu|^{\frac{d-1}{2}})$. Combining this with \eqref{est:trans}  yields 
 \begin{align*} 
 \Big|\iint\Big(\int_{k\le  |\xi|^2<k+1}\hspace{-25pt} e^{i\psi(\nu^{-\frac12}(x,y),\nu^{\frac12}\xi)}d\xi+\widetilde R_{\nu,k}(x,y)\Big)f(x)g(y)dxdy\Big|\lesssim k^{-1+\frac{d}{2}(\frac1p-\frac1q)}
 \end{align*}
 provided that $f$ and $g$ are supported in $\sqrt\nu B$ and $\|f\|_p=\|g\|_{q'}=1$. 
 Here $\widetilde R_{\nu,k}(x,y)$ satisfies $|\widetilde R_{\nu,k}|\lesssim \nu^{-\frac{d}{2}}(k\nu)^{\frac{d-1}{2}}$. From \eqref{spectral-phase}, we note that the phase function $\psi(\nu^{-\frac12}(x,y),\nu^{\frac12}\xi)\to \langle x-y,\xi\rangle$ as $\nu\to \infty$. Thus, taking $\nu\rightarrow\infty$, we obtain
 \begin{align*} 
 \Big|\iint\Big(\int_{k\le  |\xi|^2<k+1}e^{i\left\langle x-y,\xi\right\rangle}d\xi\Big)f(x)g(y)dxdy\Big|\lesssim k^{-1+\frac{d}{2}(\frac1p-\frac1q)}
 \end{align*}
 if $\|f\|_p=\|g\|_{q'}=1$.   This yields  the estimate \eqref{scaled-euclid}
 which is equivalent to \eqref{est:euclid} as seen in the above (e.g. (\ref{est:saturated})). 
\end{proof}

\subsection{Estimates away from $\sqrt \lambda \mathbb S^{d-1}$}

Let us set $B_r=B(0,r)$, $0<r<1$. 
In this subsection we are concerned with the estimate over the set  $B_r \times B_r$ 
which strengthens the result in Theorem \ref{thm-locest} for $p,q$ satisfying $1/p+1/q=1$.  In Section \ref{sec:endpoint} and Section \ref{off-d} we will deal with the case in which $x,y$ are close to $\sqrt \lambda \mathbb S^{d-1}$.   
The condition \eqref{det-con} is no longer satisfied with $B_r$, $r>1/\sqrt 2$. Appearance of points $x,y$ for which $\cD(x,y)=0$ gives a rise to significant change in the boundedness property 
of $\Pi_\lambda$ since the favorable  $O(\lambda^{-\frac12})$ bound on $L^1-L^\infty$  estimate is not available any more. In what follows we work with the scaled operator 
$\fP_\lambda$ instead of $\Pi_\lambda$ (\eqref{eq:scaled-op}) and the estimate for $\Pi_\lambda$ can be deduced by \eqref{eq:norm-scaled}.

\begin{prop}
\label{away-sphere}
Let  $\delta_\circ$ be a small positive constant, $r=1-\delc$.  
Then, there are  constants $C=C(\delc)$ and  $c=c(\delc)$\footnote{In fact, $c(\delc)\sim \sqrt{\delc}$.} such that the following hold\,$:$ 
\vspace{-3pt}
\begin{itemize}
[leftmargin=0.9cm, labelsep=0.3 cm, topsep=0pt]

\item [$(i)\,$]
If $2^{-j}\le c$,  the estimate 
\begin{equation}
\label{eq:away-sphere00}
\| \chi_{B_r} \fP_\lambda [ \psi_j^\kappa ] \chi_{B_r}\|_{p\to p'}\le  C\lambda^{\scaleppp} 2^{(\frac{d+1}{2}\dpp-1)j }, \quad \kappa=\pm, \pm\pi, 
\end{equation}
holds for $1\le p\le 2$. 

\item [$(ii)$] If $2^{-j}\ge c$,  we have the same estimates in \eqref{eq:away-sphere00} for $6/5 < p\le 2$.  
Furthermore, at the endpoint  $p=6/5$ we have restricted weak type estimate with the same bound.  
\end{itemize}
  \end{prop}
  
  Making use of the estimate \eqref{eq:away-sphere00} and Lemma \ref{tt-st}, we can  obtain  off-diagonal estimates for  $\chi_{B_r} \fP_\lambda \chi_{B_r}$. It seems natural to expect that the same  sharp bound as  \eqref{est-loc} continues to hold for $\|\chi_{B_r} \fP_\lambda \chi_{B_r}\|_{p\to q}$ (via \eqref{eq:norm-scaled}). However, due to the degeneracy of the phase function 
 the range of the sharp bounds differs from that of Theorem \ref{thm-locest}.  We summarize such  estimates in  Corollary \ref{local-r}.

 We dedicate the rest of this section to proving Proposition \ref{away-sphere}.

 \subsection{2nd-order derivative of $\cP$}  
As $\cD(x,y)$ gets close to zero
the oscillatory integral $I_j$ does not lend an estimate which  is good enough for our purpose. 
Unlike the case in which $|\cD(x,y)|$ is bounded away from zero,  we need to use the second order derivative $\partial_s^2\cP$ to get the correct  order of decay  in $\lambda$. 
However, as is to be seen later, the second order derivative alone does not simply give a favorable lower bound since $\partial_s^2\cP$ also vanishes.
A typical strategy to get around this kind of difficulty  might be to combine  the  lower bounds  from the first order and second order derivatives but it is not viable since the zeros of $\partial_s \cP$ and $\partial_s^2 \cP$ merge to a single point as $\cD(x,y)\to 0$ (see \eqref{qcos}).  This leads us to break  the integral dyadically  away from the zero of $\partial_s^2\cP$. 
To this end,  we need to take a look at the zeros of  $\partial_s^2 \mathcal P$.

A computation shows
\begin{equation}
\label{ph-dd}
 \partial_s^2 \mathcal P(x,y,s) =-\frac{\inp xy\mathcal R(x,y,\cos s) }{ \sin^3 s},
\end{equation}
where  
\begin{equation}
\label{rxy}
\mathcal R:=\mathcal R(x,y,\tau) :=  \tau^2-{\inp xy}^{-1}(|x|^2+|y|^2)\tau+ 1.
\end{equation}
As before, we need to identify the zeros of $\partial_s^2 \mathcal P(x,y,s)$. The polynomial $\mathcal R$ has two distinct  roots $\tau^\pm(x,y)$ as long as $x-y\neq 0$ and $x+y\neq 0$: 
\Be
\label{roots}  \tau^\pm(x,y) =\frac{|x|^2+|y|^2\pm |x+y||x-y|}{2\inp xy}=\frac{2\inp xy}{|x|^2+|y|^2\mp |x-y||x+y|}.
\Ee
Since $|\tau^+(x,y)|>1$ if $|x+y||x-y|\neq 0$, the root $\tau^-(x,y)$ is more relevant for our purpose. 
We define $S_c$ on $B_2\times  B_2\setminus ( \{x =y\}\cup \{x=-y\}) $  by setting 
\begin{equation}
\label{sc-}\cos S_c(x,y):= \tau^-(x,y).
\end{equation}

\begin{lem}
\label{scs}   $S_c$ defines a smooth function from $B_2\times  B_2\setminus ( \{x =y\}\cup \{x=-y\})$ to $(0,  \pi)$. We have $S_c(x,y)\sim  |x-y|^\frac12$ and
\Be \label{dist-tau} 1-\tau^-\sim  |x-y|, \quad 1+\tau^-\sim |x+y|  \Ee if  $|x+y|+|x-y|\ge c$ for some $c>0$. 
 \end{lem}

It is clear that $\pi-S_c(x,y)=S_c(x,-y)$ because $\cos(\pi-s)=-\cos s$.  Hence, $\pi-S_c(x,y)\sim |x+y|^\frac12$ provided that 
$|x+y|+|x-y|\ge c$ for some $c>0$. 

\begin{proof}[Proof of Lemma \ref{scs}]  From \eqref{roots} it is easy to see that  $\tau^-(x,y)\in (-1,1)$, and we also note that the function $\tau^-$ is smooth on   $B_2\times  B_2\setminus ( \{x =y\}\cup \{x=-y\})$. 
In fact, this is obvious when $|\inp xy|>0$ but, if $|\inp xy|$ is small, we can write $\tau^-(x,y)=2\inp xy/(|x|^2+|y|^2+|x-y||x+y|)$ making use of 
the identity 
\Be 
\label{t-id}
|x+y|^2|x-y|^2=(|x|^2+|y|^2)^2-4\inp xy^2. 
\Ee 
  Since $\cos^{-1}$ is smooth on the interval $(-1, 1)$, it follows that $S_c$ is smooth on   $B_2\times  B_2\setminus ( \{x =y\}\cup \{x=-y\})$.
Since  $1-\cos S_c(x,y)\sim S_c^2(x,y)$, to see $|x-y|^\frac12\sim  S_c(x,y)$  it suffices to observe that 
\begin{equation}
\label{sc}
\begin{aligned}
1-\cos S_c(x,y)&=\frac{(|x+y|-|x-y|)|x- y|}{2\inp xy} 
\\&= \frac{2|x- y|}{|x+y|+|x-y|}
\end{aligned}
\end{equation} 
because    $|x+y|^2-|x-y|^2=4\inp xy$.   Finally, to get \eqref{dist-tau} 
it is enough to note  $1\pm\tau^-=2|x\pm y|/(|x+y|+|x-y|)$, which follows from a simple manipulation.
\end{proof}

The next lemma shows how  the zeros of $\partial_s \cP$ and $\partial_s^2 \cP$ are related to each other in terms of $\cD(x,y)$. 
\begin{lem}
The following identity holds\,$: $
\begin{equation}  
\label{qcos}\mathcal Q(\cos S_c(x,y))=-\,\frac{2|x-y||x+y|}{|x-y||x+y|+ |x|^2+|y|^2-2\inp xy^2}  \cD(x,y).
\end{equation}
\end{lem}

\begin{proof}
Indeed, using \eqref{sc-} and \eqref{q-def},  we note that $\mathcal Q(\cos S_c(x,y))$ equals 
\begin{align*}
4^{-1} \inp xy^{-2} \Big ((|x|^2+|y|^2-|x-y||x+y|)^2-4 \inp xy^2+4\inp xy^2|x-y||x+y|\Big) .
\end{align*}
Since  \[(|x|^2+|y|^2-|x-y||x+y|)^2-4 \inp xy^2=- |x-y||x+y|(|x+y|-|x-y|)^2,\] we have 
\begin{align*}
\mathcal Q(\cos S_c(x,y))=\frac{|x-y||x+y|}{2\inp xy^{2}}\Big (  |x+y||x-y|  -|x|^2-|y|^2+2\inp xy^2 \Big).
\end{align*}
We also note  
$ |x+y|^2|x-y|^2  -(|x|^2+|y|^2-2\inp xy^2)^2=-4\inp xy^2\cD(x,y)$
using \eqref{t-id}. 
 Hence, we get \eqref{qcos}.
\end{proof}

\subsection{Proof of Proposition \ref{away-sphere}} 
\label{proof-away}
In order to show Proposition \ref{away-sphere}  we  start with partitioning the set $B_r\times B_r$ into finitely many sets $\{A_n\times A_m\}_{n,m}$  of small diameter less than $\vepc$ so that 
$\cD(A_n\times A_m)$ is contained in an interval of length $\vepc$. The constant $\varepsilon_\circ$ is to be taken small enough later.  We need to show
\[ \|\chi_{ A_n} \fP_\lambda [ \psi_j^\kappa ] \chi_{A_m}\|_{p\to q}\le  C\lambda^{\scalep} 2^{(\frac{d+1}{2}\dpq-1)j }\]
for each $n, m$.

If $|\cD(x,y)|\gtrsim \vepc $ for  $(x,y)\in A_n\times A_m$,  
from Corollary \ref{pj-det} and \eqref{symmetric}  we get  the desired estimate  \eqref{eq:away-sphere00} for $1\le p\le 2$ with $B_r$ replaced by  $A_n$ and $A_m$.   Thus, it is sufficient to show 
\Be 
\label{est-aa}
\| \chi_{A} \fP_\lambda [ \psi_j^\kappa ] \chi_{A'}\|_{p\to p'}\le  C\lambda^{\scaleppp} 2^{(\frac{d+1}{2}\dpp-1)j }
\Ee
 under the assumption that $A,$ $A'\subset B_r$ and 
 \begin{equation}
\label{det-small0}
 |\cD(x,y)|\ll  \vepc,   
  \quad (x,y)\in A\times A'.
\end{equation}  
This gives $|x|^2+|y|^2=1+\inp xy^2+O(\vepc)$. 
Using   \eqref{t-id},  we also have
\Be 
\label{eq:string}
\begin{aligned}
|x-y|^2|x+y|^2
 &= (1-\inp xy^2)^2+O(\vepc)
 \\
 &=(2-|x|^2-|y|^2)^2+O(\vepc)
 \end{aligned}
\Ee
for $(x,y)\in A\times A'$.
Note that $|x|^2+ |y|^2\le 2-3\delc$ because $(x,y)\in B_r\times B_r$. Taking  $\vepc$  small enough, we have 
\begin{align}
\label{eq:lowerxy}
|x-y||x+y|&\gtrsim \delc,  \quad (x,y)\in A\times A'.
\end{align} 
We now prove Proposition \ref{away-sphere} by considering the cases $2^{-j}\ll  \sqrt\delc$ and 
$2^{-j}\gtrsim  \sqrt\delc$,  separately. 
The latter  case is more involved and we hand it later. 

 \subsection*{{Case $2^{-j}\ll \sqrt\delc$ and  $\kappa=\pm \pi, \pm$} } Since  the conditions  \eqref{det-small0} and \eqref{eq:lowerxy} are invariant under $y\to -y$,  by \eqref{symmetric}   
 it is sufficient to show \eqref{est-aa} with $\kappa=+$. 
From \eqref{roots}, Lemma \ref{scs}, and \eqref{eq:lowerxy} we note that  $S_c(x,y)\in [c \sqrt\delc,  \pi-c \sqrt\delc]$ for some $c>0$ if $(x,y)\in A\times A'$. 
From this we observe  $|\mathcal R(x,y,\cos s)|\ge c$ for some constant $c>0$ on the support of $\psi^+_j$, so  it follows from \eqref{ph-dd}  that
$ |\partial_s^2 \mathcal P| \gtrsim  2^{3j}$ if $s\in \supp \psi_j^+$.\footnote{If  $|\langle x,y\rangle|\ge c$ for some small  $c>0$ it is clear  because $S_c(x,y)\not\in\supp \psi_j^+$. Otherwise, i.e., if $|\langle x,y\rangle|\ll 1$,  then  $|x|^2+|y|^2\sim 1$ because   $|\mathcal D|\ll \varepsilon_\circ$. Thus, $\partial_s^2 \mathcal P=\frac{(|x|^2+|y|^2)\cos s+O(\varepsilon_\circ+c)}{\sin ^3s}\sim 2^{3j}$ on $\supp \psi_j^+$ with small enough $\vepc, c>0$.} This and Lemma \ref{osc-lemma} give 
\[
  \Big|\int   e^{i\lambda \mathcal P(x,y,s)} \mathfrak a(s) \psi_j^+(s)    ds\Big|\lesssim  \lambda^{-\frac12}  2^{\frac{d-3}2 j} 
  \]
 for $(x,y)\in A\times A'$. 
Thus, we have  $\| \chi_A \fP_\lambda [ \psi_j^+] \chi_{A'}\|_{1\to \infty}\lesssim    \lambda^{-\frac12} 2^{\frac{d- 3}2 j}$. Interpolating  
this and  the $L^2$ estimate  $\| \chi_A \fP_\lambda [ \psi_j^+ ] \chi_{A'}\|_{2\to 2}\lesssim   \lambda^{-\frac d2}2^{-j}$ which follows from \eqref{easy-l2} and scaling, 
we get 
\[ \| \chi_{ A} \fP_\lambda [ \psi_j^+ ] \chi_{A'}\|_{p\to q}\lesssim  C\lambda^{\scalep} 2^{(\frac{d-1}{2}\dpq-1)j },\]
which is, in fact, better than the desired estimate \eqref{est-aa}.

We now turn to the remaining case of  $\psi_j^\kappa$ where  $2^{-j}\gtrsim  \sqrt\delc$ and $\kappa=\pm$, $\pm \pi$, $0$. 

{\bf Case  $2^{-j}\gtrsim \sqrt\delc$ and $\kappa=\pm, \pm\pi, 0$.} Again by \eqref{symmetric}  it is sufficient to show the estimate \eqref{est-aa} with $\kappa=+, 0$.  
 $\supp \psi_j^+$ and $\supp \psi^0$ may  contain the critical point  $S_c$, so this makes it difficult  to get the correct  $\lambda^{\scalep}$ decay.
Since $2^{-j}\gtrsim c$,  the estimates for $\chi_{ A} \fP_\lambda [ \psi_j^+ ] \chi_{A'}$ and $\chi_{ A} \fP_\lambda [ \psi^0 ] \chi_{A'}$ can be handled by the same argument. Therefore, we 
only prove the estimate \eqref{est-aa} with $\psi_j^+$, $2^{-j}\gtrsim \sqrt\delc$.

Decomposing the set $A\times A'$ into a finite number of disjoint sets, 
we may assume 
 that  
\begin{equation}
\label{eps-con0}
|S_c(x,y)- S_c(x',y')|\le  \vepc 
\end{equation}
whenever $(x,y)$ and $(x',y')$ are contained in $A\times A'.$ This is clearly possible since  $S_c$ has bounded derivatives because of  \eqref{eq:lowerxy}.

Let  $(x_\ast, y_\ast)\in A\times A'$ and set  
\[\psi_{\ast}(s):= \psi_\circ \Big(\frac{s-S_c(x_\ast, y_\ast)}{2\vepc}\Big)
\]  
 where $\psi_\circ\in C_c^\infty(-2,2)$ such that $\psi_\circ(s)=1$ for $s\in [-1,1]$.  
Using this, we break 
\[\chi_{ A} \fP_\lambda [ \psi_j^+ ] \chi_{A'}=  \chi_{ A} \fP_\lambda [ \psi_\ast \psi_j^+  ] \chi_{A'}+\chi_{ A} \fP_\lambda [ (1-\psi_\ast)\psi_j^+  ] \chi_{A'}.\]
The second operator is easy to deal with.  
We note from \eqref{eps-con0} that   $|\partial_s^2 \cP|\gtrsim\vepc$ on the support of  $(1-\psi_\ast)\psi_j^+$ because $|s-S_c(x,y)|\gtrsim \vepc$ if $s\in \supp (1-\psi_\ast)\psi_j^+$. Hence,  by Lemma \ref{osc-lemma} we get  $\|\chi_{ A} \fP_\lambda [ (1-\psi_\ast)\psi_j^+] \chi_{A'}\|_{1\to \infty} \lesssim \lambda^{-\frac12}$.  We also have 
$\|\chi_{ A} \fP_\lambda [ (1-\psi_\ast)\psi_j^+  ] \chi_{A'}\|_{2\to 2} \lesssim \lambda^{-\frac d2}$ by \eqref{easy-l2}.  
By interpolation between these estimates we get the estimate $\|\chi_{ A} \fP_\lambda [ (1-\psi_\ast)\psi_j^+  ] \chi_{A'}\|_{p\to p'}\lesssim \lambda^{\scaleppp}$.  Therefore, we  need only  to show 
\[ \| \chi_{ A} \fP_\lambda [ \psi_\ast \psi_j^+  ] \chi_{A'}\|_{p\to p'}\lesssim  C\lambda^{\scaleppp} \] 
for $6/5<p\le 2$. 
Further localization to an interval of length $\sim\vepc$ is important in proving the $L^2$ estimate \eqref{ppjl2} below (see Lemma \ref{nonzero-det}).

 We begin with making additional decomposition away from $S_c$ by inserting $\widetilde\psi(2^l (\cdot-S_c))$: 
 \[   \chi_A \mathfrak P_\lambda[\psi_\ast\psi_j^+]\chi_{A'}=\sum_l  \chi_A \mathfrak P_\lambda[\psi_\ast\psi_j^+\widetilde \psi(2^l(\cdot-S_c))]\chi_{A'}. \] 
 Since $\psi_\ast$ is supported in an interval of length $\lesssim \vepc$, 
we may  clearly assume $2^{-l}\lesssim \vepc$ since  $\mathfrak P_\lambda[\psi_\ast\psi_j^+\widetilde \psi(2^l(\cdot-S_c))]=0$ otherwise by the support property of $\widetilde \psi$.  Note that the kernel of $\mathfrak P_\lambda[\psi_\ast\psi_j^+\widetilde \psi(2^l(\cdot-S_c))]$ is given by 
\[  I_{j,l}=\int   e^{i\lambda \mathcal P(x,y,s)} \mathfrak a(s)\psi_\ast(s)\psi_j^+(s)   \widetilde\psi(2^l (s-S_c(x,y))) ds.\] 
Since $2^{-j}\gtrsim c$, we have $\sin s\sim 1$. Using \eqref{ph-dd}  we see 
$  |\partial_s^2 \mathcal P| \sim    2^{-l}$  on the support $\widetilde \psi(2^l (\cdot-S_c))\psi_j^+$.  
Thus, Lemma \ref{osc-lemma} gives
$ |I_{j,l}|\lesssim \lambda^{-\frac12} 2^{\frac12l}$, from which we obtain 
\begin{equation}
\label{ppjl}
 \|   \chi_A \mathfrak P_\lambda[\psi_\ast\psi_j^+\widetilde \psi(2^l(\cdot-S_c))]\chi_{A'} \|_{1\to \infty} \lesssim   \lambda^{-\frac12} 2^{\frac l2} .
 \end{equation}
On the other hand, we claim 
\begin{equation}
\label{ppjl2}
 \|  \chi_A \mathfrak P_\lambda[\psi_\ast\psi_j^+\widetilde  \psi(2^l(\cdot-S_c))]\chi_{A'} \|_{2\to 2} \lesssim    \lambda^{-\frac d2}2^{-l}.
 \end{equation}
 Assuming this for the moment, we prove  \eqref{est-aa}.
Interpolation between \eqref{ppjl} and \eqref{ppjl2} gives 
\[ \|  \chi_A \mathfrak P_\lambda[\psi_\ast\psi_j^+ \widetilde\psi(2^l(\cdot-S_c))]\chi_{A'}  f \|_{p'} \lesssim  \lambda^{\frac{d-1}2\dpp-\frac d2} 
             2^{(\frac32\dpp-1)l} \|f\|_p\]
Then, by summation  over $l$  we get the desired estimate \eqref{est-aa} for $2^{-j}\gtrsim c$  provided that $\frac1p-\frac1{p'}< \frac23$. If $\frac1p-\frac1{p'}= \frac23$, by $(c)$ in Lemma \ref{s-trick} we get the  restricted weak type bound.

To complete the proof of Proposition \ref{away-sphere} it remains to show \eqref{ppjl2}. 
Our proof of \eqref{ppjl2} is based on the $L^2$ estimate for oscillatory integral operator which is a generalization of Plancherel's theorem. 
For $a\in C_c^\infty(\mathbb R^d\times \mathbb R^d)$ and $\phi$ which is smooth on the support of $a$,  
we  define 
\Be\label{def:osc} \mathcal O_\lambda[\phi, a]f:=\int e^{i\lambda  \phi(x,y) } a(x,y) f(y) dy.\Ee
To show \eqref{ppjl2}, we use   the following. 

\begin{lem}{\cite{hor} $($also see \cite[Theorem 2.1.1]{so93}$)$}  
\label{generalized}  
Suppose $ \det( \partial_x \partial_y^\intercal \phi )\neq 0$ on the support of $a$, then we have the estimate
$\|\mathcal O_\lambda[\phi, a]f\|_2\le C\lambda^{-\frac d2}\|f\|_2.$ 
\end{lem}

\subsection{Proof of \eqref{ppjl2}: $L^2$ estimate}  
\label{mild-l2} We can not use $L^2$ isometry of the propagator $e^{it H}$ (\eqref{easy-l2}) any more   because the  inserted cutoff function $\widetilde\psi(2^l (\cdot-S_c(x,y))$  disturbs orthogonality.  
Instead we use Lemma \ref{generalized}. 
For simplicity let us set 
\Be 
\label{sxyl}
\sxyl(x,y, s):=2^{-l}s+S_c(x,y).
\Ee
After change of variables $s\to  2^{-l}s+S_c(x,y)$ 
we see that the kernel of the operator 
$\chi_A \mathfrak P_\lambda[\psi_\ast\psi_j \widetilde \psi(2^l(\cdot-S_c))]\chi_{A'}$ is given by  
\[   2^{-l}  \chi_ {A}(x) \chi_{A'}(y) \int   e^{i\lambda \mathcal P(x,y,  \sxyls )}   (\mathfrak a\psi_\ast\psi_j)(\sxyls)\widetilde\psi(s)\,ds.\] 
Clearly we may replace $\chi_ {A}$ and $\chi_{A'}$ with smooth functions $\widetilde \chi_ {A}$ and $\widetilde\chi_{A'}$ which are adapted to the sets $A$ and $A'$. More precisely, 
$\widetilde \chi_ {A}=1$ on $A$,  $\widetilde\chi_{A'}=1$ on $A'$, and supports of $\widetilde \chi_ {A},\widetilde \chi_ {A'}$  are contained in $c\varepsilon_0$-neighborhoods of 
$A$, $A'$, respectively. 
Let us set 
\begin{align}
\label{phase-s} \Phi_s(x,y)&:=\mathcal P(x,y, \sxyls),  
\\
\nonumber 
 a_s(x,y)&:=\widetilde \chi_ {A}(x) \widetilde\chi_{A'}(y)  (\mathfrak a\psi_\ast\psi_j)(\sxyls) \widetilde \psi(s).
   \end{align}
 In order to prove \eqref{ppjl2} it is sufficient to show 
\begin{equation}
\label{tl2} 
\|\mathcal O_\lambda[\Phi_s, a_s] f\|_2\lesssim \lambda^{-\frac d2} \|f\|_2
\end{equation}   
uniformly in $s\in \supp \psi$.  Since 
the phase and amplitude functions are smooth and uniformly bounded in $C^\infty$, by Lemma \ref{generalized}   it is enough for \eqref{tl2} to show the following.

\begin{lem}
\label{nonzero-det}  Let $x, y\in B_{1-\delc}$ for some $\delta_\circ>0$. 
 Suppose  \eqref{det-small0}, \eqref{eq:lowerxy}, and  $2^{-l}\lesssim \vepc$  holds. Then, if $\vepc>0$ is small enough,   we have 
\Be
\label{partialxy0}
  \det(\partial_x \partial_y^{\intercal}   \Phi_s(x,y))\neq 0, \quad (x,y)\in A\times A'.
  \Ee
\end{lem}

In fact, for  \eqref{tl2} we need to verify  \eqref{partialxy0}   for $(x,y)\in \supp \widetilde \chi_A\times  \supp \widetilde  \chi_{A'}$   but  the assertion clearly remains valid by continuity if we take $\vepc$ small enough. The following observation significantly simplifies  the proof of Lemma \ref{nonzero-det}. 

\begin{lem}  
\label{id} Let $(x,y)\in B_2\times  B_2\setminus ( \{x =y\}\cup \{x=-y\})$.   Then, the following identities hold: 
\Be
\label{xyS}
\begin{aligned}
\sin S_c(x,y) \partial_x S_c(x,y)&=\cos S_c(x,y)\big(a(x,y)x-b(x,y)y\big),\\
 \sin S_c(x,y) \partial_y^\intercal S_c(x,y)&=\cos S_c(x,y)\big(a(x,y) y^\intercal-b(x,y)x^\intercal\big),
\end{aligned}
\Ee
where 
\[a(x,y)=\frac{2}{|x-y||x+y|}, \  \ \  b(x,y)=\frac{|x|^2+|y|^2}{\inp xy |x-y||x+y|}.\] 
\end{lem}

\begin{proof}
Differentiating both sides of the first equality in \eqref{sc} and rearranging the terms, we  get 
\begin{align*}
\sin S_c(x,y) & \partial_x  S_c(x,y)= \mathrm F(x,y)      x  + \mathrm G(x,y)  y,
\end{align*}
where 
\begin{align*}
\mathrm F(x,&y)=  \frac{1}{2\inp xy} \Big( \frac{|x-y|}{|x+y|}+\frac{|x+y|}{|x-y|}-2\Big), 
\\  
\mathrm G(x,y)= \frac1{2\inp xy^2}\Big[  \inp xy & \Big( \frac{|x-y|}{|x+y|}-\frac{|x+y|}{|x-y|}+2\Big)-\Big({|x+y|}|x-y|-{|x-y|}^2\Big) \Big]. 
\end{align*}
It is easy see that  
$\mathrm F(x,y)=\frac{|x|^2+|y|^2-|x+y||x-y|}{\inp xy |x+y||x-y|}$. Thus, by  \eqref{sc-} and  \eqref{roots}  we see 
\[ \mathrm F(x,y)=  \frac{2\cos S_c(x,y)}{|x-y||x+y|}.  \]
 Using \eqref{t-id},  a computation shows
 \begin{align*}
\mathrm G(x,y)=  \frac{(|x|^2+|y|^2)(|x-y||x+y|-|x|^2-|y|^2)}{{2\inp xy^2} |x+y||x-y|}
=-  \frac{\cos S_c(x,y)(|x|^2+|y|^2)}{\inp xy  |x-y||x+y| }.
\end{align*}
For the second equality we use 
  \eqref{sc-} and  \eqref{roots}.  We thus get  the first identity in \eqref{xyS}.  Since $S_c(y,x)=S_c(x,y)$, the second identity  in \eqref{xyS} follows from the first one by interchanging the roles of $x,y$.
\end{proof}

\begin{proof}[Proof of Lemma  \ref{nonzero-det}]
Differentiating  both side of  \eqref{phase-s}, we obtain 
\begin{equation}
\label{partialxy}
\begin{aligned}
\partial_x\partial_y^{\intercal} &  \Phi_s(x,y) =  \partial_x \partial_y^{\intercal} \mathcal P(x,y, \sxyl)
+
\partial_x S_c \partial_y^{\intercal}  \partial_s \mathcal P(x,y,  \sxyl) +
\\[3pt]
\quad \partial_x&\partial_s\mathcal P(x,y,  \sxyl) \partial_y^{\intercal}  S_c 
 +\partial_s^2\mathcal P(x,y,  \sxyl) 
\partial_x S_c  \partial_y^{\intercal}  S_c +  \partial_s\mathcal P(x,y,  \sxyl) 
 \partial_x\partial_y^{\intercal}  S_c.
\end{aligned}
\end{equation}
Let us set 
\begin{equation}
\label{hhh}
\begin{aligned}
\mathbf H = \partial_x \partial_y^{\intercal} & \mathcal P(x,y, S_c) +  
\partial_x S_c \partial_y^{\intercal}  \partial_s \mathcal P(x,y, S_c)  + \partial_x\partial_s\mathcal P(x,y,S_c) 
 \partial_y^{\intercal}  S_c.
 \end{aligned}
 \end{equation}

From Lemma \ref{scs} and  \eqref{eq:lowerxy}  we see that $S_c$ is a smooth function with
bounded derivatives on $A\times A'$.  $\sin s$ is bounded below since $2^{-j}\gtrsim  c$. So, from \eqref{ph-dd}  it follows that $\partial_s^2\mathcal P(x,y,  \sxyls)=O(2^{-l})=O(\vepc)$. 
By \eqref{ph-d},  \eqref{qcos} and \eqref{det-small0} we also see that $\partial_s\mathcal P(x,y,  \sxyls)=O(2^{-l}+\vepc)=O(\vepc)$.
Since $S_c$  is smooth on $A\times A'$,  using \eqref{eps-con0} and \eqref{partialxy}  we have 
\[
\begin{aligned}
\partial_x\partial_y^{\intercal}   \Phi_s(x,y) &= \mathbf H +O(\vepc).
\end{aligned}
\]

The matter is now reduced to showing $|\det\mathbf H|\ge c$ for some $c>0$. To this end,  we set 
\Be
\label{abcd}
\begin{aligned}
\gamma:= \inp xy,  \ \  \  \alpha:=\frac\gamma {1-\gamma^2}, \\
 \beta:= \frac{1+\gamma^2}{1-\gamma^2}, \ \  \ \delta:= \frac1 {1-\gamma^2}.
 \end{aligned}
\Ee
Since $1-\inp xy^2\gtrsim \delc$ by \eqref{eq:string}, we have $|x-y||x+y|=1-\gamma^2 + O(\vepc)$. 
A computation shows 
\[\cos S_c(x,y)-\inp xy= \frac{2\inp xy \mathcal D(x,y)}{|x|^2+|y|^2-2\inp xy^2+|x-y||x+y|} .\]
Using \eqref{det-small0} and \eqref{eq:lowerxy}, 
we have  $\cos S_c(x,y)-\inp xy=O(\vepc)$, that is to say,  $\cos S_c(x,y)=\gamma  +O(\vepc)$. From \eqref{eq:string} we also have 
$|x|^2+|y|^2=1+\gamma^2+O(\vepc)$. 
We note that $1-\gamma^2\gtrsim \delc$ because of \eqref{eq:string}.  Thus,  combining all those  in the above with  Lemma \ref{id}, we see 
\begin{align*}
\sin S_c(x,y) \partial_x S_c(x,y)&=2\alpha x-\beta y + O(\vepc),
\\
\sin S_c(x,y) \partial_y^\intercal S_c(x,y)&=2\alpha y^\intercal-\beta x^\intercal + O(\vepc).
\end{align*}
On the other hand,  by differentiating both side of \eqref{ph-d} we get
\begin{align*}
( \partial_x\partial_s\mathcal P)(x,y, S_c(x,y) )
 &=\frac{\cos S_c(x,y) y-x}{\sin^2 S_c(x,y)}=
\alpha y-\delta x+ O(\vepc),  
 \\   
 (\partial_y^\intercal \partial_s\mathcal P)(x,y, S_c(x,y))
 &=\frac{\cos S_c(x,y) x^\intercal-y^\intercal}{\sin^2 S_c(x,y)}=\alpha x^\intercal-\delta y^\intercal+ O(\vepc).
   \end{align*}
Since  $ \partial_x \partial_y^{\intercal} \mathcal P(x,y, S_c(x,y))=-\frac1{\sin S_c(x,y)} \,\mathbf I_d $,  putting together  the above identities in \eqref{hhh},   we have
\[
   \sin  S_c(x,y)  \mathbf H
    = -\mathbf I_d+   (2\alpha^2 +\beta\delta)(xx^\intercal+yy^\tr)  -2\alpha\beta yx^\tr  -4\alpha \delta  xy^\tr  +O(\vepc)
  \]
and, using  \eqref{abcd}, we get  
\[  
(1-\gamma^2)^2\sin  S_c(x,y)  \mathbf H=\mathbf M +O(\vepc),
\]
where 
\Be 
\label{sc-matrix}
\mathbf M= -(1-\gamma^2)^2\mathbf I_d+   (1+3\gamma^2)(xx^\intercal+yy^\tr)  -2\gamma(1+\gamma^2) yx^\tr  -4\gamma xy^\tr.\Ee
Since $\delc\lesssim 1-\gamma^2\lesssim 1$ and $\sin  S_c(x,y)\sim 1$, taking $\vepc$ small enough  it is sufficient to show $|\det\mathbf M|\gtrsim c$ for some $c>0$. 

From \eqref{sc} and \eqref{phase-s} we see  $S_c(\mathbf U x, \mathbf Uy)=S_c(x, y)$ and 
$\Phi_s(\mathbf Ux, \mathbf Uy)=\Phi_s(x, y)$ for any $\mathbf U \in \mathrm O(d)$.  Hence, choosing $\mathbf U \in \mathrm O(d)$ such that  
$\mathbf Ux=(r,0,0,\dots, 0)$ and $\mathbf U y=(\rho,h,0,\dots, 0)$ (see \eqref{rot}),     we may assume 
\begin{equation} 
\label{coordinates}
x=(r,0,0,\dots, 0), \  \  \  y=(\rho,h,0,\dots, 0). 
\end{equation}
Then  
$\mathbf M_{i,j}=0$ if  $i\neq j$ and  $i$ or $j\ge 3$, and  $\mathbf M_{j,j}= -(1-\gamma^2)^2$ if $j\ge 3$.
 Therefore, 
\Be
\label{det-h}
 \det  \mathbf M=  (-1)^{d-2} (1-\gamma^2)^{2(d-2)} \det    \widetilde {\mathbf M}
 +  O(\vepc),
  \Ee
where $\widetilde {\mathbf M}$ is the $2\times 2$ matrix given by  $\widetilde {\mathbf M}:= (\mathbf M_{i,j})_{1\le i,j\le 2}$. 
The matrix $\widetilde {\mathbf M}$ is easy to compute.   In fact, we only need to consider the first two components of the vectors 
$x$ and $y$. Let us set $\widetilde  x=(r,0)$ and $\widetilde  y=(\rho, h)$.  From \eqref{sc-matrix} it is clear that 
$
  \widetilde {\mathbf M}= -(1-\gamma^2)^2\,\mathbf I_2+   
  (1+3\gamma^2)(\widetilde x \widetilde x^\intercal+ \widetilde y \widetilde y^\tr)  -2\gamma(1+\gamma^2) \widetilde y \widetilde x^\tr  -4\gamma \widetilde x \widetilde y^\tr.
$
Hence, a routine computation gives 
\begin{align*}  \widetilde {\mathbf M}=&\begin{pmatrix}
-(1-\gamma^2)^2& 0\\ 
0&  -(1-\gamma^2)^2
\end{pmatrix}
  + (1+3\gamma^2)
           \begin{pmatrix}
               r^2+\rho^2& \rho h\\ 
                \rho h&  h^2
             \end{pmatrix} 
\\[4pt] 
&\qquad \qquad\qquad\quad-\begin{pmatrix}
2\gamma(3+\gamma^2)r\rho & 4\gamma r h\\ 
2\gamma(1+\gamma^2)r h  & 0
\end{pmatrix}.
\end{align*}
We note $r\rho=\gamma$ and $r^2+\rho^2+h^2=\gamma^2+1+O(\vepc)$ (by \eqref{det-small0}). After rearranging the matrix $\widetilde {\mathbf M}$,  we get
\[
\widetilde {\mathbf M}=  \begin{pmatrix}
-(1+3\gamma^2)h^2    &  \big ((1+3\gamma^2)\rho-4\gamma r\big)h
\\ 
\big((1+3\gamma^2)\rho-2\gamma(1+\gamma^2)r\big)h   &  (1+3\gamma^2) h^2 -(1-\gamma^2)^2
\end{pmatrix}  +O(\vepc).
\]
Using   $r^2+\rho^2+h^2=\gamma^2+1+	O(\vepc)$,  again,  and $r\rho=\gamma$,  a computation shows 
\[ \det  \widetilde {\mathbf M}=  h^2(1-\gamma^2)^2 |x|^2+O(\vepc). \] 
Let  $\theta(x,y)$ denote the angle  between $x$ and $y$. Since $h=|y|\sin \theta(x,y)$, we have
$ \det  \widetilde {\mathbf M}=  |x|^2|y|^2\sin^2 \theta(x,y)(1-\inp xy^2)^2 +O(\vepc)$.
Note $(1-|x|^2)(1-|y|^2)\gtrsim \delc^2$ because  $x,y\in B_{1-\delc}$. So, by \eqref{det-small0} and \eqref{angle},  it follows that  $|x|^2|y|^2 \sin^2\theta(x,y)\gtrsim \delc^2$ if $\vepc\ll \delc^2$.  
Recalling $1-\gamma^2\gtrsim \delc$,  we conclude that $\det  \widetilde {\mathbf M}\gtrsim \delc^4$  if  $\vepc$ is small enough. Finally we combine this with \eqref{det-h} to finish the proof. 
\end{proof}

\begin{rem}  
As  is clear in the course of the proof, to have the argument in the proof of Lemma \ref{nonzero-det} work,  it is important that $\sxy$ and its derivatives are uniformly bounded, that is to say, $\cos \sxy$ stays away from $1$ and $-1$. In order to handle the case where  $\cos \sxy$  gets closer to $1$ or $-1$,  we need further decomposition and refined analysis,  which we carry  out in what follows. 
\end{rem}

\section{Localization on annuli and $L^2$ estimate}
 In this section we consider the estimate over  $A_\mu^\sigma\times A_{\mu'}^{\sigma'}$, $\mu'\le \mu\ll 1$ with $\sigma,\sigma'=\pm$, $\circ$ (see \eqref{def-annalus} and \eqref{coarse-decomp} below).   We recall  \eqref{ph-d} and that $\cD(x,y)$ is the discriminant of the quadratic equation $\cQ(x,y,\tau)=0$. 
 As observed in Section \ref{sec:local}, the value of $\cD$ plays an important role in determining the boundedness of $\Pi_\lambda$.   
Then, we  write
\begin{align}
 \label{angle} 
\mathcal D(x,y)=-|x|^2|y|^2 \sin^2\theta(x,y) +(1-|x|^2)(1-|y|^2), 
 \end{align}
 where   $\theta(x,y)$ denotes  the angle  between $x$ and $y$. 
 Since $|(1-|x|^2)(1-|y|^2)|\sim \mm$ for $(x,y)\in \aas$ for $\sigma, \sigma'=\pm$,   
comparative size of the angle $\theta(x,y)$  against $\mmt$ controls the value of $\cD$. Thus, we are led to make a decomposition of  $A_\mu^\sigma\times A_{\mu'}^{\sigma'}$ which is convenient for 
controlling the angle between $x$ and $y$.    

For the purpose  we use a Whitney type decomposition of   $\mathbb S^{d-1}\times \mathbb S^{d-1}$ away from its diagonal. Consequently, we  get the decomposition \eqref{decomp} below, which  provides an efficient way of localization. In particular, this allows us  to use the results in Section \ref{sec:local}.

\subsection{Sectorial decomposition of annuli}    

Following the typical dyadic decomposition process, for each $\nu\ge 0$ we partition
$\mathbb S^{d-1}$ into spherical caps $\Theta_k^\nu$  such that  $\Theta_k^\nu \subset \Theta_{k'}^{\nu'}$ for some $k'$ whenever $\nu\ge \nu'$ 
and $c_d  2^{-\nu}\le 
  \diam (\Theta_k^\nu)  \le C_d2^{-\nu}$ for some constants $c_d$, $C_d>0$. Let  
$\nu_\circ:=\nu_\circ(\mu, \mu')$ denote  the integer  $\nu_\circ$ such that 
\Be\label{nu0}  \mm/2<  2^6 C_d^2 2^{-2\nu_\circ} \le \mm.  
\Ee
Then, we may write
\[\mathbb S^{d-1}\times \mathbb S^{d-1}=\bigcup_{ \nu:  2^{-\nu_\circ} \le    2^{-\nu} \lesssim 1 }\,\, \bigcup_{k\sim_\nu k'} \Theta_k^\nu\times \Theta_{k'}^{\nu}\] 
where  $k\sim_\nu k'$ implies $\dist (\Theta_k^{\nu}, \Theta_{k'}^{\nu})\sim 2^{-\nu}$ if  $\nu> \nu_\circ$ and  $\dist (\Theta_k^{\nu}, \Theta_{k'}^{\nu})\lesssim 2^{-\nu}$ if  $\nu= \nu_\circ $.  For example, see \cite[p.971]{TVV}.  It should be noted that the sets $\Theta_k^{\nu_\circ}$ and $\Theta_{k'}^{\nu_\circ}$  may be not separated because we stop the decomposition 
procedure at $\nu=\nu_\circ$.  For fixed $\mu$ and $\mu'$, we define 
\[  A_{\mu,k}^{\sigma,\nu}=\Big\{ x\in A_{\mu}^\sigma:  \frac{x}{|x|}\in    \Theta_{k}^\nu\Big\}, 
\ \ \ \sigma\in\{+,\circ,-\}. \] 
Throughout the paper we assume that the indices $k$, $k'$ are associated to $\mu$, $\mu'$, respectively, and we drop the subscript $\mu$, $\mu'$ to simplify notation. 
Thus we may write 
\begin{equation}
\label{decomp}
A_{\mu}^\sigma \times  A_{\mu'}^{\sigma'}=   \bigcup_{ \nu:   2^{-\nu_\circ}\le     2^{-\nu} \lesssim 1 }\,\, \bigcup_{k\sim_\nu k'} A_k^{\sigma,\nu}\times A_{k'}^{\sigma',\nu}, 
\quad
\sigma, \sigma'\in\{+,\circ,-\}.
\end{equation}
For simplicity  we also set
\[  \chi_{k}^{\sigma, \nu}=\chi_{A_{k}^{\sigma,\nu}}, 
\quad \sigma\in\{+,\circ,-\}. \]

\subsection{Estimates for $I_j$ over  $\aask$} 
We  obtain further estimates for $I_j(x,y)$ defined by \eqref{S3-ptkernel} while $(x,y)\in \aask$.  We separately handle the cases 
$2^{-\nu}\gg \mu$  and $2^{-\nu}\lesssim \mu$.

\begin{lem}
\label{case-a}  Let   $0<\mu'\le \mu\ll 2^{-\nu}\le 1/100$. 
If    $(x,y)\in  A_k^{\circ, \nu}\times A_{k'}^{\circ, \nu}$ with $k\sim_\nu k'$, then, for any $N>0$, we have  
\begin{equation}
\label{a-case-kernel}
   | I_j(x, y)|    \lesssim  \begin{cases}   
                                 2^{- j} \big( \lambda2^{j}2^{-2\nu}+1\big)^{-N} ,        
                                           &  2^{-2j} 
\lesssim   2^{-\nu}                                      \\
                                  2^{- j}    \big( \lambda 2^{-3j}+1\big) ^{-N},      
                                               &   2^{-2j}  \gg 2^{-\nu}  
                          \end{cases}\,.
                       \end{equation}
\end{lem}

Taking $N=1/2$ in  \eqref{a-case-kernel} gives bounds on $\|\chi_{k}^{\circ,\nu} \mathfrak P_\lambda [ \psi_j ]  \chi_{k'}^{\circ,\nu}\|_{1\to \infty}$. 
 Then, we interpolate those bounds with $\|\chi_{k}^{\circ,\nu} \mathfrak P_\lambda [ \psi_j ]  \chi_{k'}^{\circ,\nu}\|_{2\to 2}\lesssim \lambda^{-\frac d2} 2^{-j}$ which follows from 
\eqref{easy-l2}. So,  we get
\begin{equation}
\label{eq:haha8}
\|\chi_{k}^{\circ,\nu} \mathfrak P_\lambda [ \psi_j ]  \chi_{k'}^{\circ,\nu}\|_{p\to p'} 
          \lesssim  \begin{cases}   
                                 \lambda^{\frac{d-1}2\dpp-\frac d2} 2^{(\frac{d-1}2\dpp -1)j}  2^{\nu\dpp },  
                                      & \!\! 2^{-2j} 
\lesssim   2^{-\nu}                                  
\\
                                 \lambda^{\frac{d-1}2\dpp-\frac d2} 2^{(\frac {d+3}2\dpp -1) j},    
                                 &  \!\!  2^{-2j}  \gg 2^{-\nu} 
                          \end{cases} \hspace{-.4cm}
\end{equation}
for  $k\sim_\nu k'$ and 
 $1\le p\le 2$. 

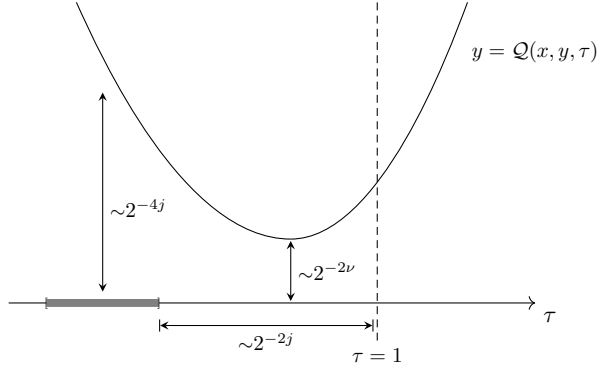
\begin{figure}[!t]
\begin{tikzpicture} [scale=0.5]
\draw[->] (1,1)--(15,1); 
\draw[ densely dashed] (10.8,0)--(10.8,9); 
\node[below right, scale=0.9] at (15,1) {$\tau$};
\node[ below, scale=0.8] at (10.8,0) {$\tau=1$};
\draw (8.5,2.7) parabola (2.8,9); 
\draw  (8.5,2.7) parabola (13.2,9);
\draw[<->,>=stealth] (8.5,1.05)--node[midway, right, scale=0.8]{{\small $\sim$} $\!\!2^{-2\nu}$}(8.5,2.65);
\node[above, scale=0.8] at (15,7.2) {$y=\mathcal Q(x, y, \tau)$};
\draw[gray, line width=3] (2,1)--(5,1);
\node[scale=0.7] at (2,1){{\tiny$[$}};\node[scale=0.7] at (5,1){{\tiny$]$}};
\draw[|<->|,>=stealth] (5,0.4)--node[midway, below, scale=0.8]{{\small $\sim$} $\!\!2^{-2j}$}(10.7,0.4);
\draw[<->,>=stealth] (3.5,1.3)--node[midway, below right, scale=0.8]{{\small $\sim$} $\!\!2^{-4j}$}(3.5,6.6);
\end{tikzpicture}
\caption{ The case $2^{-\nu}\gg \mu$, $(x,y)\in A^{\circ, \nu}_k\times A^{\circ, \nu}_{k'}$, and $2^{-2j}\gg 2^{-\nu}.$ The thick  line segment contains the set $\cos (\supp \psi_j)$. }
\label{fig:pm}

\end{figure}

\begin{proof}[Proof of Lemma \ref{case-a}]
If  $(x,y)\in \aack$, the angle $\theta(x,y)$ is $\sim 2^{-\nu}$ and $|x-y|\sim 2^{-\nu}$ since $\mu\ll 2^{-\nu}$. Thus, 
from \eqref{angle}  it follows that $-\cD(x,y)\sim 2^{-2\nu}$  and  we  also have
\[   |x-y|\sim 2^{-\nu},  \  \   -\cD(x,y)\sim  2^{-2\nu},   \   \   |1-\inp{x}{y}|\lesssim 2^{-\nu}\]
for $(x,y)\in \aack$. 
For the last one may use the identity $|x-y|^2+\cD(x,y)=(1-\inp xy)^2$. From \eqref{ph-d}  we see (see Figure \ref{fig:pm})  
\[
|\partial_s \cP(x,y,s) |
            \gtrsim   \begin{cases}   
                                 2^{-2\nu}2^{2j},           
                                         & 2^{-2j} \lesssim   2^{-\nu} 
                                         \\
                                \  2^{-2j} ,       
                                & 2^{-2j}  \gg 2^{-\nu}                           \end{cases} \, ,
                          \  \   \  \   (x,y)\in \aack 
\]
  for $s\in \supp \psi_j$. 
   On the other hand, using  \eqref{ph-d}, a computation shows  
\[  \Big| \ddk(x,y,s)\Big| \lesssim   \begin{cases}   
                                 2^{-2\nu}2^{(1+k)j},            
                                         & 2^{-2j} 
\lesssim   2^{-\nu}                                          \\
                                \  2^{-4j}  2^{(1+k)j},       
                                &  \, 2^{-2j}  \gg 2^{-\nu} 
                          \end{cases} \,   ,
                          \  \   \  \   (x,y)\in\aack \] 
on the support of $\psi_j$. Hence, using  Lemma  \ref{S3-osc} we get \eqref{a-case-kernel}.
\end{proof}

\begin{lem} \label{case-b} Let $0<\mu'\le \mu\le 1/100$, $\sigma,\sigma'\in \{+,-\}$,  and $(x,y)\in A_k^{\sigma, \nu}\times A_{k'}^{\sigma', \nu}$ with $k\sim_\nu k'$. Assume that $2^{-\nu}\lesssim \mu.$
\vspace{-6pt}
\begin{enumerate} 
[leftmargin=.5cm, labelsep=.15cm, topsep=1pt]
\item[$(a)$] Suppose that $\mu'\ll\mu$. Then we have the following\,$:$
\begin{enumerate}
[topsep=1pt]
    \item [$(a1)$] For any $N>0$, we have
    \begin{equation}
\label{b-case-kernel}
     | I_j(x,y) |      \lesssim  \begin{cases}   
                                 2^{- j} \big( \lambda 2^{j}\mu^{2}+1\big)^{-N} ,        &  2^{-j}  \ll   \smu  \ \ 
                                                                \\
                                   2^{- j} \big( \lambda 2^{-3j}+1\big)^{-N},        & 2^{-j}\,  \gg    \smu
                          \end{cases}.
\end{equation} 
    \item[$(a2)$] If $2^{-j}\sim \sqrt{\mu}$, we have  the first case of \eqref{b-case-kernel} for any $N$ for $(\sigma, \sigma') = (-,+)$, $(-,-)$.
\end{enumerate}
\item[$(b)$] 
Suppose that $\mu'\sim\mu$. Then we have the following\,$:$
\begin{enumerate}[ topsep=1.5pt]
    \item [$(b1)$] If $2^{-j}\ll\sqrt{\mu}$, then the first case of \eqref{b-case-kernel} is valid for any $N$ under the additional assumption that $\sigma = \sigma'$, $|\mathcal D(x,y)|\ll\mu^2$.
    \item[$(b2)$] If $2^{-j}\gg\sqrt{\mu}$, then the second case of  \eqref{b-case-kernel} holds  for any $N$.
\end{enumerate}
\end{enumerate}  
\end{lem}

Similarly 
 as before,  the estimate \eqref{b-case-kernel} gives  $L^1$--$L^\infty$ estimates for $\chi_{k}^{\pm,\nu} \mathfrak P_\lambda [ \psi_j ]  \chi_{k'}^{\pm,\nu}$.  
 Taking $N=1/2$ in  \eqref{b-case-kernel}  and  interpolating those consequent estimates and  the $L^2$ estimate $\|\chi_{k}^{\pm,\nu} \mathfrak P_\lambda [ \psi_j ]  \chi_{k'}^{\pm,\nu}\|_{2\to 2}\lesssim \lambda^{-\frac d2} 2^{-j}$ which  follows  from \eqref{easy-l2},   we obtain
\begin{equation}
\label{eq:sec8} 
         \|\chi_{k}^{\sigma,\nu} \mathfrak P_\lambda [ \psi_j ]  \chi_{k'}^{\sigma',\nu}\|_{p\to p'}          
                \lesssim  \begin{cases}                                   
                      \lambda^{\frac{d-1}2\dpp-\frac d2} 2^{(\frac{d-1}2\dpp -1)j}  \mu^{-\dpp} ,    & \  2^{-j}   \ll  \smu            
                                                                \\                                
                      \lambda^{\frac{d-1}2\dpp-\frac d2} 2^{(\frac{d+1}2\dpp -1)j}  2^{ \frac \nu 2 \dpp },  & \ 2^{-j}  \sim  \smu 
                                                                 \\                        
                        \lambda^{\frac{d-1}2\dpp-\frac d2} 2^{(\frac {d+3}2\dpp  -1) j} ,         & \ 2^{-j}\,   \gg   \smu     
\end{cases}\end{equation}
for $k\sim_\nu k'$,  $\mmt \ll 2^{-\nu} \lesssim \mu$, and  $\sigma,\sigma'\in \{+,-\}$  provided  $1\le p\le 2$. The $L^1$--$L^\infty$ estimate in the second case $2^{-j}  \sim  \smu$   follows from Lemma \ref{oscillatory}  since $|\cD(x,y)|\sim  2^{-2\nu}$ if $(x,y)\in A_k^{\sigma,\nu}\times A_{k'}^{\sigma',\nu}$, $k\sim_\nu k'$, and $\mmt \ll 2^{-\nu}$.

\begin{figure}[!t]
\hspace{-35pt}\begin{minipage}[t]{.35\textwidth}
\vspace{0pt}

\begin{tikzpicture} [scale=0.47]
\draw[->] (0,1)--(12,1);  
\draw[densely dashed] (11.2,0.3)--(11.2,8.2); 
\node[above, scale=0.8] at (11.2,8.2) {$\tau=1$};
\node[below , scale=0.9] at (12.2,0.8) {$\tau$};

\draw (6,2.3) parabola (0,8); 
\draw  (6,2.3) parabola (12,8);
\draw[<->,>=stealth] (6,1.05)--node[midway, right, scale=0.8]{ {\small $\sim$}  $\!\!2^{-2\nu}$}(6,2.25);
\draw[|<->|,>=stealth] (6,0.4)--node[midway, below, scale=0.8]{{\small $\sim$} $\!\!\mu$}(11.2,0.4);
\node[above, scale=0.8] at (3.1,6.4) {$y=\mathcal Q(x, y, \tau)$};
\draw[gray, line width=3] (9.5,1)--(10,1);
\node[scale=.62] at (9.5,1){{\tiny$[$}};\node[scale=.62]  at (10,1){{\tiny$]$}};
\draw[<->,>=stealth] (9.75,1.2)--node[midway,  right, scale=0.8]{{\small $\sim$} $\!\!\mu^2$}(9.75,4.3);

\node[below, scale=.9] at (6,-0.8) {($\mathbf{I}$) $2^{-j}\ll \sqrt{\mu}$ and  $A_k^{+,\nu}\times A_{k'}^{+,\nu}$};

\end{tikzpicture}
\end{minipage}
\hspace{2.3cm}
\begin{minipage}[t]{.35\textwidth}
\vspace{0pt}

\begin{tikzpicture} [scale=0.47]
\draw[->] (4,1)--(16,1);  
\draw[densely dashed] (9,0.6)--(9,8.2); 
\node[above, scale=0.8] at (9,8.2) {$\tau=1$};
\node[above left, scale=0.9] at (16, 0.0) {$\tau$};
\draw (10,2.3) parabola (4,8); 
\draw  (10,2.3) parabola (16,8);
\draw[<->,>=stealth] (10,1.05)--node[midway, right, scale=0.8]{ {\small $\sim$} $\!\!2^{-2\nu}$}(10,2.2);
\draw[|<->|,>=stealth] (10,0.4)--node[midway, below, scale=0.8]{{\ \ \small $~\sim$} $\!\!\mu$}(8,0.4);
\node[above, scale=0.8] at (12.5,6) {$y=\mathcal Q(x, y, \tau)$};
\draw[gray, line width=3] (7,1)--(8,1);
\node[scale=.7] at (7,1){{\tiny$[$}};\node[scale=.7]  at (8,1){{\tiny$]$}};
\draw[<->,>=stealth] (7.5,1.2)--node[midway,  left, scale=0.8]{{\small $\sim$} $\mu^2$}(7.5,3.2);
\node[below, scale=.9] at (10,-0.8) {($\mathbf{I\!I}$) $2^{-j}\sim \sqrt{\mu}$ and  $A_k^{-,\nu}\times A_{k'}^{-,\nu}$};

\end{tikzpicture}
\end{minipage}

\caption{Particular cases when $\sqrt{\mu\mu'} \ll 2^{-\nu}\lesssim \mu$. The thickened line segment  contain $\cos(\supp \psi_j )$. }
\label{fig:2}
\end{figure}
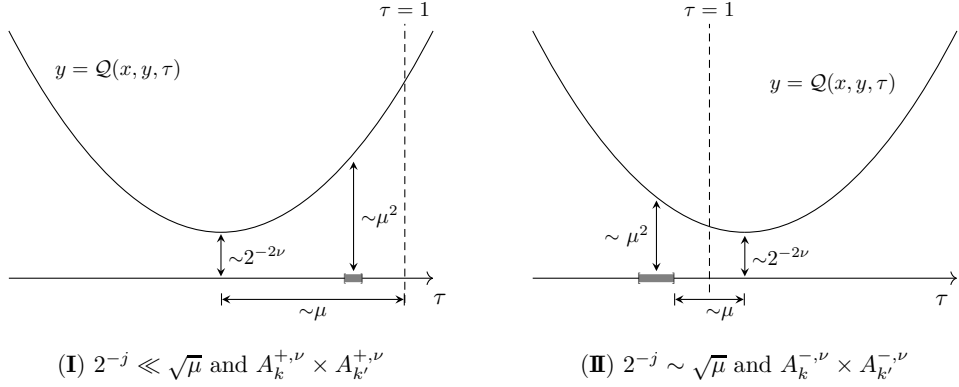

\begin{proof}[Proof of Lemma \ref{case-b}] 
To show $(a1)$, we consider the two cases 
\[ \mmt \ll 2^{-\nu}\lesssim \mu, \quad  \mmt\gtrsim  2^{-\nu}, \]
 separately. 
We  begin by  noting  that $ |1-\inp{x}{y}|\sim \mu $ for $(x,y)\in A_k^{\sigma, \nu}\times A_{k'}^{\sigma', \nu}$ since $\mu'\ll \mu$ and  $|1-\inp xy|\sim 
|1-\inp xy^2|=|1-|x|^2|y|^2|+O(2^{-2\nu}) $.

We first handle   the case $\mmt \ll 2^{-\nu}\lesssim \mu$.  Since $2^{-\nu}\lesssim \mu$ and $\mu'\le \mu\ll 1$, 
from \eqref{angle}  
 it follows that 
 $0< -\cD(x,y)\sim  2^{-2\nu}\lesssim  \mu^2,$ $(x,y)\in A_k^{\sigma, \nu}\times A_{k'}^{\sigma', \nu}.$
Since $\cQ(x,y,\cos s)= (\cos s- \inp xy)^2-\mathcal D$, $ |1-\inp{x}{y}|\sim \mu $, and $1-\cos s\sim 2^{-2j}$ on the support of $\psi_j$,  
we have 
\Be
\label{Q-}
\cQ(x,y,\cos s)
            \sim   \begin{cases}   
                                 \mu^2,            &  2^{-j}   \ll  \smu  
                                                                 \\
                                  \   2^{-4j},         &  2^{-j}\,   \gg   \smu
                          \end{cases}, \  \   \ \ (x,y)\in  A_k^{\sigma, \nu}\times A_{k'}^{\sigma', \nu}
\Ee
Thus,  from \eqref{ph-d}  we see that
\Be
\label{lower-case}
|\partial_s \cP(x,y,s) |
            \gtrsim   \begin{cases}   
                                 \mu^2 2^{2j},            &  2^{-j}   \ll  \smu  
                                                                 \\
                                  \   2^{-2j},         &  2^{-j}\,   \gg   \smu
                          \end{cases}, \  \   \ \ (x,y)\in  A_k^{\sigma, \nu}\times A_{k'}^{\sigma', \nu}
  \Ee
 if  $s\in \supp\psi_j$.
  On the other hand, using  \eqref{ph-d} and by a computation we get
\Be 
\label{lower-case2}
\Big| \ddk(x,y,s)\Big| \lesssim   \begin{cases}   
                                 \mu^2 2^{(1+k)j} ,           
                                          &  2^{-j}   \ll  \smu  
                                         \\
                                \ 2^{-4j}  2^{(1+k)j},         
                                &  2^{-j}   \gg  \smu  
                          \end{cases} \,   ,
                          \  \   \  \   (x,y)\in A_k^{\sigma, \nu}\times A_{k'}^{\sigma', \nu}  \Ee
on the support of $\psi_j$. Thus, by  Lemma \ref{S3-osc} we obtain \eqref{b-case-kernel}. 

The remaining  case  $\mmt\gtrsim  2^{-\nu}$ can be handled in the same manner.  Indeed,  from \eqref{angle} we note  that 
$|\cD(x,y)|\lesssim \mm$ for $(x,y)\in A_k^{\sigma, \nu}\times A_{k'}^{\sigma', \nu}$ since $\mu'\ll \mu$. Then, following the same argument in the above, we see 
\eqref{Q-} holds. 
Once we have those (upper and lower) bounds on $\cQ(x,y,\cos s)$, the rest of argument is similar to that of the previous case.  
So, we omit the detail.

We now show $(a2)$.  Note that $2(\inp xy-1)=|x|^2-1+|y|^2-1-|x-y|^2$. Since $|x-y| \lesssim \max(\mu,2^{-\nu})$, $2^{-\nu}\lesssim \mu$ and $\mu'\ll \mu\ll 1$,   we have $\inp xy -1\sim \mu$ for  $(x,y)\in  A_k^{-, \nu}\times A_{k'}^{\pm, \nu}$. 
Since $\cQ(x,y,\cos s)= (\cos s- \inp xy)^2-\mathcal D(x,y) \ge \cQ(x,y,1)=|x-y|^2\gtrsim \mu^2$ and $ 2^{-j} \lesssim  \smu$,   we obtain  $\cQ(x,y,\cos s)\sim \mu^2$ 
on $\supp \psi_j$ (see $(\mathbf {I\!I})$ in  Figure \ref{fig:2}). Now, by \eqref{ph-d} it follows that
$
|\partial_s \cP(x,y,s) |  \sim  \mu^2 2^{2j}
$
 for $(x,y)\in A_k^{-,\nu}\times A_{k'}^{\pm,\nu}$,  and we also have the upper bound $
 |\ddk(x,y,s)|\lesssim \mu^2 2^{(1+k)j}$  if $s\in \supp \psi_j$ because $2^{-j}\lesssim \smu$. Therefore, by Lemma \ref{S3-osc} we get the desired estimate.   

It remains to show $(b1)$, $(b2)$. They can be shown by a slight modification of the previous arguments.  
We first consider $(b2)$.  Since $2^{-\nu}\lesssim \mu$,  we note   $|\cD(x,y)|\lesssim \mu^2$ and $ |1-\inp{x}{y}|\lesssim \mu $ for $(x,y)\in A_k^{\sigma, \nu}\times A_{k'}^{\sigma', \nu}$.  We also have $\cQ(x,y,\cos s)=(\cos s- \inp xy)^2-\mathcal D\sim 2^{-4j}$ on the support of $\psi_j$ since $2^{-j}\gg   \smu$. Therefore, by the bound on  $\cQ(x,y,\cos s)$  we have the estimates of the second cases of \eqref{lower-case} and \eqref{lower-case2}. Thus, we  get the desired estimate repeating  the same argument.      

We finally deal with $(b1)$.  Since $1-\inp xy^2=1-|x|^2+1-|y|^2-\cD(x,y)$,  $(x,y)\in A_k^{\sigma, \nu}\times A_{k'}^{\sigma, \nu}$, and $|\cD(x,y)|\ll \mu^2$,  we have $ |1-\inp{x}{y}|\sim \mu $.   Thus, we have $\cQ(x,y,\cos s)\sim \mu^2$ if $s\in \supp \psi_j$ since $ 2^{-j} \ll  \smu$. 
In the same manner as before we have the estimates of the first case of \eqref{lower-case} and \eqref{lower-case2}. The remaining  is identical as before, so we omit the detail. 
\end{proof}

\subsection{Reduction to small separation}
Let $c$ be a  positive constant such that $1/(200d)\le c\le  1/(100d)$, and let $0<\mu, \mu'\ll c$.  We partition  the annuli  $A^\circ_\mu$   and $A^\circ_{\mu'}$ into  finitely many disjoint sets
of diameter $\sim c$ such that 
\Be
\begin{aligned}
\label{coarse-decomp}
A_\mu^\circ&=\{ x: | 1-|x||\le 2\mu \}=\bigcup  A_j,   
\\
   A_{\mu'}^\circ&= \{ x: | 1-|x||\le  2\mu' \}=\bigcup  A_j'
\end{aligned}
 \Ee
 and the volumes of  $ A_j$ and  $A_j'$ are $\sim c^{d-1} \mu$ and $\sim c^{d-1} \mu'$, respectively. 
 This coarse decomposition allows us to sort out the minor parts whose contributions are less  significant.  Indeed, if  $\dist (A_j, A_k')\ge c$ and $\dist (A_j, -A_k')\ge c,$  the contribution of $\chi_{A_j} \fP_\lambda\chi_{A_k'}$ is negligible.    More precisely, we have the following.

\begin{lem}
\label{large-sep}   Let $0<\mu,\mu'\ll  (100d)^{-2}$,  $A\subset \Ac_\mu$, and $A'\subset  \Ac_{\mu'}$. 
Suppose $\dist (A, A')\ge c$  and $\dist (A, -A')\ge c$  for  a constant $ (200d)^{-1}\le c\le (100d)^{-1}$. Then,   we have 
\Be\label{improved}
\| \chi_{A} \fP_\lambda[\psi_j^\kappa ] \chi_{A'}f\|_{q}\lesssim   \lambda^{-N} 2^{-Nj}  (\mm)^{\frac12(1-\dpq)} \|f\|_p, \quad \kappa=0,\, \pm, \, \pm \pi,
\Ee
for any $N>0$
provided that $1/p+1/q=1$ and $1\le p\le 2$.   Furthermore, $(i)$ if $\mu\sim \mu'$, then \eqref{improved} holds for any $p,q$ satisfying  $1\le p\le q\le \infty$.
\end{lem}

 Throughout the paper, we occasionally use the following  lemma, which follows from Young's inquality.

\begin{lem}
\label{trivial}   Let $E$ and $F$ be measurable sets of finite measure.  Suppose 
$  |Tf(x)|\le D \int  \chi_E(x)\chi_F(y) |f(y)| dy.$
Then $\|T\|_{2\to 2}\le D\sqrt{|E||F|}$. Additionally, $(a)$ if $|E|,|F|\lesssim B$, then $\|T\|_{p\to q}\lesssim D B^{1-\dpq}$  for $1\le p\le q\le \infty$. 
\end{lem}

\begin{proof}[Proof of Lemma \ref{large-sep}]   Since  $\sin^2\theta(x,y)\gtrsim c$ and $0<\mu,\mu'\ll c$,   from \eqref{angle}  we see $-\cD(x,y)\sim c$ 
and $\cQ(x,y, \tau)\gtrsim c$ for $\tau\in \mathbb R$ if $(x,y)\in A\times A'$.  So, by \eqref{ph-d}  we get  $|\partial_s \cP |\gtrsim 2^{2j}$ on the support of $\psi_j^\kappa$. Also, 
it is easy to show $|\partial_s^l \cP |\lesssim 2^{(l+1)j}$ for $l\ge 1$ if $s\in \supp\psi_j^\kappa$. 
Hence, by  Lemma \ref{S3-osc}  we have
\Be
\label{rapid}
 \Big| \int   e^{i\lambda \mathcal P(x,y,s)} \psi_j^\kappa(s) ds\Big|\lesssim 2^{- j}\big( \lambda2^{j} +1\big)^{-N}, \quad (x,y)\in A\times A'
\Ee
for any $N$. Since $|A|\lesssim \mu$ and $|A'|\lesssim \mu'$, 
by the estimate \eqref{rapid} and Lemma \ref{trivial} 
we get  $\| \chi_{A} \fP_\lambda[\psi_j^\kappa ] \chi_{A'}\|_{2\to 2}\lesssim   \lambda^{-N} 2^{-Nj}  (\mm)^{\frac12}$ for any $N$. Interpolation between this and the obvious $L^1$--$L^\infty$ estimate 
gives the desired bound \eqref{improved} with $1/p+1/q=1$.  Finally, for the statement $(i)$ we use $(a)$ in Lemma \ref{trivial} to get \eqref{improved} for $1\le p\le q\le \infty$ because 
$|A|, |A'|\lesssim \mu$. \end{proof}

We also have the following if $A$ and $A'$ have small diameters.

\begin{lem}
\label{small-sep}  Let  $(200d)^{-1} \le c \le  (100d)^{-1}$, $0<\mu,\mu'\ll  (100d)^{-2}$. Suppose that 
$A\subset  A_\mu^\circ$, $A'\subset  A_{\mu'}^\circ$ and  $A$, $A'$ are of diameter $\le c$. Then,
\vspace{-6pt}
\begin{enumerate} 
[leftmargin=0.85cm, labelsep=0.3 cm, topsep=0pt]
\item[$(a)$] if $\dist (A, A')< c$,   we have  \eqref{improved} for $\kappa=\pm \pi,\, 0$, 
\item[$(b)$] if $\dist (A, -A')< c$,  we have  \eqref{improved} for $\kappa=\pm,\, 0$, 
\end{enumerate}
\vspace{-6pt}
provided that  $1/p+1/q=1$ and $1\le p\le 2$. 
Furthermore, $(i)$ if $\mu\sim \mu'$, then the same statements  $(a)$ and $(b)$  hold for $p,q$ satisfying  $1\le p\le q\le \infty$.
\end{lem}

\begin{proof}  By  the second identity in \eqref{symmetric} it is sufficient to show $(b)$.  
Since $\dist (A, -A')$ $< c$,  we have $|1+\inp xy |\le 3c$ 
whereas $\cos s \ge - \cos 2^{-3}$ on the support of 
$\psi_j^{\kappa }$, $\kappa=\pm, 0$,  and $|\cD(x,y)|\le 2c$ by \eqref{angle}.  So, we have $|\cQ(x,y, \cos s)|\sim 1$  by \eqref{q-def}.  Thus, 
we see $|\partial_s \cP(x,y,s)|\gtrsim  2^{2j}$ and, similarly as before, $|\partial_s^l \cP |\lesssim 2^{(l+1)j}$ for $l\ge 1$  if $s\in \supp\psi_j^{\kappa}$, $\kappa=\pm,0$.  Hence, Lemma  \ref{S3-osc} gives the estimate \eqref{rapid}. 
Once we have \eqref{rapid}, the desired estimate \eqref{improved}  follows by the same argument as in the proof of Lemma \ref{large-sep}. 
The statement $(i)$ can be shown similarly as before, so we omit the detail.    
\end{proof}

\begin{rem}
\label{separated}
The bounds in Lemma \ref{large-sep} and Lemma \ref{small-sep} are much smaller  than what  we intend  to prove for 
$ \chi_{\mu}^\sigma \fP_\lambda[\psi_j^\kappa ] \chi_{\mu'}^{\sigma'}$. Thus, by a finite decomposition of $A_\mu^\sigma$ and $A_{\mu'}^{\sigma'}$ $($see \eqref{coarse-decomp}$)$ it is sufficient to show the estimate 
for $ \chi_{A} \fP_\lambda[\psi_j^\kappa ] \chi_{A'}$ while assuming that $A\subset A_\mu^\sigma$,  $A'\subset A_{\mu'}^{\sigma'}$ are of  small diameter and 
\Be
\label{distanced}
\dist (A, A')\le c \text{ and } \kappa=\pm, \text{ or } \dist (A, -A')\le c \text{ and } \kappa=\pm \pi.
\Ee
By the second identity in \eqref{symmetric} the estimate for the second case can be deduced from the one for the first case, so it is sufficient to consider the first case only. Moreover, 
the estimate for the case $\kappa=-$ can be deduced from that for the case $\kappa=+$ because of the first identity in \eqref{symmetric}. Thus we are reduced to handling the case 
$\dist (A, A')\le c \text{ and } \kappa=+$. 
\end{rem}

\subsection{Refined $L^2$ bound}
In this subsection we obtain $L^2$ estimates for $\Pi_\lambda$ over the set $A_{\lambda,\mu}^\pm\times A_{\lambda,\mu'}^\pm$ which we use  later. 

It was shown in Koch and Tataru \cite{T05} that localization near the sphere $\sqrt \lambda \mathbb S^{d-1}$ yields an improved $L^2$ estimate for $\Pi_\lambda$.   More precisely, they obtained 
$
 \| \Pi_\lambda   \chi_{\lambda, \mu}^\pm \|_{2\to 2} \lesssim   \mu^\frac14
$
(see \eqref{kt-nu}), 
which is clearly equivalent to 
\Be
\label{l2-ktt}
 \| \chi_{\lambda,\mu}^\pm \Pi_\lambda  \chi_{\lambda, \mu}^\pm \|_{2\to 2} \lesssim   \mu^\frac12.
 \Ee  
 This estimate can be  strengthened  further when $0<\mu\lesssim \lambda^{-2/3}$, see the estimate \eqref{smallmu} in 
Proposition \ref{prop:estext}.   For later use we record the following 
 which seemly looks stronger than  \eqref{l2-ktt} but it is equivalent to  \eqref{l2-ktt} as is easy to show using the $TT^*$ argument. 

\begin{lem} \label{l2-mmj}   Let $0<\mu'\le \mu\le 1$. Then,  we have     
\[\|\chi_{\lambda,\mu}^\sigma \Pi_\lambda  \chi_{\lambda, \mu'}^{\sigma'}  \|_{2\to 2} \lesssim  \mmf, \quad  \sigma, \sigma'\in \{+,-\} .     \]
\end{lem}
The estimate \eqref{l2-ktt} has its counterpart \eqref{jl2} below  for the operator $\Pi_\lambda[ \psi_j ]$ which  may be regarded as an extension of  \eqref{l2-ktt} 
since \eqref{jl2} implies \eqref{l2-ktt}. 
\begin{lem}  
\label{local-l2}
Let  $0< \mu\le 1$ and $2^{-j}\gg \smu$. Then, for  $\kappa=\pm, \pm \pi, 0$,  we have
\begin{equation}
\label{jl2}
   \| \chi_{\lambda, \mu}^\sigma \Pi_\lambda[ \psi_j^\kappa ]  \chi_{\lambda, \mu}^\spp    \|_{2\to 2} \lesssim 2^j\mu, \quad  \sigma, \sigma'\in \{+,-\} . 
   \end{equation}
\end{lem}
It is easy to see that the estimate \eqref{jl2} implies \eqref{l2-ktt}. Indeed, we have the bound $\|\Pi_\lambda[ \psi_j^\kappa ]\|_{2\to 2}\lesssim 2^{-j}$  from   \eqref{easy-l2}, thus
$ \|\sum_{2^{-j}\lesssim \smu}  \Pi_\lambda[ \psi_j^\kappa ] \|_{2\to 2}\lesssim \sum_{2^{-j}\lesssim \smu}  2^{-j}\lesssim \smu$ for $\kappa= \pm, \ \pm \pi, 0. $
On the other hand, from \eqref{jl2} we have 
\[ \Big \| \chi_{\lambda, \mu}^\sigma\Big( \sum_{2^{-j}\gg \smu}  \Pi_\lambda[ \psi_j^\kappa ]\Big) \chi_{\lambda,\mu}^\spp \Big\|_{2\to 2}\lesssim \sum_{2^{-j}\gg\smu} \mu 2^j  \lesssim \smu, 
\  \  \  \kappa= \pm, \ \pm \pi, \ 0. \]
Thus, combining the above two estimates  with  \eqref{decomp-proj}  we get  \eqref{l2-ktt}.

In order to show Lemma \ref{local-l2},  we use  the following Lemma \ref{kkpsum} and Lemma \ref{trivial1}. 

\begin{lem}\label{kkpsum}  Let $1\le p\le q\le \infty$ and let $T$ be an operator from $L^p$ to $L^q$.   Suppose we have the estimate 
$\|\spxkl   T  \spxklp\|_{p\to q}\le B$ whenever $k\sim_\nu k'$. Then, with $C$ only depending on $d$,   we have 
\[
\|\sum_{k\sim_\nu k'} \spxkl T  \spxklp\|_{p\to q}\le C B.
\]
\end{lem}

\begin{proof}  Since  $\{A_k^{\sigma, \nu}\}_k $  are disjoint,   for each $k$ there are at most as many as $O(1)$  $k'$ such that  $k\sim_\nu k'$. Thus,  we have
\begin{align*}
\|  \sum_{k\sim_\nu k'} \spxkl T  \spxklp  f\|_q\lesssim \big( \sum_{k\sim_\nu k'}  \|   \spxkl T  \spxklp  f\|_q^q \big)^{1/q} 
\lesssim  B \big( \sum_{k'}  \| \spxklp  f\|_p^q \big)^{1/q} .
\end{align*}
The last  is clearly bounded by $CB\|f\|_p$ because $p\le q$ and $\{A_k^{\sigma', \nu}\}_k$  are disjoint. 
\end{proof}

\begin{lem}
\label{trivial1}  Let $T$ be an operator from $L^p$ to $L^q$.  Let $u, v, u', v' $ are positive measurable  functions such that $u'\le u$ and $v'\le v$. Then, 
$\|u' T v'\|_{p\to q} \le  \|uT v\|_{p\to q}$. 
\end{lem}

\begin{proof}[Proof of Lemma \ref{local-l2}]  Instead of $\Pi_\lambda [\psi_j^\kappa]$ we consider the rescaled operator $\fP_\lambda[\psi_j^\kappa]$  and show
the bound
\[
     \| \chi_{\mu}^\sigma \fP_\lambda [ \psi_j^\kappa ]  \chi_{ \mu}^{\sigma'}    \|_{2\to 2} \lesssim \lambda^{-\frac d2}  2^j  \mu 
  \]
  for $\kappa= \pm, \ \pm \pi, \ 0,$   which is equivalent to \eqref{jl2} as is clear from \eqref{eq:norm-scaled}.  After  decomposing $A_\mu^\sigma$ and $A_{\mu'}^\spp$ into subsets of small diameter (e.g., \eqref{coarse-decomp}), by  Lemma \ref{large-sep} and Lemma \ref{small-sep} (also see Remark \ref{separated}) with $p=2$, 
   it is sufficient to show 
  \Be
  \label{eq:l2-scaled}
   \| \chi_{A}\fP_\lambda [ \psi_j^+]  \chi_{A'}   \|_{2\to 2} \lesssim 2^j  \mu \lambda^{-\frac d2}
  \Ee
  under the assumption that $A\subset A_\mu^\sigma$, $A'\subset A_{\mu}^\spp$ are of diameter $\sim c$, and $\dist(A,A')\le c$ for a small positive constant $c$.
     For the purpose we make use of the decomposition \eqref{decomp}.
    Using \eqref{decomp} with $\mu=\mu'$ and Lemma \ref{trivial1},  we  note
\Be
\label{eq:++}
  \|\chi_A \fP_\lambda [ \psi_j^+ ]  \chi_{A'}\|_{2\to 2}\lesssim \sum_{\nu\le \nu_\circ} \big\| \sum_{k\sim_\nu k'}\spxkl \fP_\lambda [ \psi_j^+ ]  \spxklp \big\|_{2\to 2} 
  \Ee
  where $2^{-\nu}\lesssim c$ and   $\supp \spxkl, \supp \spxklp$ are contained in  a $O(c)$-neighborhood of  $A$ and $A'$, respectively.  We show  the right hand side of the above is bounded by $C\lambda^{-\frac d2}2^j  \mu$ by splitting the sum over $\nu$ into  the two  cases 
\[ 2^{-\nu}\gg \mu,  \quad   \mu \gtrsim  2^{-\nu},\] separately. 

  We first consider the case $2^{-\nu}\gg \mu$.  Since $ 2^{-\frac d2 j}(d/ds)^k (\mathfrak a{\psi}_j)=O(2^{kj})$ and $2^{-\nu}\gg\mu$,
 from Lemma \ref{case-a} we have 
\[
| \int   e^{i\lambda \mathcal P(x,y,s)} (\mathfrak a \psi_j)(s) ds | \lesssim   2^{\frac{d-2}2j} \big( \lambda2^{j}2^{-2\nu}+1\big)^{-N} 
\]
if $(x,y)\in  \aask$, $\ksim$, and $2^{-\nu}\gtrsim\, 2^{-2j}$.   Combining  this estimate and Lemma \ref{trivial},  we see  
\[  \|\spxkl  \fP_\lambda [ \psi_j^+ ]  \spxklp\|_{2\to 2}\lesssim  2^{\frac{d-2}2j} (\lambda  2^j   2^{-2\nu})^{-N} \mu 2^{-(d-1)\nu}, \quad \ksim\]
 for any $N>0$ since $ |A_ k^{\pm,\nu}|, | A_{k'}^{\pm, \nu}|\sim \mu 2^{-(d-1)\nu}$. 
Taking particularly $N=d/2$    and  using Lemma \ref{kkpsum} we get $  \|\sum_{k\sim_\nu k'}\spxkl \fP_\lambda [ \psi_j^+ ]  \spxklp\|_{2\to 2}\lesssim  2^{-j } \mu 2^{\nu} \lambda^{-\frac d2} .$
 Hence, summation over $\nu: 2^{-\nu}\gtrsim\, 2^{-2j}$ gives
\[ \sum_{\nu: 2^{-\nu}\gtrsim\, 2^{-2j} }  \|\sum_{k\sim_\nu k'}\spxkl \fP_\lambda [ \psi_j^+ ]  \spxklp\|_{2\to 2}\lesssim  \lambda^{-\frac d2} 2^{j } \mu.\]
Thus, we need only to consider the sum over ${\nu: \mu\ll 2^{-\nu}\ll 2^{-2j} }$. Since  $2^{-\nu}\ll 2^{-2j} $, 
by  Lemma \ref{case-a} it follows that 
\Be
 \label{1st1st} \Big| \int   e^{i\lambda \mathcal P(x,y,s)} (\mathfrak a \psi_j^+)(s)ds  \Big| \lesssim   2^{\frac{d-2}2j} \big(\lambda  2^{-3j }+ 1\big)^{-N}
 \Ee
for  $(x,y)\in   \aask$, $\ksim$.
 As before,  using \eqref{1st1st}  with $N=d/2$,   Lemma \ref{trivial},  and Lemma \ref{kkpsum}, we have 
 $ \| \sum_{k\sim_\nu k'} \spxkl \fP_\lambda [ \psi_j^+ ]  \spxklp\|_{2\to 2} \lesssim \lambda^{-\frac d2} \mu  2^{(2d-1)j}  2^{-(d-1)\nu}$. Thus, 
 we obtain 
 \begin{equation}
 \label{eq:l2l2l2}   \sum_{ \mu \ll 2^{-\nu} \ll 2^{-2j}}   \| \sum_{k\sim_\nu k'} \spxkl \fP_\lambda [ \psi_j^+ ]  \spxklp\|_{2\to 2}
  \lesssim     \lambda^{-\frac d2} \mu  2^j .  \!\!\!\!\!
 \Ee
 
 {The remaining case  $2^{-\nu}\lesssim  \mu$ can be handled similarly.}  Since $2^{-2j}\gg \mu$, by $(b)$ in
 Lemma \ref{case-b} we have \eqref{1st1st}  for $(x,y)\in A_{k}^{\sigma,\nu}\times A_{k'}^{\spp,\nu}$. 
Using this estimate,   Lemma \ref{trivial}, and Lemma \ref{kkpsum} we get 
 \[
\sum_{2^{-\nu} \lesssim  \mu}   \| \sum_{k\sim_\nu k'}\spxkl \fP_\lambda [ \psi_j^+ ]  \spxklp\|_{2\to 2}
 \lesssim 2^{(2d-1)j} \lambda^{-\frac d2}\mu^{d}\lesssim     \lambda^{-\frac d2} \mu\,2^j.
 \] 
 Therefore, we obtain the desired estimate \eqref{eq:l2-scaled}. 
\end{proof}

\begin{lem} Let $0<\mu',  \mu\ll 1$ and $\sigma, \sigma'=\pm$. If $\mu'\le c\mu$ for a small enough $c>0$, then there are constants  $c_1<  c_2$ such that,   for $\kappa=\pm, \pm\pi$,   we have  
\label{l2-mmj2}  
\begin{align*}
\Big\| \sum_{j:2^{-j} \le c_1 \sqrt \mu} \,  \chi_{\mu}^\sigma \fP_\lambda [ \psi_j^\kappa ]  \chi_{\mu'}^\spp   \Big \|_{2\to 2} \lesssim    \lambda^{-\frac d2}  \mmf ,
\\
 \Big\| \sum_{j:2^{-j} \ge  c_2 \sqrt \mu} \,  \chi_{\mu}^\sigma \fP_\lambda [ \psi_j^\kappa ]  \chi_{\mu'}^\spp  \Big \|_{2\to 2} \lesssim    \lambda^{-\frac d2}   \mmf .      
  \end{align*}
\end{lem}

This lemma   plays a significant  role in obtaining the unbalanced improvement. We refer the reader forward to Lemma \ref{LL22}.

\begin{proof} 
We begin by noting that $\mu\gg \mmt$. Similarly as in Lemma \ref{local-l2}, using  a preparatory decomposition such as \eqref{coarse-decomp}, by  Lemma \ref{large-sep}, Lemma \ref{small-sep} (Remark \ref{separated}) with $p=2$, and \eqref{symmetric}, 
 it is sufficient to show that   
 \begin{align}
 \label{I-est}  
\sum_{j:2^{-j} \le c_1 \sqrt \mu} \big\|  \,  \chi_{A} \fP_\lambda [ \psi_j^+ ] \chi_{A'}   \big \|_{2\to 2} \lesssim    \lambda^{-\frac d2}  \mmf ,
\\
 \label{III-est}  
\sum_{j:2^{-j} \ge  c_2 \sqrt \mu}   \big\| \,  \chi_{A}  \fP_\lambda [ \psi_j^+ ]  \chi_{A'}    \big \|_{2\to 2} \lesssim    \lambda^{-\frac d2}   \mmf,      
  \end{align}
   while assuming $A\subset A_\mu^\sigma$, $A'\subset A_{\mu'}^\spp$ are of diameter $\le c$, and $\dist(A,A')\le c$.
Hence, using  the decomposition \eqref{decomp},  we have \eqref{eq:++} and  we may assume $2^{-\nu}\lesssim  c.$
We prove \eqref{I-est} and \eqref{III-est} by  separately considering the cases 
\[ 2^{-\nu}\gg \mu,  \quad   \mu\gtrsim  2^{-\nu}.\]
We first show \eqref{III-est}. Recalling \eqref{eq:++} in which  we may assume $A_k^{\sigma, \nu }\subset A$ and $A_{k'}^{\sigma', \nu }\subset A'$, we  take summation  over $\nu$ and consider the sum over $j$ afterward.

{\bf Sum over $\nu:2^{-\nu}\gg \mu$.}  From \eqref{a-case-kernel} with $N=d/2$ and Lemma \ref{trivial} 
it follows that 
\begin{equation}
\label{kernel:est-case-a}
\|\spxkl \mathfrak P_\lambda [ \psi_j ]  \spxklp\|_{2\to 2} 
          \lesssim  \begin{cases}   
                                \lambda^{-\frac d2} 2^{- j}  2^\nu \mmt,   
                                            &   2^{-2j} \lesssim 2^{-\nu} , 
                                 \\
                                   \lambda^{-\frac d2} 2^{(2d-1) j}  2^{-(d-1)\nu}  \mmt,            
                                           &    2^{-2j}\gg  2^{-\nu}.
                          \end{cases} 
\end{equation}
Splitting the sum $ \sum_{\nu:2^{-\nu}\gg \mu}$ into $ \sum_{\nu:2^{-\nu}   \gtrsim   2^{-2j}} +  \sum_{\nu:\mu \ll 2^{-\nu}   \ll   2^{-2j}}$,  by Lemma \ref{kkpsum} and the estimate \eqref{kernel:est-case-a} we get
\begin{align}
\label{III-a}
 \| \sum_{\nu:2^{-\nu}\gg \mu}\sum_{k\sim_\nu k'} \spxkl \mathfrak P_\lambda [ \psi_j ]  \spxklp\|_{2\to 2} 
  &  \lesssim    2^j\lambda^{-\frac d2} \mmt .         
 \end{align}
{\bf Sum over $\nu: \mathrm \mu\gtrsim  2^{-\nu}$.}    Choosing large enough $c_2$ such that $2^{-j}  \gg  \smu$,  from \eqref{b-case-kernel} with $N=d/2$ (i.e., $(b)$ in Lemma \ref{case-b})  and Lemma \ref{trivial} 
we have  the second case of \eqref{kernel:est-case-a}. Thus, 
by Lemma \ref{kkpsum} we get
\begin{align}
 \| \sum_{\nu: \mu \gtrsim 2^{-\nu}} \ \sum_{k\sim_\nu k'} \spxkl \mathfrak P_\lambda [ \psi_j ]  \spxklp\|_{2\to 2}
                   \label{III-b} &         \lesssim   \lambda^{-\frac d2} \mmt\, 2^{(2d-1) j} \mu^{d-1} .
        \end{align}
 Using the estimates \eqref{III-a} and \eqref{III-b} and taking summation along $j: 2^{-j}\gg \sqrt\mu$, we get  the desired estimate  \eqref{III-est} since $\mmt\,\sum_{2^{-j}\gg \smu}( 2^j+2^{(2d-1) j} \mu^{d-1}) \lesssim \mmf\,. $

{We now show  \eqref{I-est}}. The argument here is similar with the previous one and  hence we shall be brief.

{\bf Sum over  $\nu:2^{-\nu}\gg \mu$.}   In this case the estimate 
 \eqref{a-case-kernel} holds,  so we have \eqref{kernel:est-case-a} as before. Taking $c_1$ small enough, i.e., $2^{-j}\ll \smu$,   
we may use the first estimate in \eqref{kernel:est-case-a}. So, by Lemma \ref{kkpsum} we get
 \Be\label{I-a} 
 \begin{aligned}
 \| \sum_{\nu:2^{-\nu}\gg \mu} \sum_{k\sim_\nu k'} \spxkl \mathfrak P_\lambda [ \psi_j ]  \spxklp\|_{2\to 2} 
       & \lesssim  \lambda^{-\frac d2} 2^{- j}   \mmt  \sum_{\nu:2^{-\nu}\gtrsim \mu}    2^\nu 
       \\
       &\lesssim  \lambda^{-\frac d2} \mmt\, 2^{- j}\mu^{-1} .
        \end{aligned}
        \Ee
{\bf Sum over  $\nu: \mu\gtrsim  2^{-\nu}$.}  We  have the first case in \eqref{b-case-kernel}  because $2^{-j}\ll \smu$.  
Thus, taking  $N=d/2$  and  using Lemma \ref{trivial}, 
we have 
$\|\pxkl \mathfrak P_\lambda [ \psi_j ]  \pxklp\|_{2\to 2} 
          \lesssim   \lambda^{-\frac d2}    2^{-j}  \mu^{-d}  \mmt  2^{-(d-1)\nu} 
$
as long as  $ \mu\gtrsim  2^{-\nu}.$   Similarly as before, by this estimate and Lemma \ref{kkpsum}  we obtain 
\begin{align}
\label{I-b} \|  \sum_{\nu:\, \mu \gtrsim 2^{-\nu}}\, \sum_{k\sim_\nu k'} \pxkl \mathfrak P_\lambda [ \psi_j ]  \pxklp\|_{2\to 2} 
              \lesssim \lambda^{-\frac d2} \mmt   2^{-j}  \mu^{-1} .     
                \end{align}  
Finally, combining the estimates  \eqref{I-a} and  \eqref{I-b} and taking summation along $j$, we get 
\[
\sum_{j:2^{-j} \le c_1 \sqrt \mu} \big\|  \,  \chi_{A} \fP_\lambda [ \psi_j^+ ] \chi_{A'}   \big \|_{2\to 2} 
\lesssim   \lambda^{-\frac d2} \mmt\,   \sum_{2^{-j}\ll \smu }2^{-j}  \mu^{-1} \lesssim \lambda^{-\frac d2} \mmf, 
\]
which gives the desired estimate \eqref{I-est}.  
\end{proof}

\section{Unbalanced improvement: Proof of Theorem  \ref{endpoint} } 
\label{sec:endpoint}
In this section we obtain  the estimate   for  $\chi_\mu^\sigma \fP_\lambda  \chi_{\mu'}^{\sigma'}$   with $ \sigma, \sigma'\in \{+, -\}$ and $\mu'\le \mu\ll 1$. 
When $\mu'\ll  \mu$, as is already mentioned in the introduction,  there is an improvement in bound which can not be recovered  by the estimate such as \eqref{kt-nu}. 
This improvement in bound is crucial in obtaining the endpoint estimate  \eqref{intro-conj}.  For simplicity we set  
\[ \mathrm B_{p,q}(\lambda,\mu, \mu'):= \lambda^{\frac{d-1}2\dpq-\frac d2}  (\mm)^{(\frac14-\frac{d+3}8\dpq)}    .\]

\begin{thm} 
\label{improvedbound}  Let $d\ge 5$ and $0<\mu'\le\mu\ll 1$.  Suppose $0<\dpp<2/(d+1)$.  Then, for some $c>0$, 
we have the estimate 
 \Be
 \label{eq:main} \|  \chi_\mu^\sigma\fP_\lambda  \chi^{\sigma'}_{\mu'}  \|_{p\to p'} \lesssim  
\bpqp  \mom^{c}, \quad \sigma, \sigma^\prime\in \{+, -\}. 
\Ee
\end{thm}

Without the extra factor  $\mom^{c}$  
the estimate \eqref{eq:main} follows from  \eqref{kt-nu}. Once we have  \eqref{eq:main}, the proof of Theorem \ref{endpoint} is 
rather straightforward.

\begin{proof}[Proof of Theorem \ref{endpoint}] Let us set \[ p_\ast=\frac{2(d+3)}{(d+5)}.\] We begin with noting from \eqref{eq:norm-scaled}  that 
 the estimate \eqref{intro-conj}  is equivalent to 
\[\|\fP_\lambda \|_{ \ps\to 2}\lesssim \lambda^{\frac{d-1}2 \delta(p_\ast,{2})-\frac d2}.\]
Let $0<\delc\ll 1$ be a small  enough constant. To show the  endpoint estimate 
\eqref{intro-conj}  we break the operator $\fP_\lambda$ as follows:
\[
 \fP_\lambda =  \fP_\lambda \chi_{\mathbb B}  + \fP_\lambda \chi_{\delc}^\circ  +  \fP_\lambda \chi_{\mathbb  B'} , 
\]
where  $\mathbb B=\{x: |x|<  1-2\delc \}$ and $\mathbb  B'=\{x:  1+2\delc< |x|\}$. 
The desired estimate for $\fP_\lambda \chi_{\mathbb B} $  follows from the local estimate  in Proposition  \ref{away-sphere}. Indeed, 
taking summation over $j$ and $\kappa$,  we get  $ \| \chi_{\mathbb B} \fP_\lambda \chi_{\mathbb B}\|_{p\to p'}\lesssim \lambda^{\scaleppp}$ for 
$0\le \dpp< 2/(d+1)$ (also see Corollary \ref{local-r}).  This with $p=p_\ast$  gives the estimate 
$\| \chi_{\mathbb B}\fP_\lambda \chi_{\mathbb B}\|_{\ps\to \ps'}\lesssim  \lambda^{\frac{d-1}2 \delta(p_\ast, \ps')-\frac d2}$. 
Thus, we have the desired estimate   \[\|\fP_\lambda \chi_{\mathbb B}\|_{\ps\to 2}\lesssim \lambda^{\frac{d-1}2 \delta(p_\ast,{2})-\frac d2}\] via $TT^\ast$ argument.
\footnote{Here, we use the idenity $\inp{\fP_\lambda\chi_E f}{\fP_\lambda\chi_E f}=\lambda^{-d/2}\inp{\chi_E\fP_\lambda\chi_E f}{f}$.} 
The estimate for $\fP_\lambda \chi_{\mathbb B'} $ is easy. Since  the kernel of $\Pi_\lambda$ decays rapidly when $|x|,$ or $ |y|\ge (1+2\delc)\sqrt \lambda$ (for example, see \cite[Theorem 4]{K94} and Proposition \ref{prop:estext})  we have $ \| \chi_{\mathbb B'} \fP_\lambda \chi_{\mathbb B'}\|_{p\to q}\lesssim \lambda^{-N}$ for any $N>0$ if $1\le p\le 2\le q\le \infty$. 

The matter is now reduced to  showing 
$\| \fP_\lambda \chi_{\delc}^\circ  \|_{\ps\to 2}\lesssim \lambda^{\frac{d-1}2 \delta(\ps , 2 )-\frac d2}.
$
Via  $TT^\ast$ argument  this estimate is clearly equivalent 
to  
\[\|\chi_{\delc}^\circ \fP_\lambda \chi_{\delc}^\circ  \|_{p_\ast\to p_\ast'}\lesssim \lambda^{\scalepq{p_\ast}{p_\ast'}},\] which we show in what follows. Recalling $\mu, \mu'$ are positive  dyadic numbers,  we  may write 
 \[\chi_{\delc}^\circ \fP_\lambda \chi_{\delc}^\circ f=\sum_{\sigma, \sigma'\in \{+, -\}} \sum_{k} \sum_{\mu, \mu': \frac\mu{\mu'} \sim 2^k}  \chi_\mu^\sigma \fP_\lambda  \chi^{\sigma'}_{\mu'}  f.\]
   Since  $A_\mu^\pm$ are almost disjoint, 
   we have
\[
\|\sum_{\mu, \mu': \frac\mu{\mu'} \sim 2^k}   \chi_\mu^\sigma  \fP_\lambda  \chi^{\sigma'}_{\mu'}   f\|_{p_\ast'}
        \lesssim  \big(\sum_{\mu, \mu': \frac\mu{\mu'} \sim 2^k} \|   \chi_\mu^\sigma  \fP_\lambda  \chi^{\sigma'}_{\mu'}    f\|_{p_\ast'}^{p_\ast'}\big)^{1/{p_\ast'}} .
        \]
 Note that $  \chi^{\sigma'}_{\mu'}  \fP_\lambda \chi_\mu^\sigma$ is the adjoint operator of $\chi_\mu^\sigma \fP_\lambda  \chi^{\sigma'}_{\mu'} $. Thus, by 
 Theorem \ref{improvedbound} and duality  it follows that  $\|\chi_\mu^\sigma \fP_\lambda  \chi^{\sigma'}_{\mu'}\|_{p_\ast\to p'_\ast}\lesssim   \min(2^{ck}, 2^{-ck})  
                                   \lambda^{\scalepq{p_\ast}{p_\ast'}} .$
 Hence, we now get 
\begin{align*}
\|\sum_{\mu, \mu': \frac\mu{\mu'} \sim 2^k}  &  \chi_\mu^\sigma  \fP_\lambda  \chi^{\sigma'}_{\mu'}   f\|_{p_\ast'}
       \lesssim \min(2^{ck}, 2^{-ck})  
                                    \lambda^{\scalepq{p_\ast}{p_\ast'}}\big( \sum_{ \mu'} \|  \chi^{\sigma'}_{\mu'}  f\|_{p_\ast}^{p_\ast'}  \big)^{1/{p_\ast'}} . 
                \end{align*}
Since ${p_\ast'}\ge p_\ast$,  $(\sum_{ \mu'} \| \chi^{\sigma'}_{\mu'}   f\|_{p_\ast}^{p_\ast'} \big)^{1/p_\ast'}\lesssim \|f\|_{\ps}$. Therefore, taking  summation  over $k$ and then over  $\sigma,$ $\sigma'$ gives 
the desired estimate.  
\end{proof}

We now provide a brief overview on the proof of Theorem \ref{improvedbound}. 
Starting with \eqref{decomp-proj},  we obtain bounds on the operator norm of  $ \chi_{\lambda,\mu}^\pm \Pi_\lambda [\psi_j^\kappa] \chi_{\lambda,\mu'}^\pm$. In order to do so, we invoke the sectorial Whitney type decomposition \eqref{decomp} of $A_{\lambda,\mu}^\pm\times A_{\lambda,\mu'}^\pm$. This allows us to control efficiently the angular interactions and sort out the critical scale  of $j$ which is given by $2^{-j}\sim \smu$ when $\mu\ge \mu'$. For the other case $2^{-j}\not \sim \smu$,  the estimates for  $ \chi_{\lambda,\mu}^\pm \Pi_\lambda [\psi_j^\kappa] \chi_{\lambda,\mu'}^\pm$  are relatively  easier to show.   
If  $2^{-j}\sim \smu$,   the critical angular separation is given by  $2^{-\nu_\circ}<2^{-\nu}\sim \mmt$ and the required $L^1$--$L^\infty$ bound of $O(\lambda^{-\frac12})$ is no available. So,  we need to make additional decomposition away from the point $S_c$ as in Section \ref{mild-l2}. Though the decomposition yields the desirable  $L^1$--$L^\infty$ bound in $\lambda$, the $L^2$ estimate becomes highly nontrivial (e.g. Proposition \ref{critical-muj}). In order to prove the sharp $L^2$ estimate,  we need to analyze 
the associated oscillatory integral operator in detail. This shall be separately done in Section \ref{proof-end}. 

The rest of this section concerns the proof of Theorem \ref{improvedbound}. 

\subsection{Reduction and decomposition} 
 If $\mu\sim \mu'$,  
the estimate \eqref{eq:main}  is relatively easier to show.  In particular,  if $\mu\gtrsim \lambda^{-\frac 23}$ and $\sigma=\sigma'=+$,  for $0\le \dpp\le 2/(d+1)$  the estimate   \eqref{eq:main} follows from  
 \eqref{kt-nu}  and scaling  (also see \eqref{kt-nu-})  since  $ \|  \chi_\mu^\sigma\fP_\lambda  \chi^{\sigma'}_{\mu'}  \|_{p\to p'}\le  \lambda^{-d/2}\|\fP_\lambda \chi^{\sigma'}_{\mu'}  \|_{p\to 2}  \|  \chi_\mu^\sigma\fP_\lambda\|_{2\to p'}$.  In Theorem  \ref{thm-annest} (Proposition \ref{mm-comparable})   all these estimates are included.   The estimate \eqref{eq:main} for  the remaining case  $\mu\lesssim \lambda^{-\frac 23}$  follows from 
\eqref{smallmu} in Proposition \ref{prop:estext}.   Thus, we may assume that 
\[c\gg \mu\gg \mu'\] 
for a small constant $c>0$ and $\mu'/\mu$ is small enough. 
The other case  either  can be handled by Proposition  \ref{away-sphere} or by the rapid decay of the kernel of $\Pi_\lambda$ on the region
$\{(x,y): |x|>\sqrt\lambda(1+c) \text{ or } |y|>\sqrt\lambda(1+c)\}$  with some small $c>0$.  

We now recall \eqref{decomp-proj} and note that the same identity holds with $\Pi_\lambda$ replaced by $\fP_\lambda$.  After partitioning  the sets $A_\mu^\sigma$, $A_{\mu'}^{\sigma'}$ into finite number of subsets  $A_j$, $A_j'$ as in \eqref{coarse-decomp}, we apply Lemma \ref{large-sep} and  Lemma \ref{small-sep} (Remark \ref{separated}).  Therefore,  discarding the minor parts whose contributions are not significant,  we only need to consider  $ \chi_{A} \fP_\lambda \chi_{A'}$ while assuming that $A\subset A_\mu^\sigma$,  $A'\subset A_{\mu'}^{\sigma'}$ are of small diameter $\lesssim c$ and \eqref{distanced} holds.   Additionally, by the second identity in \eqref{symmetric} it is enough to consider the first case in \eqref{distanced}.  Therefore, the matter is now reduced to showing 
 \[
\|  \sum_{\kappa=\pm}\sum_{j\ge 4}\chi_{A} \fP_\lambda [\psi_j^\kappa]  \chi_{A'} \|_{p\to p'} \lesssim  
\bpqp \mom^c
\]
while $\dist (A, A')\le c$. Indeed, by Lemma \ref{small-sep} we have 
\Be
\label{eq:la2}
 \sum_{\kappa=\pm\pi, 0}\sum_j \|  \chi_{A } \fP_\lambda [\psi_j^\kappa]  \chi_{A'} \|_{p\to p'} \lesssim  
\lambda^{-N} (\mm)^{\frac12(1-\dpp)}
\Ee
for $1\le p\le 2$ because  $\dist (A, A')\le c$. So, we may  disregard $ \sum_{\kappa=\pm\pi, 0} \chi_{A } \fP_\lambda [\psi_j^\kappa]  \chi_{A'}$.

 Since  $\phi_\lambda(\mU x,\mU y,t)=\phi_\lambda(x,y,t)$, by rotation we may  additionally assume 
where
\[
A, A'\subset \mathbb  A:= \Big\{ x :   \Big  |\frac{x}{|x|}-e_1\Big|\le \frac 1{25d}\,\Big\}\,.
\]
We use  \eqref{decomp} to make a  decomposition. So, in order to show the desired estimate   it  suffices  to obtain
the estimate 
\Be
\label{eq:la}
\|  \sum_{\kappa=\pm}\sum_{j\ge 4}\sum_{ \nu\le \nu_\circ}\sum_{k\sim_\nu k'} \chi_{k}^{\sigma,\nu} \fP_\lambda [\psi_j^\kappa] \chi_{k'}^{\sigma',\nu}  \|_{p\to p'} \lesssim  
\bpqp\mom^c
\Ee
under the assumption that 
\Be 
\label{AA}
A_k^{\sigma,\nu}\times A_{k'}^{\sigma',\nu}\subset (\fA\cap A_\mu^\circ)\times (\fA\cap A_{\mu'}^\circ), \quad  2^{-\nu}\le 1/(50d).
\Ee

We prove Theorem  \ref{improvedbound} by establishing \eqref{eq:la}. We first show the estimate \eqref{eq:la} for $(\sigma, \sigma')=(+,+)$. The other cases 
$(\sigma, \sigma')\not=(+,+)$ are much easier  and we handle them at the end of this section.

\subsection{Proof of  \eqref{eq:la} with $(\sigma, \sigma')=(+,+)$}
The  {proof of Proposition \ref{away-sphere}} 
(Section \ref{proof-away})  heavily relies on the fact that  $S_c$ (and  other related quantities)  is a smooth function with uniformly bounded derivatives on 
the considered regions.  
However,  as  $\mu$, $\mu'$ get small,  this fact no longer remains valid on $A^{\pm}_{\mu}\times A^{\pm}_{\mu'}$, so  we can not obtain  the sharp estimates in a similar way. 
To get around this issue,  we use the decomposition \eqref{decomp} and this leads us to  handle the decomposed operators  by considering various  different cases, separately.  
In the course of doing so,   we specify the main part where $2^{-j}\sim \sqrt \mu$,  $2^{-\nu}\sim \mmt$,  and $|\cD(x,y)|<\vepc \mm$.  
The estimate for this core case  is much more involved, so we postpone its proof until Section \ref{proof-end}.

 We  show \eqref{eq:la}  by separately  considering the cases 
\[   2^{-\nu}\gg \mmt, \quad    \mmt \gtrsim  2^{-\nu} > 2^{-\nu_\circ}, \quad   \nu=\nu_\circ. \] 

\subsubsection{Sum over  $\nu:2^{-\nu}\gg\mmt$}  
\label{sec:large}
The desired estimate is relatively simpler to show and we  obtain it using   Lemma \ref{case-a} and Lemma \ref{case-b}. 
\begin{lem} 
\label{lala}  
Let $d\ge 2$ and $\mu'\ll \mu$.   If  ${2}/{(d+3)}<\dpp<  {2}/{(d-1)}$, then for some $c>0$  we have 
\Be 
\label{eq:lala}
\|  \sum_{\kappa=\pm}\sum_{j\ge 4} \sum_{ \mmt \ll 2^{-\nu}}\sum_{k\sim_\nu k'} \chi_{k}^{+,\nu} \fP_\lambda [\psi_j^\kappa] \chi_{k'}^{+,\nu} \|_{p\to p'} \lesssim  
\bpqp
\mom^{c}. 
 \Ee
\end{lem}

\begin{proof}  In order to make use of the  lower bounds on $\partial_s\cP$,  we  further divide  the sum $\sum_{ \mmt \ll 2^{-\nu}}$  considering the two cases:
\[ 2^{-\nu}  \gg   \mu , \quad   \mmt \ll  2^{-\nu}\lesssim \mu.   \] 
{\bf Case  $2^{-\nu}\gg   \mu$.}  
Using the estimate \eqref{eq:haha8} and Lemma \ref{trivial1}, we have 
\[ 
\|\chi_{k}^{+,\nu} \mathfrak P_\lambda [ \psi_j^\pm ]  \chi_{k'}^{+,\nu}\|_{p\to p'} 
          \lesssim  \begin{cases}   
                                 \lambda^{\frac{d-1}2\dpp-\frac d2} 2^{(\frac{d-1}2\dpp -1)j}  2^{\nu\dpp },  
                                      & \!\! 2^{-2j} 
\lesssim   2^{-\nu}                                  
\\
                                 \lambda^{\frac{d-1}2\dpp-\frac d2} 2^{(\frac {d+3}2\dpp -1) j},    
                                 &  \!\!  2^{-2j}  \gg 2^{-\nu} 
                          \end{cases} \hspace{-.4cm}
\]
for  $k\sim_\nu k'$.  Splitting  the sum over $j$ into two cases $j: 2^{-\nu}   \gtrsim   2^{-2j}$ and $j:2^{-\nu}   \ll   2^{-2j}$, 
we see $\sum_j \|\chi_{k}^{+,\nu} \mathfrak P_\lambda [ \psi_j^\pm ]  \chi_{k'}^{+,\nu}\|_{p\to p'}\le C\lambda^{\frac{d-1}2\dpp-\frac d2} 2^{(\frac{d+3}4\dpp-\frac12)\nu}$ if 
$\ksim$ since $\frac 2{d+3} <\dpp<  \frac{2}{d-1}$. 
Combining this estimate  with   Lemma \ref{kkpsum} gives  
\[
\|\sum_{k\sim_\nu k'} \sum_{j} \chi_{k}^{+,\nu} \mathfrak P_\lambda [ \psi_j ^\pm]  \chi_{k'}^{+,\nu}\|_{p\to p'}
\le C  \lambda^{\frac{d-1}2\dpp-\frac d2}    2^{(\frac{d+3}4\dpp-\frac12)\nu}. 
\]  Thus,   if ${2}/{(d+3)}<\dpp<  {2}/{(d-1)}$, summation over  $\nu: 2^{-\nu}\gg \mu$ yields 
\[
 \Big\| \sum_{\kappa=\pm} \sum_j\sum_{2^{-\nu} \gg \mu} \sum_{k\sim_\nu k'}\chi_{k}^{+,\nu} \mathfrak P_\lambda [ \psi_j^\kappa]  \chi_{k'}^{+,\nu} \Big\|_{p\to p'} 
  \lesssim    \bpqp \mom^{\frac{d+3}8\dpp -\frac14}.
 \]

{\bf Case $ \mmt \ll 2^{-\nu}\lesssim \mu$.} Let $c_1$, $c_2$ be the constants in Lemma \ref{l2-mmj2}  and  split 
\Be \label{three-sums}
\begin{aligned} 
& \,\,\,\sum_j  \sum_{k\sim_\nu k'}\chi_{k}^{+,\nu} \mathfrak P_\lambda [ \psi_j^\kappa]  \chi_{k'}^{+,\nu}\\
            =&\,\,\Big(\sum_{2^{-j}   \le c_1  \smu} +\sum_{   c_1  \smu< 2^{-j}   <c_2  \smu} +\sum_{ 2^{-j}  \ge c_2  \smu}  \Big) \sum_{k\sim_\nu k'}\chi_{k}^{+,\nu} \mathfrak P_\lambda [ \psi_j^\kappa]  \chi_{k'}^{+,\nu}
                       \\
            =&:\,\cI^\kappa(\nu)+\cI\!\cI^\kappa(\nu)+\cI\!\cI\!\cI^\kappa(\nu).
            \end{aligned}
            \Ee
Using the first case estimate in \eqref{eq:sec8} and Lemma \ref{kkpsum} and taking  sum  over $j$, we get  
\Be
\label{I-mid}
\begin{aligned}
\|\cI^\kappa(\nu)\|_{p\to p'}  &\lesssim   \lambda^{\frac{d-1}2\dpp-\frac d2}   \mu^{\frac12-\frac{d+3}4\dpp} 
\\
& = 
\bpqp \mom^{\frac{d+3}8\dpp-\frac14} 
\end{aligned}
\Ee
for $\kappa=\pm$ provided that  $\dpp< {2}/({d-1})$. If $\dpp> {2}{(d+3)}$, similarly using the third estimate in \eqref{eq:sec8} and  taking sum along $j$, we obtain  
\begin{align}
\label{III-mid}
\|\cI\!\cI\!\cI^\kappa(\nu)\|_{p\to p'}&\lesssim   
         \bpqp \mom^{\frac{d+3}8\dpp-\frac14}, \quad \kappa=\pm .
        \end{align}
Since $ \mmt \ll 2^{-\nu}\lesssim \mu$, using the second inequality in \eqref{eq:sec8},    
we get 
\Be
\label{II-mid}
\|\cI\!\cI^\kappa(\nu)\|_{p\to p'}
\lesssim  
\bpqp \mom^{\frac{d+1}8\dpp-\frac14}, \quad \kappa=\pm. 
\Ee
The current estimate \eqref{II-mid} for $\cI\!\cI^\kappa(\nu)$ is not desirable in that improvement of the bound  with the  factor $\mom^{c} $  is possible only if $\dpp>\frac2{d+1}$. 
However, we can enlarge the range of $p$  using  the following, which is a consequence of Lemma \ref{l2-mmj2}.

\begin{lem}  Let  $A\subset A_\mu^\sigma$,  $A'\subset A_{\mu'}^{\sigma'}$ be of diameter $\lesssim c$ and  satisfy $\dist (A, A')\le c$.  Then, for $\sigma, \sigma'=\pm,$
\label{LL22} 
we have 
\Be
\label{L2L2} 
\big\|  \sum_{\kappa=\pm}  \  \sum_{ j:  c_1  \smu< 2^{-j}   <c_2  \smu} \chi_{A} \fP_\lambda [\psi_j^\kappa]  \chi_{A'} \big\|_{2\to 2}\lesssim  \lambda^{-\frac d2} \mmf.
\Ee
\end{lem}
 
However, it is also possible to obtain the same bound on $\| \chi_{\mu}^\pm \fP_\lambda [ \psi_j^\kappa]  \chi_{\mu'}^\pm\|_{2\to2} $ additionally assuming that $\lambda^{-2/3}\lesssim\mu'\le \mu\le 1/4$, see Lemma \ref{l2-mmj23}. 

\begin{proof}[Proof of Lemma \ref{LL22}] From Lemma \ref{l2-mmj}, \eqref{eq:norm-scaled}, and Lemma \ref{trivial1}  it follows that  $\|  \chi_{A} \fP_\lambda \chi_{A'} \|_{2\to 2}\lesssim  \lambda^{-\frac d2} \mmf$. 
Thus, recalling \eqref{decomp-proj},  by \eqref{eq:la2}  we have 
\Be\label{L2222}
\|  \sum_{\kappa=\pm}\sum_{j\ge 4}\chi_{A} \fP_\lambda [\psi_j^\kappa]  \chi_{A'} \|_{2\to 2}\lesssim  \lambda^{-\frac d2} \mmf.
\Ee
From Lemma   \ref{l2-mmj2}  and \ref{trivial1} we see that both $\|\sum_{2^{-j}\ge c_2\sqrt \mu }\chi_{A} \fP_\lambda [\psi_j^\pm]  \chi_{A'} \|_{2\to 2}$ and  
$\|\sum_{2^{-j}\le  c_1\sqrt \mu }\chi_{A} \fP_\lambda [\psi_j^\pm]  \chi_{A'} \|_{2\to 2}$ are bounded by $C  \lambda^{-\frac d2} \mmf$. 
Therefore, the estimate \eqref{L2L2} follows. 
\end{proof} 

By Lemma \ref{LL22},  we particularly have $\|  \sum_{\kappa=\pm}\cI\!\cI^\kappa(\nu)\|_{2\to 2}\lesssim  \lambda^{-\frac d2} \mmf$. 
On the other hand, by \eqref{II-mid} we have  $\|  \sum_{\kappa=\pm}\cI\!\cI^\kappa(\nu)\|_{p\to p'}
\lesssim   \bpqp \mom^{c}$ for some $c>0$  if $\dpp>\frac2{d+1}$. We interpolate  these two estimates   to get  
\[
\|\sum_{\kappa=\pm}\!\cI\!\cI^\kappa(\nu)\|_{p\to p'}
\lesssim \bpqp \mom^{c}
\]
 for $\dpp >0$. Since $ \mmt \ll 2^{-\nu}\lesssim \mu $, there are only as many as  $O(\log (\mu/\mu'))$ $\nu$. 
 Combining the above  estimate and   \eqref{I-mid}, \eqref{III-mid} and taking sum along $\nu$, we get 
\[ 
 \|  \sum_{\kappa=\pm} \sum_j \sum_{\mmt \ll 2^{-\nu}\lesssim \mu}  \sum_{k\sim_\nu k'} \chi_{k}^{+,\nu} \mathfrak \fP_\lambda [\psi_j^\kappa] \chi_{k'}^{+,\nu} \|_{p\to p'} \lesssim  
 \bpqp   \mom^{c}   
\]
for some $c>0$ 
provided that  
${2}/{(d+3)}<\dpp<  {2}/{(d-1)}$. Therefore, we obtain  \eqref{eq:lala} combining this  and the estimate for the case  $2^{-\nu}\gg   \mu$. 
\end{proof}

We have dealt with the case $2^{-\nu} \gg \mmt$. We turn to the remaining cases  $\mmt \gtrsim  2^{-\nu} > 2^{-\nu_\circ}$ and $\nu=\nu_\circ$.

\subsubsection{Case $\nu=\nu_\circ$} 
\label{sec:simmm} We note  from \eqref{angle}  and \eqref{nu0}  that  $\cD(x,y)\sim \mm$ for $(x,y)\in A_{k}^{+,\nuc}\times A_{k'}^{+,\nuc}$ of  $k\sim_\nuc k'$. 

\begin{lem}
\label{mmu-det-large} Let $d\ge 2$. If 
$  {2}/{(d+3)}<\dpp<  {2}/{(d-1)},  then $ for some $c>0$ 
 we have 
 \Be 
\label{bbbbb}
\| \sum_{\kappa=\pm}\sum_j  \sum_{k\sim_\nuc k'} \chi_{k}^{+,\nuc} \fP_\lambda[\psi_j^\kappa] \chi_{k'}^{+,\nuc} \|_{p\to p'} \lesssim   \bpqp  \mom^{c}.  \Ee
 \end{lem} 
\begin{proof} 
We first note that the estimates in \eqref{eq:sec8} hold with $\sigma, \sigma'=+$, $\nu=\nuc$, and $k\sim_\nuc k'$.
Indeed, from \eqref{angle} and \eqref{nu0} we have  $\cD(x,y)\sim \mm$ for $(x,y)\in A_{k}^{+,\nuc}\times A_{k'}^{+,\nuc}$, $k\sim_\nuc k'$. 
By  Lemma \ref{case-b},  Lemma \ref{oscillatory}, and \eqref{symmetric}  the estimates in \eqref{eq:sec8} hold with  $\nu=\nu_\circ$,  $p=1$ and $q=\infty$. In fact,  
 the second case  in \eqref{eq:sec8} follows from  Lemma \ref{oscillatory}  since $\mm\sim  \cD(x,y)$.  As before, we interpolate those  estimates 
with the easy $L^2$ estimate from \eqref{easy-l2} to get  the estimates in \eqref{eq:sec8} with $\nu=\nuc$ and $k\sim_\nuc k'$.

Recalling \eqref{three-sums}, we need to show that  the norms $\|\cI^\kappa(\nuc) \|_{p\to p'}$, $\|\cI\!\cI^\kappa(\nuc) \|_{p\to p'}$, and $\|\cI\!\cI\!\cI^\kappa(\nuc) \|_{p\to p'}$ are bounded 
by    $C\bpqp  \mom^{c}$. 
Since we have  the estimates in \eqref{eq:sec8} with $\sigma, \sigma'=+$,  $\nu=\nuc$, and $k\sim_\nuc k'$, by the same argument as before  we have the estimates \eqref{I-mid}, \eqref{III-mid}, and \eqref{II-mid}. Hence,   we get 
\Be
\label{3+1}
\sum_{\kappa=\pm} \big(\|\cI^\kappa(\nuc) \|_{p\to p'}+\| \cI\!\cI\!\cI^\kappa(\nuc)\|_{p\to p'}\big)\lesssim \bpqp  \mom^{c}
\Ee
for some $c>0$ provided that  ${2}/{(d+3)}<\dpp<  {2}/{(d-1)}$.   
On the other hand, using the  second case estimate in \eqref{eq:sec8} with $\sigma, \sigma'=+$, $\nu=\nuc$,  we have
\[
\| \sum_{\kappa=\pm} \cI\!\cI^\kappa(\nuc) \|_{1\to \infty} 
          \lesssim      \mathrm B_{1,\infty}(\lambda, \mu, \mu')      \mom^{\frac{d+1}8\delta(1,\infty)-\frac14}.   
                  \]
                  Since we are assuming \eqref{AA},  using \eqref{L2L2}, 
 Lemma \ref{trivial1}, and Lemma \ref{kkpsum} successively, we obtain 
\[ \| \sum_{\kappa=\pm} \cI\!\cI^\kappa(\nuc) \|_{2\to 2}\lesssim  \mathrm B_{2,2}(\lambda, \mu, \mu').\]
Indeed, by Lemma \ref{kkpsum} it is sufficient to show   \[\|  \sum_{\kappa=\pm} \sum_{   c_1  \smu< 2^{-j}   <c_2  \smu} 
 \chi_{k}^{+,\nuc} \fP_\lambda[\psi_j^\kappa] \chi_{k'}^{+,\nuc}\|_{2\to 2}  \lesssim \mathrm B_{2,2}(\lambda, \mu, \mu')\] for $ {k\sim_\nuc k'}$. This in turn follows from 
 \eqref{L2L2} and Lemma \ref{trivial1}.  Interpolation between the above two estimates  for  $\sum_{\kappa=\pm} \cI\!\cI^\kappa(\nuc)$ 
gives
           \[  \| \sum_{\kappa=\pm} \cI\!\cI^\kappa(\nuc) \|_{p\to p'} 
          \lesssim      \bpqp     \mom^c\]  if $\dpp>0$.  This and  \eqref{3+1} yield  \eqref{bbbbb} for  some $c>0$ provided that  ${2}/{(d+3)}<\dpp<  {2}/{(d-1)}$. 
      \end{proof}

\subsubsection{Sum over $\nu: \mmt \gtrsim  2^{-\nu} > 2^{-\nu_\circ}$}  Since there are only as many as $O(1)$ $\nu$ we only have to consider a single $\nu$. 
\begin{prop}
\label{l0} 
Let $d\ge 5$ and $\mmt \gtrsim  2^{-\nu} > 2^{-\nu_\circ}$.  If $0<\dpp<2/(d-1)$, then we have  
\[
\| \sum_{\kappa=\pm}\sum_j  \sum_{k\sim_\nu k'} \chi_{k}^{+,\nu} \fP_\lambda[\psi_j^\kappa] \chi_{k'}^{+,\nu} \|_{p\to p'} \lesssim   \bpqp  \mom^{c}
\] for some $c>0$.
\end{prop}

Once we have Proposition \ref{l0}, we may  complete the proof of  the estimate \eqref{eq:la} with $\sigma, \sigma'=+$.

\begin{proof}[Proof of \eqref{eq:la} with $\sigma, \sigma'=+$]  
Putting together the estimates in Lemma \ref{lala}, Lemma \ref{mmu-det-large}, and Proposition  \ref{l0},
we have \eqref{eq:la} for some $c>0$  provided that   ${2}/{(d+3)}<\dpp<  {2}/{(d-1)}$. To extend the range of $p$ we interpolate the consequent  estimate with 
\[ \| \sum_{\kappa=\pm}\sum_j \sum_{ \nu\le \nu_\circ} \sum_{k\sim_\nu k'} \chi_{k}^{+,\nu} \fP_\lambda[\psi_j^\kappa] \chi_{k'}^{+,\nu} \|_{2\to 2} \lesssim   \mathrm B_{2,2}(\lambda, \mu, \mu'),\] 
which follows from \eqref{L2222},  Lemma \ref{trivial1}, and Lemma \ref{kkpsum}. 
 Therefore, we get    \eqref{eq:la} if $ 0<  \dpp <  \frac2{d-1}$.  
\end{proof}

In order to show Proposition \ref{l0} we first recall \eqref{three-sums}. From \eqref{angle} and \eqref{nu0} it is clear that 
$|\cD(x,y)|\lesssim\mm$ for $(x,y)\in A_{k}^{+,\nu}\times A_{k'}^{+,\nu}$ if  $k\sim_\nu k'$. Since $\mu'\ll\mu$,  $|\cD(x,y)|\ll \mu^2$, so 
by  Lemma \ref{case-b} we have  the estimate \eqref{b-case-kernel}. Thus,  the first and the third case estimates in \eqref{eq:sec8} continue to be valid.  
Taking sum over $j$  gives the estimates \eqref{I-mid} and \eqref{III-mid} for ${2}{(d+3)}<\dpp< {2}/({d-1})$.  On the other hand, by Lemma  
\ref{l2-mmj2}
we also have  
\[ \|\cI^\kappa(\nu)\|_{2\to 2} \lesssim \mathrm B_{2,2}(\lambda, \mu, \mu'), \quad  \|\cI\!\cI\!\cI^\kappa(\nu)\|_{2\to2}\lesssim \mathrm B_{2,2}(\lambda, \mu, \mu'). \]
Interpolating those estimates and  \eqref{I-mid} and \eqref{III-mid}  for ${2}{(d+3)}<\dpp< {2}/({d-1})$, respectively,  we get 
\[  \sum_{\kappa=\pm}  \|\cI^\kappa(\nu)\|_{p\to p'} + \sum_{\kappa=\pm}  \| \cI\!\cI\!\cI^\kappa(\nu)\|_{p\to p'} \lesssim \bpqp \mom^c \]
for some $c>0$ if $0<\dpp<2/(d-1)$.  Thus, it remains to show the estimate for $\cI\!\cI^\kappa(\nu)$. By Lemma \ref{kkpsum}, 
we only need to consider the  bound on the operator
$\chi_{k}^{+,\nu} \fP_\lambda  [ \psi_j^\kappa]  \chi_{k'}^{+,\nu}$ with $k\sim_\nu k'$. Thus, the proof of Proposition \ref{l0} is completed if we show  
\begin{equation}
\label{bbb011}
\| \sum_{\kappa=\pm} \sum_{   c_1  \smu< 2^{-j}   <c_2  \smu}  \chi_{k}^{+,\nu} \mathfrak P_\lambda [ \psi_j^\kappa]  \chi_{k'} \|_{p\to p'} \lesssim   \bpqp  \mom^{c}\
\end{equation} 
for $k\sim_\nu k'$. 

 Since the estimate is invariant under the transformation $(x,y)\to (\mathbf U x,\mathbf U y)$, $\mathbf U\in \mathrm O(d)$, by rotation  we may assume
\Be
\label{position0}
\begin{aligned}
A_{k}^{+,\nu}\subset \Big\{x: |x_1-1|\sim \mu, \,|\ol x|\lesssim  \mmt  \Big \},
\\
A_{k'}^{+,\nu} \subset \Big\{y: |y_1-1|\sim \mu', \, |\ol y|\lesssim  \mmt  \Big\}.
\end{aligned}
\Ee
Here $x=(x_1, \ol x)\in \mathbb R\times \mathbb R^{d-1}$.

For $x^*\in \mathbb R^d$ and $r>0$,  we denote 
\[\mathfrak R_{\mu'}^{\mu}(x^\ast, r)=\Big\{  (x_1, \ol x): |x_1-x_1^\ast|\le r \mu, \,, |\ol x-\ol{x}^\ast|\le r\mmt   \Big \}. \] 
If $(x,y)\in  A_{k}^{+,\nu}\times A_{k'}^{+,\nu}$,  unlike the previous cases,  $\cD(x,y)$ may  vanish.  
To handle this case, we  further decompose  $A_{k}^{+,\nu}\times A_{k'}^{+,\nu}$ to localize $\mathcal D$ tightly. 
To this end,  we break  $A_{k}^{+,\nu}\times A_{k'}^{+,\nu}$ into finitely many rectangles $\mathfrak R_{\mu'}^{\mu}(x_k, \epsilon)\times \mathfrak R_{\mu}^{\mu'}(y_l, \epsilon)$  so that  $ |\mathcal D(x,y)|< \vepc\mu\mu'$, or $ |\mathcal D(x,y)|\ge \vepc\mu\mu'$ holds on each of those rectangles.  
This particular form of the rectangles shall be important for our argument in Section \ref{proof-end}. Fortunately, this can be achieved by the following Lemma 
 if we take a sufficiently small $\epsilon>0$.

\begin{lem}\label{dychotomy} Let  $A_{k}^{+,\nu}$ and $A_{k'}^{+,\nu}$ satisfy \eqref{position0}.   Suppose 
$(x,y), (x^\ast, y^\ast) \in A_{k}^{+,\nu}\times A_{k'}^{+,\nu}$ such that  
$x \in \mathfrak R_{\mu'}^{\mu}(x^\ast, \epsilon)$ and $y\in \mathfrak R_{\mu}^{\mu'}(y^\ast, \epsilon)$. 
Then,  
\Be
\label{det-control} |\mathcal D(x,y)-\mathcal D(x^\ast,y^\ast)|\lesssim \epsilon \mm.
\Ee
\end{lem}

\begin{proof} 
We denote $z=(z_1, \dots, z_{2d})=(x,y)$ and  $z^*=(z_1^*, \dots, z_{2d}^*)=(x^\ast, y^\ast)$. Then, set 
$\mathcal D_0=\mathcal D(z)$ and 
\[ \mathcal D_k=\mathcal D(z_1^*,\dots, z_k^*, z_{k+1}, \dots, z_{2d}), \quad 1\le k\le 2d.\]
So $\mathcal D_{2d}=\mathcal D(z^*).$ Since $\D(x,y)-\D(x^\ast,y^\ast)=\mathcal D_0-\mathcal D_{2d}$, we have
$
 \D(x,y)-\D(x^\ast,y^\ast)=\sum_{k=0}^{2d-1}\mathcal D_k-\mathcal D_{k+1}.
$
Thus, it is sufficient to show 
\[ \mathcal D_k-\mathcal D_{k+1}= O(\epsilon \mu'\mu), \quad   1\le k\le 2d. \]

We note   $|x_1-x_1^\ast|\le \epsilon \mu$, $|y_1-y_1^\ast|\le \epsilon \mu'$, 
and $|x_j-x_j^*|, |y_j-y_j^*|\le \epsilon \mmt$, $j\ge 2$, and   
 \begin{align*}
 \mathcal D(x,y)=1+\sum_{j,l=1}^d (x_jy_j)(x_ly_l) -x_1^2-\dots-x_d^2-y_1^2-\dots-y_d^2\,.  
 \end{align*} 
Thus, $\mathcal D_0-\mathcal D_1= (x_1^\ast-x_1)(x_1^\ast+x_1)(1-y_1^2)+ 2(x_1-x_1^*)y_1\sum_{l\neq 1}x_ly_l$. So, we get $\mathcal D_0-\mathcal D_1  =O(\epsilon \mu'\mu)
$ 
 because $|1-y_1|\lesssim \mu'$ and $|x_l|,|y_l|\lesssim \sqrt{\mu\mu'}$ for $l=2,\cdots, d.$ Similarly, we have 
$ \mathcal D_d-\mathcal D_{d+1}= (y_1-y_1^\ast)(y_1+y_1^\ast)(1-x_1^2)+ 2x_1^\ast(y_1-y_1^\ast)\sum_{l\neq1} x_l y_l=O(\epsilon \mu'\mu).$
 When $1\le k\le d-1$, we note that    $ \mathcal D_k-\mathcal D_{k+1}$ is equal to 
\[
(x_j^2-(x_j^\ast)^2)y_j^2-(x_j^2-(x_j^\ast)^2)+2(x_j-x_j^\ast)y_j^\ast\sum_{l\neq j} x_l y_l.
\]
Similarly, when $d+1\le k\le 2d-1$,  $\mathcal D_k-\mathcal D_{k+1}$ equals 
$(y_j^2-(y_j^\ast)^2)x_j^2-(y_j^2-(y_j^\ast)^2)+2x_j^\ast(y_j-y_j^\ast)\sum_{l\neq j} x_l y_l.
$ 
Thus, we see  $ \mathcal D_k-\mathcal D_{k+1}=O(\epsilon \mu'\mu)$ for $1\le k\le d-1$ and $d+1\le k\le 2d-1$.  
 Therefore, we get \eqref{det-control}. 
\end{proof}

Let $\epsilon=c\vepc$ with a sufficiently small $c>0$. 
Let $\{B\}=\{\mathfrak R_{\mu'}^{\mu}(x_k, \epsilon): 1\le k\le K\}$, $\{B'\}=\{\mathfrak R_{\mu}^{\mu'}(y_l, \epsilon): 1\le l\le L\}$ be collections of  almost disjoint rectangles  
which  cover  $A_{k}^{+,\nu}$, $ A_{k'}^{+,\nu}$, respectively.  Then, taking $c$ small enough, by Lemma \ref{dychotomy} 
we may assume that one of the following holds:
\begin{align}
\label{large-d}
  &|\mathcal D(x,y)|\gtrsim \vepc \mm, \  \  \  \forall (x,y)\in   B\times B',  
   \\ 
   \label{small-det}
     \ & |\mathcal D(x,y)|\ll \vepc \mm, \  \  \  \forall (x,y)\in   B\times B'.     
  \end{align} 
  In order to show \eqref{bbb011} (Lemma \ref{trivial1}),  we may clearly assume 
\Be 
\label{bb-decomp} A_{k}^{+,\nu}\times A_{k'}^{+,\nu}=   \bigcup B\times B'.  \Ee
Thus, the matter is reduced to obtaining the estimate, for $p$ satisfying $0<\dpp<2/(d-1)$, 
\begin{equation}
\label{bbb01}
\| \sum_{\kappa=\pm} \sum_{   c_1  \smu< 2^{-j}   <c_2  \smu}  \chi_{B} \fP_\lambda[\psi_j^\kappa] \chi_{B'} \|_{p\to p'} \lesssim   \bpqp  \mom^{c}\
\end{equation} with some $c>0$ 
while $B$, $B'$ satisfy either \eqref{large-d} or \eqref{small-det}. 
We can handle the first case in the same manner as in the proof of Lemma \ref{mmu-det-large}. 
Since $|\mathcal D(x,y)|\sim  \mm$, by Lemma \ref{oscillatory} and \eqref{symmetric}  we have 
\[
 \| \sum_{\kappa=\pm} \sum_{   c_1  \smu< 2^{-j}   <c_2  \smu}  \chi_{B} \fP_\lambda[\psi_j^\kappa] \chi_{B'} \|_{1\to \infty}\lesssim  \mathrm B_{1,\infty}(\lambda, \mu, \mu')      \mom^{\frac{d+1}8\delta(1,\infty)-\frac14}.   
 \]  
 On the other hand, by \eqref{L2L2}  and Lemma \ref{trivial1}  we have
\[ \| \sum_{\kappa=\pm} \sum_{   c_1  \smu< 2^{-j}   <c_2  \smu}  \chi_{B} \fP_\lambda[\psi_j^\kappa] \chi_{B'} \|_{2\to 2}\lesssim \mathrm B_{2,2}(\lambda, \mu, \mu').\]
Interpolation between these two estimates gives \eqref{bbb01} if $0<\dpp<2/(d-1)$.  
The estimate \eqref{bbb01} for the other case \eqref{small-det} is immediate from the next proposition.

\begin{prop} \label{critical-muj}   Let $d\ge 5$. 
Let $2^{-j}\sim \smu$ and $2^{-\nu}\sim \mmt $.  Suppose 
\eqref{position0} and \eqref{small-det} holds for $(x,y)\in B\times B' $. Then, if $\frac{4}{3d+4}<\dpp<  \frac23$, there is a constant $c>0$ such that 
\[
  \|\chi_{B} \mathfrak P_\lambda [ \psi_j ]  \chi_{B'}\|_{p\to p'} \lesssim  \bpqp \mom^c .
\]
Moreover, $(i)$ we have  restricted weak type estimate  for  $\chi_{B} \mathfrak P_\lambda [ \psi_j ]  \chi_{B'}$ if $\dpp=2/3$ and $\mu=\mu'$ for $d\ge 2$.
\end{prop}

 The  statement ($i$) is not to be in use  in this section but we make use of  it later to prove Theorem \ref{thm-annest}.

 Proposition \ref{critical-muj} implies \eqref{bbb01} while \eqref{small-det} holds for $(x,y)\in B\times B' $ provided that $\frac{4}{3d+4}<\dpp<  \frac23$. 
 Since  we have  \eqref{L2L2}, the range is extended via interpolation in the same manner as before so that  \eqref{bbb01}  holds if  $0<\dpp<  \frac23$. 
 To complete the proof of  Proposition \ref{l0} it now remains to show Proposition \ref{critical-muj}.  

However, the proof of Proposition \ref{critical-muj}  is more involved. We postpone it to the next section. Instead, 
we close this section by providing the proofs of \eqref{eq:la} with $(\sigma,\sigma')\neq(+,+)$.

\begin{rem} When $d=1$ our argument fails.  Especially, the estimate in Proposition \ref{critical-muj} is not viable.  
 Under the assumption of Proposition \ref{critical-muj}  we have   
$\cD(x,y)=(1-x^2)(1-y^2)\sim \mu\mu'$ for $(x,y)\in B\times B'$ because $d=1$.  So,  from Lemma \ref{oscillatory} it follows
$\|\chi_{B} \mathfrak P_\lambda [ \psi_j ]  \chi_{B'}\|_{1\to \infty} \lesssim  \lambda^{-\frac12}(\mm)^{-1/4}$. On the other hand, we have $\|\chi_{B} \mathfrak P_\lambda [ \psi_j ]  \chi_{B'}\|_{2\to 2} \lesssim (\mm)^{1/4}$ by Lemma \ref{l2-mmj23}.  Interpolation between these two estimates only gives
$
 \|\chi_{B} \mathfrak P_\lambda [ \psi_j ]  \chi_{B'}\|_{p\to p'} \lesssim   \bpqp
$ 
for $1\le p\le 2$.
 However, the additional factor  $\mom^c$ is not  available.  \end{rem}

\subsection{Proof of \eqref{eq:la} with $(\sigma,\sigma')\neq(+,+)$}   
In what follows we show the estimate \eqref{eq:la} for the remaining cases but thanks to a favorable bound on the oscillatory integral we don't need the sophisticated argument  used  for the case  $\sigma,\sigma'=+$.   We show the bounds considering  the cases  
\[ 2^{-2\nu}\gg \mm,  \quad 2^{-2\nu}\lesssim  \mm,\]
 separately.

The case $2^{-2\nu}\gg \mm$  can be handled without changing the argument used in Section \ref{sec:large} since the estimates in Lemma \ref{case-a}, and Lemma \ref{case-b} remain valid for the case  $(\sigma,\sigma')=(\pm, -), (-, +).$ 
Indeed, we may use the estimates \eqref{eq:haha8}, \eqref{eq:sec8}, and \eqref{L2L2}.   The rest of argument is exactly same as in the proof of \eqref{eq:lala} with  $\sigma,\sigma'=+$, so we omit the details.

For the remaining case $2^{-2\nu}\lesssim  \mm$,  we need only to show 
\Be
\label{eq:rem}
\|  \sum_{\kappa=\pm}\sum_{j\ge 4}  \sum_{k\sim_\nu k'} \chi_{k}^{\sigma,\nu} \fP_\lambda \chi_{k'}^{\sigma',\nu} \|_{p\to p'} \lesssim  
\bpqp \mom^c
\Ee
because there are only finitely many $\nu$. We first consider the estimate \eqref{eq:rem} with $\sigma=-, \sigma'=\pm$. 
Since $\mm\ll \mu^2$, we  have   the second case of  \eqref{b-case-kernel} if $2^{-j}\gg\smu$,  and otherwise  we may use $(a2)$ in Lemma \ref{case-b}. 
Also, see $(\mathbf I\!\mathbf I)$ in Figure \ref{fig:2}.
Taking $N=1/2$ in these estimates gives  $L^1$--$L^\infty$ estimates. Interpolating them with  $ \|\chi_{k}^{-,\nu} \mathfrak P_\lambda [ \psi_j^\kappa ]  \chi_{k'}^{\pm,\nu}\|_{2\to 2}\lesssim \lambda^{-\frac d2}2^{-j}$ which follows from \eqref{easy-l2} and scaling \eqref{eq:norm-scaled},    
we have,   for $k\sim_\nu k'$ and  $\kappa=\pm$,  
\[       \|\chi_{k}^{-,\nu} \mathfrak P_\lambda [ \psi_j^\kappa ]  \chi_{k'}^{\pm,\nu}\|_{p\to p'}          
                \lesssim  \begin{cases}                                   
                      \lambda^{\frac{d-1}2\dpp-\frac d2} 2^{(\frac{d-1}2\dpp -1)j}  \mu^{-\dpp} ,    & 2^{-j}   \lesssim   \smu,         
                                                                                                                                 \\                        
                        \lambda^{\frac{d-1}2\dpp-\frac d2} 2^{(\frac {d+3}2\dpp  -1) j} ,         & 2^{-j}\,   \gg   \smu,   
\end{cases} 
\]
for  $1\le p\le 2$.  Splitting the cases  $2^{-j}\lesssim \smu$ and $2^{-j}\gg \smu$ and using these estimates,  we get 
\[\sum_j \max_{\ksim} \|\chi_{k}^{-,\nu} \mathfrak P_\lambda [ \psi_j^\kappa ]  \chi_{k'}^{\pm,\nu}\|_{p\to p'} \lesssim \bpqp \mom^{c}, \quad \kappa=\pm, \] 
 if ${2}/{(d+3)}<\dpq<  {2}/{(d-1)}$.  Hence, by Lemma \ref{kkpsum} and \eqref{symmetric}  we have \eqref{eq:rem}. This estimate can be interpolated with \eqref{L2222} to give 
\eqref{eq:rem}  with $\sigma=-, \sigma'=\pm $ for $0<\dpq<2/(d-1)$.

Finally, the estimate \eqref{eq:rem} with $(\sigma, \sigma')=(+, -)$ is easier  because $-\cD(x,y)\sim \mm$ for $(x,y)\in A_k^{+,\nu}
\times A_{k'}^{-, \nu}$. In fact, one may repeat the same lines of argument in Section \ref{sec:simmm} to show the desired estimate 
for  $0<\dpp<2/(d-1)$. We omit the details.

\section{Proof of Proposition \ref{critical-muj}} 
\label{proof-end}

In this section we prove Proposition \ref{critical-muj}. We recall \eqref{position0} and \eqref{bb-decomp}.  
After additional decomposition, rotation, and reflection \footnote{To make it hold that $ |x_2|\ll \mmt$ and $ y_2\sim   \mmt,$ we may use  $(x_2, y_2)\to -(x_2, y_2)$ 
for the kernel of the operator $\chi_{B} \mathfrak P_\lambda [ \psi_j ]  \chi_{B'}$. }
we may assume 
\Be
\label{positionor}
\begin{aligned} 
B \subset & \big\{x: |x_1-1|\sim \mu,  \ \  |x_2|\ll \mmt,  \ \  |\wt x|\lesssim  \mmt  \, \big\},
\\ 
B'\subset &\big\{x: |y_1-1|\sim \mu',   \ \ \ y_2\sim   \mmt, \ \ \ |\wt y|\lesssim  \mmt \,  \big\}
\end{aligned}
\Ee
because the assumption \eqref{small-det} and \eqref{angle} implies $\theta(x,y)\sim \mmt$ for $(x,y)\in B\times B'$.

\subsection{Additional decomposition of $B$}

\begin{figure}
\begin{tikzpicture}[scale= .9]
\draw (-2,0)--(0,0); 
\draw (0.1,0)--(1.9,0);
\node  at (0.1,0) {${\mathlarger{ \mathlarger \wr}}$};
\node  at (0.0,0) {${\mathlarger{ \mathlarger \wr}}$};

            \draw[->] (-2,3)--(-2,3.5) node [above] at (-2,3.5) {$\bar x$};
            \draw[->] (-2,3)--(-1.5,3) node [right] at (-1.5,3) {$x_1$};

\draw [->] (6.3,0)--(11,0); \node[below] at (11,-0.05) {$x_1$};

\draw (10.5,-0.05)--(10.5,0.05); \node[below, scale=0.8] at (10.5,-0.05) {$1$};
\draw [dotted,line width=0.6, domain=-5.45:13] plot ({17.57*cos(\x)-8}, {17.57*sin(\x)});
\draw [dotted,line width=0.6, domain=-5.5:13] plot ({17.15*cos(\x)-8}, {17.15*sin(\x)});
\draw [dotted,line width=0.6, domain=-6:13] plot ({14.27*cos(\x)-8}, {14.27*sin(\x)});
\draw [dotted,line width=0.6, domain=- 7:13] plot ({9.95*cos(\x)-8}, {9.95*sin(\x)});
\draw[thick] (1.9,-0.4)--(6.3,-0.4); \draw [thick](1.9,0.4)--(6.3,0.4);  
\draw[<->,>=stealth] (2,0.7) -- node[midway, above] {{\tiny $\sim$}$\mu $} (6.1,0.7);
\draw[thick] (1.9,-0.4) --(1.9,0.4); \node at (2.3,0) {$Q_1$};
\draw[thick] (2.7,-0.4) --(2.7,0.4); \node at (3.1,0) {$Q_2$}; 
\draw[thick] (3.5,-0.4)--(3.5,0.4); \node at (4.3,0) {$\cdots$};
\draw[thick] (5.5,-0.4)--(5.5,0.4); \node at (5.9,0) {$Q_n$};
\draw[thick] (6.3,-0.4)--(6.3,0.4);
\draw[|<->|,>=stealth] (2.7,-0.6) -- node[midway, below, scale=0.9] {$\varepsilon_\circ \sqrt{\mu\mu'} \,$} (3.5,-0.6);
\draw[|<->|,>=stealth] (6.5,-0.4) -- node[midway, above right, scale=0.9] {$\varepsilon_\circ \sqrt{\mu\mu'} $} (6.5,0.4);
\draw[thick] (9,1.4) -- (9.5,1.4); \draw[thick] (9,2.2) -- (9.5,2.2); 
\draw[thick] (9,1.4)--(9,2.2); \draw[thick] (9.5,1.4)--(9.5,2.2);
\node at (9.25,1.78) {$B'$};
\draw[|<->|,>=stealth] (9,1.2) -- node[midway, below, scale=.9] { \ \  {\tiny $\sim$}$\mu'$} (9.5, 1.2);
\draw[|<->|,>=stealth] (9.7,1.4) -- node[midway,  right, scale=0.9] { $\varepsilon_\circ \sqrt{\mu\mu'} $} (9.7,2.2);
\end{tikzpicture}
\caption{The cubes $Q_1, \dots, Q_n$,  and $B'$.}
\label{fig:cubes}
\end{figure}
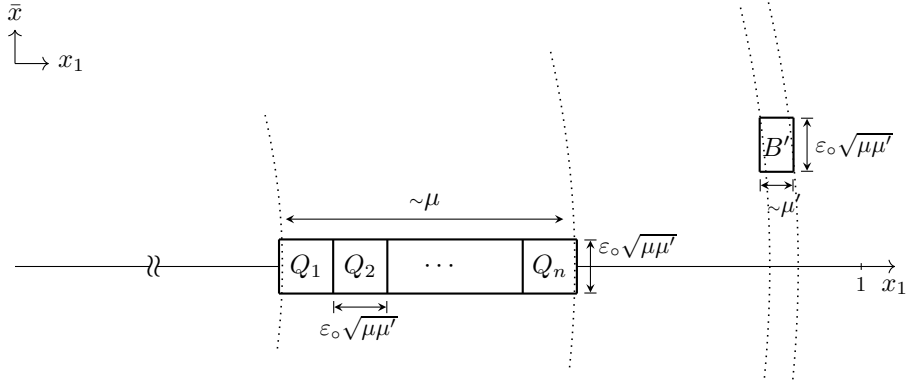

We further decompose $B$ into smaller cubes 
$Q_1,\,\dots\,, Q_n$,  $n\lesssim  \sqrt\frac{\mu}{\mu'}$ of  side length $\vepc \mmt$ so that 
\[ B=\bigcup_{k=1}^n Q_k .\] 
(See Figure \ref{fig:cubes}.)
From \eqref{positionor} and  \eqref{sc}, 
we see  that 
\[
 S_c(Q_1\times B'),\,  \dots \,, S_c(Q_n\times B') 
\subset [c_2 \sqrt {\mu}, c_1 \sqrt {\mu}]
\]
for some constants $c_1>c_2>0$ and each of them is again contained in boundedly overlapping  intervals  $I_k$ of length about $\vepc \smp$.    Indeed, 
since  $(x,y)\to \cos \sxy(x,y)$  is Lipschitz continuous as is clear from \eqref{sc}, 
$  |\cos S_c(x,y)- \cos S_c(x^*,y^*)|\lesssim  \vepc\sqrt{\mu\mu'}$ for all $(x,y),  (x^*,y^*) \in  Q_k\times B'. $
From Lemma \ref{scs}  
we note that $S_c(x,y)\sim  \sqrt{ \mu}$ because $|x+y|+|x-y|\gtrsim l$ if $(x,y)\in B\times B'$.
 Since  $s\sim \smu$ on the support of $\psi_j$, by  the above inequality  and the mean value theorem  we have 
\begin{equation}
\label{eps-con}
|S_c(x,y)- S_c(x^*,y^*)|\lesssim  \vepc\smp
\end{equation}
whenever $(x,y)$ and $(x^*,y^*)$ are in $Q_k\times B'$.

We decompose $\chi_{B} \mathfrak P_\lambda  \chi_{B'}=\sum_k\chi_{Q_k} \mathfrak P_\lambda  \chi_{B'} $. To control each term 
we may  simply apply the triangle inequality. However,  this leads to large loss in bound, so we exploit a weak localization property which we obtain 
using estimate for the  associated oscillatory integral.

Let $\fc(Q_k\times B')$ be the center of the set  $Q_k\times B'$, $k=1,\dots, n$, and for a large constant $C>0$ let us set 
\begin{align*}
  \psi_{Q_k}^0&= \widetilde \psi_0\Big(\frac{s-S_c(\fc(Q_k\times B'))}{C\vepc  \smp}\Big),
  \\ 
\qquad \psi_{Q_k}^m&=  \widetilde \psi\Big(\frac{s-S_c(\fc(Q_k\times B'))}{C\vepc  2^m\smp}\Big), \quad  m\ge 1,  
\end{align*} 
where   $\widetilde \psi_0=\sum_{j\le 0} \widetilde\psi(2^{-j}\cdot)$ and $\widetilde \psi=\psi(|\cdot|)$. 
We  clearly have
\begin{align*}  
\chi_B
\fP_\lambda [ \psi_j ]  \chi_{B'} = \sum_{m\ge 0}\sum_{k}  \mathfrak P^m_k 
                        := \sum_{m\ge 0} \sum_k  \chi_{Q_k}\fP_\lambda[\psi_j \psi_{Q_k}^m]\chi_{B'}. 
                                            \end{align*}
In what follows we may assume  that $2^m \sqrt {\mu'} \lesssim \smu  $ because $ \psi_j\psi_{Q_k}^m=0$ otherwise due to the support property of $\psi$. 
The next lemma control the contribution from  $\sum_{m\ge 1}\sum_{k}  \mathfrak P^m_k$.
\begin{lem}  
\label{m-m}  Let  $d\ge 2$. If $\dpp> \frac{4}{3d+4} $, for some $c>0$ we have 
\Be \label{m1} \|\sum_{m\ge 1} \sum_k \mathfrak P^m_k\|_{p\to p'} \lesssim   \bpqp \mom^c.
\Ee
\end{lem}

Refining the argument further, it is also possible to extend the range $p$ where the estimate  \eqref{m-m} holds. However, we do not pursue the matter since it does not give any improvement of the consequent result.  

\begin{proof}  We note that $|s-S_c(\fc(Q_k\times B'))|\sim 2^m \smp $ implies $|\cos s-\cos S_c(\fc(Q_k\times B'))|\sim 2^m \mmt $ for $m\ge 1$. 
Since $\cos S_c(x,y)$ is a root of  the equation $ \mathcal{Q}(x,y, \tau)=0$, $ | \inp xy-\cos S_c(x,y)|^2= {|\cD(x,y)|}.$ Combining these observations with  \eqref{q-def} and $\eqref{small-det}$, 
 we have $\cQ(x,y,\cos s)\sim 2^{2m} \mm$  if $(x,y)\in Q_k\times B'$ and $s\in  \supp \psi_{Q_k}^m$. 
Thus, from    \eqref{ph-d}   we see 
\[
|\partial_s \mathcal P| \gtrsim   (2^m\mmt)^2  \mu^{-1}\sim   2^{2m} \mu'
\]
on the support of $\psi_{Q_k}^m$. Since $2^m \sqrt {\mu'} \lesssim \smu  $,  we also have $|\partial_s^k \mathcal P|\lesssim  2^{2m} \mu' (2^{m} \sqrt{\mu'} )^{1-k}$  on $\supp\psi_{Q_k}^m$and 
$ |(d/ds)^k(\psi_j \psi_{Q_k}^m)|\lesssim  (2^{m} \sqrt{\mu'} )^{-k}$.  
Integration by parts  (Lemma \ref{S3-osc}) gives the estimate 
\begin{equation}
\label{haha}
|(\chi_{Q_k}\fP_\lambda[\psi_j \psi_{Q_k}^m] \chi_{B'})(x,y)|
\lesssim  \chi_{Q_k}(x)\chi_{B'}(y) 2^m\smp  \mu^{-\frac d4}\, \big(1+ \lambda 2^{3m} \mup^\frac 32\big)^{-N}
  \end{equation}
  for any $N$.  Since  $Q_k$ are essentially disjoint,    taking $N=1/2$  we get 
 \begin{equation}
 \label{onetoinfty} 
  \| \sum_k \mathfrak P^m_k 
 \|_{1\to \infty}
         \lesssim   \lambda^{-\frac12} 2^{-\frac12 m} (\mm)^{-\frac{d+3}8+\frac14}\mom^{\frac{d-1}8}
         \end{equation}
          for $m\ge 1$. 
Since $|Q_k|\sim (\mm)^\frac d2$ and $|B'|\sim  \mu' (\mm)^\frac {d-1}2 $, taking $N=d/2$ in \eqref{haha}  and using Lemma \ref{trivial} yield 
\begin{align*}
 \|  \chi_{Q_k}\fP_\lambda[\psi_j \psi_{Q_k}^m] \chi_{B'} \|_{2\to 2} 
                             &\lesssim   \lambda^{-\frac d2}  2^{m(1-\frac{3d}{2})} \mup^{\frac34- \frac{d}{4}}\mu^{\frac d4-\frac14}. 
\end{align*}

On the other hand, since the cutoff function $\psi_j \psi_{Q_k}^m$ is independent of  $x,y$ and is supported in an interval of length $2^m\smp$, by the isometry of the propagator  
(i.e., \eqref{easy-l2}) and scaling we have 
$\|  \chi_{Q_k}\fP_\lambda[\psi_j \psi_{Q_k}^m] \chi_{B'} \|_{2\to 2}\lesssim  \lambda^{-\frac d2} 2^m\smp .$  
Combining these two estimates via interpolation  to remove $2^m$,
we get 
\[
\|  \chi_{Q_k}\fP_\lambda[\psi_j \psi_{Q_k}^m] \chi_{B'} \|_{2\to 2}\lesssim  \lambda^{-\frac d2} \mup^{\frac{2d+1}{6d}}\mu^\frac{d-1}{6d} . 
\] 
Since $\{Q_k\}$ are essentially disjoint and  there are as many as $O(\sqrt{\mu/\mu'})$ of $\{Q_k\}$, it is clear that  $\| \sum_k \mathfrak P^m_k \|_{2\to2} \lesssim  \mom^{-\frac14} \max_{k} \| \chi_{Q_k}\fP_\lambda(\psi_j \psi_{Q_k}^m) \chi_{B'} \|_{2\to 2} $.  Therefore, we get  
\begin{align*}
\| \sum_k &\mathfrak P^m_k \|_{2\to2} \lesssim    \lambda^{-\frac d2} \mup^{\frac{2d+1}{6d}}\mu^\frac{d-1}{6d} \mom^{-\frac14} = \lambda^{-\frac d2}( \mm)^\frac14  \mom^{-\frac16+\frac1{6d}}.  
\end{align*}
We now interpolate this with \eqref{onetoinfty} to get 
\[\| 
 \sum_k \mathfrak P^m_k \|_{p\to p'}
           \lesssim  2^{-\frac m2\dpp}    \bpqp
            \mom^{\frac{d-1}{6d}(\frac{3d+4}{4}\dpp-1)} .
\] 
So, summation along $m\ge 1$  gives  the desired estimate \eqref{m1} with some $c>0$ if $ \dpp> \frac{4}{3d+4}.$ 
\end{proof}

\subsection{Estimate for $\sum_{k}\mathfrak P^0_k $}    This section and subsequent sections are  devoted to showing the following, 
 which completes the proof of Proposition \ref{critical-muj} being combined with Lemma \ref{m-m}.

\begin{prop}
\label{mm2}  Let $d\ge 5$.   If  $6/(d+5)<\dpp<  2/3$, then we have, for some $c>0$,
\begin{equation}
\label{m0}
\| \sum_{k}\mathfrak P^0_k \|_{p\to p'} \lesssim   \bpqp \mom^c.  
\end{equation}
We also have restricted weak type estimate for $\sum_{k}\mathfrak P^0_k$ if $\dpp=2/3$ and $\mu=\mu'$ for $d\ge 2$.
\end{prop}

\begin{proof}[Proof of Proposition \ref{mm2}]
As before, in order to get the desired bound $O(\lambda^{-\frac12})$ for $L^1$--$L^\infty$ estimate we need additionally to decompose $\mathfrak P^0_k$ by inserting  the cutoff function  $\widetilde \psi(2^{l}(s-S_c(x,y))$. That is to say, 
\begin{equation}
\label{mmm}
 \mathfrak P^0_k=  \sum_l \mathfrak P^{0}_{k,l}:=  \sum_l   \chi_{Q_k}\fP_\lambda[\psi_j \psi_{Q_k}^0\widetilde \psi(2^{l}(\cdot-S_c)) ] \chi_{B'}.
 \end{equation}
This type of decomposition away from the critical point  is unavoidable since  the decay of order $-\frac12$  is only available when we have the same order of decay estimate for the oscillatory integral. Clearly, for each $l$ we have 
\Be 
\label{interval}
2^{-l}\le \vepc \smp
\Ee
since the summand becomes zero otherwise. The condition \eqref{interval} plays  a crucial role in what follows. It allows us to suppress the errors to manageable size by taking $\vepc$ small enough. 

We note that 
\Be
\label{p2-lower}
|\partial_s^2 \mathcal P(x,y,s)|\gtrsim 2^{-l}, \quad s\in \supp   \widetilde\psi(2^{l}(\cdot-S_c(x,y))) \cap \supp \psi_j.
\Ee
Indeed, from \eqref{rxy} and \eqref{roots}  we have $\mathcal R(x,y,\cos s)=(\cos S_c(x,y)- \cos s)( \tau_+(x,y)-\cos s)$.  On the other hand,   we have $\tau_+(x,y)-1\sim  |x-y|\sim \mu$ by \eqref{roots} if $(x,y)\in B\times B'$. Since 
 $\tau_+(x,y)-\cos s\ge \tau_+(x,y)-1\sim \mu $ for $(x,y)\in B\times B'$ and $s\in \supp \psi_j$, using  the mean value theorem we have  
$|\mathcal R(x,y,\cos s)|\sim  \mu^\frac32 |S_c(x,y)-s|\sim   \mu^\frac322^{-l}$. So, \eqref{p2-lower} follows by \eqref{ph-dd}.  

Using \eqref{S3-ptkernel}, 
we have 
$
 \|\mathfrak P^{0}_{k,l} \|_{1\to \infty}\lesssim\lambda^{-\frac12}2^{\frac l2} \mu^{-\frac d4} =\lambda^{-\frac12}2^{\frac l2}  ( \mm)^{-\frac d8} \mom^{\frac d8}
$
by \eqref{p2-lower} and Lemma \ref{osc-lemma}. 
Thus, it follows that 
 \begin{equation}
\label{mkl-est}
 \|\sum_{k} \mathfrak P^{0}_{k,l} \|_{1\to \infty}
 \lesssim\lambda^{-\frac12}2^{\frac l2}  ( \mm)^{-\frac d8} \mom^{\frac d8}
 \end{equation}
 since  $Q_k$ are  essentially disjoint. 
  On the other hand, for any $\epsilon>0$ we claim 
 \Be
 \label{eq:l2ml2}
\|\mathfrak P^{0}_{k,l} \|_{2\to 2}\lesssim  \lambda^{-\frac d2} 2^{-l}   \mom^{-\frac34-\epsilon}.
\Ee
 This is to be shown later (see Lemma \ref{key-osc}). Since there are as many as $O(\mom^{-\frac12})$ essentially disjoint $Q_k$,  
 \eqref{eq:l2ml2} gives 
 $
\|\sum_k\mathfrak P^{0}_{k,l} \|_{2\to 2}\lesssim  \lambda^{-\frac d2}  2^{-l} \mom^{-1-\epsilon}.
$
  Interpolating this estimate and  \eqref{mkl-est},  we get
  \[ 
  \| \sum_k \mathfrak P^{0}_{k,l}\|_{p\to p'} \lesssim  \lambda^{\scaleppp}  2^{(\frac32\dpp-1)l}  (\mm)^{-\frac d8\dpp}  \mom^{-1+\frac{d+8}{8}\dpp-\epsilon }
  \]
  for $1\le p\le 2.$ 
Summation  over $l$ gives 
\[\| \sum_k \mathfrak P^{0}_k\|_{p\to p'}  \lesssim   \bpqp   \mom^{-\frac 34+\frac{d+5}{8}\dpp-\epsilon}\] for 
$p$ satisfying  $\dpp< 2/3$ because of  \eqref{interval}. 
Since $\epsilon$ can be taken arbitrarily  small, we have 
\eqref{m0} if $6/(d+5)<\dpp <2/3$.

For $d\ge 2$  the restricted weak type bound for $\dpp=2/3$ and $\mu=\mu'$ follows from the above estimate and Lemma \ref{s-trick}.
\end{proof}

For completion of  the proof of Proposition \ref{critical-muj} it remains to show \eqref{eq:l2ml2}.
 
 \subsection{Proof of \eqref{eq:l2ml2}}  Now we proceed to show the $L^2$ estimate  \eqref{eq:l2ml2}.   
 Concerning the operator  $\mathfrak P^{0}_{k,l}$, 
 we may replace $\chi_{Q_k}$ in \eqref{mmm} with a smooth function  $\widetilde \chi_ {Q_k}$  adapted to $Q_k$ such that $ \widetilde\chi_ {Q_k}=1 $ on $Q_k$ 
 and $\widetilde \chi_ {Q_k}$ is supported in 
 a $\varepsilon_\circ \sqrt{\mu\mu'}$ neighborhood of $Q_k$.  
 
 We recall the phase function $\Phi_s$  given by \eqref{phase-s}
and  set 
\Be
\label{as-as}
  A_s(x,y):= \widetilde \chi_ {Q_k}(x) \chi_{B'}(y)   \big ( 2^{-\frac d2j} \mathfrak a\psi_j\psi_{Q_k}^0  \big)(\sxyls) .
  \Ee
Then, changing variables $s\to 2^{-l}s+S_c(x,y)$  (see Section \ref{proof-away}) and recalling the definition \eqref{def:osc},  we see
\[  \mathfrak P^{0}_{k,l} f= 2^{\frac d2j} 2^{-l} 
\int  \widetilde \psi(s)  \mathcal O_\lambda  [\Phi_s, A_s] f\,ds.\] 
Since $2^{-j}\sim \smu$, for \eqref{eq:l2ml2} it is sufficient to show the following.

\begin{lem}\label{key-osc} Let $d\ge 2$ and $s\in \supp \psi$. Suppose \eqref{interval} holds. Then, for any $\epsilon>0$ we have the estimate
\[\|\mathcal O_\lambda [\Phi_s, A_s]\|_{2\to 2}\le C  \mu^\frac d4 \lambda^{-\frac d2}  \mom^{-\frac34-\epsilon}\]
with $C$ independent of $s$. 
\end{lem}

As is to be seen below,  the above oscillatory integral estimate displays a different nature.  The phase $\Phi_s$ is basically contained in the category of 
the oscillatory integral operators considered  in Lemma \ref{generalized}. However, as  $\mu'/\mu$ gets smaller,  it behaves badly  along the variables $x_1, y_1$.

 In order to prove Lemma 
\ref{key-osc} we  need the following several lemmas which  provide control over derivatives of the phase $\Phi_s$ and the amplitude $A_s$.

\begin{lem}  Let $2^{-l}\le \vepc \smp
$, and $s\in \supp \psi$. 
\label{sc-derivative} If $(x,y)\in B\times B'$, then 
\begin{align}
\label{est-partial}
\partial_{x}^\alpha\partial_y^\beta S_c&=O\Big( \mu^{\frac12-|\alpha|-|\beta|}\Big),  
\\
\label{est-amp}
\partial_{x}^\alpha A_s&=O( (\mm)^{-\frac{|\alpha|}{2}}).
\end{align} 
\end{lem}

\begin{proof} Let us set $G(x,y):=|x-y|^{1/2}$. 
To show \eqref{est-partial}, from \eqref{sc} we first note  that $S_c(x,y)$ is close to $G(x,y)$, and one may expect  
$S_c(x,y)$ behaves like $G(x,y)$, which satisfy \eqref{est-partial} and \eqref{est-amp} instead of $S_c$. More rigorously, we can prove \eqref{est-partial} inductively  by making use of \eqref{sc}.

We first note from Lemma \ref{scs}  that \eqref{est-partial} holds with $\alpha=\beta=0$  since $\dist(B,B')\sim \mu$. 
Then, we assume  \eqref{est-partial} is true for $|\alpha|+|\beta|\le N$. Applying $\partial_{x}^\alpha\partial_y^\beta $ on both side of \eqref{sc}  while $|\alpha|+|\beta|=N+1$,  it is not difficult to see that 
\Be
\label{derivative2}
\sin \sxy(x,y) \partial_{x}^\alpha\partial_y^\beta S_c(x,y) +O\big(\mu^{-N}\big)= O\big(\mu^{-N}\big)
\Ee
for $(x,y)\in B\times B'.$ Indeed the right hand side of \eqref{sc} behaves  as if  it were $|x-y|$ and the terms other than the first one in the  left hand side of  \eqref{est-partial} are given by a linear combination of 
products of  either $\sin \sxy(x,y)$, or $\cos \sxy(x,y)$,   $\prod_{i=1}^l\partial_{x}^{\ol\alpha_i}\partial_y^{\ol\beta_i} S_c(x,y)$, $l\ge 2$ and $\sum_{i=1}^l ( |\ol\alpha_i|+|\ol\beta_i|)\le N+1$ . Thus, our induction assumption 
shows these terms  are $O(\mu^{-N})$ and we get \eqref{derivative2}. 
Since  $\mu'\le \mu$, $|x-y|\sim \mu$ and $\sxy(x,y)\sim \smu$ for $(x,y)\in B\times B'$, 
\eqref{derivative2} gives  $\partial_{x}^\alpha\partial_y^\beta S_c(x,y)=O(\mu^{-N-1/2})$ as desired. 

Once we have 
\eqref{est-partial}, the proof of \eqref{est-amp} is easier. Since $\frac{d^k}{ds^k}(2^{\frac d2j} \mathfrak a\psi_j)=O(\mu^{-\frac k2})$,  from \eqref{as-as} it is sufficient to show that 
\begin{align*}
\partial_{x}^\alpha\partial_y^\beta \Big( \rho \Big(\frac{2^{-l}s+S_c(x,y)-s_0}{C\vepc\smp}\Big)\Big)&=O\big((\mmt)^{-|\alpha|-|\beta|}\big), 
\\
\partial_{x}^\alpha\partial_y^\beta \Big ( \rho \Big(\frac{2^{-l}s+S_c(x,y)-s_0}{C\vepc\smu}\Big)\Big)&=O(\mu^{-|\alpha|-|\beta|})
\end{align*}
for any $s_0$ provided that $\rho\in C_c^\infty(-1,1)$ with $|\frac{d^k}{ds^k}\rho| \le C$  for $k=0,1,\dots$. 
The other factors can be handled easily.  Since we have \eqref{est-partial}, both of the bounds  follow from straightforward computation. 
\end{proof}

\begin{lem}
\label{lem:ph-derivative}   Let $2^{-l}\le \vepc \smp$.  If $(x,y)\in Q_k\times B'$  and $0<\mu'\le \mu\ll1$,  then for $|\alpha|, |\beta|\ge 1$ we have 
\begin{equation}
\label{ph-derivative}
 \partial_x^\alpha \partial_y^\beta  \Phi_s(x,y)=O(\mu^{\frac32-|\alpha|-|\beta|}).
 \end{equation}
\end{lem}

\begin{proof} 
We recall \eqref{phase-s}  and \eqref{p-def}. 
 If $|\alpha|=1$, we write 
 \begin{align*}
 \ \partial_{x}^\alpha \Phi_s(x,y,s)
=I+I\!I+ I\!I\!I,
\end{align*}
where 
\begin{align*}
& I=-\frac{  \mathcal Q(x,y,\cos \sxyl (x,y))}{\sin^2 \sxyl (x,y)} \partial_{x}^\alpha\sxy, \qquad I\!I=  x^\alpha  \frac{\cos  \sxyl (x,y) -1}{\sin  \sxyl (x,y) }, 
 \\[3pt]
 & \qquad\qquad\qquad \ I\!I\!I=(x^\alpha-y^\alpha) \csc  \sxyl (x,y).
 \end{align*}
 Since $|\mathcal D|\le \vepc\mu\mu'$ for $(x,y)\in Q_k\times B'$, by \eqref{interval} and   \eqref{qcos}  we note 
\Be 
\label{Q-size} \mathcal Q(x,y,\cos \sxyl(x,y) )=O(\vepc \mu\mu' ), \quad (x,y)\in Q_k\times B'.
\footnote{It is clear from \eqref{eps-con}. Also one may use the mean value theorem to see $\mathcal Q(\cos ( 2^{-l}s+S_c(x,y)))
=\partial_\tau \mathcal Q(\cos(  2^{-l}s^*+S_c(x,y)))\sin (  2^{-l}s^*+S_c(x,y)) 2^{-l}s +O(\vepc \mu\mu')$ for some $s^*$. }
 \Ee
We also have 
$x^\alpha =O ({(\mm)}^{|\ol\alpha| /2})$ and $x^\alpha-y^\alpha=O (\mu^{\alpha_1}{(\mm)}^{|\ol\alpha| /2}) $ for  $(x,y)\in Q_k\times B'$. 
Here we denote $\alpha=(\alpha_1, \ol \alpha)\in \mathbb N_0\times \mathbb N_0^{d-1}.$
Since  $\sxyl\sim \smu$, and $\sxy(x,y)\sim \smu$, by \eqref{est-partial} we see  
$I=O(\mu'\mu^{\frac12-|\alpha|})$ and $I\!I=O(\smu {(\mm)}^{|\ol\alpha| /2})$. 
Similarly, since $|x-y|\lesssim \mu$, we get $I\!I\!I=O(\mu^{\alpha_1-\frac12}  {(\mm)}^{|\ol\alpha| /2})$.  We also see  the same estimate  holds for $\partial_{y}^\alpha \Phi_s(x,y,s)$.  

The estimates for higher order 
derivatives can be shown in the  similar manner by direct differentiation. So, we omit the detail.  
\end{proof}

For a given matrix $\mathbf M$ we denote by $\mathbf M_{i,j}$ the $(i,j)$-th entry of $\mathbf M$ and 
define  $\widetilde{\mathbf M}$ to be  the $(d-1)\times (d-1)$ submatrix of $\mathbf M$  given by $\widetilde{\mathbf M}_{i,j}=\mathbf M_{i+1,j+1}$  
for $1\le i,j\le d-1$.  The following lemma shows that the matrix $\partial_x \partial_y^{\intercal}   \Phi_s(x,y)$ takes a particular form depending on $\mu$ and $\mu'$. 
This fact plays an important role   in proving Lemma \ref{key-osc}.  
\begin{lem}  
\label{Hmatrix-lower}
Let $2^{-l}\le \vepc \smp
$ and set \[\mathfrak M=-  \sin S_c(x,y) \partial_x \partial_y^{\intercal}   \Phi_s(x,y).\]  
If $(x,y)\in Q_k\times B'$ and $0<\mu'\le\mu\ll1$, then
\Be
\label{dxyPhi}
\mathfrak  M=:\begin{pmatrix} 
 \mathfrak  M_{1,1} & v^\intercal 
 \\
 w & \widetilde{\mathfrak M}
\end{pmatrix}
=
\begin{pmatrix}  \sim \mu'/\mu &  O\big(\sqrt{ {\mu'}/{\mu}}\,\big)\,
\\  
\,O\big(\sqrt{ {\mu'}/{\mu}}\,\big)  & O (1)
\end{pmatrix},
\Ee
and 
\begin{align}
\label{det-value}
\det{\mathfrak M} &\sim \mom.
\end{align}
Additionally, if $\mu'\ll \mu$, then we have 
\begin{align} 
\label{det-value1}
\det{\widetilde{\mathfrak M}}  &\sim  1.
\end{align}
\end{lem}

To prove Lemma \ref{key-osc} we combine Lemma \ref{Hmatrix-lower} and  Lemma \ref{generalized-small2} below. 
The oscillatory integral considered in Lemma \ref{generalized-small2}  differs from  the typical oscillatory integral (cf. Lemma \ref{generalized}) in that 
the phase and amplitude functions do not behave well in certain directions.

\begin{lem}
\label{generalized-small2} Let $0<\omega\le 1$.  Let $a$ be a smooth function supported in $S_\omega:=[-\omega, \omega]\times B_{d-1}(0,1)\times [-\omega^2, \omega^2]\times B_{d-1}(0,1)$ such that 
\Be
\label{a-con}
 |\partial_{x}^{\alpha}  a|\le C\omega^{-\alpha_1}.
 \Ee  
Suppose that  $(\partial_x\partial_y^\intercal \phi)_{i,j}=O(1)$, $|\!\det \partial_x\partial_y^\intercal \phi|\sim 1$, $|\!\det \partial_{\bar x}\partial_{\bar y}^\intercal \phi|\sim 1$, and 
\Be
\label{p-con}
| \partial_{x}^{\alpha}\partial_{y}^{\beta}  \phi| \lesssim  \omega^{-2 +|\ol \alpha|+|\ol \beta|},  \quad 
|\alpha|, |\beta|\ge 1 
\Ee 
 on $S_\omega$. Then,   for any $\epsilon>0$
we have 
\Be \label{l2-oa} 
\|\mathcal O_\lambda [\phi, a]  f\|_2\lesssim \lambda^{-\frac d2}  \omega^{-\frac12-\epsilon} \|f\|_2. 
\Ee 
\end{lem}

In the above, the factor $\omega^{-\frac12-\epsilon}$ is important for our purpose but there is no reason to believe that it  is optimal.   It seems to be likely that the bound can be improved by further refinement of our argument. However, this does not  look simple, especially,  in a general form such as in Lemma \ref{generalized-small2}. Instead, to obtain an improved bound, one may try to exploit the particular structure of the phase $\Phi_s$  and amplitude $A_s$. Unfortunately, this leads to considerable increase of complexity.  We hope to pursue the matter elsewhere near future.

 \begin{proof}[Proof of Lemma \ref{key-osc}]
Let $(x_0, y_0)$ be the center of the cube  $Q_k\times B'$. 
We make change of variables  to balance the size of entries in the matrix $\mathfrak M$. 
Let us set 
\begin{align*}
L(x,y)&= \big( \mu  x_1,\mmt\,\,\ol x, \mu y_1, \mmt\,\, \ol y\big)+  (x_0,  y_0),
\\ 
\nonumber \widetilde \Phi(x,y)&=(\smu\mu')^{-1}\Phi_s(L(x,y)) ,\\
 \widetilde{A}(x,y)&=  A_s(L(x,y)) .
\end{align*}
By the change of variables $(x,y)\to L(x,y)$, it is easy to see 
\Be 
\label{op-norm2}
\|\mathcal O_\lambda [\Phi_s, A_s]\|_{2\to 2}  =\mom^{-\frac12}(\mm)^{\frac d2}
\Big\|\mathcal  O_{\smu\mu'\lambda} \big[ \wt\Phi,  \wt A\,\big]\Big\|_{2\to 2}.
\Ee
We recall that  $Q_k$ is a cube of side length $\vepc \mmt$ and $ B'$ is of dimensions $\vepc\mu'\times\times\vepc\mmt\times \cdots \vepc\mmt$, and  hence  $\wt A$ is supported in 
$ B_1(0, \sqrt{\mu'/\mu})\times B_{d-1}(0, 1) \times  B_1(0, \mu'/\mu)\times B_{d-1}(0,1).$ 
From \eqref{est-amp} it follows that 
\begin{equation}
\label{scaled-amp2}
\partial_{x_1}^{\alpha_1} \partial_{\ol x}^{\ol \alpha} \wt A=O\big(\mom^{-\frac{\alpha_1}2}\big).
\end{equation}
Note that $ \sin S_c(x,y)\sim \smu$. Using \eqref{dxyPhi}, \eqref{det-value},  and \eqref{ph-derivative}, successively,  we get
\begin{align} 
 \det \partial_x  \partial_y^\intercal   \wt \Phi  &\sim 1,  
 \nonumber
  \\
\label{scaled-det2}
 (\partial_x  \partial_y^\intercal   \wt \Phi)_{i,j}&=O(1),    \quad 1\le i, j\le d
 \\
\label{scaled-ph2} 
\partial_{x_1}^{\alpha_1} \partial_{\ol x}^{\ol \alpha}\partial_{y_1}^{\beta_1} \partial_{\ol y}^{\ol \beta} 
  \wt \Phi &=  O \big(\mom^{\frac{|\ol \alpha |+ |\ol \beta | }2 -1}\big), \quad |\alpha|, |\beta|\ge 1 
\end{align}
 for $(x,y)\in \supp \wt A$. Therefore,   if $\mu\sim \mu'$, the desired estimate follows from Lemma \ref{generalized}. 
 If $\mu\gg \mu'$, then we have
\[ \det \partial_{\ol x} \partial_{\ol y}^\intercal   \wt \Phi \sim 1,\]
which follows from \eqref{det-value1}. 
Hence,  combining this with  \eqref{scaled-amp2}, \eqref{scaled-ph2}, and \eqref{scaled-det2}, we may apply Lemma \ref{generalized-small2} with $\omega=(\mu'/\mu)^{1/2}$ to the oscillatory integral operator 
$ \mathcal  O_{{\smu\mu'\lambda}}[ \wt \Phi,  \wt A\,]$. Thus, we  get 
\[
\Big\|\mathcal  O_{\smu\mu'\lambda}\big[ \wt\Phi,  \wt A\,\big]\Big\|_{2\to 2}
\lesssim 
\lambda^{-\frac d2}  (\mu')^{-\frac d2}\mu^{-\frac d4} \mom^{-\frac14}\mom^{-\epsilon} .
\]
By \eqref{op-norm2} we obtain the desired estimate. 
\end{proof}

Regarding the proof of Lemma \ref{generalized-small2},  the typical argument (for example, see \cite[pp. 377--379]{St93})  which deals with the $L^2$ oscillatory integral estimate does not work well here since the integration by parts along two different directions  does not necessarily split  in  a nice manner. We get around this  by using a simple localization  argument.

\begin{proof}[Proof of Lemma \ref{generalized-small2}]  
By finite decomposition and translation we may replace $S_\omega$ with 
$[-\vepc\omega, \vepc\omega]\times B_{d-1}(0,\vepc)\times [-\vepc\omega^2, \vepc\omega^2]\times B_{d-1}(0,\vepc)$ with a small enough $\vepc$. Let us set 
\[ I_\lambda(y,y')=\int e^{i\lambda ( \phi(x,y)-\phi (x,y'))}  a(x,y)\overline {a}(x,y')dx,\]
which is the kernel of $\mathcal O_\lambda [\phi, a]^*\mathcal O_\lambda [\phi, a]$.   
We first claim that 
\begin{equation}
\label{iii}
| I_\lambda(y,y')| \lesssim C\omega (1+\lambda\omega^2|y-y'|)^{-N}.
 \end{equation}
To show this, let us set 
\begin{align*}
\Psi(x)&=  \phi(\omega x_1, \ol x,y)-\phi (\omega x_1, \ol x, y'), 
\\
\quad A(x) &=   a(\omega x_1, \ol x,y)\overline {a}(\omega x_1, \ol x,y'). 
\end{align*}
Then we write 
\[ I_\lambda(y,y')=\omega \int e^{i\lambda \Psi(x) }  A(x)dx.\]
Clearly, by \eqref{a-con} we have $\partial^\alpha_x A=O(1)$ for all $\alpha$. We claim that 
\begin{align}
|\nabla  \Psi(x)|&\gtrsim \omega |y-y'|, 
\label{lowerbdd}
\\
|\partial_x^\alpha   \Psi(x)|&\lesssim |y-y'|,   \quad \forall |\alpha|\ge 2.
 \label{upperbdd}
\end{align}
Since $\big(\frac{\nabla  \Psi(x)}{i\lambda |\nabla   \Psi(x)|^2}  \cdot \nabla\big) e^{i\lambda \Psi(x) }=e^{i\lambda \Psi(x) }$, 
using the operator $\frac{\nabla  \Psi(x)}{\lambda |\nabla   \Psi(x)|^2}  \cdot \nabla$, we perform integration by parts $N$ times to get 
\[ |I_\lambda(y,y')|\lesssim   \omega \sum_{\ell =0}^N \ \  \sum_{\sum |\alpha_i| =N+\ell, \ |\alpha_i|\ge 2  } \ 
\int \frac{\prod_{i=1}^\ell |\partial_x^{\alpha_i} \Psi(x)|}{\lambda^N|\nabla  \Psi(x)|^{N+\ell}}  \mathcal A(x) dx,    \] 
where $\mathcal A$ is a bounded function supported in $B_d(0,1)$.  By \eqref{lowerbdd}  and \eqref{upperbdd} we get 
the estimate \eqref{iii}.  

Now, we show \eqref{lowerbdd} and \eqref{upperbdd}. We first note that 
$\nabla \Psi(x)= ( \omega \partial_{x_1}\phi(\omega x_1, \ol x,y)-\omega\partial_{x_1}\phi (\omega x_1, \ol x, y'), 
\partial_{\ol x}\phi(\omega x_1, \ol x,y)-\partial_{\ol x}\phi (\omega x_1, \ol x, y'))$. Hence,  we have
$|\nabla \Psi(x)|\ge \omega |\nabla \phi(\omega x_1, \ol x,y)- \nabla \phi(\omega x_1, \ol x,y')|$. For \eqref{lowerbdd}
it is sufficient to show that \[ |\nabla \phi(\omega x_1, \ol x,y)- \nabla \phi(\omega x_1, \ol x,y')|
\gtrsim |y-y'|.\]  Using Taylor's expansion, $\nabla \phi(\omega x_1, \ol x,y)- \nabla \phi(\omega x_1, \ol x,y')$ is equal to
\[\partial_x\partial_y^\intercal \phi (\omega x_1, \ol x,y')(y-y')
+\!\!\!\sum_{i, \alpha: |\alpha|=2} O(E_{i,\alpha}|(y-y')^\alpha|), \] 
where $E_{i,\alpha}=\sup_{(x,y)\in S_\omega} |\partial_{x_i}\partial_{y}^\alpha  \phi(x,y)|$. Since 
$|y_1-y'_1|\le \vepc \omega^2$, by \eqref{p-con}  
 it is clear that $E_{i,\alpha}|(y-y')^\alpha|=O(\vepc|y-y'| )$.  We also have $|\partial_x\partial_y^\intercal \phi (\omega x_1, \ol x,y')(y-y')|\sim |y-y'|$ because $(\partial_x\partial_y^\intercal \phi)_{i,j}=O(1)$ and $|\!\det \partial_x\partial_y^\intercal \phi|\sim 1$. So,  we get  $ |\nabla \phi(\omega x_1, \ol x,y)- \nabla \phi(\omega x_1, \ol x,y')|
\gtrsim |y-y'|.$  For \eqref{upperbdd} it is sufficient to show that $\partial_x^\alpha \partial_{y_j} (\phi(\omega x_1, \ol x,y))=O(1)$, $|\alpha|\ge 2$ by the mean value theorem and this is clear since  $\partial_x^\alpha \partial_{y_j} (\phi(\omega x_1, \ol x,y))=O( \omega^{-2+\alpha_1+|\ol\alpha|+|\ol\beta|})=O(1)$ if $|\alpha|\ge 2$. 

By virtue of  \eqref{iii}, it is enough to show \eqref{l2-oa} while assuming that 
$f$ is supported in a cube  of side length $\lambda^{-1}\omega^{-2-\epsilon}$. Indeed, let $\{\mathfrak q\}$ 
be  essentially disjoint closed cubes of side length $\lambda^{-1}\omega^{-2-\epsilon}$ which 
covers $B_d(0,1)$.  So, we may write $f=\sum_\fq f_{\fq}$ where $f_\fq=\chi_\fq f$ and then we have 
\begin{align*}
\|\mathcal O_\lambda [\phi, a] f\|_2^2 
   &=\big(\!\!\!\sum_{\dist(\fq, \fq')\ge  \lambda^{-1}\omega^{-2-\epsilon} }  \! + \!\!\!\sum_{\dist(\fq, \fq')< \lambda^{-1}\omega^{-2-\epsilon} } \big)
\int  I_\lambda(y,y') f_{\fq'}(y') \ol f_{\fq} (y)  dy dy'\\
           &=:A+B.
           \end{align*}
   Note that      $\lambda\omega^2 |y-y'|\gtrsim  \omega^{-\epsilon}$ for $(y,y')\in \fq\times\fq'$ if $\dist(\fq, \fq')\ge  \lambda^{-1}\omega^{-2-\epsilon}$. 
   So,  taking large $N=N(\epsilon)$  we  see from \eqref{iii} that 
\[  |A|\lesssim \int  \omega^M (1+\lambda\omega^2|y-y'|)^{-N/2} |f(y)||f(y')| dydy' . \] 
Thus $ |A|\le C\lambda^{-d} \omega^{M-2d} \|f\|_2^2$. Therefore, the matter reduces  to showing 
\[|B|\lesssim \lambda^{-d}  \omega^{-1-\epsilon} \|f\|_2^2.\]
By disjointness of $\{\fq\}$, for this  it suffices  to show
$|\int  I_\lambda(y,y') f_{\fq'}(y') \ol f_{\fq} (y)  dy dy'|\lesssim  \lambda^{-d}  \omega^{-1-\epsilon}\|f_{\fq'}\|_2\|f_\fq\|_2$ with  $\dist(\fq, \fq')
\lesssim  \lambda^{-1}\omega^{-2-\epsilon}$. In turn, the estimate clearly follows from 
\Be
\label{single-cube}
  \|\mathcal O_\lambda [\phi, a] g\chi_\fq\|_2\lesssim   \lambda^{-\frac d2}  \omega^{-\frac12-\epsilon}\|g\|_2
  \Ee
  by the Cauchy-Schwarz's inequality. 
  
   To show \eqref{single-cube} we further divide $\fq$ into almost disjoint rectangular  slabs  $S_1, \dots, S_L$ of dimensions $\lambda^{-1}\times \lambda^{-1}\omega^{-2-\epsilon}\times\cdots \times \lambda^{-1}\omega^{-2-\epsilon}$ such that $\fq=\bigcup_{k=1}^L
S_k$. We claim
\Be
\label{single-slab}\|\mathcal O_\lambda [\phi, a] g\chi_{S_k} \|_2\lesssim \lambda^{-\frac d2}\omega^\frac12   \|g\|_2.
\Ee
By the triangle inequality and \eqref{single-slab} we  have 
\[\|\mathcal O_\lambda [\phi, a] g\chi_\fq\|_2\le \sum_{k=1}^L \| \mathcal O_\lambda [\phi, a] g  \chi_{S_k} \|_2 
\lesssim  \lambda^{-\frac d2}  \omega^\frac12 \sum_{k=1}^L \|g  \chi_{S_k} \|_2 .\]
Hence, the estimate \eqref{single-cube} follows because $  \sum_{k=1}^L \|g  \chi_{S_k} \|_2\lesssim \omega^{-1-\frac\epsilon 2}\|g\|_2$. 
This  is clear by the Cauchy-Schwarz inequality because 
 $L=O(\omega^{-2-\epsilon})$.

In order to show \eqref{single-slab}, it is sufficient to prove that, for $1\le k\le L$,
\begin{equation}
\label{decay-est}
| I_\lambda(y,y')|\lesssim C\omega (1+\lambda|\ol y-\ol y'|)^{-N},  \quad  y, y'\in S_k
\end{equation}
because $\int_{S_k} \omega (1+\lambda|\ol y-\ol y'|)^{-N} dy\lesssim \omega\lambda^{-d}$. 
Since $|y_1-y'_1|\le \lambda^{-1}$ for $y, y'\in S_k$ and $\partial_{\ol x}\partial_{y_1} \phi=O(1)$,  
$ \partial_{\ol x}\phi(x,y)-\partial_{\ol x} \phi (x,y')=\partial_{\ol x}\phi(x,y)-\partial_{\ol x} \phi (x,y_1, \ol {y'})+O(\lambda^{-1}).$
Recalling $|\det \partial_{\ol x}\partial_{\ol y}\phi|\sim 1$, we see 
\[ |\lambda( \partial_{\ol x}\phi(x,y)-\partial_{\ol x} \phi (x,y'))|\gtrsim \lambda |\ol y-\ol {y'}|+O(1).\] 
Therefore,  integration by parts in $\ol x$  gives the desired bound \eqref{decay-est}.  \end{proof}

For the rest of this section we prove  Lemma \ref{Hmatrix-lower}. 

\subsection{Proof of Lemma \ref{Hmatrix-lower}}  As before,  it is enough to show Lemma \ref{Hmatrix-lower} under the assumption \eqref{coordinates}, i.e., 
$x=(r,0,0,\dots, 0)$ and $y=(\rho,h,0,\dots, 0)$, which we assume for the rest of this section.   We can use a  rotation $\mathbf U$ to  cover  the general case. 
Indeed, for $(x,y)\in B\times B'$,  there is a rotation matrix $\mathbf U$ such that 
$\mathbf Ux=(r,0,0,\dots, 0)$, $\mathbf Uy=(\rho,h,0,\dots, 0)$,  and  
\Be
\label{rotation-matrix}
\mathbf U
=
\begin{pmatrix}  
      \mathbf U_{1,1} &  v
     \\  
w &  \widetilde{\mathbf U}
\end{pmatrix}
  =
   \begin{pmatrix}  1+ O(\mm) &  O(\mmt)
               \\  
           O(\mmt)  &  O(1)
   \end{pmatrix}.
\Ee 
To see this, 
first  choose  a rotation matrix  $\mathbf R$  such that $\mathbf R x=(r,0,\dots, 0)$. Since  $x\in B$, it is clear that $\mathbf R_{1,1}= e_1\cdot\mathbf R e_1 =1+O(\mm)$ and 
$\mathbf R_{1,j}= O(\mmt),$ 
$ \mathbf R_{j,1} = O(\mmt)$, $2\le j\le d$. Next, we may choose 
another rotation $\mathbf S$ such that  $\mathbf Se_1=e_1$ and $\mathbf S\mathbf R y=(\rho, h,0,\dots, 0)$. Then, it is easy to see  that $\mathbf S\mathbf R$ takes the form  in \eqref{rotation-matrix} because $y/|y|=e_1+O(\mmt)$ if $y\in B'$.   
Similarly,  we see $\mathbf U^\intercal$  is also of the same form as in \eqref{rotation-matrix}. 

By \eqref{sc} and \eqref{phase-s},  $S_c(\mathbf U x, \mathbf Uy)=S_c(x, y)$ and 
$\Phi_s(\mathbf Ux, \mathbf Uy)=\Phi_s(x, y)$ for $\mathbf U\in \mathrm O(d)$. Using the second identity, a computation shows 
$ \mathbf U^\intercal \partial_x \partial_y^\intercal \Phi_s(\mathbf Ux, \mathbf Uy)  \mathbf U= \partial_x \partial_y^\intercal \Phi_s(x, y).$
Thus, we have 
\Be 
\label{identity} 
\mathfrak M(\mathbf Ux, \mathbf Uy)=  \mathbf U  \mathfrak M(x,y) \mathbf U^\intercal . 
\Ee
Now we claim that \emph{if  the matrix $\mathfrak M(\mathbf Ux, \mathbf Uy)$ satisfies \eqref{dxyPhi},  \eqref{det-value}, and  \eqref{det-value1} in place of $\mathfrak M$, 
then the same is true for $\mathfrak M(x,y)$.} It is clear from \eqref{identity}  
that $\mathfrak M(x,y)$ satisfies 
\eqref{det-value}.  Therefore, we may assume   \eqref{coordinates} to show  Lemma \ref{Hmatrix-lower}.

To verify the claim,   we  need only to show
 \eqref{dxyPhi} and  \eqref{det-value1} for  $\mathfrak M(x,y)$. Since $\mathbf U$ and $\mathbf U^\intercal$ take the form in \eqref{rotation-matrix}, 
using  block matrix multiplication, it is easy to  see that 
\[  \mathfrak M(x,y)= \mathbf U^\intercal \mathfrak M(\mathbf Ux, \mathbf Uy) \mathbf U =\begin{pmatrix}  
     a  +O(\mu') &  O(\sqrt{\mu'/\mu})
     \\  
O(\sqrt{\mu'/\mu})  & \quad \widetilde{\mathbf U}^\intercal\,\widetilde{\mathfrak M}(\mathbf Ux, \mathbf Uy)  \, \widetilde{\mathbf U}   +O(\mu')
\end{pmatrix}\]
with $a\sim   \mu'/\mu$. 
Since $ \mu'\ll \mu'/\mu$,  $\mathfrak M(x,y)$  takes the same form as in  \eqref{dxyPhi}. Also note that 
$|\det \widetilde{\mathbf U}|\sim 1$ and  $|\det \widetilde{\mathbf U}^\intercal|\sim 1$ because $\mathbf U$ and $\mathbf U^\intercal$ are of the same form in \eqref{rotation-matrix}.
Also,  $\det \widetilde{\mathbf U}$ and $\det \widetilde{\mathbf U}^\intercal$  have the same sign, so  $\det (\widetilde{\mathbf U}^\intercal\,\widetilde{\mathfrak M}(\mathbf Ux, \mathbf Uy)  \, \widetilde{\mathbf U})  \sim 1$. This shows  
\eqref{det-value1} is true for $\mathfrak M(x,y)$ because $\mu'\ll 1$ and completes the proof of our claim.

Recalling \eqref{partialxy} we begin with sorting out  harmlessly small parts.   Let us set
\begin{align*}
\mathbf G(x,y)&:=
    -\sin\sxyl\big(\partial_x S_c  \partial_y^{\intercal}  \partial_s \mathcal P(x,y,  \sxyl) 
+ \partial_x\partial_s\mathcal P(x,y,  \sxyl) \partial_y^{\intercal}  S_c  \big),
 \\
\mathbf F(x,y)&:= -\sin\sxyl\partial^2_s\mathcal P(x,y,\sxyl)\partial_x S_c\partial_y^{\intercal} S_c,
\\ 
 \mathbf M(x,y)&:= \  \mathbf I_d+ \mathbf G(x,y) + \mathbf F(x,y).
\end{align*}
Since $\sxyl\sim 2^{-j}\sim \smu$,  by \eqref{ph-d},  \eqref{Q-size}, and Lemma \ref{sc-derivative}  we see 
\[
| \sin\sxyl(x,y,s) \partial_s\mathcal P(x,y,  \sxyls) 
 \partial_x\partial_y^{\intercal}  S_c(x,y)|\lesssim  \vepc  \mu'/\mu, \   \   \ (x,y)\in Q_k\times B'.      
\]
Since $ -\sin s\,\partial_x \partial_y^{\intercal}  \mathcal P(x,y,s)= \mathbf I_d$,  from \eqref{partialxy} we have 
\[
    -\sin \sxyl \partial_x\partial_y^\intercal\Phi_s(x,y) = \mathbf M(x,y)+  O(\vepc  {\mu'}/\mu).
  \]
To show Lemma \ref{Hmatrix-lower}, taking small enough $\vepc>0$,  it is sufficient to show that 
the matrix $\mathbf M$ satisfies \eqref{dxyPhi}, \eqref{det-value} and \eqref{det-value1} in place of $\mathfrak M$.   

Let us set 
\begin{align}
\label{def-a-c}
\ba:= \frac{|x|^2+|y|^2}{\inp xy}, 
     \quad\quad \bc:=\frac{\cu}{\su|x+y||x-y|}. 
\end{align}
Then,  from Lemma \ref{id} we have
\[
    \partial_x S_c(x,y)=\bc(2x-\ba y), \quad \partial_y^\intercal S_c(x,y)=\bc(2y^\intercal-\ba x^\intercal).
\]
On the other hand,  a simple computation shows 
\begin{align*}
    \partial_x\partial_s\mathcal P(x,y, \sxyl)&=\frac{1}{\sin^2  \sxyl }(\cos  \sxyl\, y-x), \\ 
     \partial_y^\intercal \partial_s\mathcal P(x,y, \sxyl)&=\frac{1}{\sin^2 \sxyl}(\cos  \sxyl \, x^\intercal-y^\intercal).
\end{align*}
Using these  identities, we note that 
\Be 
\label{FG}
\begin{aligned}
    \mathbf G 
    &= -\frac \bc\sul  \big((2\cul+\ba)(xx^\intercal+yy^\intercal)-4xy^\intercal-2\ba\cul yx^\intercal \big),\\
    \mathbf F     &=\sul \, \partial_s^2\mathcal P (x,y, \sxyl) \, \bc^2\big(2\ba(xx^\intercal+yy^\intercal)-4 xy^\intercal-\ba^2 yx^\intercal\big).
\end{aligned}
\Ee
Since we are assuming  \eqref{coordinates} it is clear that $\mathbf F_{i,j}=0$ and $\mathbf G_{i,j}=0$ if $i$ or $j>2$ and, hence,  $\mathbf M_{i,j}=\delta_{ij}$ if  $i$ or $j> 2$ where 
$\delta_{ij}$ is the Kronecker delta.  Thus, to  show Lemma \ref{Hmatrix-lower} it is sufficient to consider the entries $\mathbf M_{1,1}$, $\mathbf  M_{1,2}$, $\mathbf  M_{2,1}$, and $\mathbf  M_{2,2}$.  For the purpose we need the following.

\begin{lem}\label{matlem}  Let $0<\mu'\le\mu\le c_\circ$ and let $x=(r,0,0,\dots, 0)$ and $y=(\rho,h,0,\dots, 0)$. 
Assume $2^{-l}\le\vepc\sqrt{\mu'}$, $1-r\sim \mu$, $1-\rho^2-h^2\sim \mu'$ and $h\sim \mmt$. Then, if $c_\circ>0$ and $\vepc$ are small enough, then  for $(x,y)\in Q_k\times B'$ and $\sxyl(x,y)\in \supp \psi_{Q_k}^0$ we have
\begin{align}
   & \mathbf M_{1,1} (x,y)=\frac{h^2}{|x-y|^2}+ O(\vepc{\mu'}/{\mu}),
   \label{m11}
   \\
   &  \mathbf M_{1,2}(x,y)=\frac{-h}{2|x-y|^2}\big((1+\cos\sxy)\rho-2 r\big)+O(\vepc\mom^\frac12), 
     \label{m12}
    \\
    &\mathbf M_{2,1}(x,y)=\frac{-h}{2|x-y|^2}\big((1+\cos S_c)\rho-2 r\cos\sxy  \big)+O(\vepc\mom^\frac12), 
     \label{m21}
    \\
     &\mathbf M_{2,2}(x,y)=  \frac{(r-\rho)^2}{|x-y|^2}+O(\vepc{\mu'}/{\mu}).
     \label{m22} 
\end{align}
\end{lem}

It should be noted that this lemma remains to be valid under the assumption that $\mu'\le \mu$ instead of $\mu'\ll \mu$. The case $\mu\sim \mu'$  can be proved in a much  simpler way  since the involved quantities are relatively large.

\begin{proof}[Proof of Lemma \ref{Hmatrix-lower}] 
Recall that it is enough to show  the matrix $\mathbf M$ satisfies \eqref{dxyPhi}, \eqref{det-value} and \eqref{det-value1} in place of $\mathfrak M$.   Once we have Lemma \ref{matlem}, the proof of Lemma \ref{Hmatrix-lower} is rather straightforward. 
Since $\mathbf M_{i,j}=\delta_{ij}$ if  $i$ or $j> 2$, for \eqref{dxyPhi} it is sufficient to note that 
\begin{align*}
    \mathbf M_{1,1}& \sim \mom, \qquad\qquad \mathbf M_{1,2}=O(\mom^\frac12), 
    \\ 
    \mathbf M_{2,1}&=O(\mom^\frac12), \qquad  \mathbf M_{2,2}=1-O\mom.
\end{align*}
This follows because $h\sim \mmt$, $|x-y|\sim \mu$, $S_c\sim \smu$, $|x-y|^2=(\rho-r)^2+h^2$, and $\rho-r=O(\mu)$. Indeed, the first and the last are clear. 
For the second and the third we simply write 
$(1+\cos\sxy)\rho-2r\cos\sxy =(1+\cos{S_c})(\rho-r)+(1-\cos{S_c})r$ and $(1+\cos S_c)\rho-2r=2(\rho-r)-(1-\cos{S_c})\rho$. Then, we observe that both are $O(\mu)$.  
As in the above, for \eqref{det-value} it is sufficient to show
\[ \det \begin{pmatrix}  \mathbf  M_{1,1}  & \mathbf  M_{1,2}
\\  
\mathbf  M_{2,1}& \mathbf  M_{2,2}
\end{pmatrix}\sim \frac {\mu'}{\mu}\,.
  \] 
  Using \eqref{m11}--\eqref{m22}, we see that   the determinant  $\det ((\mathbf M_{i,j} )_{1\le i,j\le 2})$ 
   equals
  \begin{align*}
      \frac{h^2}{4|x-y|^4} \big((1-\cos S_c)^2r^2+(3+\cos S_c)(1-\cos S_c)(r-\rho)^2\big)+ O(\frac{\vepc\mu'}{\mu}).
\end{align*}
Thus,  we have $\det ((\mathbf M_{i,j} )_{1\le i,j\le 2})  \sim \frac{\mu'}\mu$ 
if $\vepc$ is small enough.  Finally,  for \eqref{det-value1} 
we only have to note that   $|\rho-r|\sim \mu$ if $\mu'\ll \mu$ since  \eqref{m22} holds and $\mathbf M_{i,j}=0$ if $i$ or $j>2$.
 \end{proof}

We now proceed to prove Lemma \ref{matlem} by showing \eqref{m11}--\eqref{m22}, separately.   Thanks to \eqref{coordinates} and \eqref{FG} we can obtain explicit expressions of the matrices $\mathbf G(x,y)$ and $\mathbf F(x,y)$.  A simple computation gives 
\[
  \!\!\! \begin{pmatrix} \mathbf G_{1,1}& \!\!\! \mathbf G_{1,2} 
   \\
   \mathbf G_{2,1}&  \!\!\!  \mathbf G_{2,2}
                      \end{pmatrix}      =\frac {-\,\bc}\sul   \begin{pmatrix}
   (\ba+2)(r-\rho)^2+2(1-\cul)h^2& \!\!\big((\ba+2\cul)\rho-4r\big)h \\
   \!\! \big((\ba+2\cul)\rho-2\ba\cul r\big)h & (\ba+2\cul)h^2
    \end{pmatrix}.
\]
Indeed, for $\mathbf G_{1,1}$ we note $\mathbf G_{1,1}=\frac{-\, \bc}\sul \big( (\ba+2\cul)(r^2+\rho^2)-2(2+\ba\cul)r\rho \big)$  from \eqref{FG} and then we combine
this with the identity $\ba r\rho=|x|^2+|y|^2=r^2+\rho^2+h^2$.    The others are clear from \eqref{FG}.  Similarly,  we also have 
\[
 \begin{pmatrix} \mathbf F_{1,1}& \mathbf F_{1,2} 
   \\
   \mathbf F_{2,1}& \mathbf F_{2,2}
                      \end{pmatrix} = \sul \, \partial_s^2\mathcal P(x,y,\sxyl) \, \bc^2
    \begin{pmatrix}
    2\ba(r^2+\rho^2)-(4+\ba^2)r\rho & (2\ba\rho-4r)h \\
    (2\rho-\ba r)\ba h & 2\ba h^2
    \end{pmatrix}.
\]

\begin{proof}[Proof of   \eqref{m11}]
The proof of \eqref{m11} is more involved  than the others.  Since $ \bc/\sul $ $\sim\mu^{-2}$, $(1-\cul)\sim\mu$, $h^2\sim\mu\mu'$, and $\frac \bc\sul  (1-\cul)h^2
=O(\mu')$, taking $c_\circ$ small enough, we may assume $\mu'\le\vepc{\mu'}/{\mu}$. Thus, we get
\Be
\label{g11}
\begin{aligned}
    \mathbf G_{1,1}&=-\frac \bc\sul (\ba+2)(r-\rho)^2+O(\frac{\vepc\mu'}{\mu}) \\
    &=\frac{-\cos{S_c}(r-\rho)^2}{\sin{S_c}\sin{\sxyl}|x-y|}\frac{|x+y|}{\inp xy}+O(\frac{\vepc\mu'}{\mu}).
\end{aligned}
\Ee

We now consider  $\mathbf F_{1,1}$. Using $|x|^2+|y|^2=r^2+\rho^2+h^2$, we  write 
\[
    \mathbf F_{1,1}=\sul\,\partial_s^2\mathcal P(x,y,\sxyl)\,\bc^2\big(2\ba(|x|^2+|y|^2)-(4+\ba^2)r\rho\big)
    -2\ba h^2\,\sul\,\partial_s^2\mathcal P(x,y,\sxyl)\,\bc^2.
\]
Observe that the second term in the right hand is $O(\vepc\mom^{3/2}).$ Indeed, 
since $S_c^l\sim \mu$  and  $\mathcal R(x,y,\cos  \sxy(x,y) )=0$, by the mean value theorem we have
\begin{equation}
\label{r-size} \mathcal R(x,y,\cos  \sxyl )=O(2^{-l}\mu^\frac32),  \ \  (x,y)\in Q_k\times B'. \end{equation}
Hence, recalling \eqref{interval} and  \eqref{ph-dd},   we see that $ \sin\sxyl\partial^2_s\mathcal P(x,y,\sxyl)=O(\vepc\mmt)$. 
Since $\bc\sim\mu^{-\frac 32}$, it follows $h^2\,\sul\,\partial_s^2\mathcal P(x,y,\sxyl)\,\bc^2=O(\vepc\mom^{3/2})$. 
Using $r\rho = \inp xy$ and \eqref{def-a-c},  we  have $2\ba(|x|^2+|y|^2)-(4+\ba^2)r\rho=|x-y|^2|x+y|^2/\inp xy= \cos^2 \sxy/(\bc^2 \sin^2 \sxy \inp xy)$. So, combining  this with \eqref{ph-dd}, we have
\begin{align*}
  \mathbf F_{1,1}&=\frac{\sul\,\partial_s^2\mathcal P(x,y,\sxyl)\cos^2\sxy}{\inp xy  \sin^2\sxy}+O(\vepc\frac{\mu'}{\mu})     
  \\&=
-\frac{\mathcal R(x,y,\cul)\cos^2\sxy}{\sult\sut}+O(\vepc\frac{\mu'}{\mu}).
\end{align*}
We put together this and \eqref{g11} for  which we  use $(r-\rho)^2=|x-y|^2-h^2$ to get 
    \begin{align*}
    \mathbf M_{1,1}& =1+\mathbf F_{1,1}+\mathbf G_{1,1}= A+ B +O(\vepc\mom),
    \end{align*}
    where
    \begin{align*} 
   A&:= 1-\frac{|x+y||x-y| \cos \sxy}{\inp xy \sul\su}-\frac{\mathcal R(x,y,\cul)\cos^2\sxy}{\sult\sut}, \\
                        B&:=  \frac{h^2}{|x-y|}  \frac{|x+y|\cu}{\sul\su \inp xy}.
\end{align*}
We first deal with $B$. Since $\sul\sim \su\sim \smu$ and $\sul-\su=\oo {2^{-l}}$, 
\begin{align*}
 B &=  \frac{h^2 |x+y| \cu}{ |x-y| \inp xy \sut} 
                                                   +O(\vepc\big(\frac{\mu'}{\mu}\big)^{\frac32}) 
                                                   =  \frac{h^2}{|x-y|^2}  +O(\vepc\big(\frac{\mu'}{\mu}\big)^{\frac32}).
    \end{align*}
For the second equality  we use the identity 
\begin{align}\label{cossin}
\frac\sut{\cos S_c}= 
  \frac{|x+y||x-y|}{\inp xy}=\tau^+-\cos\sxy,
\end{align} 
which follows from \eqref{sc-} and \eqref{sc}. Here $\tau^+$ is given by \eqref{roots}. 
Thus, to complete the proof of \eqref{m11},  we only have to show 
\begin{align}\label{m11error}
   A=O(\vepc\mom).
\end{align}
By the definition of $S_c$ (see \eqref{sc-}) $\mathcal R(x,y,\cul)=(\cul-\cu)(\cul-\tau^+)$.  Using  this and \eqref{cossin}, we may write $A $ as follows: 
\begin{align*}
     \frac{\cu(\tau^+-\cu)\sul (\sul-\su)}{\sut\sult}
                         -\frac{\cos^2 S_c(\cul-\cu)
                       (\cul-\tau^+)}{\sut\sult}.
                       \end{align*}
                       Since $\cul-\cu=O(\smu  2^{-l})$ and $\sul-\su=O(2^{-l})$, we may 
                       replace the second appearing $\cu$ and the first appearing $\sul$ with
                       $\cul$ and $\su$, respectively,  allowing $O(\vepc\mu'/\mu)$ error.   Thus,  
\begin{align*}
     A                                       &=               
   \frac{\cu(\tau^+-\cul)}{\sut\sult}  \Big(\su (\sul-\su)
                        + \cu(\cul-\cu)\Big) +O(\vepc\frac{\mu'}{\mu}).
                        \end{align*}
                        Using  elementary trigonometric identities, we have 
\[             A        =          
 -  \frac{\cu(\tau^+-\cul)}{\sut\sult}  2\sin^2{\Big(\frac{\sxyl-S_c}{2}\Big)}+O(\vepc\frac{\mu'}{\mu}) .
                       \]
Since $S_c^l-S_c=O(2^{-l})$, by  \eqref{interval}  the estimate \eqref{m11error} follows. 
\end{proof}

\begin{proof}[Proof of \eqref{m12}--\eqref{m22}] We start with observing that  the contributions from $\mathbf F$ are negligible. Indeed, 
\[ \mathbf F_{1,2}=O(\vepc \mu'/\mu), \quad \mathbf F_{2,1}=O(\vepc \mu'/\mu),  \quad \mathbf F_{2,2}=O(\vepc \mom^\frac32). \] 
For the first we  use $2-\ba=O(|x-y|^2)=O(\mu^2)$ and   $  \sul \, \partial_s^2\mathcal P(x,y,\sxyl) \, \bc^2=O(\vepc \mu^{-\frac52}(\mu')^\frac12)$ which follows from \eqref{r-size}. The others can be handled similarly.  We may discard  $\mathbf{F}_{1,2}$, $\mathbf{F}_{2,1}$, and $\mathbf{F}_{2,2}$, which are harmless. 
Thus, in order to show \eqref{m12}--\eqref{m22} it is sufficient to show $\mathbf{G}_{1,2}$, $\mathbf{G}_{2,1}$, and $1+\mathbf{G}_{2,2}$ satisfy \eqref{m12}--\eqref{m22}, respectively, 
in place of $\mathbf{M}_{1,2}$, $\mathbf{M}_{2,1}$, and $\mathbf{M}_{2,2}$.

We first consider the term $\mathbf{G}_{1,2}$. We write
\begin{align*}
    \mathbf G_{1,2}=-2\frac{h\bc}\sul\big((1+\cu)\rho-2r\big)-\frac{\rho h\bc}\sul (\ba-2)-\frac{2\rho h\bc}\sul(\cul-\cu).
\end{align*}
Since $\ba-2=O(\mu^2),\,\cul-\cu=O(\vepc\sqrt{\mu'\mu})$,  the second and third terms on the right-hand side  are bounded 
by $O(\vepc\mom^\frac12)$. So, we only need to consider the first term. 
Using \eqref{def-a-c} and \eqref{cossin},  we see 
\[
\frac{4\bc}\sul =\frac{4\inp xy}{|x-y|^2|x+y|^2}\frac{\su}{\sul}  =\frac{1}{|x-y|^2}\frac{\su}{\sul}+O(1).
\]
Since $\sul\sim \su\sim \smu$ and $\sul-\su=\oo {2^{-l}}$,  we get 
\begin{equation}
\label{a-com}
\frac{4\bc}\sul=\frac{1}{|x-y|^2}+O(2^{-l}\mu^{-\frac52})+O(1).
\end{equation}
Using  this, the identity $(1+\cu)\rho-2r=O(\mu)$,  and \eqref{interval}, we get
\begin{align*}
    \mathbf G_{1,2}=\frac{-h}{2|x-y|^2}\big((1+\cu)\rho-2 r\big)+O(\vepc\mom^\frac12)
\end{align*}
and hence \eqref{m12}.  

We can handle $\mathbf G_{2,1}$ in the  similar manner, so we shall be brief.  Since $\ba-2=O(\mu^2)$ and $(1+\cul)\rho-2\cul r=O(\mu)$,  
using \eqref{a-com} and $\cul-\cu=O(\vepc\mmt)$, we have 
\begin{align*}
    \mathbf{G}_{2,1}   =\frac{-h}{2|x-y|^2}\big((1+\cu)\rho-2\cu r\big)+O(\vepc\mom^\frac12).
   \end{align*} 
This gives \eqref{m21}. 

Finally, we consider $1+ \mathbf G_{2,2}$. Since $\ba+2\cul-4=O(\mu)$ and $ {h^2\bc}/\sul=O(\mu'/
\mu)$, 
\begin{align*}
    1+\mathbf G_{2,2}&=1-\frac {4h^2\bc}\sul +O(\mu')
   =1-\frac{h^2}{|x-y|^2}+O(\vepc\frac{\mu'}{\mu})  \\
   &=\frac{(r-\rho)^2}{|x-y|^2}+O(\vepc\frac{\mu'}{\mu}).
\end{align*}
For the second equality we use  \eqref{a-com} and \eqref{interval}. Thus, we get \eqref{m22}.
\end{proof}

\section{Estimates with $p,q$ off the line  of duality over  $A_\mu^\pm$}
\label{off-d}
In this section we study the estimate  \eqref{est-ann}. We expect that the estimate \eqref{est-ann} should be true for $1\le p\le 2\le q\le \infty$ with possible exceptions of some endpoint cases. While complete characterization of such estimates  remains open, we verify  that the  estimate  \eqref{est-ann} holds true  on a large set of exponents.

As seen in the previous sections, compared with $\chi_{\lambda,\mu}^+ \Pi_\lambda\chi_{\lambda,\mu}^+$,  the estimate for $\chi_{\lambda,\mu}^- \Pi_\lambda\chi_{\lambda,\mu}^-$ is easier to prove.  In fact, the bound on $\|\chi_{\lambda,\mu}^- \Pi_\lambda\chi_{\lambda,\mu}^-\|_{p\to q} $ is much smaller than that on 
\eqref{est-ann}. In  Proposition \ref{prop:estext}  we obtain the sharp bound on  $\|\chi_{\lambda,\mu}^- \Pi_\lambda\chi_{\lambda,\mu}^-\|_{p\to q} $ for $1\le p\le 2\le q\le \infty$.

In order to state our result we need additional notations. Let 
\begin{align*}
    \mathfrak G=\left(\frac{2d^2+7d-7}{2(2d-1)(d+1)},\frac{2d-3}{2(2d-1)}\right).
\end{align*}
Note that the line segment $[\mathfrak A, (5/6,1/6)]$ and the line $x-y=2/{(d+1)}$ meet each other at $\mathfrak G$ if $d\ge 3$, and 
$\mathfrak G=\mathfrak G'=\left( 5/6,1/6\right)$ if $d=2$.  See Figure \ref{xyx} and Figure \ref{xyxz}.

Recalling Definition \ref{points-abcd}, we define regions $\mathfrak L_1$, $\mathfrak L_2 $, and $\mathfrak L_3 \subset [1/2,1]\times[0,1/2]$.  For $d\ge 3$, 
let $\mathfrak L_1$ be the closed pentagon with vertices $(1/2, 1/2),  \mathfrak A,  \mathfrak G,  \mathfrak G', \mathfrak A'$ from which  two points $\mathfrak G$ and $\mathfrak G'$ are removed,   $\mathfrak L_2$  be the closed triangle with vertices $\mathfrak A , (1,1/2),  \mathfrak D$ 
from which the point $\mathfrak D$ is deleted,   and $\mathfrak L_3$ be  the closed pentagon given  by vertices $(1,0),  \mathfrak D,  \mathfrak G,  \mathfrak G', \mathfrak D'$  excluding 
  $\mathfrak G$, $\mathfrak D$, $\mathfrak G'$, and $\mathfrak D'$ (Figure \ref{xyxz}). If $d=2$, $\mathfrak L_1$ becomes the closed square with vertices $(1/2, 1/2),  \mathfrak A,  \mathfrak G=\mathfrak G', \mathfrak A'$,  and  $\mathfrak L_3$ is  the closed quadrangle with vertices $(1,0),  \mathfrak D,  \mathfrak G=\mathfrak G', \mathfrak D'$ from which the points  $\mathfrak D$ and $\mathfrak D'$ are deleted (Figure \ref{xyx}).

\begin{figure}
\centering
\begin{minipage}{.4\textwidth}
  \centering
 \begin{tikzpicture}[scale=0.6]
\draw[<->] (0,7) node[left]{$\frac1q$}--(0,0)--(6.8,0)node[right]{$\frac1p$};
\draw (0,0) rectangle (6,6);
\draw (0,6)--(0,2)--(4,2)--(4,6); 
\draw (0,2)--(3,0); 
\draw (4,6)--(6,3); 
\draw (4,2)--(3,0); \draw (4,2)--(6,3); 
\node[left] at (0,2) {$\mathfrak A'$}; \node[above] at (4,6) {$ \mathfrak A$}; 
\draw (3,0) circle [radius=0.07]; \node[below] at (3,0) {$\mathfrak D'$};
\draw (6,3) circle [radius=0.07]; \node[right] at (6,3) {$\mathfrak D$}; 
\node[] at (2,4) {$\mathfrak L_1$};
\node[] at (16/3,5) {$\mathfrak L_2$}; \node[] at (1,2/3) {$\mathfrak L_2'$};
\node[] at (5,1) {$\mathfrak L_3$};
\node[anchor=north west] at (3.9,2.1) {$\mathfrak G$}; 
\node[left] at (0,0) {$0$};\node[below] at (0,0) {$\frac12$};
\node[left] at (0,6) {$\frac12$};\node[below] at (6,0) {$1$};
\end{tikzpicture}
\caption{when $d=2$. }
\label{xyx}
\end{minipage}
\hspace{30pt}
\begin{minipage}{.4\textwidth}
  \centering
 \begin{tikzpicture}[scale=0.6]
 
\draw[<->] (0,7) node[left]{$\frac1q$}--(0,0)--(6.8,0)node[right]{$\frac1p$};
\draw (0,0) rectangle (6,6);
\draw (0,6)--(0,3)--(12/5,12/5) --(18/5,18/5)--(3,6); 
\draw (0,3)--(2,0); 
\draw (3,6)--(6,4); 
\draw (12/5,12/5)--(2,0); \draw(18/5,18/5)--(6,4); 
\node[left] at (0,3) {$\mathfrak A'$}; \node[above] at (3,6) {$ \mathfrak A$}; 
\draw (2,0) circle [radius=0.07]; \node[below] at (2,0) {$\mathfrak D'$};
\draw (6,4) circle [radius=0.07]; \node[right] at (6,4) {$\mathfrak D$}; 
\node[] at (1.5,4.5) {$\mathfrak L_1$};
\node[] at (5,16/3) {$\mathfrak L_2$}; \node[] at (2/3,1) {$\mathfrak L_2'$};
\node[] at (4.2,1.8) {$\mathfrak L_3$};
\draw (12/5,12/5) circle [radius=0.07]; \node[right] at (12/5,12/5) {$\mathfrak G'$};
\draw (18/5,18/5) circle [radius=0.07]; \node[below] at (18/5,18/5) {$\mathfrak G$}; 
\node[left] at (0,0) {$0$};\node[below] at (0,0) {$\frac12$};
\node[left] at (0,6) {$\frac12$};\node[below] at (6,0) {$1$};
\end{tikzpicture}
\caption{when $d\ge 3$. }
\label{xyxz}
\end{minipage}
\end{figure}

\begin{thm}\label{thm-annest}  Let $d\ge 2$. 
Let $(1/p,1/q)\in(\mathfrak L_1\cup\mathfrak L_2\cup\mathfrak L_2'\cup\mathfrak L_3)$ and $\lambda^{-\frac23}\le\mu\le 1/4$.
Then, there is a constant $C$, independent of  $\lambda,\mu$, such that
\eqref{est-ann} holds.  Additionally,  $(i)$ we have  
\[
    \|\chi_{\lambda,\mu}^+ \Pi_\lambda\chi_{\lambda,\mu}^+ f\|_{L^{\frac{2d}{d-1},\infty}}\lesssim (\lambda\mu)^{\frac{d-3}{4}}\|f\|_{1}.\]
    \end{thm}

As is to be shown below, the bounds in Theorem \ref{thm-annest} are sharp, see Proposition \ref{lower-mu}. 
It is possible to obtain the estimate for $(1/p,1/q)$ contained in the triangle $\Delta \mathfrak A\mathfrak D\mathfrak G$ and $(\Delta \mathfrak A\mathfrak D\mathfrak G)'$ by interpolating 
the  sharp estimate  \eqref{est-ann}  with $\ppq\not\in \Delta \mathfrak A\mathfrak D\mathfrak G\cup (\Delta \mathfrak A\mathfrak D\mathfrak G)'$.  
However, those estimates  obtained  by interpolation  are not sharp, so they are not included in the statement of Theorem \ref{thm-annest}.  
The exponents $\beta(p,q)$ and $\gamma(p,q)$ have the same regimes of  change.   
We suspect that the bounds on $\|\chi_{\lambda,\mu}^+ \Pi_\lambda\chi_{\lambda,\mu}^+\|_{p\to q}$ may similarly behave as the local estimates in Theorem \ref{thm-locest}.

\subsection{Off-diagonal estimate for $\chi_\mu^+ \fP_\lambda [ \psi_j ]  \chi_{\mu}^+$.}   In order to prove Theorem \ref{thm-annest},  we make use 
of the $TT^*$ argument in Section 2. For the purpose we first consider estimate for $\chi_\mu^+ \fP_\lambda [ \psi_j ]  \chi_{\mu}^+$ with $q=p'$.

\begin{prop} 
\label{prop:mm-comparable} Let $0<\mu \le 1/4$ and $\sigma=+$. If $0\le \dpp< \frac23$, we have 
\begin{equation}
\label{diagonal}
 \|  \chi_\mu^\sigma \fP_\lambda [ \psi_j ]  \chi_{\mu}^\sigma \|_{p\to p'}\lesssim  \lambda^{\scaleppp} 2^{(\frac{d+1}2\dpp-1)j} \mu^{-\frac12\dpp}.
\end{equation}
\end{prop}

\begin{proof}  
We set,  for simplicity, 
\[\cB_{p,j}:= \lambda^{\scaleppp} 2^{(\frac{d+1}2\dpp-1)j} \mu^{-\frac12\dpp}.\]
To show \eqref{diagonal}, we  make use of the decomposition \eqref{decomp} with $\mu=\mu'$. By Lemma \ref{large-sep} and Lemma \ref{small-sep}, we may assume 
$2^{-\nu}\le 1/100$ and  it is sufficient to show 
\Be
\label{bbjj}
\| \sum_{ \nu\le \nu_\circ}\sum_{k\sim_\nu k'} \chi_{k}^{\sigma,\nu} \fP_\lambda [\psi_j] \chi_{k'}^{\sigma,\nu}  \|_{p\to p'} \lesssim  
\bp
\Ee
while assuming  \eqref{AA}. The contribution from the other cases is minor.   
As before we show \eqref{bbjj} considering the cases  $2^{-\nu}\gg   \mu$ and $2^{-\nu}\lesssim \mu$, separately.

We first consider the sum over $\nu: 2^{-\nu}\gg \mu$.  If $\mu\gtrsim 2^{-2j} $,  by the first case estimate in \eqref{eq:haha8} and  Lemma \ref{kkpsum} we get
\[ \|\sum_{2^{-\nu}\gg \mu} \sum_{k\sim_\nu k'} \chi_{k}^{\sigma,\nu} \fP_\lambda [\psi_j] \chi_{k'}^{\sigma,\nu} \|_{p\to p'} 
          \lesssim   
                                 \lambda^{\frac{d-1}2\dpp-\frac d2} 2^{(\frac{d-1}2\dpp -1)j}\mu^{-\dpp}\lesssim \bp.
                                 \]
When $\mu\ll  2^{-2j} $, splitting the sum $\sum_{\nu: 2^{-\nu}\gg \mu} =\sum_{\nu: 2^{-\nu}\gtrsim 2^{-2j}}+\sum_{\nu: 2^{-2j}\gg 2^{-\nu}\gg \mu} $, the first sum is readily  handled by using the estimate  \eqref{eq:haha8} as in the above.  So we only need to show 
 \Be
 \label{eq:temp} \Big \|\sum_{\nu: 2^{-2j}\gg 2^{-\nu}\gg \mu} \sum_{k\sim_\nu k'} \chi_{k}^{\sigma,\nu} \fP_\lambda [\psi_j] \chi_{k'}^{\sigma,\nu}\Big \|_{p\to p'}\lesssim \bp.\Ee
However,   direct summation along $\nu$ using the second case  in \eqref{eq:haha8} produces unwanted additionally $\log (1/\mu)$ factor. 
We get around this  by using \eqref{eq:l2l2l2} and the following  observation, which is an easy consequence of 
disjointness of the sets $\{A_{k}^{\sigma,\nu} \times A_{k'}^{\sigma,\nu} \}_{k\sim_\nu k', \nu}$ and the second case estimate in  \eqref{eq:haha8}: 
\[\|\sum_{\nu: 2^{-2j}\gg 2^{-\nu}\gg \mu} \sum_{k\sim_\nu k'}\chi_{k}^{\sigma,\nu} \fP_\lambda [\psi_j] \chi_{k'}^{\sigma,\nu} \|_{1\to \infty} 
          \lesssim  \lambda^{-\frac12} 2^{\frac {d+1}2 j}. \] 
Interpolating this and the $L^2$--$L^2$ bound which follows from \eqref{eq:l2l2l2},  we see the left hand side of 
          \eqref{eq:temp} is bounded by 
\[C \lambda^{\frac{d-1}2\dpp-\frac d2} 2^{(\frac{d-1}2\dpp +1)j}  \mu^{1-\dpp}\] for $1\le p\le 2$. 
Since this is bounded by  $C \bp$,   the desired bound  \eqref{eq:temp}   therefore follows.

To show the estimate for the sum over $\nu: 2^{-\nu_\circ}\le 2^{-\nu}\lesssim \mu$, by Lemma \ref{kkpsum} it is sufficient to show 
\Be
\label{comcom} \| \chi_{k}^{\sigma,\nu} \fP_\lambda [\psi_j] \chi_{k'}^{\sigma,\nu} \|_{p\to p'} \lesssim     \bp, \quad k\sim_\nu k' 
\Ee 
for each $\nu$ satisfying  $2^{-\nu_\circ}\le 2^{-\nu}\lesssim \mu$ since there are only $O(1)$ as many as  $\nu$.  If $\nu=\nu_\circ$ we have   $|\cD(x,y)|\sim  \mu^2$ as is clear from \eqref{angle}. Hence,  we get  \eqref{comcom} for $1\le p\le 2$ interpolating  
the $L^1$--$L^\infty$ estimate \eqref{apjb} (Corollary  \ref{pj-det}) and  $\| \chi_{k}^{\sigma,\nu} \fP_\lambda [\psi_j] \chi_{k'}^{\sigma,\nu} \|_{2\to 2}\lesssim \lambda^{-\frac d2}2^{-j}$ which follows from \eqref{easy-l2}.  So, it remains to show \eqref{comcom} for the case 
 $2^{-\nu_\circ}<2^{-\nu}\lesssim \mu$.

 We show the estimate \eqref{comcom}    by considering   the three cases  $2^{-j} \gg  \smu $,  $2^{-j}\sim \sqrt \mu$, and  $2^{-j} \ll  \smu$,  separately. 
For the first case, using $(b2)$ in  Lemma \ref{case-b},  we have the same bound  as in the third case of  \eqref{eq:sec8}. Consequently, the desired bound \eqref{comcom} follows
since $ \lambda^{\frac{d-1}2\dpp-\frac d2} 2^{(\frac {d+3}2\dpp  -1) j}\lesssim \bp$.   
We need to handle the remaining two cases $2^{-j}\sim \sqrt \mu$ and  $2^{-j} \ll  \smu$. 

Noting that  the sets $A_{k}^{+, \nu}$ and $ A_{k'}^{+,\nu}$, $\ksim$  are  contained in cubes of side length $\sim \mu$ 
and distanced about $\sim \mu$ from each other, 
 we  break   $A_{k}^{+,\nu}$ and $ A_{k'}^{+,\nu}$ into finitely many (essentially disjoint) cubes 
$B, B'$ of side length $\vepc \mu$. Indeed, we may assume 
\[
A_{k}^{+,\nu}=\bigcup B, \quad  A_{k'}^{+,\nu}=\bigcup B' 
\] (for example see \eqref{bb-decomp}), and 
 for \eqref{comcom} it is enough to show
\Be
\label{bbbb} \|\chi_B\mathfrak P_\lambda [ \psi_j ]  \chi_{B'}\|_{p\to p'} \lesssim     \bp. \Ee 
By Lemma \ref{dychotomy} we may assume either \eqref{large-d} or \eqref{small-det} with $\mu=\mu'$. 
If \eqref{large-d} holds, i.e., $\cD(x,y)\gtrsim \vepc \mu^2$ for $(x,y)\in B\times B'$, the estimate \eqref{bbbb} follows from Corollary  \ref{pj-det} and \eqref{easy-l2} as before.  
Hence, we may assume $\cD(x,y)\ll \vepc \mu^2$ for $(x,y)\in B\times B'$.  Since $2^{-\nu_\circ}< 2^{-\nu}\lesssim \mu$ and $\cD(x,y)\ll \vepc \mu^2$, by $(b1)$ in Lemma \ref{case-b} 
we have 
\[
 \|\chi_B\mathfrak P_\lambda [ \psi_j ]  \chi_{B'}\|_{1\to \infty} \lesssim     2^{\frac{d-1}{2}j}(\lambda\mu)^{-\frac12},  \quad 2^{-j}\ll \smu. 
\]
Interpolating this and  $\|\chi_B\mathfrak P_\lambda [ \psi_j ]  \chi_{B'}\|_{2\to 2}\lesssim \lambda^{-\frac d2}2^{-j}$ gives the desired estimate  \eqref{bbbb} if $2^{-j}\ll \smu$.

Finally,  we deal with 
the case  $2^{-j} \sim  \smu$. The desired bound \eqref{bbbb} follows from 
Proposition \ref{critical-muj} if  $4/(3d+4)<\dpp<2/3$.  We interpolate the estimate with  
  $ \|\chi_B\mathfrak P_\lambda [ \psi_j ]  \chi_{B'}\|_{2\to 2} \lesssim  \lambda^{-\frac d2} 2^{-j}$ to get  \eqref{bbbb}  for   $0\le \dpp<2/3$. This completes the proof. 
\end{proof}

\begin{rem} 
\label{about-bbjj} In the above we have actually shown \eqref{bbjj}  for $0\le \dpp\le 1$ when  $2^{-j}\gg \smu$, or  $2^{-j}\ll \smu$. However, 
we have  the estimate \eqref{bbjj}  on the restricted range $0\le \dpp< 2/3$ when   $2^{-j}\sim \smu$. 
\end{rem}
By means of  Proposition 
\ref{prop:mm-comparable} we obtain the following, which we use for the proof of Theorem \ref{thm-annest}.

\begin{prop} 
\label{mm-comparable} Let $d\ge 3$,  $\mu \le 1$ and $\sigma=+$. If $\ppq\in \mathfrak L_1$, 
 we have
\begin{equation}
\label{mmm0}
 \|  \chi_{\lambda,\mu}^\sigma \Pi_\lambda \chi_{\lambda,\mu}^\sigma  \|_{p\to q}
     \lesssim   
                                 \lambda^{-\frac12\dpq} \mu^{-\frac{d+3}4\dpq+\frac12} .
             \end{equation}
Moreover, we have restricted weak type $(p,q)$ estimate for $\chi_{\lambda,\mu}^\sigma \Pi_\lambda \chi_{\lambda,\mu}^\sigma$ if $(1/p, 1/q)=\mathfrak G, \mathfrak G'$ for $d\ge 3$. \end{prop}

\begin{proof}
We  basically follow  the same argument in the {proof of Lemma  \ref{tt-st}}, so we shall be brief.  Since $\|\Pi_\lambda[\psi_j] \chi_{\lambda,\mu}^\sigma f\|_{2}^2=\inp{\Pi_\lambda[\psi_j]^*\Pi_\lambda[\psi_j] \chi_{\lambda,\mu}^\sigma f}{\chi_{\lambda,\mu}^\sigma f}$, 
we have
\begin{align*}
\|\Pi_\lambda[\psi_j] \chi_{\lambda,\mu}^\sigma f\|_{2}^2\le \sum_{k: 2^{-k} \lesssim 2^{-j}}
\| \chi_{\lambda,\mu}^\sigma (\Pi_\lambda[\psi_j]^*\Pi_\lambda[\psi_j])_k \chi_{\lambda,\mu}^\sigma f\|_{p'}\|f\|_p,
\end{align*}
where $(\Pi_\lambda[\psi_j]^*\Pi_\lambda[\psi_j])_k$ is defined by \eqref{pk}.
We now recall  \eqref{bk}. By Proposition \ref{prop:mm-comparable} (via scaling) and Minkowski's inequality
we have 
\[\| \chi_{\lambda,\mu}^\sigma  (\Pi_\lambda[\psi_j]^*\Pi_\lambda[\psi_j])_k  \chi_{\lambda,\mu}^\sigma\|_{p\to p'} \lesssim  2^{-j} \lambda^{-\frac12\dpp}   
                                   2^{(\frac{d+1}{2}\dpp-1)k } \mu^{-\frac12\dpp} 
                               \] for $0\le \dpp< 2/3$.
Thus, taking sum over $k$, we get   
$\|\Pi_\lambda[\psi_j]  \chi_{\lambda,\mu}^\sigma   \|_{p\to 2}\lesssim \lambda^{-\frac12\delta(p,2)}   
2^{(\frac{d+1}{2}\delta(p,2)-1)j } \mu^{-\frac12\delta(p,2)}$. This clearly implies 
\[
\|  \chi_{\lambda,\mu}^\sigma \Pi_\lambda[\psi_j]  \chi_{\lambda,\mu}^\sigma  \|_{p\to 2}\lesssim \lambda^{-\frac12\delta(p,2)}   
2^{(\frac{d+1}{2}\delta(p,2)-1)j } \mu^{-\frac12\delta(p,2)}
\]
if $\frac{d+1}{2}\dpp< 1$. As before,  the missing  case $p=2(d+1)/(d+3)$ can be recovered 
 thanks to the bilinear interpolation argument (e.g., Keel and Tao \cite{keel-tao}), also see the proof of Lemma \ref{tt-st}. 
Duality and interpolation with \eqref{diagonal} (after rescaling)  yield 
\begin{equation}
\label{interpol}
\| \chi_{\lambda,\mu}^\sigma \Pi_\lambda[\psi_j] \chi_{\lambda,\mu}^\sigma\|_{p\to q}\lesssim \lambda^{-\frac12\dpq}   2^{(\frac{d+1}{2}\dpq-1)j }  \mu^{-\frac12\dpq}
\end{equation}           
provided that $\ppq$ is contained in the closed quadrangle with vertices  $(1/2, 1/2)$,  $\mathfrak A$,  $\mathfrak A'$,  and $(5/6, 1/6)$ from which the point  $(5/6, 1/6)$ is excluded. 
Thus,  by  $(c)$ in  Lemma \ref{s-trick} we get  
\[ \|\sum_{j\ge 4}\chi_{\lambda,\mu}^\sigma \Pi_\lambda[\psi_j] \chi_{\lambda,\mu}^\sigma f \|_{q,\infty}\lesssim \lambda^{-\frac12\dpq}     \mu^{-\frac12\dpq} \|f\|_{p,1}\]
 for $(1/p, 1/q)=\mathfrak G, \mathfrak G'$.\footnote{Note that the line segments $[\mathfrak A, (5/6, 1/6)]$, $[\mathfrak A', (5/6, 1/6)]$ meet
the line $x-y=\frac 2{d+1}$ at $\mathfrak G$, $\mathfrak G'$, respectively.}  (See Figure \ref{xyxz}.)  Thus, by \eqref{symmetric} we get the same estimate for  the operator $\sum_{j}\chi_{\lambda,\mu}^\sigma \Pi_\lambda[\psi_j^\kappa] \chi_{\lambda,\mu}^\sigma$, $\kappa=-, \pm\pi$.  Since the estimate \eqref{interpol} remains valid with $\kappa=0$,  using \eqref{decomp-proj}  and combining all these estimates,  we get 
the restricted weak type $(p,q)$ estimate  for $ \chi_{\lambda,\mu}^\sigma \Pi_\lambda \chi_{\lambda,\mu}^\sigma$ with $(1/p, 1/q)=\mathfrak G, \mathfrak G'$.  Here we also use  the fact that 
$-\frac12 \dpq=-\frac{d+3}4 \dpq+\frac12$ if $(1/p, 1/q)=\mathfrak G, \mathfrak G'$. 
Real interpolation among these estimates  gives \eqref{mmm0}  for $\ppq\in (\mathfrak G, \mathfrak G')$.  

We may further interpolate those consequent estimates with \eqref{l2-ktt} to get 
\eqref{mmm0} on an extended range, particularly, for $p,q$ satisfying $1/p+1/q=1$ and $0\le \dpq\le 2/(d+1)$.  Thus,  we get 
$\|\Pi_\lambda \chi_{\lambda,\mu}^\sigma\|_{p\to 2}\lesssim    \lambda^{-\frac12\delta(p,2)} \mu^{-\frac{d+3}4\delta(p,2)+\frac14}$  for $\frac{2(d+1)}{d+3} \le p\le 2$.
Since  $\| \chi_{\lambda,\mu}^\sigma \Pi_\lambda \chi_{\lambda,\mu}^\sigma\|_{p\to 2}\le \| \chi_{\lambda,\mu}^\sigma \Pi_\lambda\|_{2\to 2} \|\Pi_\lambda \chi_{\lambda,\mu}^\sigma\|_{p\to 2}$, 
using $
 \|   \chi_{\lambda, \mu}^\sigma \Pi_\lambda   \|_{2\to 2} \lesssim   \mu^\frac14
$(see \eqref{l2-ktt}), we obtain 
\[
\|\chi_{\lambda,\mu}^\sigma\Pi_\lambda \chi_{\lambda,\mu}^\sigma\|_{p\to 2}\lesssim    \lambda^{-\frac12\delta(p,2)} \mu^{-\frac{d+3}4\delta(p,2)+\frac12}
\]
 for $\frac{2(d+1)}{d+3} \le p\le 2$. This gives \eqref{mmm0} for $\ppq\in [(1/2, 1/2),  \mathfrak A ]$. Thus, by duality we also get 
 \eqref{mmm0} for $\ppq\in [(1/2, 1/2),  \mathfrak A' ]$. Interpolation between  these estimates and
  \eqref{mmm0} for $\ppq\in (\mathfrak G,  \mathfrak G' )$, which we have already shown, 
 we get \eqref{mmm0}  for  $\ppq\in \mathfrak L_1$.
 \end{proof}

\subsection{Weak type estimate for $\mu$-local case} 
In point of view of the $L^1$--$L^{\frac{2d}{d-1},\infty}$ estimate in Theorem \ref{thm-locest} the estimates in the next proposition  seem to be natural. 
Also, it is not difficult to see that the estimates are optimal, see  Proposition \ref{lower-mu}.

\begin{prop}  \label{prop:weakmu} Let $d\ge 2$. Let $r\le 1-\delta_\circ$  for some $\delta_\circ>0$. Then, we have  
\Be
\label{est:brweak} \|\chi_{\sqrt\lambda B_r}  \Pi_{\lambda}\chi_{\sqrt \lambda B_r}  \|_{L^{\frac{2d}{d+1},1} \rightarrow L^\infty}\lesssim \lambda^{\frac{d-3}{4}}.
\Ee
If  $ \lambda^{-\frac 23}\lesssim \mu \ll 1$,  we  have
\Be
\label{est:muweak}\|\chi_{\lambda,\mu}^\sigma\Pi_\lambda \chi_{\lambda,\mu}^\sigma\|_{L^{\frac{2d}{d+1},1} \rightarrow L^\infty}\lesssim (\lambda\mu)^{\frac{d-3}{4}}, \quad 
\sigma=+. 
\Ee
\end{prop}

By virtue of \eqref{est:muweak}  and Proposition \ref{mm-comparable},  we now can prove Theorem \ref{thm-annest}. 

 \begin{proof}[Proof of Theorem \ref{thm-annest}]  By duality the estimate \eqref{est:muweak}    proves the statement $(i)$ in Theorem \ref{thm-annest}, so we need only to prove 
the $L^p$--$L^q$ estimate for $\chi_{\lambda,\mu}^\sigma\Pi_\lambda \chi_{\lambda,\mu}^\sigma$.  When $d\ge 3$,  by 
Proposition \ref{mm-comparable}  we have \eqref{est-ann}  for $\ppq\in \mathfrak L_1$.

We observe
\Be
\label{eq:2to00}
\|\chi_{\lambda,\mu}^+ \Pi_\lambda \chi_{\lambda,\mu}^+\|_{2\to\infty}\lesssim\mu^{\frac14}(\lambda\mu)^{\frac{d-2}{4}}. 
\Ee
Since $\|\chi_{\lambda,\mu}^+ \Pi_\lambda \chi_{\lambda,\mu}^+\|_{2\to\infty} \le  \| \Pi_\lambda \chi_{\lambda,\mu}^+\|_{2\to 2}\|\chi_{\lambda,\mu}^+ \Pi_\lambda\|_{2\to \infty}$,
this follows by \eqref{kt-nu} and  $ \| \Pi_\lambda   \chi_{\lambda, \mu}^\pm \|_{2\to 2} \lesssim   \mu^\frac14$ (equivalently, \eqref{l2-ktt}). 
Then, by (real) interpolation between the estimate with $\ppq=\mathfrak A'$,  \eqref{est:muweak},  and  \eqref{eq:2to00},  
    we get \eqref{est-ann}  for $\ppq\in \mathfrak L_2'$. 
The estimates for $\ppq\in \mathfrak L_2$ follow from duality.

 For the estimate with  $\ppq\in \mathfrak L_3$  we interpolate the estimates for $\ppq=(1,0), \mathfrak D, \mathfrak G, \mathfrak G', \mathfrak D'. $
See Figure \ref{xyxz}. For $d=2$, we use the bounds in \eqref{kt-nu} together with duality. Interpolating  those estimates and  the estimates with $\ppq=(1,0), \mathfrak D, \mathfrak G= \mathfrak G', \mathfrak D'$ we get all the desired estimates. 
\end{proof}

Let us now consider the estimate for $\chi_{\sqrt \lambda B_r} \Pi_\lambda \chi_{\sqrt \lambda B_r}$, $r=1-\delc$,  which we have dealt with in Section 3.  
Since we have the estimate \eqref{eq:away-sphere00} in Proposition \ref{away-sphere} which takes the place of  \eqref{diagonal}, following the same way as  in the proof of Proposition 
\ref{mm-comparable},  one obtains the corresponding estimate for 
$\chi_{\sqrt \lambda B_r} \Pi_\lambda \chi_{\sqrt \lambda B_r}$ with bound  $C \lambda^{-\frac12\dpq}$.  Using these estimates together with  \eqref{est:brweak} (and its dual estimate)  we obtain the following, which we state without proof.

 \begin{cor}\label{local-r}  Let  $\delta_\circ>0$, $r=1-\delc$, and let $(1/p,1/q)\in(\mathfrak L_1\cup\mathfrak L_2\cup\mathfrak L_2'\cup\mathfrak L_3)$.  
Then, we have 
\[\| \chi_{\sqrt \lambda B_r} \Pi_\lambda \chi_{\sqrt \lambda B_r}\|_{p\to q}\lesssim \lambda^{\beta(p,q)}.\]   \end{cor}

\subsection{Proof of Proposition \ref{prop:weakmu}}
To show Proposition \ref{prop:weakmu}, we use the next  Lemma \ref{ker10} which  can  be shown by making use of Lemma \ref{S4-asym} and following the proof of \cite[Lemma 3.2.2]{Th93} (also, see  \cite{T05}). 

\begin{lem}
\label{ker10}  Let $ \lambda^{-\frac 23}\lesssim \mu \le 1/4$. Then, we have the following\,$:$
\vspace{-6pt}
\begin{enumerate} 
[leftmargin=0.85cm, labelsep=0.3 cm, topsep=0pt]
\item[$(i)\,$]    $|\Pi_{\lambda}(x,y)|\lesssim 
(\lambda\mu)^{\frac{d-2}{2}}$ for $x,y \in A_{\lambda, \mu}^\pm$. 
\item[$(ii)$] $|\Pi_{\lambda}(\sqrt\lambda x,\sqrt \lambda y)|\lesssim \lambda^{\frac{d-2}{6}}$ if $||x|-1|,  ||y|-1|\le \lambda^{-\frac23}$.
\end{enumerate}
\end{lem}

Making use of Lemma \ref{ker10},  we obtain the following which is  useful for getting estimate for the kernel given by spectral decomposition. 

\begin{lem}
\label{PI-est}  Let $ \lambda^{-\frac 23}\lesssim \mu \le 1/4$ and $ 2^{-j}\gtrsim  (\lambda\mu)^{-1} $. Then, there is a constant $C$, independent of  $x,$ $y$, $\lambda$, and $\mu$,  such that 
\[\sum_{\lambda'} \big(1+2^{-j}|\lambda-\lambda'|\big)^{-N}|\Pi_{\lambda'}(x,y)| \le C2^j(\lambda\mu)^{\frac{d-2}{2}}, \quad x,y \in A_{\lambda, \mu}^+.\]
\end{lem}

 Since  we also need to deal with the case $\lambda\neq \lambda'$, we can not directly apply Lemma \ref{ker10}.
To overcome this, we use an argument based on monotonicity.  To this end, let us set 
\begin{align*}
 A_{\lambda,\mu}^e &=\{x:|x|\ge\lambda^{1/2}(1-\mu)\}, 
 \\ 
 \chi_{\lambda,\mu}^e&=\chi_{A_{\lambda,\mu}^e}.
 \end{align*}
Using Lemma \ref{ker10}, \eqref{smallmu} and  \eqref{est:ext--} with $p=1$ and $q=\infty$,   it is easy to see 
\Be
\label{PI-00} 
\|\chi_{\lambda,\mu}^e \Pi_\lambda(\cdot,\cdot) \chi_{\lambda,\mu}^e\|_{L^\infty_xL^\infty_y} \lesssim    (\lambda\mu)^{\frac{d-2}{2}}, \quad  \lambda^{-\frac 23}\lesssim \mu \ll 1.
\Ee
Indeed,  the estimate \eqref{PI-00} is equivalent to $\|  \chi_{\lambda,\mu}^e \Pi_\lambda\|_{2\to \infty}\lesssim  (\lambda\mu)^{\frac{d-2}{4}}$. Since the kernel of $\Pi_\lambda$ decays rapidly if 
$|x|$ or $ |y|>1+c$ for some $c>0$, this estimate  can be shown by combining the three estimates $\|  \chi_{\lambda,\mu}^+ \Pi_\lambda\|_{2\to \infty}\lesssim  (\lambda\mu)^{\frac{d-2}{4}}$, $ \lambda^{-\frac 23}\lesssim \mu \lesssim 1$, $ \|\chi_{\lambda, \mu}^\circ
\Pi_\lambda \|_{2\to \infty }\lesssim \lambda^{\frac{d-2}{12}}$, $ \mu \sim \lambda^{-\frac 23}$, and  $ \|\chi_{\lambda, \mu}^-\Pi_\lambda \|_{2\to \infty}\lesssim \lambda^{\frac{d-2}{12}}(\lambda\mu^{3/2})^{-N}$, $ \lambda^{-\frac 23}\lesssim \mu \lesssim  1$.  These  estimates follow  from Lemma \ref{ker10}, \eqref{smallmu},  and  \eqref{est:ext--}, respectively, by taking $p=1$ and $q=\infty$.

\begin{proof}[Proof of Lemma \ref{PI-est}]
In order to show Lemma \ref{PI-est}, it is sufficient to show 
\Be\label{pi-sum}
\sum_{\lambda'} \big(1+2^{-j}|\lambda-\lambda'|\big)^{-N}\|\chi_{\lambda,\mu}^e \Pi_{\lambda'}(\cdot,\cdot) \chi_{\lambda,\mu}^e\|_{L^\infty_xL^\infty_y} \le C2^j(\lambda\mu)^{\frac{d-2}{2}}.
\Ee
To deal with the sum, for each $\lambda, \lambda'$  let us set 
\Be
\label{ell}
\ell(\rho)=
\begin{cases}
   \  \  \   (\lambda/\lambda')^{\frac23}\rho, \quad &\lambda \ge \lambda',
    \\
   \,   (\lambda'-\lambda)/\lambda+\rho, \quad &\lambda < \lambda'.
    \end{cases} \Ee
 Since  $\lambda^{-\frac23}\lesssim\mu$,  we note that $\ell(\mu) \gtrsim(\lambda')^{-\frac23}$. Since  $(\lambda')^\frac12(1-\ell(\mu))\le \lambda^\frac12(1-\mu)$, it follows  
 \Be 
\label{inclusion} 
A_{\lambda,\mu}^e\subset A_{\lambda',\ell(\mu)}^e
 \Ee 
for  $\lambda^{-\frac23}\lesssim\mu$. 
 So, we have  $  \|\chi_{\lambda,\mu}^e \Pi_{\lambda'}(\cdot,\cdot) \chi_{\lambda,\mu}^e\|_{L^\infty_xL^\infty_y}\le \|\chi_{\lambda', \ell(\mu)}^e \Pi_{\lambda'}(\cdot,\cdot) \chi_{\lambda', \ell(\mu)}^e\|_{L^\infty_xL^\infty_y}$. Thus, by \eqref{PI-00} we see that   the left hand side of \eqref{pi-sum} is bounded above by 
 \begin{align*}
 C  \sum_{\lambda'} \big(1+2^{-j}|\lambda-\lambda'|\big)^{-N} (\lambda'\ell(\mu))^{\frac{d-2}{2}}.
 \end{align*}
 To prove \eqref{pi-sum} we show the above sum is bounded by $C2^j(\lambda\mu)^{\frac{d-2}{2}}$. 
 Recalling \eqref{ell} and considering separately the cases $\lambda\ge \lambda'$ and $\lambda<\lambda'$,  it is sufficient to show
  \begin{align*}
 (\lambda\mu)^{\frac{d-2}{2}}  \sum_{\lambda'\le \lambda} \big(1+2^{-j}(\lambda-&\lambda')\big)^{-N} (\lambda'/\lambda)^{\frac{d-2}{6}}
 \lesssim  2^j(\lambda\mu)^{\frac{d-2}{2}},\qquad
   \\
     (\lambda\mu)^{\frac{d-2}{2}} \sum_{\lambda'>\lambda} \big(1+2^{-j}(\lambda'-\lambda)\big)^{-N}& (\lambda'/\lambda)^\frac{d-2}2 \big((\lambda'-\lambda)/(\lambda\mu)+1\big)^{\frac{d-2}{2}}
      \lesssim  2^j(\lambda\mu)^{\frac{d-2}{2}}.
 \end{align*}
Both follow from a simple computation. For the second inequality we use  the condition $ 2^{-j}\gtrsim  (\lambda\mu)^{-1}$.
\end{proof}

The proofs of \eqref{est:muweak} and \eqref{est:brweak} are similar. 
We first prove the estimate \eqref{est:muweak}.   The proof of  \eqref{est:brweak} is later given by modifying that of \eqref{est:muweak}.

 \begin{proof}[Proof of \eqref{est:muweak}] 
 In order to prove \eqref{est:muweak},  we first reduce the matter to showing \eqref{eq:core} by obtaining various estimates. Then, we dedicate the rest of the proof to proving \eqref{eq:core}.

 Recalling \eqref{decomp-proj}  we need to show that  \eqref{est:muweak} holds if $\Pi_\lambda$ is replaced by  $\sum_j \Pi_\lambda[\psi_j^++\psi_j^-]$,  $\sum_j \Pi_\lambda[\psi_j^{\pi}+\psi_j^{-\pi}]$, and  $\Pi_\lambda [\psi^0]$.   
  Concerning the bounds for $\sum_j \Pi_\lambda[\psi_j^++\psi_j^-]$ and $\sum_j \Pi_\lambda[\psi_j^{\pi}+\psi_j^{-\pi}]$, it is sufficient to consider only  $\sum_j \Pi_\lambda[\psi_j^++\psi_j^-]$ since $\Pi_\lambda[\psi_j^++\psi_j^-](x,y)=C\,\Pi_\lambda[\psi_j^{\pi}+\psi_j^{-\pi}](x,-y)$ with  $|C|=1$ 
   as can  be easily checked by Lemma \ref{kernel} (cf. \eqref{symmetric}).   
 We begin with dividing the sum  
 \[  
 \Pi_{m}+ \Pi_{e} :=\sum_{j : 2^{-j}\ge  (\lambda\mu)^{-1}  } \Pi_\lambda[\psi_j^++\psi_j^-] + \sum_{j : 2^{-j}\le   (\lambda\mu)^{-1}  } \Pi_\lambda[\psi_j^++\psi_j^-] .
 \]  
 The operator $\Pi_e$ is easier to deal with.  
 Since $\sum_{j : 2^{-j}\le   (\lambda\mu)^{-1}  } (\psi_j^++\psi_j^-)=\psi_\circ( (\lambda\mu) \cdot)$ for some $\psi_\circ\in C_c^\infty(-2,2)$, we note  $\Pi_{e}=\Pi_\lambda[\psi_\circ( (\lambda\mu) \cdot)]$. 
 
We handle $\Pi_e$ via spectral decomposition. In fact, using the decomposition  \eqref{spectral-proj} and Lemma \ref{PI-est}, we have
 \begin{align}
 \label{Pie100}
|\Pi_e(x,y)|\lesssim \sum_{\lambda'} (\lambda\mu)^{-1}\big(1+(\lambda\mu)^{-1}|\lambda-\lambda'|\big)^{-N}|\Pi_{\lambda'}(x,y)| \lesssim (\lambda\mu)^{\frac{d-2}{2}}
\end{align}
for $x, y \in A_{\lambda, \mu}^\sigma$.
Similarly, we recall $\Pi_\lambda(x,y)=\sum_{|\alpha|=\lambda} \Phi_\alpha(x)  \Phi_\alpha(y) $ and  use  \eqref{spectral-proj}. Then, by  orthonormality  between the Hermite functions and Lemma \ref{PI-est}  we have 
\Be
\label{Pie200}
\int |\Pi_e(x,y)|^2 dy\lesssim  \sum_{\lambda'} (\lambda\mu)^{-2}\big(1+(\lambda\mu)^{-1}|\lambda-\lambda'|\big)^{-N}|\Pi_{\lambda'}(x,x)| \lesssim (\lambda\mu)^{\frac{d-4}{2}}
\Ee
 for $x \in A_{\lambda, \mu}^\sigma$.  
 The above two estimates  \eqref{Pie100}  and  \eqref{Pie200}  give  $\| \chi_{\lambda,\mu}^\sigma\Pi_e  \chi_{\lambda,\mu}^\sigma\|_{1\to \infty}\lesssim (\lambda\mu)^{\frac{d-2}{2}}$ and 
$\|  \chi_{\lambda,\mu}^\sigma\Pi_e \chi_{\lambda,\mu}^\sigma\|_{2\to \infty}\lesssim (\lambda\mu)^{\frac{d-4}{4}}$, respectively. The second  follows from  \eqref{Pie200}  and the Cauchy-Schwarz inequality since 
\Be
\label{2-infty}
\|T\|_{2\to \infty}=\|T \|_{L_x^\infty(L_y^2)}
\Ee
 if $T(x,y)$ is the kernel of an operator $T$. 
 Interpolation between those two estimates  gives the bound
\Be
\label{hhhh}
\| \chi_{\lambda,\mu}^\sigma \Pi_e  \chi_{\lambda,\mu}^\sigma\|_{\frac{2d}{d+1}\rightarrow\infty}\lesssim(\lambda\mu)^{\frac{d-3}{4}}, 
\Ee
which is acceptable.

 Thus, in order to show \eqref{est:muweak}  we need to obtain the bound on the operator norm of  $\chi_{\lambda,\mu}^\sigma  \Pi_m \chi_{\lambda,\mu}^\sigma$ and $\chi_{\lambda,\mu}^\sigma\Pi_\lambda[\psi^0]\chi_{\lambda,\mu}^\sigma$ from $L^{\frac{2d}{d+1},1}$ to $L^\infty$.  
Thanks to \eqref{symmetric} we only have to show 
\begin{align}
 \|\sum_{j : 2^{-j}\ge  (\lambda\mu)^{-1}  }\chi_{\lambda,\mu}^\sigma  \Pi_\lambda[\psi_j^+] \chi_{\lambda,\mu}^\sigma   \|_{L^{\frac{2d}{d+1},1}\rightarrow L^\infty}
&\lesssim(\lambda\mu)^{\frac{d-3}{4}},\label{mmm-0}\\
 \|\chi_{\lambda,\mu}^\sigma  \Pi_\lambda[\psi^0] \chi_{\lambda,\mu}^\sigma   \|_{L^{\frac{2d}{d+1},1}\rightarrow L^\infty}\label{000}
&\lesssim(\lambda\mu)^{\frac{d-3}{4}}.
\end{align}
For \eqref{mmm-0} we separately consider  the cases $2^{-j}\not\sim \smu$ and $2^{-j}\sim \smu$. Here, by $2^{-j}\not\sim \smu$ we mean that  $2^{-j}\gg \smu$ or $2^{-j}\ll \smu$.  The cases $2^{-j}\not\sim \smu$ and $2^{-j}\sim \smu$ need to be handled differently.   The first case is easier.  

Indeed, we claim 
\[
\|\sum_{j: 2^{-j}\not \sim \smu,\,   2^{-j}\ge  (\lambda\mu)^{-1} }\chi_{\lambda,\mu}^\sigma\Pi_{\lambda}[\psi_j^+]\chi_{\lambda,\mu}^\sigma\|_{L^{\frac{2d}{d+1},1}\rightarrow L^\infty}\lesssim (\lambda\mu)^{\frac{d-3}{4}}.
\] 
Note $\frac{d-1}{d} \frac12+\frac{1}{d}1 =\frac{d+1}{2d}$. Then, by $(b)$  in Lemma \ref{s-trick}  the estimate is a straightforward consequence of the following estimates:  
\begin{align}
\label{10-2} 
\|\chi_{\lambda,\mu}^\sigma\Pi_{\lambda}[\psi_j^+]\|_{2\rightarrow\infty}&\lesssim 2^{-\frac{j}{2}}(\lambda\mu)^{\frac{d-2}{4}},  \quad\, 2^{-j}\gtrsim  (\lambda\mu)^{-1}, 
\\
\label{10-3} \| \chi_{\lambda,\mu}^\sigma \Pi_{\lambda}[\psi_j^+]  \chi_{\lambda,\mu}^\sigma\|_{1\rightarrow\infty}&\lesssim 2^{\frac{d-1}{2}j}(\lambda\mu)^{-\frac12}, \quad 2^{-j}\not\sim \smu. 
\end{align}
It should be noticed that the first estimate \eqref{10-2} is valid without  the assumption $2^{-j}\not\sim \smu$ and the estimate \eqref{10-2} can be shown similarly as before. Indeed, 
using \eqref{spectral-proj},  orthonormality, and Lemma \ref{PI-est},  we get 
\Be
\label{2to00}
\begin{aligned}
\int_{\mathbb R^d} |\Pi_\lambda[\psi^+_j](x,y)|^2 dy   &\lesssim 2^{-2j}  \sum_{\lambda'}\big(1+2^{-j}|\lambda-\lambda'|\big)^{-N} |\Pi_{\lambda'}(x,x)| \\
                                                                                     &\lesssim 2^{-j}(\lambda\mu)^{\frac{d-2}{2}}
\end{aligned}
\Ee
for $x\in A_{\lambda, \mu}^\sigma$.  Thus, by \eqref{2-infty} the estimate \eqref{10-2} follows.   To show \eqref{10-3}, we now work with the rescaled operator $\fP_\lambda$. 
We use the decomposition \eqref{decomp} with $\mu=\mu'$. Then,  by   Lemma \ref{large-sep} and Lemma \ref{small-sep} we may assume  
$2^{-\nu}\le 1/100$ and it is sufficient for \eqref{10-3} to show  \eqref{bbjj} with $p=1$
while we assume \eqref{AA}. In fact, the desired estimate is actually  shown in the course of the proof of Proposition \ref{prop:mm-comparable}, see Remark \ref{about-bbjj}.

To complete the proof of  \eqref{mmm-0} it remains to show
\Be
\label{final-weak}
\|\chi_{\lambda,\mu}^\sigma\Pi_{\lambda}[\psi_j^+]\chi_{\lambda,\mu}^\sigma\|_{L^{\frac{2d}{d+1},1}\rightarrow L^\infty}\lesssim (\lambda\mu)^{\frac{d-3}{4}}, \quad {2^{-j} \sim \smu}.
\Ee
Before we proceed to prove \eqref{final-weak},  we show the estimate \eqref{000} which is much easier. Indeed,   $\Pi_\lambda[\psi^0]$ corresponds to  $\Pi_\lambda[\psi^+_j]$, $2^{-j}\sim1$. Thus, repeating the same argument as before, we have \eqref{10-2} and \eqref{10-3} with $\psi_j^+$  replaced by  $\psi^0$, so \eqref{000} follows from interpolation between these two estimates.

In order to prove \eqref{final-weak}, we first note that the estimate  \eqref{final-weak} is equivalent to  
 $
\|\chi_{\mu}^\sigma\fP_{\lambda}[\psi_j^+]\chi_{\mu}^\sigma\|_{L^{\frac{2d}{d+1},1}\rightarrow L^\infty}\lesssim \lambda^{-\frac12}  \mu^{\frac{d-3}{4}}
$ (cf. \eqref{eq:norm-scaled}). 
By  \eqref{decomp},  Lemma \ref{large-sep},  and Lemma \ref{small-sep},  it is enough to show that 
\[     \| \sum_{\nu} \sum_{k\sim_\nu k'} \chi_k^{\sigma,\nu} \fP_\lambda [ \psi_j^+ ]   \chi_{k'}^{\sigma,\nu}   \|_{L^{\frac{2d}{d+1},1}\rightarrow L^\infty}\lesssim \lambda^{-\frac12}  \mu^{\frac{d-3}{4}},  \quad 2^{-j}\sim \smu
\]
under the assumption  \eqref{AA}. 
Bounds on the sum over  the case $2^{-\nu}\gg \mu$ and $ \nu=\nu_\circ$ are  easy to obtain.    By Lemma \ref{case-a} we have 
$ \sum_{2^{-\nu}\not\sim\smu} \| \chi_k^{\sigma,\nu} \fP_\lambda [ \psi_j^+ ]   \chi_{k'}^{\sigma,\nu}\|_{1\to \infty}\lesssim   \lambda^{-\frac12} \mu^{-\frac {d+1}4}$
because the sets $\{A_{k}^{\sigma,\nu} \times A_{k'}^{\sigma,\nu} \}_{k\sim_\nu k', \nu}$ are disjoint.  The estimate \eqref{10-2} and scaling  (together with Lemma  \ref{trivial1} and Lemma \ref{kkpsum})  
yield 
\[\|{\sum_{2^{-\nu}\not\sim \mu}} \sum_{k\sim_\nu k'} \chi_k^{\sigma,\nu} \fP_\lambda [ \psi_j^+ ]   \chi_{k'}^{\sigma,\nu} \|_{2\to \infty} \lesssim \lambda^{-\frac12}  \mu^{\frac{d-1}{4}}.\] 
Hence, by interpolation  we obtain 
\[{
 \| \sum_{2^{-\nu}\not\sim \mu}}\sum_{k\sim_\nu k'} \chi_k^{\sigma,\nu} \fP_\lambda [ \psi_j^+ ]   \chi_{k'}^{\sigma,\nu} \|_{{\frac{2d}{d+1}}\rightarrow \infty} \lesssim \lambda^{-\frac12}  \mu^{\frac{d-3}{4}}.
\]

Therefore, we are reduced to showing, for $k\sim_\nu k',$ $ 2^{-j}\sim \smu,$ and $2^{-\nu}\sim \mu$, 
\Be
\label{eq:core}
\|  \chi_k^{\sigma,\nu} \fP_\lambda [ \psi_j^+ ]   \chi_{k'}^{\sigma,\nu} \|_{L^{\frac{2d}{d+1},1}\rightarrow L^\infty}\lesssim \lambda^{-\frac12}  \mu^{\frac{d-3}{4}}.
\Ee
 However, the strategy we have used in the above does not work since we can not get the desirable 
 $L^1$--$L^\infty$ estimate for $\pnchi \fP_\lambda [ \psi_j^+]  \chi_{k'}^{+,\nu} $.

 We use the decomposition \eqref{bb-decomp} with $\mu=\mu'$. 
 So, we may assume that $B$, $B'$ are almost disjoint cubes of sidelength {$\sim c\vepc \mu$}. Similarly as in Lemma \ref{dychotomy}  we may assume either $|\cD(x,y)|\ll \vepc \mu^2$ or $|\cD(x,y)|\gtrsim \vepc \mu^2$ if $(x,y)\in B\times B'$. Regarding the latter case  we have $ \|\chi_B \fP_\lambda[\psi_j]  \chi_{B'}\|_{1\to \infty}\lesssim \lambda^{-\frac12}  \mu^{-\frac{d+1}{4}} $ by Lemma \ref{oscillatory}. Interpolation with $ \|\chi_B \fP_\lambda[\psi_j]  \chi_{B'}\|_{2\to \infty}\lesssim \lambda^{-\frac12}  \mu^{\frac{d-1}{4}} $, which 
 follows form  \eqref{10-2},  yields the desired bound. 

Therefore, it is enough to show 
\Be 
\label{sim}
\| \chi_B \fP_\lambda [ \psi_j ]  \chi_{B'}  \|_{L^{\frac{2d}{d+1},1}\rightarrow L^\infty}\lesssim \lambda^{-\frac12}  \mu^{\frac{d-3}{4}}, \quad 2^{-j}\sim \smu
\Ee
under the additional assumption that  $|\cD(x,y)|\ll \vepc \mu^2$ for $(x,y)\in B\times B'.$
We  may also assume that $B$ and $B'$ satisfy \eqref{positionor} with $\mu=\mu'$.   Then, by (additional) finite decomposition and  the argument which  yields  \eqref{eps-con} we may assume 
\[
|S_c(x,y)- S_c(x',y')|\le  \vepc\smu
\]
whenever $(x,y)$ and $(x',y')$ are in $B\times B'$. 

Fix $(x^\ast, y^\ast)\in B\times B'$.  To localize the operator we need the similar preparatory decomposition  which is used in Section  \ref{proof-away}.  As before, let us set  
\[\psi^{\ast}(s):= \psi_\circ \Big(\frac{s-S_c(x^\ast, y^\ast)}{2\vepc\smu}\Big)
\]  
where $\psi_\circ\in C_c^\infty(-2,2)$ and $\psi(s)=1$ if $|s|\le 1$. We break  
\[\chi_{ B} \fP_\lambda [ \psi_j ] \chi_{B'}=  \chi_{ B} \fP_\lambda [ \psi^\ast \psi_j  ] \chi_{B'}+\chi_{ B} \fP_\lambda [ (1-\psi^\ast)\psi_j  ] \chi_{B'}.\]
The second operator $\chi_{ B} \fP_\lambda [ (1-\psi^\ast)\psi_j  ] \chi_{B'}$ can be easily handled.  
Since $|\partial_s^2 \cP|\gtrsim \vepc \smu$ on   $\supp (1-\psi^\ast)\psi_j$, by Lemma \ref{osc-lemma} we have
\[ 
\|\chi_{ B} \fP_\lambda [ (1-\psi^\ast)\psi_j  ] \chi_{B'}\|_{1\rightarrow\infty}\lesssim \lambda^{-\frac12}\mu^{-\frac {d+1}4}.
\]
On the other hand, since $\mathcal F((1-\psi^\ast)\psi_j)=O(\smu \big(1+\smu|\tau|\big)^{-N} )$ and $2^{-j}\sim \smu\gtrsim (\lambda\mu)^{-1}$,  
we also have
\[
\|\chi_{ B} \fP_\lambda [ (1-\psi^\ast)\psi_j  ] \chi_{B'}\|_{2\rightarrow\infty}\lesssim  \lambda^{-\frac12} \mu^{\frac{d-1}{4}}.
\]
In fact, this can be shown in the same manner as \eqref{10-2}.  
Then, interpolation between these estimates gives the estimate 
$ \|\chi_{ B} \fP_\lambda [ (1-\psi^\ast)\psi_j  ] \chi_{B'}  \|_{2d/(d+1)\rightarrow \infty}\lesssim \lambda^{-\frac12}  \mu^{\frac{d-3}{4}}$, which is acceptable. 
Thus, to prove \eqref{sim} it is sufficient to consider $\chi_{ B} \fP_\lambda [ \psi^\ast \psi_j  ] \chi_{B'}.$

Let $\widetilde \chi_ {B}$ and  $\widetilde \chi_ {B'}$ be smooth functions adapted to $B$ and $B'$, respectively, such that they are respectively supported in
$\widetilde B=\{x: \dist(x, B)\le c \vepc\mu \}$ and $\widetilde B'=\{x: \dist(x, B')\le c\vepc \mu \}$ for some small $c>0$ and $\widetilde \chi_ {B}(x)=1$ if $x\in B$ and  $\widetilde \chi_ {B'}(y)=1$ if $y\in B'$.
The desired estimate for $\chi_{ B} \fP_\lambda [ \psi^\ast \psi_j  ] \chi_{B'}$  follows if we show 
\Be 
\label{sim-ast}
\| \widetilde  \chi_B \fP_\lambda[\psi^*\psi_j] \widetilde \chi_{B'}  \|_{L^{\frac{2d}{d+1},1}\rightarrow L^\infty}\lesssim \lambda^{-\frac12}  \mu^{\frac{d-3}{4}}, \quad 2^{-j}\sim \smu
\Ee
while we are assuming  
\[ |\cD(x,y)|\ll \vepc \mu^2, \quad (x,y)\in  \widetilde B\times \widetilde B'.\]
Actually, we obtain a little better estimate than we need.

As in Section \ref{mild-l2}, we break $\widetilde  \chi_B \fP_\lambda[\psi^*\psi_j] \widetilde \chi_{B'}$ away from the critical point $S_c$: 
\[  \widetilde  \chi_B \fP_\lambda[\psi^*\psi_j] \widetilde \chi_{B'}=\sum_l \fP_l^\ast :=\sum_{l}\widetilde \chi_B \fP_\lambda[\psi^*\psi_j \widetilde\psi (2^l(\cdot-S_c))] \widetilde \chi_{B'}. \] 
Here, we may assume   $2^{-l}\lesssim \vepc \smu$ because  $\fP^\ast_l=0$ otherwise (cf. \eqref{interval}). 
 Since $|\partial_s^2 \mathcal P|\gtrsim 2^{-l}$ on the support of $\widetilde\psi(2^{l}(\cdot-S_c))\psi_j$ (see \eqref{p2-lower}), by Lemma \ref{osc-lemma} we have
$
\|\fP^\ast_l\|_{1\rightarrow\infty}\lesssim 2^{\frac l2}\lambda^{-\frac12}\mu^{-\frac d4}.
$
Thus, for  the estimate \eqref{eq:core}  it is sufficient to show   
\Be
\label{Pi2}
\| \fP^\ast_l\|_{2\rightarrow\infty}\lesssim 2^{-\frac l2}\lambda^{-\frac12} \mu^{\frac{d-2}{4}}.
\Ee
Indeed, interpolation gives 
\[\| \fP^\ast_l\|_{p\rightarrow\infty} \lesssim  2^{l(\frac2p-\frac32)}\lambda^{-\frac12}\mu^{-\frac {d-1}p+\frac{3d-4}4}\]
 for $1\le p\le 2$. 
Since $2^{-l}\lesssim \vepc\smu$, summation along $l$ yields
$  \|\sum_{l} \fP^\ast_l\|_{p\rightarrow\infty} 
\lesssim \lambda^{-\frac12}\mu^{-\frac dp+\frac{3d-1}4}$
for $p>4/3$. In particular, we have $ \|\sum_{l} \fP^\ast_l\|_{{2d}/(d+1) \rightarrow\infty} \lesssim  \lambda^{-\frac12}  \mu^{\frac{d-3}{4}}$ for $d\ge 3$.
If $d=2$, one can use  $(b)$ in Lemma \ref{s-trick}  to get   $\|\sum_{l} \fP^\ast_l\|_{L^{\frac 43,1} \rightarrow L^\infty} \lesssim  \lambda^{-\frac12}\mu^{-\frac{1}4}$. 
Thus, we get the desired \eqref{sim-ast} for $d\ge 2$.

In order to prove \eqref{Pi2},  after changing of variables  we note 
\[\fP^\ast_l(x,y) 
 = 2^{\frac d2j} 2^{-l} 
\int e^{i\lambda \Phi_s(x,y)}A(x,y,s) ds, \]
where $\Phi_s$ is given by \eqref{phase-s} and 
\[
  A(x,y,s)=  \widetilde \chi_ {B}(x) \widetilde\chi_{B'}(y)  (2^{-\frac d2j} \mathfrak a\psi^*\psi_j)(\sxyls) \widetilde \psi(s)\] 
  (cf. \eqref{as-as}).
We recall $\sxyl$ is given by \eqref{sxyl}.  By \eqref{2-infty}, it is sufficient for \eqref{Pi2} to show 
\Be  
\label{Pi4} \int  |\fP^\ast_l(x,y)|^2 dy\lesssim 2^{-l}\lambda^{-1} \mu^{\frac{d-2}{2}}.
\Ee
Let us write   
\Be
\label{Pi3}
\int  |\fP^\ast_l(x,y)|^2 dy= 2^{-2l} 2^{dj}\iint\bigg(\int  A(x,y,s) \overline{A(x,y,t) }e^{i\lambda\Psi(x,y,s,t)}dy\bigg) dtds,
\Ee
where 
\[ \Psi(x,y,s,t):=\cP(x,y,\sxyls )-\cP(x,y,\sxylt). 
\] 

We now claim
\Be 
\label{lower-deriv}
|\partial_y \Psi(x,y,s,t)|\gtrsim 2^{-l} |s-t| 
\Ee
for $(x,y)\in \widetilde B\times \widetilde B'$ if we take  $\vepc$ sufficiently small.  In the same manner as in the proofs of  Lemma \ref{sc-derivative} and Lemma \ref{lem:ph-derivative} in Section \ref{proof-end},   it is not difficult to see 
\[ \partial_y^\alpha \Psi=O(\mu^{-|\alpha|+1}2^{-l} |s-t| ),\quad \partial_y^\alpha  A=O(\mu^{-|\alpha|})\]  
for $(x,y)\in \widetilde B\times \widetilde B'$. In fact, for the first estimate
we use the mean value theorem together with Lemma \ref{lem:ph-derivative}.  Thus, routine integration by parts (similarly as in the first part of  the proof Lemma \ref{generalized-small2})   gives \footnote{One may rescale as before, that is to say, $y\to \mu y+y_{B'}$ where $ y_{B'}$ is the center of $B'$. } 
\[
\Big| \int  A(x,y,s) A(x,y,t) e^{i\lambda\Psi(x,y,s,t)}dy \Big| \lesssim (1+\lambda \mu 2^{-l}|s-t|)^{-N} \mu^d. 
\]
Combining this  with \eqref{Pi3} and integrating in $s$,  we get \eqref{Pi4}.

It remains to show \eqref{lower-deriv}.  Let us set 
\begin{align*}
E&=
\big(\partial_s \cP( x,y,\sxyls)-\partial_s \cP(x,y,\sxylt)\big)\partial_y S_c(x,y),
\\
F&=
\partial_y \cP(x,y,\sxyls)-\partial_y \cP(x,y,\sxylt).
\end{align*}
Then, we have 
\[\partial_y \Psi(x,y,s,t)=E+F.\]
Using  \eqref{sxyl}, by the mean value theorem   
 we have 
 \[E=
\partial_s^2 \cP( x,y, S_c(x,y)+ 2^{-l}s^*)2^{-l}(s-t)\partial_y S_c(x,y)\] for some $s^*$, $|s^*|\lesssim 1$. By \eqref{rxy} and \eqref{roots},  we note  $\mathcal R(x,y,\cos (S_c(x,y)+ 2^{-l}s^*))$ 
is equal to 
\[ (\cos S_c(x,y)- \cos (S_c(x,y)+ 2^{-l}s^*))( \tau_+(x,y)-\cos (S_c(x,y)+ 2^{-l}s^*)).\] 
Since $S_c(x,y)+ 2^{-l}s^*\sim \smu$ and $S_c(x,y)\sim \smu$,  we see that  $\cos S_c(x,y)- \cos (S_c(x,y)+ 2^{-l}s^*)=O(\mu^{-\frac12})$  and  $ \tau_+(x,y)-\cos (S_c(x,y)+ 2^{-l}s^*)=O(\mu)$.
Thus, we have 
\[ \mathcal R(x,y,\cos (S_c(x,y)+ 2^{-l}s^*))=O(\mu^\frac32 2^{-l} ).\] 
 We now  recall $2^{-l}\lesssim \vepc \smu$ and  $\partial_y S_c(x,y)$$=O(\mu^{-\frac12})$ for $(x,y)\in \widetilde B\times \widetilde B'$ (\eqref{est-partial}). Thus,  by \eqref{ph-dd} 
it follows that 
$
E=O(\vepc 2^{-l} |s-t|) .$  To show \eqref{lower-deriv}  it is enough to verify 
\Be\label{eq:F} |F|\gtrsim 2^{-l} |s-t| \Ee taking $\vepc>0$ small enough. We exploit the particular form of $\cP$.   
For simplicity, 
 let us set 
$\sxyl(s):=2^{-l}s+S_c(x,y)$ fixing $x,y$.  By  a direct computation we get
\begin{align}
\label{FF}
F  
=\mathbf a y  -\mathbf b (x-y ).
\end{align}
where 
\[\mathbf a= \Big(\frac{\cos \sxyl(s) -1 }{\sin \sxyl(s) } -\frac{\cos \sxyl(t)-1}{\sin \sxyl(t) }\Big), \qquad \mathbf b= \Big(\frac{1 }{\sin \sxyl(s) } -\frac1{\sin \sxyl(t) }\Big). \] 
It is clear that $|\mathbf a|\sim  2^{-l} |s-t| $ by the mean value theorem and, similarly, $|\mathbf b|\sim   2^{-l} |s-t| \mu^{-1}$ because $\sxyl\sim \smu$. Since  $B$, $B'$ satisfy 
\eqref{positionor} with $\mu=\mu'$, the angle between 
two vectors  $y$ and $\mu^{-1}(x-y)$ are bounded below by a constant  $\sim 1$ because $x\in B$ and $y\in B'$. Since $|y|\sim 1$ and  $\mu^{-1}|x-y|\sim 1$,  we conclude \eqref{eq:F}. 
This completes the proof of Proposition \ref{prop:weakmu}.
\end{proof}

We now prove \eqref{est:brweak},  which essentially  corresponds to  the estimate  \eqref{est:muweak} with $\mu\sim 1$. 
So, the proof is similar with that of  \eqref{est:muweak} but we need to modify it slightly since $|x|$ and $|y|$ are not bounded away from the zero.

\begin{proof}[Proof of  \eqref{est:brweak}]  
By scaling (\eqref{eq:norm-scaled}), the estimate  \eqref{est:brweak} is equivalent to 
\[\|\chi_{B_r}  \fP_{\lambda}\chi_{B_r} f  \|_{\infty}\lesssim \lambda^{-\frac12} \|f\|_{\frac{2d}{d+1},1}.\] 
As in Section \ref{proof-away} we partition the set $B_r\times B_r$ into finitely many sets $\{A_n\times A_m\}_{n,m}$  of small diameter $\le c\vepc$ such that either $|\cD(x,y)|\gtrsim \vepc$ or  $|\cD(x, y)|\ll \vepc$ holds for $(x,y)\in A_n\times A_m$.  If $|\cD(x, y)|\gtrsim \vepc$ holds for $(x,y)\in A_n\times A_m$,  by Theorem \ref{thm-locest}  we have $\| \chi_{A_n}\fP_\lambda  \chi_{A_m}\|_{L^{\frac{2d}{d+1},1} \to L^\infty} \lesssim \lambda^{-1/2}$.  Hence, we only need to show 
\[
 \| \chi_{A} \fP_\lambda   \chi_{A'}  f\|_{\infty}\lesssim \lambda^{-\frac12} \|f\|_{\frac{2d}{d+1},1}
\]
 while assuming that $A, A'\subset B_r$ are of diameter $c\epsilon_\circ$  and  $|\cD(x,y)|\ll  \vepc$ for $(x,y)\in A\times A'$ (\eqref{det-small0}). 
Similarly as in the proof of  \eqref{est:muweak},  the matter reduces to showing 
\[ \|\sum_{2^{-j}\ge \lambda^{-1}} \chi_{A} \fP_\lambda [\psi_j]   \chi_{A'} f\|_{\infty}\lesssim \lambda^{-\frac{1}2} \|f\|_{\frac{2d}{d+1},1}.\]
Indeed, it is easy to show $\|\sum_{2^{-j}< \lambda^{-1}} \fP_\lambda [\psi_j]  f\|_{2d/(d+1)\to \infty}\lesssim \lambda^{-\frac{1}2} $ in the same way  as 
\eqref{hhhh} is shown (see \eqref{Pie100} and \eqref{Pie200}) 
because $|\Pi_\lambda(x,y)|\lesssim \lambda^\frac{d-2}2$.   Furthermore, using the same bound, \eqref{spectral-proj} and orthogonality between the Hermite functions as before, via rescaling (\eqref{eq:norm-scaled}) we have $\|  \fP_\lambda[\psi_j] \|_{2\to \infty} \lesssim  2^{-\frac{1}{2}j}\lambda^{-\frac12}$(for instance, see \eqref{2to00}) for $2^{-j}\ge \lambda^{-1}$, which trivially gives 
\Be \label{eq:2200}\|  \chi_{A} \fP_\lambda [\psi_j]   \chi_{A'}\|_{2\to \infty} \lesssim  2^{-\frac{1}{2}j}\lambda^{-\frac12},  \quad 2^{-j}\ge \lambda^{-1}.  \Ee

  If $(x,y)\in A\times A'$, taking sufficiently small $\vepc>0$,  we either have  $|x|\ge 1/2$  or $|y|\ge 1/2$
 because we are assuming \eqref{det-small0}. For the moment  we assume $|y|\ge 1/2$ for $y\in A'$. We consider the case $|x|\ge 1/2$  later, which can be handled similarly by exchanging the roles of $x,y$. 
 
 Since $A\subset B_{1-\delta_\circ}$, taking $\vepc\ll \delta^2_\circ $  we also have  
 \[ 
 |x-y|\gtrsim  \delta_\circ^2
 \] because $\cD(x,y)=(1-|y|^2)^2+O(|x-y|)$.  Since $\cQ(x,y,1)=|x-y|^2 \gtrsim  \delc^2$  and   \eqref{det-small0} holds, from \eqref{ph-d} we note that 
 $|\partial_s\cP|\gtrsim 2^{2j}$ on the support of $\psi_j$ provided that $2^{-j}\ll \delta_\circ^2$.  
 By Lemma \ref{osc-lemma} we have  $\| \chi_{A} \fP_\lambda [\psi_j]  \chi_{A'} \|_{1\to\infty}\lesssim  2^{\frac{d-4}{2}j}\lambda^{-1} $. 
 Combining this with the estimate \eqref{eq:2200} via interpolation, we get 
 \[
 \|\sum_{\lambda^{-1}\le 2^{-j}\ll \delta_\circ^2}  \chi_{A} \fP_\lambda [\psi_j]  \chi_{A'} f\|_\infty\lesssim  \sum_{2^{-j}\ll \delta_\circ^2} 
 2^{-\frac{3}{2d}j} \lambda^{-\frac{d+1}{2d}}  \|f\|_{\frac{2d}{d+1}}
 \lesssim  \lambda^{-\frac12}  \|f\|_{\frac{2d}{d+1}}.
  \]
So, we may assume $2^{-j}\gtrsim \delta_\circ^2$ and there are only finitely many $\fP_\lambda [\psi_j] $.  It suffices to show the estimate for each $ \chi_{A} \fP_\lambda [\psi_j]  \chi_{A'}$. 
Additionally,  with a small $c>0$ we may assume 
\Be
\label{scc-lower}
c\delc^2\le S_c(x,y) \le \pi/2
\Ee for $(x,y)\in A\times A'$.  
Otherwise, we  have $|\partial_s^2\cP|\gtrsim c>0$ on the support of $\psi_j$, so we get  $\|\chi_{A} \fP_\lambda [\psi_j]  \chi_{A'}\|_{1\to \infty}\lesssim \lambda^{-\frac12}$.
This and \eqref{eq:2200} give $ \|\chi_{A} \fP_\lambda [\psi_j]  \chi_{A'} f\|_\infty\lesssim \lambda^{-\frac12}  \|f\|_{\frac{2d}{d+1}}$.

 We now decompose 
 $\widetilde  \chi_A \fP_\lambda[\psi_j] \widetilde \chi_{A'}$  away from the critical point $\sxy$.  Since  $  |\partial_s^2 \mathcal P| \sim    2^{-l}$  on $\supp(\widetilde \psi(2^l (\cdot-S_c))\psi_j)$,  we have  the estimate  $\|  \widetilde\chi_A \fP_\lambda[\psi_j \widetilde\psi (2^l(\cdot-S_c))] \widetilde \chi_{A'}\|_{1 \to \infty}\lesssim 2^{\frac l2}\lambda^{-\frac12}$.
 By interpolation and summation over $l$ we  need only to show
  \[\|  \widetilde\chi_A \fP_\lambda[\psi_j \widetilde\psi (2^l(\cdot-S_c))] \widetilde \chi_{A'}\|_{2 \to \infty}\lesssim 2^{-\frac l2}\lambda^{-\frac12}.\] 
  We proceed as  in the proof of \eqref{Pi4}. Clearly, we may assume $2^{-l}\lesssim \vepc$. Recalling \eqref{Pi3},  it is sufficient to show \eqref{lower-deriv}  whenever $(x,y)\in A\times A'$.  To do this,  it is enough to verify \eqref{eq:F} if $(x,y)\in A\times A'$ and $\vepc$ is small enough. 
  In turn, \eqref{eq:F} follows if  we show
  \Be 
  \label{eq:dist}
  \dist\big(\frac{x-y}{|x-y|},  \pm \frac y{|y|}\big) \ge c\delta_\circ^2, \quad \forall (x,y)\in A\times A'
    \Ee 
    for some positive constant $c>0$. 
  Indeed, one may use an elementary fact that  $|a\mathbf v+ b\mathbf u|\gtrsim  \min\{|\mathbf v-\mathbf u|,|\mathbf v+\mathbf u|\}( |a|+|b|)$ whenever $\mathbf v, \mathbf u\in \mathbb S^{d-1}$ and $a,b\in \mathbb R$. By \eqref{scc-lower} it follows that   $|\mathbf a|\sim  2^{-l} |s-t| $ and $|\mathbf b|\sim   2^{-l} |s-t|$. Thus, using  \eqref{FF} and \eqref{eq:dist},  we get \eqref{eq:F} with $\vepc$ small enough.

   We now show \eqref{eq:dist}. 
     Suppose \eqref{eq:dist} fails, then we have either $\dist\big(\frac{x-y}{|x-y|},   \frac y{|y|}\big) \ll \delta_\circ^2$  or $\dist\big(\frac{x-y}{|x-y|}, - \frac y{|y|}\big) \ll \delta_\circ^2$ for some $(x,y)\in A\times A'$.   Using \eqref{angle},  we  have $\cD(x,y)\gtrsim  \delta_\circ^2$ in each case because  $x,y\in B_{1-\delta_\circ}$. This is a contradiction if we take $\epsilon_\circ\ll  \delta_\circ^2$ because we are assuming \eqref{det-small0}. 
     
     Finally, if  $|x|\ge 1/2$ for $x \in A$, we write $F=(\mathbf a+\mathbf b) (y-x)  +\mathbf  ax$ and routine adaptation of the previous argument interchanging the roles of $x,y$  
     proves \eqref{lower-deriv}.\end{proof}

\subsection{Estimates over $A_{\lambda,\mu}^-\times A_{\lambda,\mu}^-$. } 
Compared with the estimate for $\Pi_\lambda$ over the set $A_{\lambda, \mu}^{+} \times   A_{\lambda, \mu}^{+}$ which we have considered in the previous sections,  it is much easier to obtain
the estimate  over the set
$A_{\lambda,\mu}^-\times A_{\lambda,\mu}^-$. The estimate in \cite[Theorem 3]{T05} implies 
\Be
\label{kt-nu-}
 \|\chi_{\lambda, \mu}^-\Pi_\lambda \chi_{\lambda, \mu}^- \|_{q'\to q}\lesssim 
                                      \begin{cases}
\lambda^{-\frac53+\frac d2\delta(q',q)}\mu^{-2+\frac d2\delta(q',q)}, &   2\le q\le \frac{2(d+1)}{d-1},
\\
\lambda^{-\frac13+\frac d6\delta(q',q)}(\lambda^{\frac23}\mu)^{-N}, &     \frac{2(d+1)}{d-1} \le q\le \infty, 
                                      \end{cases}
\Ee
for   $\lambda^{-2/3}\le\mu\le1/4$  and  $
\|\chi_{A_{\lambda,\mu_\circ}^\circ}\Pi_\lambda \chi_{A_{\lambda,\mu_\circ}^\circ} \|_{q'\to q}\le C \lambda^{-\frac13+\frac d6\delta(q',q)}$ for $2\le q\le \infty$ where 
$\mu_\circ=\lambda^{-2/3}$.  Our approach based on  the decomposition \eqref{decomp} allows us to obtain better bounds than those  on the extended range $1\le p\le 2\le q\le \infty$.

\begin{prop}\label{prop:estext}
Let $(p,q)\in [1,2]\times [2,\infty]$. If  $\lambda^{-2/3}\le\mu\le 1/4$,  for any $N$ we have 
\begin{align}\label{est:ext--}
    \|\chi_{\lambda, \mu}^-\Pi_\lambda\chi_{\lambda, \mu}^-\|_{p\rightarrow q}\lesssim \lambda^{-\frac13+\frac d6\delta(p,q)}(\lambda\mu^{3/2})^{-N}.
\end{align}
If  $0<\mu\lesssim \lambda^{-2/3}$ and $d\ge 2$,  we have the estimate 
\begin{align}\label{smallmu}
\|\chi_{\lambda, \mu}^\circ
\Pi_\lambda \chi_{\lambda, \mu}^\circ\|_{p\rightarrow q}\lesssim  \lambda^{\frac 13+\frac{d-4}{6}\delta(p,q)}\mu^{1-\delta(p,q)}.
\end{align}
\end{prop}

Being combined with \eqref{eq:main}, those estimates give such weighted estimates as  in \cite[Theorem 3]{T05}. 
Concerning sharpness of  \eqref{smallmu},  one can show  the lower bound $ \lambda^{\frac 13+\frac{d-4}{6}\delta(p,q)}\mu^{1-\delta(p,q)} \lesssim
\|\chpl
\Pi_\lambda \chpl\|_{p\rightarrow q}$, $0<\mu\lesssim \lambda^{-2/3}$  by considering a suitable input function and the asymptotic expansion of the Hermite functions (for example, see Olver  \cite{olv}).

\begin{proof} In order to show \eqref{est:ext--}, 
by scaling \eqref{eq:norm-scaled} it is sufficient to obtain the equivalent estimate 
\begin{align}
\label{outer-}
    \|\chi_{A_{\mu}^-}\fP_\lambda\chi_{A_{\mu}^-}\|_{p\rightarrow q}\lesssim \lambda^{-\frac d2(1-\delta(p,q))}
    \lambda^{-\frac13+\frac d6\delta(p,q)}(\lambda\mu^{\frac32})^{-N}.
\end{align}
By \eqref{decomp} with $\sigma, \sigma'=-$, $(i)$ in  Lemma \ref{large-sep},  Lemma \ref{small-sep},  and \eqref{symmetric} 
we need only to deal with the cases $\dist(A_k^{-,\nu}, A_{k'}^{-,\nu}) \le 1/100$ and $\dist(A_k^{-,\nu}, -A_{k'}^{-,\nu}) \le 1/100$. By reflection  $y\to -y$, it is enough to consider 
the first case.  To show \eqref{outer-} we use the  estimates 
\Be
\label{est:--1}
\begin{aligned}
&\|\chi_k^{-,\nu}\fP_\lambda\chi_{k'}^{-,\nu}\|_{1\rightarrow \infty}
\lesssim
   2^{\frac{d-2}{4}\nu}(\lambda 2^{-\frac 32\nu})^{-N},     \quad  2^{-\nu}\gg\mu, \
    \\ 
    &\|\chi_k^{-,\nu}\fP_\lambda\chi_{k'}^{-,\nu}\|_{1\rightarrow \infty}
\lesssim
      \lambda^{\frac{d-2}{6}}(\lambda\mu^{\frac 32})^{-N},
     \qquad  2^{-\nu}\lesssim\mu,
   \end{aligned}
   \Ee
 with $k\sim_\nu k'$.  
To obtain the first estimate in  \eqref{est:--1},  by Lemma \ref{small-sep} it is sufficient to show  
\[\sum_{j} \| \chi_k^{-,\nu}\fP_\lambda[\psi_j^+] \chi_{k'}^{-,\nu}\|_{1\to \infty}\lesssim  2^{\frac{d-2}{4}\nu}(\lambda 2^{-\frac 32\nu})^{-N}\] since  we are assuming that $\dist(A_k^{-,\nu}, A_{k'}^{-,\nu}) \le 1/100$.  Splitting the cases $2^{-2j}\lesssim 2^{-\nu}$ and $2^{-2j}\gg 2^{-\nu}$,    this estimate follows from \eqref{a-case-kernel}.  The second estimate in \eqref{est:--1} is clear  from \eqref{kt-nu-} with $q=\infty$ because $\lambda^{-2/3}\le\mu\lesssim 1$.   
 
Once we have the estimate  \eqref{est:--1}, the proof of \eqref{outer-}  is straightforward. Indeed, using Lemma \ref{trivial} and  the estimate \eqref{est:--1}, we have
\[
\big\| \chi_k^{-,\nu}\fP_\lambda\chi_{k'}^{-,\nu}\big\|_{p\rightarrow q}
\lesssim 
\begin{cases}  \big(\mu 2^{-(d-1)\nu}\big)^{1-\delta(p,q)}
    2^{\frac{d-2}{4}\nu}(\lambda 2^{-\frac 32\nu})^{-N}, \quad & \mu\ll2^{-\nu},
    \\
\big(\mu 2^{-(d-1)\nu}\big)^{1-\delta(p,q)}    \lambda^{\frac{d-2}{6}}(\lambda\mu^{\frac 32})^{-N},  \quad & \mu\gtrsim 2^{-\nu}.
\end{cases}
\]
 Since $\lambda^{-2/3}\le\mu$,  summation over $\nu$ followed by a simple manipulation yields  
\begin{align*}
      \underset{\nu}{\sum}
  \max_{{k\sim_\nu k'}} \big\|
    \chi_k^{-,\nu}\fP_\lambda\chi_{k'}^{-,\nu}\big\|_{p\rightarrow q}
    \lesssim\lambda^{-\frac d2(1-\delta(p,q))}
    \lambda^{-\frac13+\frac d6\delta(p,q)}(\lambda\mu^{\frac32})^{-N}.
\end{align*}
 Thus, by Lemma \ref{kkpsum} we get the desired estimate \eqref{outer-}. 

We show \eqref{smallmu} in a similar manner. As before, we may assume $\dist(A_k^{\circ,\nu}, A_{k'}^{\circ,\nu}) \le 1/100$. It suffices  to show  
\[
\big\|\chi_k^{\circ,\nu}\fP_\lambda\chi_{k'}^{\circ,\nu}\big\|_{1\rightarrow \infty}
\lesssim
\begin{cases}
   2^{\frac{d-2}{4}\nu}(\lambda 2^{-\frac 32\nu})^{-N},      \, &2^{-\nu}\gg\lambda^{-\frac23}, \
    \\ 
      \qquad \lambda^{\frac{d-2}{6}},
     \, &2^{-\nu}\lesssim\lambda^{-\frac23},
\end{cases} 
  \]
  for $k\sim_\nu k'$. Once we have the estimate, one can obtain \eqref{smallmu} in the same way as  above. Since $\lambda^{-2/3}\gtrsim \mu$, the first case estimate is clear from  the first one in \eqref{est:--1}  and 
  the case second estimate  follows from  $(ii)$ in  Lemma \ref{ker10}.
  \end{proof}

Using the estimates \eqref{est:ext--} and \eqref{smallmu}, we can strengthen the estimate in Lemma \ref{LL22} under a mild assumption so that 
we don't have to make use of  cancellation between the operators.  These estimates can be used to provide a slightly different argument to prove the endpoint estimate. 

\begin{lem}
 \label{l2-mmj23}   Let $\lambda^{-\frac23}\lesssim\mu'\le \mu\le 1/4$ and $2^{-j}\sim \mu^{\frac12}$.  Then we have
\begin{equation} \label{jl222-}    \|   \chi_{\lambda, \mu}^\circ \Pi_\lambda [ \psi_j ]  \chi_{\lambda, \mu'}^\circ   \|_{2\to 2} \lesssim   \mmf .       
 \end{equation}
\end{lem}

\begin{proof}  We keep using the same notations in the proof of Lemma \ref{PI-est}. As before, in order to prove  \eqref{jl222-} it is convenient to consider a little stronger estimate. 
 Indeed, it  clearly suffices to show
\Be\label{stronger}   
\|\mmfpj\|_{2\to2} \lesssim  (\mu\mu')^\frac14 .
\Ee
To this end, we claim 
\Be \label{stronger1}  \|  \chi_{\lambda, \mu}^e \Pi_\lambda  \chi_{\lambda, \mu'}^e \|_{2\to 2} \lesssim (\mu\mu')^\frac14, \quad \mu\gtrsim \lambda^{-2/3},\Ee
which strengthens  the estimate in Lemma  \ref{l2-mmj}. 

 Assuming the estimate \eqref{stronger1} for the moment, we prove \eqref{stronger}. 
Since $\widehat {\psi_j}(\tau)=O(2^{-j} (1+2^{-j} |\tau|)^{-N})$,
using \eqref{spectral-proj} we have 
\[
\|\mmfpj\|_{2\to 2}\le   2^{-j} \sum_{\lambda' } \left(1+2^{-j}(\lambda-\lambda')\right)^{-N}\|\mmfp\|_{2\to 2} .
\]  
Recalling  \eqref{ell} and \eqref{inclusion}, we have $\|\mmfp\|_{2\to 2}\le\|\chi^e_{\lambda',\ell(\mu)}\Pi_{\lambda'}\chi^e_{\lambda',\ell(\mu')}\|_{2\to 2}$ 
by Lemma \ref{trivial1}.
 Thus, the estimate \eqref{stronger1} gives 
$
\|\mmfp\|_{2\to 2} 
 \lesssim 
         ( \ell(\mu)\ell(\mu'))^{\frac14}.
      $    
      Therefore,  \eqref{stronger} follows from  the inequality 
\[  2^{-j} \sum_{\lambda' } \left(1+2^{-j}(\lambda-\lambda')\right)^{-N}  ( \ell(\mu)\ell(\mu'))^{\frac14}\lesssim  (\mu\mu')^{\frac14}.\]       
As in the proof of Lemma \ref{PI-est}, a computation shows this since  
 $2^{-j}\lambda \gtrsim 1$, $2^{-j}\sim\mu^{\frac12}$, and $\mu,\mu'\gtrsim\lambda^{-2/3}$. We omit the detail. 
 
We now prove \eqref{stronger1}.   Since $ \|  \chi_{\lambda, \mu}^e \Pi_\lambda  \chi_{\lambda, \mu'}^e \|_{2\to 2}\le 
\|  \chi_{\lambda, \mu}^e \Pi_\lambda \|_{2\to 2} \|   \Pi_\lambda \chi_{\lambda, \mu'}^e \|_{2\to 2}$ and $\| \Pi_\lambda \chi_{\lambda, \mu'}^e \|_{2\to 2}=\|  \chi_{\lambda, \mu'}^e \Pi_\lambda  \|_{2\to 2}$ by duality,  it is enough to show 
\[\|\chi_{\lambda, \mu}^e \Pi_\lambda \|_{2\to 2}\le \mu^\frac14.\]
We decompose 
 \[ \chi_{\lambda, \mu}^e = \sum_{\lambda^{-2/3}\le\mu' \le \mu} 
                \chi^+_{\lambda, \mu'}+\chi_{\lambda,\lambda^{-2/3}}^\circ+  \sum_{\lambda^{-2/3}\le \mu'< 1 } \chi^-_{\lambda, \mu'} +  \sum_{1\le \mu'  }  \chi^-_{\lambda, \mu'}.\]
Since  the kernel  $ \Pi_\lambda(x,y)$ decays rapidly for $|x|, |y|\ge 2\sqrt\lambda$, by the triangle inequality we have  $ \sum_{1\le \mu'  } \|     \chi^-_{\lambda, \mu'} \Pi_\lambda \|_{2\to 2}\lesssim \lambda^{-N}$. 
Using  the estimates \eqref{l2-ktt},  \eqref{smallmu}, and  \eqref{est:ext--} successively, we see that $\|  \chi_{\lambda, \mu}^e \Pi_\lambda \|_{2\to 2}$ is bounded by  a constant times 
\begin{align*}
  \sum_{\lambda^{-2/3}\le \mu' \le \mu}  (\mu')^\frac14+\lambda^{-\frac 16}
 +\sum_{\lambda^{-2/3}\le \mu' \le 1}\lambda^{-\frac16}(\lambda(\mu')^\frac32)^{-N} 
 +\lambda^{-N}\lesssim \mu^{\frac14}
\end{align*}
as desired.   
\end{proof}

\subsection{Sharpness of the bound \eqref{est-ann}} Before closing this section we show that  the bound \eqref{est-ann} can not be improved, that is to say, 
\Be \label{lowerbd}  \|\chi_{\lambda,\mu}^+ \Pi_\lambda\chi_{\lambda,\mu}^+\|_{p\to q}\gtrsim \lambda^{\beta(p,q)}\mu^{\gamma(p,q)}. \Ee
 For this it is sufficient to show the  lower bounds \eqref{ann1}--\eqref{ann3} in Proposition \ref{lower-mu} below. 
Comparing those lower bounds one can easily see \eqref{lowerbd}.

\begin{prop}
\label{lower-mu}
Let $d\ge 2$, $\lambda\gg 1$, and $\lambda^{-\frac23}\le\mu\ll 1$. Then, we have 
\begin{align}\label{ann1}
    &\|\chi_{\lambda,\mu}^+ \Pi_\lambda \chi_{\lambda,\mu}^+ \|_{p\to q}\gtrsim(\lambda\mu)^{-1+\frac d2(\frac1p-\frac1q)}, 
   \\
    \label{ann2}
    &\|\chi_{\lambda,\mu}^+ \Pi_\lambda \chi_{\lambda,\mu}^+ \|_{p\to q}\gtrsim\lambda^{-\frac12(\frac1p-\frac1q)}\mu^{\frac12-\frac{d+3}{4}(\frac1p-\frac1q)}, 
  \\
    \label{ann3}
    &\|\chi_{\lambda,\mu}^+ \Pi_\lambda \chi_{\lambda,\mu}^+ \|_{p\to q}\gtrsim\lambda^{-\frac12+\frac d2(1-\frac1p-\frac1q)}\mu^{-\frac14+d(\frac34-\frac1p-\frac{1}{2q})}.
\end{align}
\end{prop}

The  lower bounds $\eqref{ann1}$, $\eqref{ann2}$ which show 
\eqref{lowerbd} for $\ppq\in \mathcal{R}_1\cup \overline{\mathcal{R}}_3$ can be obtained  by modifying the proofs of \eqref{local2}, \eqref{local1}. Instead of repeating the previous argument, we show  that the upper bound $C\lambda^{\beta(p,q)}\mu^{\gamma(p,q)}$ in \eqref{est-ann}  can not be improved to any better one, i.e.,  $C\lambda^{\beta(p,q)-\epsilon_1}\mu^{\gamma(p,q)+\epsilon_2}$ with any $\epsilon_1, \epsilon_2>0$. 
This can be done by making use of the known lower bound  from $L^p$ to $L^{p'}$. Indeed,  
we recall \eqref{kt-nu} of which sharpness was shown in \cite{T05}. Hence, 
by the  $TT^*$ argument we have 
\begin{align}\label{est:ann1,2}
\|\chi_{\lambda,\mu}^+ \Pi_\lambda \chi_{\lambda,\mu}^+ \|_{p\to p'}\ge C\lambda^{\beta(p,p')}\mu^{\gamma(p,p')}
\end{align}
 for any $1\le p\le 2$ with  $C$ depending only  on $d$.  
Suppose
\[   \|\chi_{\lambda,\mu}^+ \Pi_\lambda \chi_{\lambda,\mu}^+ \|_{p_0\to q_0}\lesssim \lambda^{\beta(p_0,q_0)-\epsilon_1}\mu^{\gamma(p_0,q_0)+\epsilon_2}  \] 
for some $\epsilon_1, \epsilon_2>0$, and $(1/p_0,1/q_0) \in\mathcal{R}_1\cup \overline{\mathcal{R}}_3$ with $p_0'\neq q_0$. By duality it follows that  $  \|\chi_{\lambda,\mu}^+ \Pi_\lambda \chi_{\lambda,\mu}^+ \|_{q_0'\to p_0'}\lesssim \lambda^{\beta(p_0,q_0)-\epsilon_1}\mu^{\gamma(p_0,q_0)+\epsilon_2}.$  Interpolation between these two estimates particularly gives  $\|\chi_{\lambda,\mu}^+ \Pi_\lambda \chi_{\lambda,\mu}^+ \|_{p_\ast \to p_\ast'}\lesssim \lambda^{\beta(p_\ast ,p_\ast')-\epsilon_1}\mu^{\gamma(p_\ast ,p_\ast')+\epsilon_2}$ where $1/p_\ast-1/p_\ast'=1/p_0-1/q_0$ because $\beta(p_0, q_0)=\beta(p_\ast ,p_\ast')$ and $\gamma(p_0, q_0)=\gamma(p_\ast ,p_\ast')$. 
This contradicts to the lower bound  \eqref{est:ann1,2} if we let $\lambda\to \infty$ and $\mu\to 0$. 

The rest of this section is devoted to proving  \eqref{ann3}.  To do so, we need the following two lemmas.

\begin{lem}\label{lem:annn3}  Let $1\le p\le 2$ and  $x_0\in A_{\lambda,\mu}^+$. Then, if $\lambda^{-\frac23}\lesssim \mu\le 2^{-2}$,  we have
\[
\|\Pi_\lambda(x_0,x)\|_{L^p(A_{\lambda,\mu}^+)}\lesssim
\lambda^{-\frac12+\frac{d}{2p}}\mu^{-\frac14+d(\frac1p-\frac14)}.
\]
\end{lem}
\begin{proof}
Since $\Pi_\lambda (\mU x, \mU y)=\Pi_\lambda(x,y)$, by rotation we may assume $x_0=|x_0| e_1$. 
Let $S_{0}=\{x\in A_{\lambda,\mu}^+ : \lambda^{-\frac12}|x-x_0|< 100\mu\}$
and 
\[S_j =\{x\in A_{\lambda,\mu}^+ : \lambda^{-\frac12}|x-x_0|\in[ 100\mu 2^{j-1}, 100\mu 2^{j})\}, \quad j\ge 1.\]
To obtain the inequality, it is enough to show the following: 
\begin{align}
\label{est1:ann3}
    & \|\Pi_\lambda(x_0,\cdot)\|_{L^p(S_0)}
    \lesssim
    \lambda^{-\frac12+\frac{d}{2p}}\mu^{-\frac14+d(\frac1p-\frac14)}, \\
 \label{est2:ann3}
    & \|\Pi_\lambda(x_0,\cdot)\|_{L^p(S_j)}
    \lesssim
    (\mu 2^j) ^{-\frac{d-2}{4}}(\lambda(\mu 2^j)^{\frac32})^{-N}
    (\lambda^{\frac d2}\mu^d 2^{(d-1)j})^\frac1p,
    \quad j \ge 1.
\end{align}
Once we have these estimates,  summation over $j$ proves Lemma \ref{lem:annn3} because $\lambda\mu^\frac32\gtrsim 1$.
The first estimate follows by using H\"older's inequality and the estimate \eqref{eq:2to00}. 
Indeed, since $x_0\in A_{\lambda,\mu}^+$ and \[
\|\chi_{\lambda,\mu}^+\Pi_\lambda\chi_{\lambda,\mu}^+\|_{2\rightarrow\infty}^2 
=\Big\|\int_{A_{\lambda,\mu}^+} \Pi_\lambda(\cdot,y)^2 \,dy
\Big\|_{L^\infty(A_{\lambda,\mu}^+)},
\]
 the left hand side of \eqref{est1:ann3} is bounded by 
\[|S_0|^{\frac1p-\frac 12} \|\Pi_\lambda(x_0,\cdot)\|_{L^2(A_{\lambda, \mu}^+)}\lesssim   |S_0|^{\frac1p-\frac 12}\|\chi_{\lambda,\mu}^+ \Pi_\lambda \chi_{\lambda,\mu}^+\|_{2\to\infty}.
\] Thus, \eqref{est1:ann3}  follows from    \eqref{eq:2to00}. 

For the second estimate \eqref{est2:ann3}  we observe 
\[\|\Pi_\lambda(x_0,\cdot)\|_{L^p(S_{j})}\lesssim \text{sup}_{k\sim_{\nu}k'}\|\chi_k^{+,\nu}\mathfrak P_\lambda \chi_{k'}^{+,\nu}\|_{1\to\infty}\,|S_{j}|^\frac1p,\]
where  $\nu$ is a number such that $2^{-\nu}\sim 100\mu 2^{j}$. 
Combining this with \eqref{est:--1}, we get the desired estimate  \eqref{est2:ann3}.
\end{proof}

\begin{lem}
\label{lem:efderiv}
    Let $\lambda\in 2\N_0+d$ and $\mu\in[\lambda^{-\frac23},\frac14]$. 
    Suppose that $h\in \mathcal S(\mathbb R^d)$ is an  eigenfunction of $H$ with eigenvalue $\lambda$, {\it i.e.,} $H h(x) = \lambda h(x)$. If  
    $
    y_0\in A_{\lambda,\mu}^+, 
    $
    then  for any $\alpha\in \N_0^d$ we have 
    \Be
    \label{eq:deriv} |\partial_y^\alpha h(y_0)|\le C (\lambda\mu)^{\frac{|\alpha|}{2}}\|h\|_{L^\infty(B(y_0,2(\lambda\mu)^{-\frac12}))}\Ee 
    with $C$ independent of $\lambda$, $\mu$ and $h$.
\end{lem}

\begin{proof}
 To show \eqref{eq:deriv} we may assume  
 \Be \label{normal} \|h\|_{L^\infty(B(y_0,2(\lambda\mu)^{-\frac12}))} = 1.\Ee
    Let $\varphi$ be a smooth cutoff function such that $\supp \varphi \subset B(0,2)$ and $\varphi\equiv 1$ on $B(0,1)$. 
    Let us set $\phi_{y_0}(y)= \varphi(\sqrt{\lambda\mu} (y-y_0))$. 
     Then,  by Fourier inversion  we write
    \begin{align*}
    \partial_x^\alpha h(y_0) = \partial_x^\alpha(\varphi_{y_0}h)(y_0)
    = (2\pi)^{-d} \iint (i\xi)^{\alpha} &\phi_{y_0}(y)e^{i(y_0-y)\cdot \xi} h(y) dy d\xi.
    \end{align*}
  Setting $\varphi_{K}(\xi)=  \varphi({\xi}/{K\sqrt{\lambda\mu}})$ with a large positive constant $K$, we split the integral 
  \begin{align*}
    \partial_x^\alpha h(y_0) = \mathrm I+ \mathrm {I\!I},
    \end{align*}
    where 
    \begin{align*}
   \mathrm I  &:=  (2\pi)^{-d} \iint (i\xi)^{\alpha} {\varphi_K}(\xi)  \phi_{y_0}(y) e^{i(y_0-y)\cdot \xi} h(y) dy d\xi ,
    \\
     \mathrm {I\!I}&:= (2\pi)^{-d} \iint (i\xi)^{\alpha} (1-\varphi_K(\xi) ) \phi_{y_0}(y)   e^{i(y_0-y)\cdot \xi} h(y) dy d\xi.
    \end{align*}
Concerning  $\mathrm I$,  we immediately  see
    $
    |{\mathrm I}|\le  C(\lambda\mu)^{\frac{|\alpha|}{2}}$ by our choice of $\varphi_K$ and $\phi_{y_0}$.    So, we need only to consider $\mathrm {I\!I}$.     
    
      Since  $(|\xi|^2+\Delta_y)h(y)=(|\xi|^2+{|y|^2}-\lambda )h(y)$, we may  write
    \[
    \mathrm {I\!I} 
    = (2\pi)^{-d} \iint {(i\xi)^\alpha}\frac{(1- \varphi_K(\xi) ) \phi_{y_0}(y) }{|\xi|^2+|y|^2-\lambda} e^{i(y_0-y)\cdot \xi} (|\xi|^2+\Delta_y)h(y) dy d\xi.
    \]
    By integration by parts, we see 
    \[
   \mathrm {I\!I} = C\iint \xi^{\alpha} (1-\varphi_K(\xi)) A_1(y,\xi) e^{i(y_0-y)\cdot \xi} h(y) dy d\xi
    \]
for a constant $C$     where
    \[
    A_1(y,\xi):=\Delta_y\bigg(\frac{ \phi_{y_0}(y) }{|\xi|^2+|y|^2-\lambda}\bigg)-2i\nabla_y\bigg(\frac{ \phi_{y_0}(y)}{|\xi|^2+|y|^2-\lambda}\bigg)\cdot\xi.
    \]
    We now note that $|\xi|^2+|y|^2-\lambda > \frac{K}{2}\lambda\mu+\frac{|\xi|^2}{2}$ for $(y,\xi)\in{\supp( \phi_{y_0})\times\supp(\varphi_K)}$ with a sufficiently large $K$ and that $|\partial_y^\alpha \phi_{y_0}(y)|\lesssim (\lambda\mu)^{\frac{|\alpha|}{2}}$ for any $\alpha\in\N_0^d$. Thus, we have  
   $
    |A_1(y,\xi)|\lesssim (\lambda\mu)^{\frac12}{|\xi|^{-1}}.
    $
    Repeating this $N$ times yields 
    \[
     \mathrm {I\!I} = C_N\iint  \xi^{\alpha} (1-{\varphi_K}(\xi)) A_N(y,\xi) e^{i(y_0-y)\cdot \xi} h(y) dy d\xi
    \]
    for a constant $C_N$ 
    where $A_N$ is a smooth function such that ${\supp} A_N \subset \supp (\phi_{y_0}\otimes\varphi_K)$ and $|A_N(y,\xi)|\lesssim |\xi|^{-N}(\lambda\mu)^{\frac N2}$.
   Using  \eqref{normal}, we consequently have 
    \begin{align*}
    | \mathrm {I\!I} |
   \lesssim (\lambda\mu)^{\frac N2}\int_{K\sqrt{\lambda\mu}\le |\xi|}  |\xi|^{|\alpha|-N}  d\xi  \int_{|y-y_0|\le 2/\sqrt{\lambda\mu}} dy  \lesssim  (\lambda\mu)^{\frac{|\alpha|}{2}}.
    \end{align*}
    Therefore, we obtain \eqref{eq:deriv}.
\end{proof}

\begin{proof}[Proof of \eqref{ann3}]  We claim that  there is a point  $x_0\in A_{\lambda,\mu}^+$ such that   
\Be
\label{eq:x0x0}
\int_{A_{\lambda,\mu}^+} \Pi_\lambda(x_0,y)^2 \,dy
\gtrsim  \mu^{\frac12}(\lambda\mu)^{\frac{d-2}{2}}.
\Ee
This in turn shows the sharpness of \eqref{eq:2to00} because the kernel $\Pi_\lambda(x,y)$ is smooth.
Assuming \eqref{eq:x0x0} for the moment we prove \eqref{ann3}.  
Let us set 
\[f(x)=\Pi_\lambda(x_0,x)\chi_{A_{\lambda,\mu}^+}(x).\] 
Then, \eqref{eq:x0x0} gives $\Pi_\lambda f(x_0)\gtrsim \mu^{\frac12}(\lambda\mu)^{\frac{d-2}{2}}$. By Lemma \ref{ker10} and \ref{lem:efderiv} we also have 
$\|\nabla\,\Pi_\lambda f\|_{L^\infty(A_{\lambda,\mu}^+)}\lesssim \mu^{\frac12}(\lambda\mu)^{\frac{d-1}{2}}$. 
Thus, by the mean value theorem 
we see that there exists a constant $c>0$ small enough such that $\Pi_\lambda f(x)\gtrsim\mu^{\frac12}(\lambda\mu)^{\frac{d-2}{2}}$ for any $x\in B(x_0,c(\lambda\mu)^{-\frac12})$.
Thus we get  
\[
\|\chi_{\lambda,\mu}^+\Pi_\lambda\chi_{\lambda,\mu}^+ f\|_q\ge
\|\Pi_\lambda f\|_{L^q\big(B(x_0,c(\lambda\mu)^{-\frac12})\cap A_{\lambda, \mu}^+\big)}
\gtrsim \mu^{\frac12}(\lambda\mu)^{\frac{d-2}{2}-\frac{d}{2q}}.
\]
Combining this and the estimate $\|f\|_p\lesssim \lambda^{-\frac12+\frac{d}{2p}}\mu^{-\frac14+d(\frac1p-\frac14)}$, which follows from Lemma \ref{lem:annn3}, we obtain \eqref{ann3}.

It remains to show \eqref{eq:x0x0} of which proof is similar to that of \eqref{local2}. 
Let $\lambda=2N +d$ as in  {\it Proof of \eqref{local2}}. We set
$$ J=\big\{\alpha\in\mathbb N^d \,: \,|\alpha|=N,\; {N\mu}/(16d) \le\alpha_j\le {N\mu}/(8d),  \ \  2\le j\le d\big\},$$
$\ell={(2\sqrt{d})^{-1}}{\sqrt{\lambda\mu}}$, and $Q_{\ell}=[ \sqrt{\lambda}(1-2\mu),\sqrt{\lambda}(1-3\mu/2)]\times  [-\ell,\ell\,]^{d-1}$.
Noting that $Q_{\ell}\subset A_{\lambda,\mu}^+$ and $|J|\sim (\lambda\mu)^{d-1}$,  we have
$$\sum_{\alpha\in J} \int_{B_{d-1}(0,(\lambda\mu)^{-\frac12})}\int_{\sqrt{\lambda}(1-2\mu)}^{\sqrt{\lambda}(1-3\mu/2)} \left|\Phi_\alpha (x_1,\bar x)\right|\, dx_1d{\bar x}\gtrsim \lambda^{\frac d4}\mu^{\frac{d+2}{4}}.$$
This is an easy consequence of Lemma \ref{S4-asym}. 
Thus, there exists $x_0\in [\sqrt{\lambda}(1-2\mu),\sqrt{\lambda}(1-3\mu/2)]\times B_{d-1}(0,(\lambda\mu)^{-\frac12})$ such that $\sum_{\alpha\in J}|\Phi_\alpha(x_0)|\gtrsim (\lambda\mu)^{\frac{3}{4}d-1}$.
Now we define  the function $g$ by setting
\[g(x)=\chi_{Q_{\ell}}(x)\sum_{\alpha\in J}c_\alpha \Phi_\alpha(x).\]
Here $c_\alpha\in \{-1,1\}$, $\alpha\in J$ such that $c_\alpha\Phi_\alpha(x_0)=|\Phi_\alpha(x_0)|$.
Using $\|\chi_{\lambda,\mu}^+\Pi_\lambda\|_{2\to 2} \lesssim \mu^{1/4}$ which is equivalent to  \eqref{l2-ktt}, we have $\|g\|_2\le\|\sum_{\alpha\in J}c_\alpha \Phi_\alpha(x)\|_{L^2(A_{\lambda,\mu}^+)}\lesssim\mu^{\frac14}\|\sum_{\alpha\in J}c_\alpha \Phi_\alpha(x)\|_2\lesssim \mu^{\frac14}(\lambda\mu)^{\frac{d-1}{2}}$.
Recalling the definition of $A_{\alpha,\beta}$ (in the proof of \eqref{local2}), we write
\[
    \Pi_{\lambda}g(x) 
      ={\mathrm I} (x)+\mathrm {I\!I}(x) :
      =  \sum_{\alpha\in J}c_{\alpha}A_{\alpha, \alpha}\Phi_{\alpha}(x)+ \sum_{\alpha\in J} \  \ \sum_{\beta: |\beta|=N, \alpha\not= \beta} c_{\alpha} A_{\alpha, \beta} \Phi_{\beta}(x) .
\]
Following the same argument used to treat the second sum in \eqref{sumsum},  one can easily show $\mathrm{I\!I}(x)=O\big((\lambda\mu)^a e^{-b\lambda\mu}\big)$ for some $a,b>0$. Now we deal with the first sum, ${\mathrm I}(x)$.
By Lemma \ref{S4-asym} we have $A_{\alpha,\alpha}\sim\mu^{\frac12}$.
Hence, our choices of $c_\alpha$ and $x_0$ ensures that there exists $C>0$ such that ${\mathrm I}(x_0)=\sum_{\alpha\in J}c_{\alpha}A_{\alpha, \alpha}\Phi_{\alpha}(x_0)\ge C\mu^{\frac12}(\lambda\mu)^{\frac34 d-1} $. 
Since $Q_{\ell}\subset A_{\lambda, \mu}^+$ and $x_0\in A_{\lambda, \mu}^+$, 
\begin{align*}
\inp{\chi_{\lambda,\mu}^+(x_0)\Pi_\lambda(x_0,\cdot) \chi_{\lambda,\mu}^+}{g}=\Pi_\lambda g(x_0)\ge
C\mu^{\frac12}(\lambda\mu)^{\frac34 d-1}-O((\lambda\mu)^a e^{-b\lambda\mu}).
\end{align*}
Since $\|g\|_2\lesssim \mu^{\frac14}(\lambda\mu)^{\frac{d-1}{2}}$, by duality we get \eqref{eq:x0x0} as desired. 
\end{proof}

\subsection*{Acknowledgements}
This work was  supported by the POSCO Science Fellowship and the NFR (Republic of Korea) grants  no. NRF-2020R1F1A1A01048520 (E. Jeong) and NRF-2018R1A2B2006298 (S. Lee and J. Ryu).

\end{document}